\documentclass[11pt]{report} % use larger type; default would be 10pt

\usepackage[utf8]{inputenc} % set input encoding (not needed with XeLaTeX)
\usepackage[papersize={8.5 in, 11 in}, nohead, includeheadfoot, left=1.5 in, right = 1 in, vmargin= 1 in]{geometry}
\usepackage[doublespacing]{setspace}
\usepackage{bm}
\usepackage{amsthm, amsfonts, amssymb, amsmath,enumitem, bbm, mathabx}
\newtheorem{thm}{Theorem}
\newtheorem{example}{Example}

\newtheorem{lem}{Lemma}
\newtheorem{prop}{Proposition}
\newtheorem{cor}{Corollary}
\newtheorem{obs}{Observation}
\newtheorem{definition}{Definition}
\newcounter{myalgctr}
\newenvironment{rem}{%      define a custom environment
   \vskip1mm\indent%         create a vertical offset to previous material
   \refstepcounter{myalgctr}% increment the environment's counter
   \textbf{Remark \themyalgctr}% or \textbf, \textit, ...
   }{\hfill$\diamond$\par}  %          create a vertical offset to following material
\numberwithin{myalgctr}{section}
\providecommand{\norm}[1]{\left\lVert#1\right\rVert}

\usepackage{graphicx} % support the \includegraphics command and options
\graphicspath{{Images/}}
\usepackage{color, float}
 \usepackage[parfill]{parskip} % Activate to begin paragraphs with an empty line rather than an indent
\usepackage{longtable}
\usepackage{dcolumn}
\usepackage{tabularx}
\usepackage{rotating}
\usepackage{xtab}
\usepackage{mathrsfs}
\usepackage{booktabs} % for much better looking tables
\usepackage{array} % for better arrays (eg matrices) in maths
\usepackage{verbatim} % adds environment for commenting out blocks of text & for better verbatim
\usepackage{subfig,tcolorbox} % make it possible to include more than one captioned figure/table in a single float
\usepackage{siunitx} % rounding number in table usign type S
\usepackage[colorlinks=true, a4paper=true, pdfstartview=FitV,
linkcolor=blue, citecolor=blue, urlcolor=blue]{hyperref}
\makeatletter
\def\namedlabel#1#2{\begingroup
    #2%
    \def\@currentlabel{#2}%
    \phantomsection\label{#1}\endgroup
}
\makeatother
\newcommand{\vertiii}[1]{{\left\vert\kern-0.25ex\left\vert\kern-0.25ex\left\vert #1 
    \right\vert\kern-0.25ex\right\vert\kern-0.25ex\right\vert}}

\newcommand{\dto}{\stackrel{D}{\to}}
\newcommand{\cE}{\mathcal{E}}
\renewcommand{\leq}{\leqslant}
\renewcommand{\geq}{\geqslant}

\newcommand{\bh}{{\boldsymbol{\theta}}}
\newcommand{\bg}{{\boldsymbol{\gamma}}}
\newcommand{\eps}{{\varepsilon}}
\newcommand{\vp}{{\phi}}
\newcommand{\cM}{\mathcal{M}}
\newcommand{\sa}{{\mathcal{C}}}
\newcommand{\e}{{\mathbb{E}}}
\newcommand{\cs}{{\bar{X}_N}}
\newcommand{\bs}{{\boldsymbol{X}}}

\newcommand{\ba}{\bm A}
\newcommand{\bt}{{\boldsymbol{\tau}}}

\newcommand{\os}{{\overline{X}_N}}
\newcommand{\cp}{{\mathscr{C}_p}}

\newcommand{\sC}{\mathscr{C}}
\newcommand{\E}{\mathbb{E}}
\newcommand{\cC}{\mathcal{C}}
\newcommand{\ma}{{\mathcal{A}}}
 
\newcommand{\be}{\mathbb{E}}
\newcommand{\R}{{\mathbb{R}}} 
\newcommand{\mysquare}[1][black]{\small\textcolor{#1}{\ensuremath\blacksquare}}

\newcommand{\cR}{{\mathcal{R}}}
\newcommand{\mf}{{\bm F}}
\newcommand{\bp}{\mathbb{P}}
\newcommand{\p}{{\mathbb{P}}} 
\newcommand{\cA}{\mathcal{A}}
\newcommand{\cI}{\mathcal{I}}
\newcommand{\pto}{\stackrel{P}{\to}}

\usepackage[subfigure]{tocloft}
\usepackage{import}
\usepackage[compact]{titlesec}
\usepackage{subfig}
\usepackage{tikz}
\usepackage{xcolor}

\newcommand{\cH}{\mathcal{H}}
\newcommand{\bw}{{\beta_{\mathrm{CW}}^*(p)}}
\newcommand{\cB}{\mathcal{B}}
\usepackage{mathtools}
\mathtoolsset{showonlyrefs}
 \usepackage{appendix}
 \usepackage{etoolbox}
 \usepackage{chngcntr}
 \AtBeginEnvironment{subappendices}{%
 \begin{center}
 \textbf{APPENDIX}
 \end{center}
 }
\newenvironment{preliminary}{}{}
\titlespacing*{\chapter}{0pt}{-33 pt}{6 pt} % The key value here is the -33 pts, I got to it by old fashioned measuring with a ruler....

\begin{document}
% Defining variables to be used throughout the document for personalization
\def\mytitle{LIMIT THEOREMS FOR DEPENDENT COMBINATORIAL DATA, WITH APPLICATIONS IN STATISTICAL INFERENCE} % Make sure this is in all caps
\def\myauthor{Somabha Mukherjee}
\def\myauthorfull{Somabha Mukherjee}
\def\mysupervisorname{Bhaswar B. Bhattacharya}
\def\mysupervisortitle{Assistant Professor of Statistics}
\newlength{\superlen}   % a "scratch" length
\settowidth{\superlen}{\mysupervisorname, \mysupervisortitle} % Width of signature line for supervisor
\def\gradchairname{Dylan Small}
\def\gradchairtitle{Class of 1965 Wharton Professor of Statistics}
\newlength{\chairlen}   % a "scratch" length
\settowidth{\chairlen}{\gradchairname, \gradchairtitle} % Width of signature line for supervisor
\newlength{\maxlen}
\setlength{\maxlen}{\maxof{\superlen}{\chairlen}}
\def\mydepartment{Statistics}
\def\myyear{2021}
\def\signatures{46 pt} % Space to accommodate the signatures, you can fiddle with this as you like

%% PRELIMINARY PAGES

\pagenumbering{roman}
\pagestyle{plain}

% TITLE PAGE

\begin{titlepage}
	\thispagestyle{empty} % No page numbers on title page, as per Manual page 8
	\begin{center}

 \onehalfspacing

\mytitle

\myauthor

A DISSERTATION

in 

\mydepartment 

For the Graduate Group in Managerial Science and Applied Economics 

Presented to the Faculties of the University of Pennsylvania

in 

Partial Fulfillment of the Requirements for the

Degree of Doctor of Philosophy

\myyear

\end{center}

\vfill % Here to make sure the page is filled

\begin{flushleft}

Supervisor of Dissertation\\[\signatures] % Space for signature, you can fiddle with this as you like

\renewcommand{\tabcolsep}{0 pt}
\begin{table}[h]
\begin{tabularx}{\maxlen}{l}
\toprule
\mysupervisorname, \mysupervisortitle\\ %Space between advisor and graduate chair, you can fiddle with this as you like
\end{tabularx}
\end{table}

Graduate Group Chairperson\\[\signatures] % Space for signature, you can fiddle with this as you likee

\begin{table}[h]
\begin{tabularx}{\maxlen}{l}
\toprule
\gradchairname, \gradchairtitle\\ %Space between advisor and graduate chair, you can fiddle with this as you like
\end{tabularx}
\end{table}
\singlespacing

Dissertation Committee % No signature necessary

Bhaswar B. Bhattacharya, Assistant Professor of Statistics

Robin Pemantle, Professor of Mathematics

Nancy R. Zhang, Ge Li and Ning Zhao Professor, Professor of Statistics

%Dylan Small, Class of 1965 Wharton Professor of Statistics

\end{flushleft}

\end{titlepage}

% COPYRIGHT NOTICE (optional)

\doublespacing

\thispagestyle{empty} % No page number as per Manual, p. 11

\vspace*{\fill}

\begin{flushleft}
\mytitle

 \copyright \space COPYRIGHT
 
\myyear

\myauthorfull\\[24 pt] % If traditional copyright then delete everything below here, but keep \end{flushleft}

This work is licensed under the \\
Creative Commons Attribution \\
NonCommercial-ShareAlike 3.0 \\
License

To view a copy of this license, visit

\url{http://creativecommons.org/licenses/by-nc-sa/3.0/}

\end{flushleft}
\pagebreak 

\begin{preliminary}

% DEDICATION (OPTIONAL)
\setcounter{page}{3}  %Makes this page Roman numeral 3. Thanks to Albert Zevelev.
\begin{center}
\textit{Dedicated to my parents, Tapas Mukherjee and Ruby Mukherjee}
\end{center}

% ACKNOWLEDGEMENT (OPTIONAL)

\clearpage
% \chapter*{ACKNOWLEDGEMENT}
\begin{center}
\bf\large ACKNOWLEDGEMENT
\end{center}
\addcontentsline{toc}{chapter}{ACKNOWLEDGEMENT} % This is to include this section in the Table of Contents

First of all, I would like to thank my parents (Tapas Mukherjee and Ruby Mukherjee) for their unmatched effort in bringing me up to where I am today. Without their support and enthusiasm, it would not have been possible for me to pursue higher studies outside my home country, India, for the past five years. I would also like to thank my elder sister (Pausali Mukherjee) for always being my role model.

Higher studies and research are not all about talent and passion for the subject; it requires a thorough understanding of the academic system, too. In this context, I was extremely fortunate to have an advisor like Bhaswar, who not only gave me an extremely comprehensive knowledge of the research world, but also taught me how actual research should be done, by training me to handle all sorts of difficulties that arise while working on real problems. Bhaswar is more than an advisor to me. At countless times he seemed like an elder brother to me, who would give me incomparable mental and moral support to uplift my spirit, whenever I was going through difficult times during my PhD career. Doing research successfully throughout the past five years would not have been possible without his priceless help and advice. 

I was very fortunate to collaborate with/learn from many wonderful Professors at Penn, such as Nancy, Robin, Dylan, Zongming and Jian. My independent research on statistical physics under Robin was the main factor that ignited my passion in this subject. The knowledge that I gathered from that course helped me immensely during my later works on inference in statistical physics models, which form the main part of this thesis. I gained my first experience in working on real statistical problems with concrete applications through a collaboration with Nancy. I would also like to thank Dylan, Zongming and Jian for sparing a lot of their time on academic discussions with me and giving me advice on various issues. I am also grateful to Professor Bodhisattva Sen from Columbia University and Professor Ayanendranath Basu from the Indian Statistical Institute for giving me research experience during my undergraduate and master years.  

I would specially like to thank two of my collaborators, Jaesung Son and Divyansh Agarwal for their wonderful contributions in my joint works with them. Both Jaesung and Divyansh are extremely efficient researchers, and it was a priviledge for me to work with them. In fact, Chapters \ref{curiech} and \ref{ch:generalmple} in this thesis are outcomes of my joint work with Jaesung and Bhaswar. I was also fortunate to have collaborators like Marcus Michelen, Stephen Melczer and Arun Kuchibhotla who were my seniors at Penn, and are now assistant professors in reputed institutes. I would specially like to thank Arun for teaching me a lot of things in mathematical statistics. He was a person whom I could ask for help at any time, and in spite of being an extremely busy researcher himself, he would always promptly agree to help me. I would also like to thank my collaborators Sagnik Halder, Debaprartim Banerjee, Professor Sumit Mukherjee, Professor George Michailidis, Professor Rohit Patra and Professor Andrew Johnson for their crucial roles in my research career. Chapter \ref{chap:covariateising} is created out of my ongoing joint work with Sagnik, George and Bhaswar.

Finally, I would like to mention three persons who influenced me the most in my academic life. One of them was my father's co-worker Udayaditya Bhattacharya who introduced me to olympiad level mathematics for the first time when I was in high-school. Without him, I would not even be pursuing my career in mathematics and statistics. The second person is Professor S.M. Srivastava, who opened my eyes to the divine beauty of Mathematics by instilling within me a deep passion for logic and axiomatic set theory. At times, when he was teaching, it felt like I was learning from the great Georg Cantor or Kurt G\H{o}del. He is the only reason behind my sublime love for mathematics. And the third person is Professor Alok Goswami who was the only person from whom I learned everything in Probability during my bachelor and master days, starting from counting balls in boxes to It\^{o} calculus and general Markov processes. I owe my research career in probability to him.

\vspace{0.2in}
\hfill Thank you,

\hfill Somabha Mukherjee.

% ABSTRACT

\clearpage
% \chapter*{ABSTRACT}
\begin{center}
ABSTRACT
\end{center}
\addcontentsline{toc}{chapter}{ABSTRACT} % This is to include this section in the Table of Contents
\begin{center}
\mytitle

\myauthor

\mysupervisorname

\end{center}

The Ising model is a celebrated example of a Markov random field, which was introduced in statistical physics to model ferromagnetism. More recently, it has emerged as a useful model for understanding dependent binary data with an underlying network structure. This is a discrete exponential family with binary outcomes, where the sufficient statistic involves a quadratic term designed to capture correlations arising from pairwise interactions. However, in many situations the dependencies in a network arise not just from pairs, but from peer-group effects. A convenient mathematical framework for capturing higher-order dependencies, is the $p$-tensor Ising model, which is a discrete exponential family where the sufficient statistic consists of a multilinear polynomial of degree $p$. This thesis develops a framework for statistical inference of the natural parameters in $p$-tensor Ising models. We begin with the Curie-Weiss Ising model, where every $p$-tuple of nodes interact with equal strengths, where we unearth various non-standard phenomena in the asymptotics of the maximum-likelihood (ML) estimates of the parameters, such as the presence of a critical curve in the interior of the parameter space on which these estimates have a limiting mixture distribution, and a surprising superefficiency phenomenon at the boundary point(s) of this curve. However, ML estimation fails in more general $p$-tensor Ising models due to the presence of a computationally intractable normalizing constant. To overcome this issue, we use the popular maximum pseudo-likelihood (MPL) method, which avoids computing the inexplicit normalizing constant based on conditional distributions. We derive general conditions under which the  MPL estimate is $\sqrt{N}$-consistent, where $N$ is the size of the underlying network. Our conditions are robust enough to handle a variety of commonly used tensor Ising models, including spin glass models with random interactions and the hypergraph stochastic block model. Finally, we consider a more general Ising model, which incorporates high-dimensional covariates at the nodes of the network, that can also be viewed as a logistic regression model with dependent observations. In this model, we show that the parameters can be estimated consistently under sparsity assumptions on the true covariate vector.

% TABLE OF CONTENTS

\clearpage
\tableofcontents

% LIST OF TABLES

%\clearpage
%\listoftables
%\addcontentsline{toc}{chapter}{LIST OF TABLES}

% LIST OF ILUSTRATIONS

\clearpage
\listoffigures
\addcontentsline{toc}{chapter}{LIST OF ILLUSTRATIONS}

% PREFACE (OPTIONAL)

% \clearpage
% \chapter*{PREFACE (optional)}
% \addcontentsline{toc}{chapter}{PREFACE} % This is to include this section in the Table of Contents
\end{preliminary}
% \chapter*{Abstract}
% \import{Chapters/}{abstract}
% \chapter*{Dedication}
% \chapter*{Acknowledgements}
% \tableofcontents
\newpage
\pagenumbering{arabic}
\pagestyle{plain} % This has to be repeated here because the lists change the style
\chapter{Introduction}\label{chap:introduction}
\begin{figure}%
	\centering
	\subfloat[\centering]{{\includegraphics[width=7.1cm,height=6.5cm]{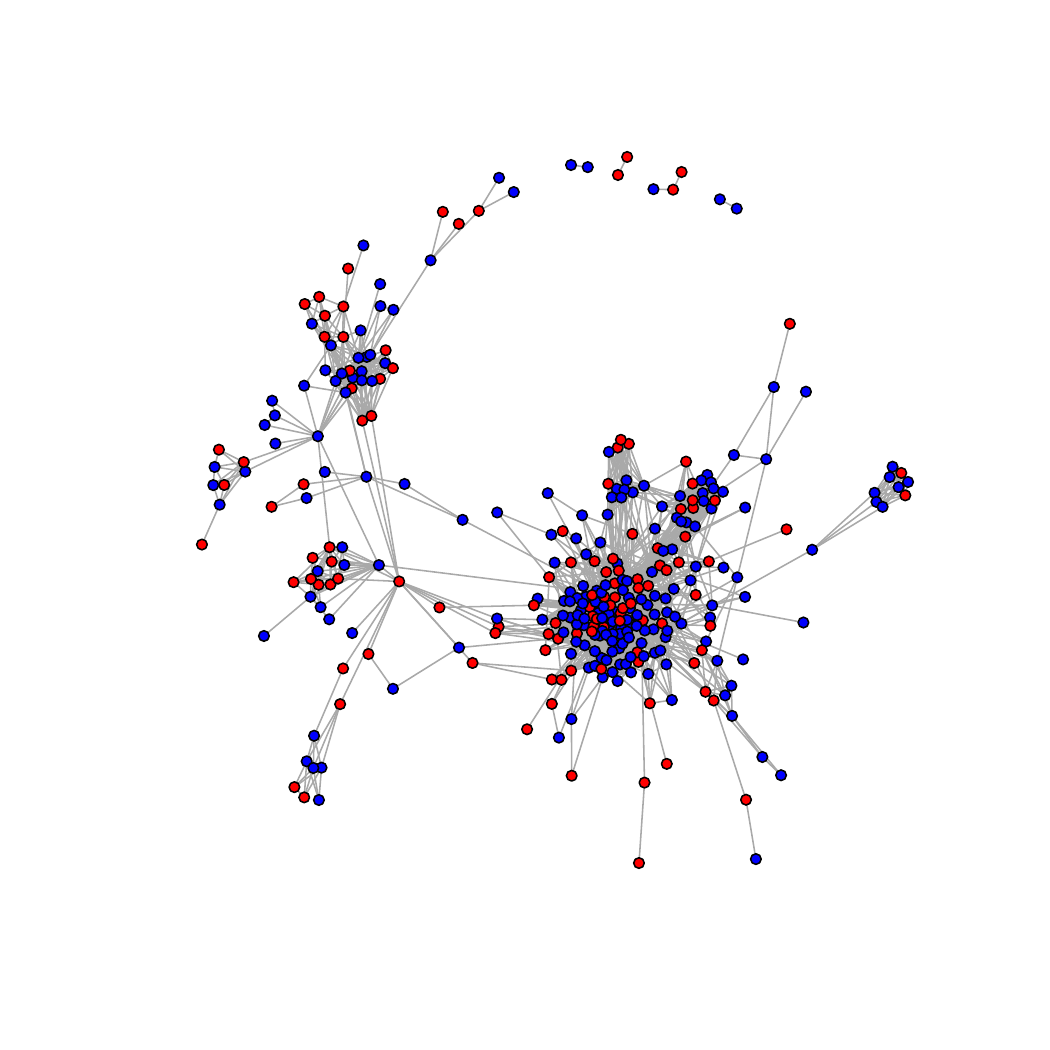} }}%
	\qquad
	\subfloat[\centering]{{\includegraphics[width=7.1cm,height=6.5cm]{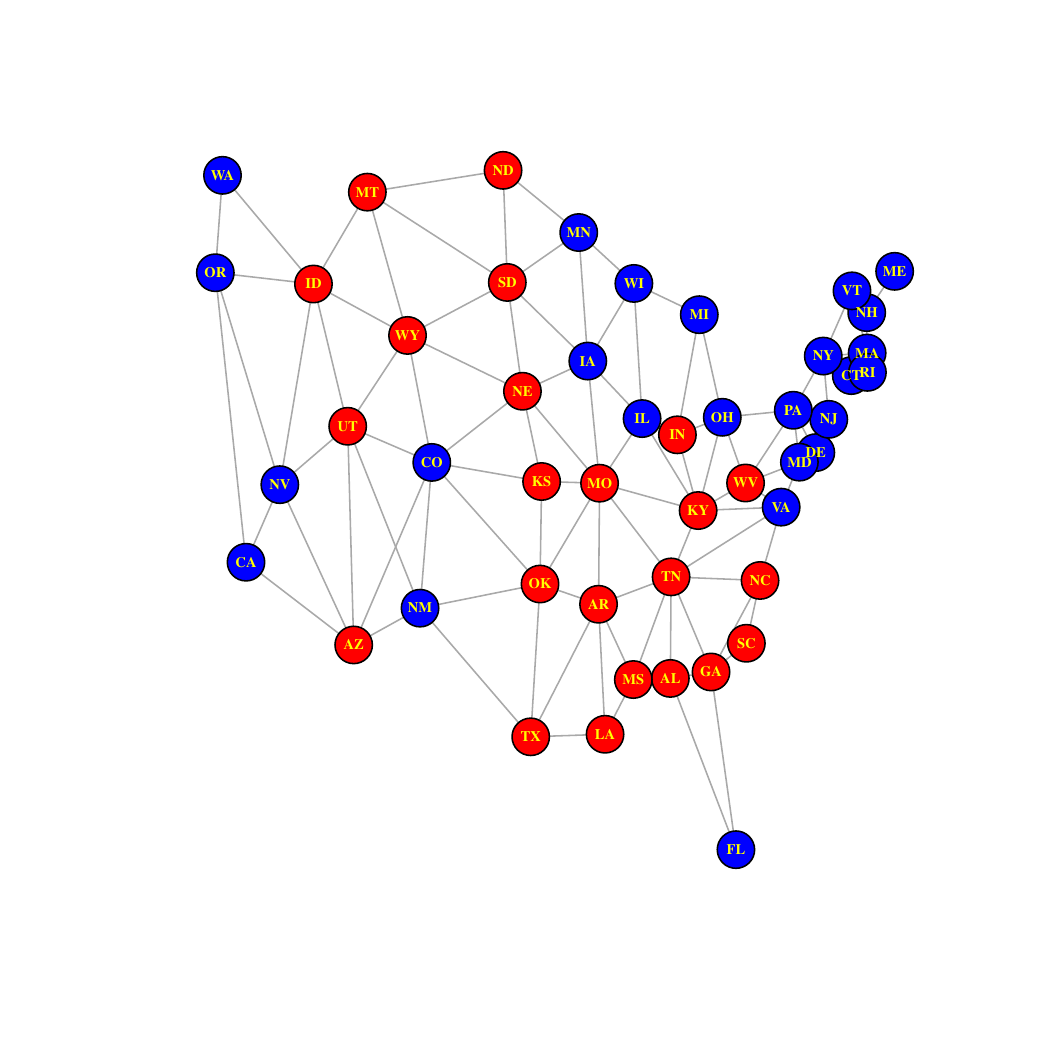} }}%
	\caption{Examples of spatially correlated network data: (a) Facebook data (red nodes are female, blue nodes are male) \cite{BM16}, (b) 2012 US election data (red states are republican, blue states are democratic)}%
	\label{fig:example}%
\end{figure}

%\import{Chapters/}{motiv_example}
The recent accumulation of dependent network data in modern statistics has made it increasingly important to develop realistic and mathematically tractable methods for modeling structure and dependence in high-dimensional distributions. Dependent data commonly arise in social and epidemic networks, spatial statistics, image databases, neural networks and computational biology. For example, in a social network like the facebook, the attributes of the users are dependent random variables, conditional on the underlying friendship network (Fig \ref{fig:example} (a)). This is because the probability that two people are friends, often depends on the similarity/dissimilarity between their attributes. The health status of individuals in an epidemic network is another example of highly correlated data. Another example of spatially correlated data is presidential election pattern across neighboring states. Figure \ref{fig:example} (b) shows the US neighborhood graph and outcome of the presidential election in the year 2012. Both the pictures in Figure \ref{fig:example} were kindly provided to me by Bhaswar Bhattacharya.

There is a massive amount of literature on how to analyze independent data in probability and statistics. However, there are many theoretical as well as practical scenarios, which demand analogous techniques for analyzing dependent data. For example, in an Erd\H{o}s-R\'enyi random graph, where each edge is present with some fixed probability, independent of the other edges, the indicators corresponding to the occurrence of triangles and higher order motif counts, are dependent random variables. Hence, in order to describe the asymptotic behavior of the total number of triangles in an Erd\H{o}s-R\'enyi random graph, one needs an asymptotic theory for the sum of dependent random variables. A more practical scenario arises in the estimation of the number of edges $|E(N)|$ in a large, inaccessible network $N$. A common strategy to do this is to sample some vertices at random from $N$ with probability $p$, and count the number of edges $T$ in the graph induced by $N$ on the sampled vertices. One can then show that the statistic $T/p^2$ is an unbiased estimator of $|E(N)|$. Now, one can write
$$T = \sum_{(i,j)\in E(N)} X_{ij}$$
where $X_{ij}$ denotes the indicator that both the nodes $i$ and $j$ in $N$ are sampled. Clearly, the collection $\{X_{ij}\}_{(i,j)\in E(N)}$ is not independent, since for any three distinct nodes $i,j,k$, $$\mathrm{Cov}(X_{ij}, X_{ik}) = p^3(1-p)\neq 0~.$$ So, once again, in order to derive theoretical properties of the estimator $T/p^2$, one needs an asymptotic theory for the sum of dependent random variables.

Another common example of spatially correlated data are the pixels of an image. In order to get a smooth image, one would require adjacent pixels to be strongly correlated. In other words, any reasonable probability model on the pixels should favor a configuration of pixels, where neighboring pixels have similar or identical states \cite{imageising}. An appropriate probability model for this setting can thus be described as
\begin{equation}\label{simple}
\mathbb{P}(\bm X) ~\propto~ \exp\left(-\beta \sum_{i\sim j} |X_i-X_j|\right)
\end{equation}
where $\bm X := (X_1,\ldots,X_N)$ is a configuration of pixel values, $i\sim j$ denotes that the pixels $i$ ands $j$ are neighbors, and the constant $\beta > 0$ is a parameter depending on the image.  

The model \eqref{simple} is a special case of the more general (second order) Markov random field, given by:

\begin{equation}\label{mrf}
\mathbb{P}(\bm X) ~\propto~ \exp\left(\sum_{i=1}^N h_{i} B_1(X_i) + \sum_{1\le i,j \le N} \theta_{ij} B_2(X_i,X_j) \right) 
\end{equation}
 
 where $\bm h := (h_1,\ldots,h_N)\in \mathbb{R}^N$, $\boldsymbol{\theta} = ((\theta_{ij}))_{1\le i,j\le N} \in \mathbb{R}^{N\times N}$, $\bm X \in \mathcal{X}^N$ for some finite set $\mathcal{X}$, $B_1:\mathcal{X}\mapsto \mathbb{R}$ is a nonzero function, and $B_2: \mathcal{X}^2 \mapsto \mathbb{R}$ is a nonzero symmetric function. For the model \eqref{simple}, $\bm h = \boldsymbol{0}$, $B_2(x,y) := |x-y|$, and $\boldsymbol{\theta} = -(\beta/2) \bm A$, where $A_{ij} =1$ if $i \sim j$ and $A_{ij}=0$ otherwise. 
 
\begin{figure}
 	\centering
 	\includegraphics[width=.4\linewidth]{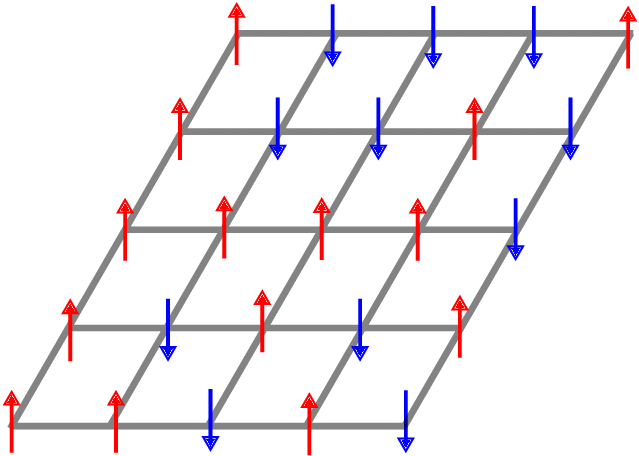}
 	\caption{Magnetic spins of particles sitting on a lattice; Picture Courtesy \cite{therm}.}
 	\label{magnetic}
 \end{figure}

 A special case of the Markov random field model \eqref{mrf} is the Ising model, which was initially developed in statistical physics to model ferromagnetism \cite{ising}. Although it was used initially as a framework for modeling interactions between particles sitting on the nodes of a network (Fig \ref{magnetic}), recently the Ising model has turned out to be particularly useful for modeling various statistical datasets with an underlying network structure (cf.~\cite{spatial,cd_trees,disease,neural,geman_graffinge,innovations} and the references therein). This is a discrete exponential family with binary outcomes, where the sufficient statistic involves a quadratic term designed to capture correlations arising from pairwise interactions and a linear term measuring the overall individual effect, and can be obtained from the model \eqref{mrf} by taking $\mathcal{X} = \{-1,1\}$, $B_1(x) = x$, $B_2(x,y) = xy$, $h_i \equiv h$ and $\boldsymbol{\theta} = \beta \bm J$ for some known symmetric matrix $\bm J$, known as the interaction matrix. $\beta \ge 0$ and $h \in \mathbb{R}$ are treated as parameters of the Ising model, with $\beta$ acting as a measure of correlation between the variables $X_1,\ldots, X_N$, and $h$ acting as an overall signal strength. In the language of statistical physics, $\beta$ and $h$ are called the inverse temperature and the external magnetic field, respectively. 
 
The $p$-tensor (spin) Ising model, a specific instance of the more general higher-order Markov random fields,  is a discrete exponential family where the sufficient statistic consists of a multilinear polynomial of degree $p\geq 2$ and a linear term, which provides an effective and mathematically  tractable way to simultaneously model both peer-group effects, between $p$-tuples of friends, and individual effects. More precisely, the $p$-tensor Ising model is given by:

\begin{equation}\label{pspinmd}
\mathbb{P}(\bm X)~\propto~ \exp\left(\beta\sum_{1\le i_1,\ldots,i_p \le N} J_{i_1\ldots i_p} X_{i_1}\ldots X_{i_p} + h \sum_{i=1}^N X_i \right)
\end{equation}

for all $\bm X \in \{-1,1\}^N$ and some known symmetric tensor $\bm J := ((J_{i_1\ldots i_p}))_{1\le i_1,\ldots,i_p \le N}$, known as the interaction tensor. For various examples and applications of this and related models in statistical physics, see  \cite{ab_ferromagnetic_pspin,pspinref1,ising_general,ferromagnetic_mean_field,turban,pspinref2} and the references therein.

%\begin{figure}%
%	\centering
%	\subfloat[\centering ]{{\includegraphics[width=7cm,height=6cm]{ad2} }}%
%	\qquad
%	\subfloat[\centering ]{{\includegraphics[width=7cm,height=6cm]{adatom} }}%
%	\caption{Multi-atom interactions on a crystal surface; Picture Courtesy: Fig \ref{fig:test111} (a): \cite{adat}, credit: Alexander Antropov, Vladimir Stegailov/Journal of Nuclear Materials}%
%	\label{fig:test111}%
%\end{figure}

 \begin{figure}
	\centering
	\includegraphics[width=7cm,height=6cm]{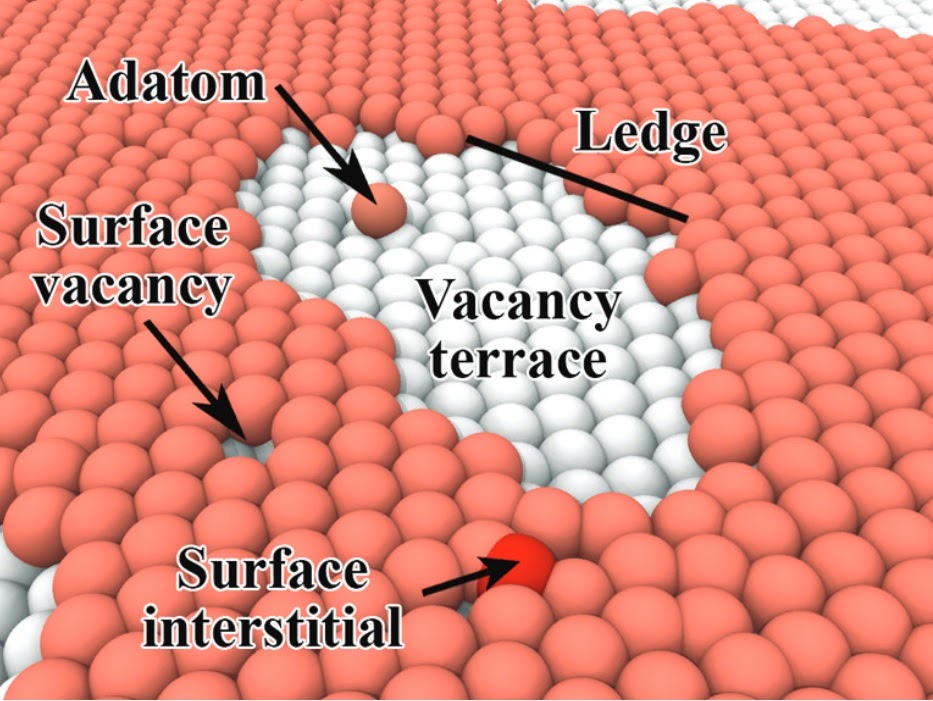}
	\caption{Atoms on a crystal surface; Picture Courtesy \cite{adat}, Credit: Alexander Antropov, Vladimir Stegailov/Journal of Nuclear Materials \cite{alexander}.}
	\label{fig:test111}
\end{figure}

The $p$-tensor Ising model acts as a useful framework in situations where the dependencies in a network arise not just from pairs, but from peer-group effects. For example, it is more likely for an individual to choose a binary attribute if many groups of friends have also chosen the same attribute. In fact, we demonstrate this phenomenon on a music recommender system data in Chapter \ref{ch:generalmple}, where we show that users' preference for a particular artist does not depend only on their pairwise interactions, but rather on higher order peer effects. A similar phenomenon also arises in various models of crystals, where the atoms on a crystal surface (see Fig \ref{fig:test111}) interact not just in pairs, but in triangles and higher order tuples.

\section{Outline of the Thesis}

In this thesis, we consider the problem of estimating the parameters $\beta$ and $h$ of the $p$-tensor Ising model \eqref{pspinmd} given a single sample from the model. This problem has been extensively studied for the $p=2$ (matrix) case, which includes, among others, the classical results on consistency and optimality of the maximum likelihood (ML) estimates for lattice models \cite{comets_exp,gidas,guyon,discrete_mrf_pickard}, and the seminal paper of Chatterjee \cite{chatterjee}, where general conditions for  $N^{\frac{1}{2}}$-consistency of the maximum pseudolikelihood estimate  (MPLE) \cite{besag_lattice,besag_nl} were derived. Various extensions and applications of the techniques in \cite{chatterjee}, in the contexts of estimation of parameters in matrix Ising models on general weighted graphs, logistic regression models with dependent observations, more general Ising models where the outcomes are influenced by various underlying networks, joint estimation of parameters, and related problems in hypothesis testing, can be found in \cite{BM16,cd_ising_I,cd_ising_II,cd_ising_estimation,pg_sm,bresler_II,rm_sm}.

However, none of these results  say anything about the limiting distribution of the estimates, and hence,  cannot be used for inferential tasks, such as constructing confidence intervals and hypothesis testing. In fact, proving general limit theorems in these models is often extremely difficult, if not impossible, because of the presence of an unknown normalizing constant (partition function) in the estimation objective function, which is both computationally and mathematically intractable for Ising models on general graphs. As a consequence, it is natural to assume certain special structures on the underlying network interactions if one desires to obtain precise results such as central limit theorems. A particularly useful structural condition which preserves several interesting properties of general systems, is to assume that all pairwise interactions between the nodes of the network are present. This is the well-known 2-tensor (matrix) Curie-Weiss model \cite{dm_gibbs_notes,ellis_book,ellis,glauber_dynamics}, which has been extensively studied in physics, probability, and statistics, and provides the foundations of our understanding of mean-field models with pairwise interactions. In particular, Comets and Gidas \cite{comets} provided a complete description of the limiting distribution of the ML estimates of the parameters in the matrix  Curie-Weiss model.

The matrix  Curie-Weiss model naturally extends to the $p$-tensor Curie-Weiss model, for any $p \geq 2$, in which the underlying tensor has all the possible $p$-tuples of interactions. The $p$-tensor Curie-Weiss model is a special case of \eqref{pspinmd}, obtained by taking all entries of the interaction tensor $\bm J$ to be $N^{1-p}$. In Chapter \ref{curiech}, we establish highly non-standard asymptotics of the ML estimates in the $p$-tensor Curie-Weiss model. In particular, we demonstrate the existence of a critical curve in the interior of the parameter space, on which the ML estimates have a limiting mixture distribution comprising of normals, half-normals and point masses. More surprisingly, we show that at the \textit{boundary points} of this curve, the ML estimates are superefficient, converging at a rate faster than the parametric rate $N^{-1/2}$, to limiting non-Gaussian distributions. The geometry of this curve also depends on the parity of the \textit{interaction factor} $p$.   

In more general Ising models, ML estimation is not possible due to the presence of an inexplicit and intractable normalizing constant in the expression of the likelihood function. To be precise, the normalizing constant for the measure \eqref{pspinmd} is given by:
$$Z_N(\beta,h) := \sum_{\bm x \in \{-1,1\}^N} \exp\left(\beta\sum_{1\le i_1,\ldots,i_p \le N} J_{i_1\ldots i_p} x_{i_1}\ldots x_{i_p} + h \sum_{i=1}^N x_i \right)~.$$ Since $Z_N(\beta,h)$ is a sum of $2^N$ many terms, for even moderately large values of $N$, it is incomputable in general (although in Chapter \ref{curiech} we will see that this can be computed in $O(n)$ time using the special structure of the Curie-Weiss network). This computational issue was circumvented by Chatterjee \cite{chatterjee}, who proposed using the maximum pseudolikelihood (MPL) estimator, which maximizes an approximation of the likelihood, obtained by taking product of the conditional distribution of each entry of $\bm X$ given the rest, over all the entries of $\bm X$. In Chapter \ref{ch:generalmple}, we use the maximum pseudo-likelihood (MPL) method to provide a  computationally efficient algorithm for parameter estimation that avoids computing the intractable partition function. We establish general conditions under which the  MPL estimate is $\sqrt N$-consistent, that is, it converges to the true parameter at rate $1/\sqrt N$. Our conditions are robust enough to handle a variety of commonly used tensor Ising models, including spin glass models with random interactions and models where the rate of estimation undergoes a   phase  transition. In particular, this includes results on $\sqrt N$-consistency of the MPL estimate in the well-known $p$-spin Sherrington-Kirkpatrick (SK) model, spin systems on general $p$-uniform hypergraphs, and Ising models on the hypergraph stochastic block model (HSBM). In fact, for the HSBM we pin down the exact location of the phase transition threshold, which is determined by the positivity of a certain mean-field variational problem, such that above this threshold the MPL estimate is $\sqrt N$-consistent, while below the threshold no estimator is consistent. Finally, we derive the precise limiting distribution of the MPL estimate in the special case of the Curie-Weiss model, which is the Ising model on the complete $p$-uniform hypergraph, at all points above its estimation threshold. Interestingly, in this case, the MPL estimate saturates the Cramer-Rao lower bound, showing that even though the MPL estimate is obtained by minimizing only an approximation of the true likelihood function for computational convenience, there is no loss in its asymptotic statistical efficiency.

In Chapter \ref{chap:covariateising}, we consider a more general model which incorporates the covariate information of the individual nodes. More precisely, it can be viewed as a generalization of the classical logistic regression model, with dependent observations. The model considered there, allows varying signal strength (external magnetic field) terms, where each of these signals is assumed to be the linear projection of some high-dimensional covariate along a fixed parameter vector. With only the signal terms present, we recover the classical high-dimensional logistic regression model, but the presence of a quadratic interaction term leads us to view this model both as an Ising model with covariates (a variant of the standard Ising model) and as a logistic regression with dependent observations (a variant of the vanilla logistic regression). This framework can be used to model the health status of individuals in an epidemic network, which are binary outcomes (healthy or ill) depending not only on the health status of neighboring individuals in that network, but also on personal health attributes like age, weight, diet, and immunity. In Chapter \ref{chap:covariateising}, we propose an $L^1$-penalized maximum pseudolikelihood approach to estimate the high-dimensional parameter in the Ising model with covariates, and show that under a sparsity assumption on the true parameter vector, our algorithm recovers it at rate $\sqrt{\log(d)/N}$, where $d$ is the dimension of the covariates and $N$ is the number of observations (size of the network).

\chapter{Maximum Likelihood Estimation in the Tensor Curie-Weiss Model}\footnotetext{This chapter is a joint work with Jaesung Son and Bhaswar B. Bhattacharya}\label{curiech}
%\import{Chapters/}{intro}

%\subsection{Organization}
In this chapter, we study the problem of parameter estimation in the $p$-tensor Curie-Weiss model, i.e. the model \eqref{pspinmd} where all $p$-tuples of interactions are present, and have equal strength. Given natural parameters $\beta \geq 0$ and $h \in \mathbb{R}$, the $p$-tensor Curie-Weiss model is a discrete exponential family on $\sa_N := \{-1,1\}^N$, defined as:
\begin{align}\label{cwwithmg}
\mathbb{P}_{\beta,h,p}(\bs)  = \frac{ \exp \left\{ \frac{\beta}{N^{p-1}} \sum_{1 \leq i_1, i_2, \ldots, i_p \leq N} X_{i_1} X_{i_2} \cdots X_{i_p} + h \sum_{i=1}^N X_i \right \} }{2^{N} Z_N(\beta,h,p)} ,
\end{align}
for $\bs :=(X_1,\ldots,X_N) \in \sa_N$.  The normalizing constant, also referred to as the partition function, $Z_N(\beta,h,p)$ is determined by the condition $\sum_{\bs \in \sa_N}\mathbb{P}_{\beta,h,p}(\bs)
\newblock=1$, that is, 
\begin{align}\label{ptnepn}
Z_N (\beta,h,p)  = \frac{1}{2^{N}} \sum_{\bs \in \sa_N}  \exp \left\{ \frac{\beta}{N^{p-1}} \sum_{1 \leq i_1, i_2, \ldots, i_p \leq N} X_{i_1} X_{i_2} \cdots X_{i_p} + h \sum_{i=1}^N X_i \right \} .
\end{align}
Denote by $F_N(\beta,h,p) := \log Z_N(\beta,h,p)$ the log-partition function of the model. Hereafter, we will often abbreviate $\mathbb{P}_{\beta,h,p}, Z_N(\beta,h,p)$, and $F_N(\beta,h,p)$, by $\mathbb{P}, Z_N$, and $F_N$, respectively, when there is no scope of confusion. For discussions on the various thermodynamic properties of this model, which in the statistical physics literature is more commonly known as the ferromagnetic $p$-spin model, refer to \cite{ab_ferromagnetic_pspin,pspinref1,ferromagnetic_mean_field,pspinref2}.

In this chapter, we consider the problem of estimating the natural parameters $\beta$ and $h$ given a single sample $\bm X \sim \mathbb{P}_{\beta, h, p}$ from the $p$-tensor Curie-Weiss model \eqref{cwwithmg}. One interesting feature of the Curie-Weiss model is that the partition function here can be computed in linear time, which is evident from the following alternative expression of $Z_N(\beta,h,p)$:

$$Z_N(\beta,h,p) = \frac{1}{2^N} \sum_{m \in \left\{-1,-1+\frac{2}{N},-1+\frac{4}{N},\ldots,1-\frac{2}{N},1\right\}} \binom{N}{\frac{N(1-m)}{2}} e^{\beta N m^p + hNm}~.$$
This makes the likelihood function easily computable, and hence, maximum likelihood (ML) estimation is possible. It is well-known, since the model \eqref{cwwithmg} has only one sufficient statistic (the sample mean $\os$), that joint estimation of the parameters $(\beta, h)$ in this model is, in general, impossible. This motivates the study of individual (marginal) estimation, that is, estimating $h$ when $\beta$ is assumed to be known and  estimating $\beta$ when $h$ is assumed to be known. As mentioned before, for the matrix $(p=2)$ Curie-Weiss model, this problem was studied in \cite{comets},  where the limiting properties of the individual ML estimates were derived. In this chapter, we consider the analogous problem for the $p$-tensor Curie-Weiss model, for $p \geq 3$. In particular, we derive precise limit theorems for the individual ML estimates of $\beta$ and $h$, hereafter, denoted by $\hat{\beta}_N$ and $\hat{h}_N$,  at all the parameter points. In addition to providing a complete description of the asymptotic properties of the ML estimates, our results unearth several remarkable new phenomenon, which we briefly summarize below. 

\begin{itemize}
	
	\item For `most' points in the parameter space, the ML estimates $\hat{\beta}_N$ and $\hat{h}_N$ are $N^{\frac{1}{2}}$-consistent and asymptotically normal (Theorem \ref{cltintr3_1} and Theorem \ref{thmmle1}). Here, the limiting variance equals the limiting inverse Fisher information, which implies that the ML estimates are, in fact, asymptotically efficient at these points (Remark \ref{remark1}). The variance of the limiting normal distribution can be easily estimated as well, hence, this result also provides a way to construct asymptotically valid confidence intervals for the model parameters (Section \ref{sec:applications}).

	\item More interestingly, there are certain `critical' points, which form a 1-dimensional curve in the parameter space, where $\hat{\beta}_N$ and $\hat{h}_N$ are still $N^{\frac{1}{2}}$-consistent, but the limiting distribution is a mixture with both continuous and discrete components. The number of mixture components is either two or three, depending on, among other things, the sign of one of the parameters and the parity of $p$. In particular, at the points where the critical curve intersects the region $h\ne 0$, the scaled ML estimates $N^{\frac{1}{2}}(\hat{\beta}_N - \beta )$ and $N^{\frac{1}{2}}(\hat{h}_N-h)$ have a surprising  three component mixture distribution, where two of the components are folded (half) normal distributions and the other is a point mass at zero (Theorem \ref{cltintr3_III} and Theorem \ref{thmmle_III}). This new phenomenon, which is absent in the matrix case, is an example of the many intricacies of the tensor model. 
	
	\item Finally, there are one or two `special' points in the parameter space, depending on whether $p \geq 3$ is odd or even, respectively, where both the individual ML estimates are superefficient, with fluctuations of order $N^{\frac{3}{4}}$ and non-Gaussian limiting distributions (Theorem \ref{cltintr3_2} and Theorem \ref{thmmle2}). 
\end{itemize} 
Our results also reveal various other interesting phenomena, such as, inconsistency of $\hat{\beta}_N$ in a region of the parameter space, and an additional (strongly) critical point, where $\hat{h}_N$ is $N^{\frac{1}{2}}$-consistent, but $\hat{\beta}_N$ is not.   These results, which are formally stated in Section \ref{sec:samplemean_mle}, together provide a complete characterization of the limiting properties of the ML estimates in the $p$-tensor Curie-Weiss model. 

An important byproduct of our analysis is a precise description of the asymptotic distribution of the sample mean (magnetization) $\os$ (Theorem \ref{cltintr3_1}), a problem which is of independent interest in statistical physics. While this has been extensively studied for the $p=2$ case, to the best of our knowledge this is the first such result for the higher order ($p \geq 3$) Curie-Weiss model. The proofs require very precise approximations of the partition function $Z_N$ and a careful understanding of the maximizers of a certain function at all points in the parameter space. One of the technical bottlenecks in dealing with tensor models is the absence of the `Gaussian transform', which allows one to relate  the partition function with certain Gaussian integrals, in models with quadratic sufficient statistics, as in the matrix  Curie-Weiss model. This method, unfortunately, does not apply when $p \geq 3$, hence, to estimate the partition function we use a more bare-hands Riemann-sum approximation (see Section \ref{approx} for details).

The rest of this chapter is organized as follows. We state our main results on the limiting distribution of the sample mean and the ML estimates in Section \ref{sec:samplemean_mle}. The proof of the limiting distribution of the sample mean is described in Section \ref{cltmagnetization}. A proof overview for the asymptotic distributions of the ML estimates is given in Section \ref{sec:pfsketch_mle}. In Section \ref{sec:applications} we describe how these limiting results can be used to construct confidence intervals for the model parameters. Various details of the proofs and other technical lemmas are given in the Appendix.

\section{Statements of the Main Results}
\label{sec:samplemean_mle}

In this section we state our main results on the limiting properties of the sample mean and the ML estimates in the $p$-tensor Curie-Weiss model. The asymptotics of the sample mean are described in Section \ref{sec:samplemean}. The limiting distributions of the ML estimates are presented in Section \ref{sec:mle_beta_h_I}. Finally, in Section \ref{sec:mle_beta_h_II} we summarize our results in a phase diagram. 

\subsection{Limiting Distribution of the Sample Mean}
\label{sec:samplemean}

The fundamental quantity of interest in understanding the asymptotic behavior of the $p$-tensor Curie-Weiss model is the sample mean $\os=\frac{1}{N} \sum_{i=1}^N X_i$.  As alluded to before, the limiting properties of $\os$ has been carefully studied for the case $p=2$ \cite{comets,ellis_book}. Here, we will consider the case $p \geq 3$, where, as discussed below, many surprises and interesting new phase transitions emerge. 

In order to state the results we need a few definitions: For $p \geq 2$ and $(\beta,h) \in \Theta:=[0, \infty) \times \R$, define the function $H = H_{\beta,h,p}:[-1,1]\rightarrow \mathbb{R}$ as
\begin{align}\label{eq:H}
H(x) := \beta x^p + hx - I(x),
\end{align} where $I(x) := \frac{1}{2}\left\{(1+x)\log(1+x) + (1-x)\log(1-x) \right\}$, for $x \in  [-1, 1]$, is the binary entropy function. The points of maxima of this function will determine the typical values of $\bar X_N$ and, hence, play a crucial role in our results. A careful analysis of the function $H$ (see Section \ref{sec:propH}) reveals that it can have one, two, or three global maximizers in the open interval $(-1, 1)$, which leads to the following definition:\footnote{For a smooth function $f: [-1, 1] \rightarrow \R$ and $x \in (-1, 1)$, the first and second derivatives of $f$ at the point $x$ will be denoted by $f'(x)$ and $f''(x)$, respectively. More generally, for $s \geq 3$, the $s$-th order derivative of $f$ at the point $x$ will be denoted by $f^{(s)}(x)$. } 

\begin{definition}\label{punique} Fix $p \geq 2$ and $(\beta,h) \in \Theta$, and let $H$ be as defined above in \eqref{eq:H}. %Then we have the following definitions: 
	\begin{enumerate} 
		
		\item The point $(\beta,h)$ is said to be $p$-{\it regular}, if the function $H_{\beta,h,p}$ has a unique global maximizer $m_* = m_*(\beta,h,p) \in (-1,1)$ and $H_{\beta,h,p}''(m_*) < 0$.\footnote{A point $m \in (-1, 1)$ is a global maximizer of $H$ if $H(m) > H(x)$, for all $x\in [-1,1]\setminus \{m\}$.} Denote the set of all $p$-regular points in $\Theta$ by $\cR_p$.
		
		\item The point $(\beta,h)$ is said to be $p$-{\it special}, if $H_{\beta,h,p}$ has a unique global maximizer $m_* = m_*(\beta,h,p) \in (-1,1)$ and $H_{\beta,h,p}''(m_*) = 0$.  
		
		\item The point $(\beta,h)$ is said to be $p$-{\it critical}, if $H_{\beta,h,p}$ has more than one global maximizer. 
		
	\end{enumerate}
\end{definition}

Note that the three cases above form a disjoint partition of the parameter space $\Theta$. Hereafter, we denote the set of $p$-critical points by $\cp$, and the set of points $(\beta, h)$ where $H_{\beta, h, p}$ has exactly two global maximizers by $\cp^+$.  We show in  Lemma \ref{derh33} that the set of points in $\cp$ form a continuous $1$-dimensional curve in the parameter space $\Theta$ (see also Figure \ref{figure:ordering1} and Figure \ref{figure:ordering2}). Next, we consider points with three global maximizers, that is $\cp \backslash \cp^+$. To this end, define 
\begin{align}\label{eq:betatilde}	
\tilde{\beta}_p := \sup\left\{\beta \geq 0: \sup_{x\in [-1,1]}H_{\beta,0,p}(x) = 0  \right\}.  
\end{align} 
Alternatively, Lemma \ref{derh33} shows that $\tilde{\beta}_p$ is the smallest $\beta \geq 0$ for which the point $(\beta, 0)$ is $p$-critical. Now, depending on whether $p$ is odd or even we have the following two cases:

%\iffalse
\begin{itemize}
	
	\item $p \geq 3$ odd: In this case Lemma \ref{derh11} shows that, for all points $(\beta, h) \in \cp$, the function $H_{\beta, h, p}$ has exactly two global  maximizers, that is, $\cp=\cp^+$. 
	
	\item $p \geq 4$ even:  Here, Lemma \ref{derh11} shows that there is a unique point 
	$\lambda_p:=(\tilde{\beta}_p,0) \in \cp$, with $\tilde \beta_p$ as defined in \eqref{eq:betatilde}, at which the function $H_{\tilde{\beta}_p,0,p}$ has exactly three global maximizers. For all other points in $(\beta, h)  \in \cp$, $H_{\beta, h, p}$ has exactly two global maximizers, that is,  $\cp=\cp^+ \cup \{\lambda_p \}$. In the case, $p\geq 4$ is even, we will refer to the point $\lambda_p$, or, equivalently, the point $\tilde{\beta}_p$, as the $p$-{\it strongly critical} point.\footnote{Note that the point $\tilde \beta_p$ is defined for all $p \geq 2$ (even or odd) as in \eqref{eq:betatilde}. However, for $p \geq 3$ odd, this point is  $p$-critical, but not $p$-strongly critical (that means it belongs to $\sC_p^+$). On the other hand, for $p=2$ this point is 2-special (see discussion in Remark \ref{remark:theta}).} Hereafter, when the need while arise to distinguish strongly critical points from other critical points, we will refer to a point which is $p$-critical but not $p$-strongly critical, as $p$-{\it weakly critical}. Note that the collection of all $p$-weakly critical points is precisely the set $\cp^+$.
\end{itemize} 
%%%%%%%%%%%%%%%%%%%%%%%%%%%%%%%%%%%%%%%%%%%%%%%%%%%%%%%%%%%%%%%%%%%%%%%%%%%%%%%%
\begin{figure}[h]
	\centering
	\begin{minipage}[c]{0.48\textwidth}
		\centering
		\includegraphics[width=3.25in]
		{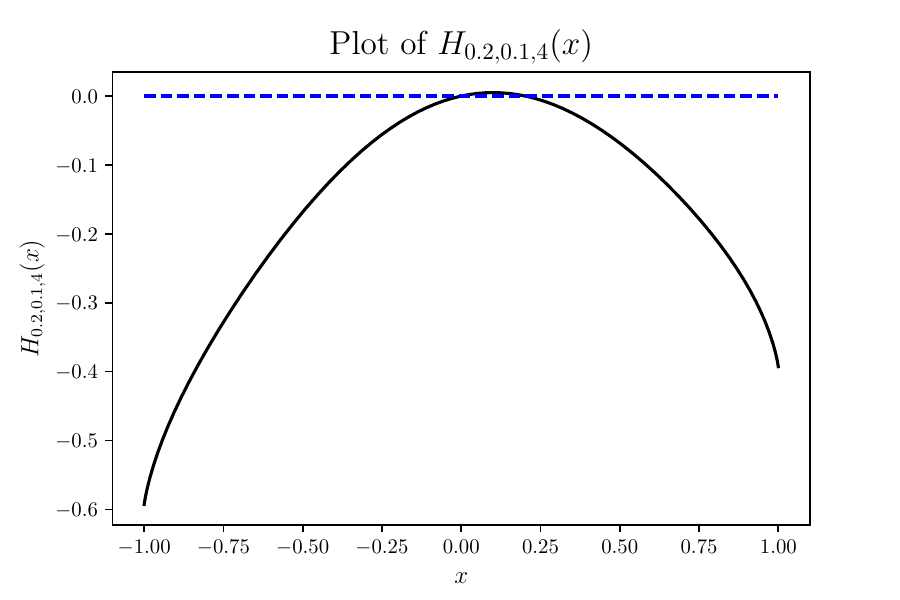}\\
		\small{(a)}
	\end{minipage}
	\begin{minipage}[c]{0.48\textwidth}
		\centering
		\includegraphics[width=3.25in]
		{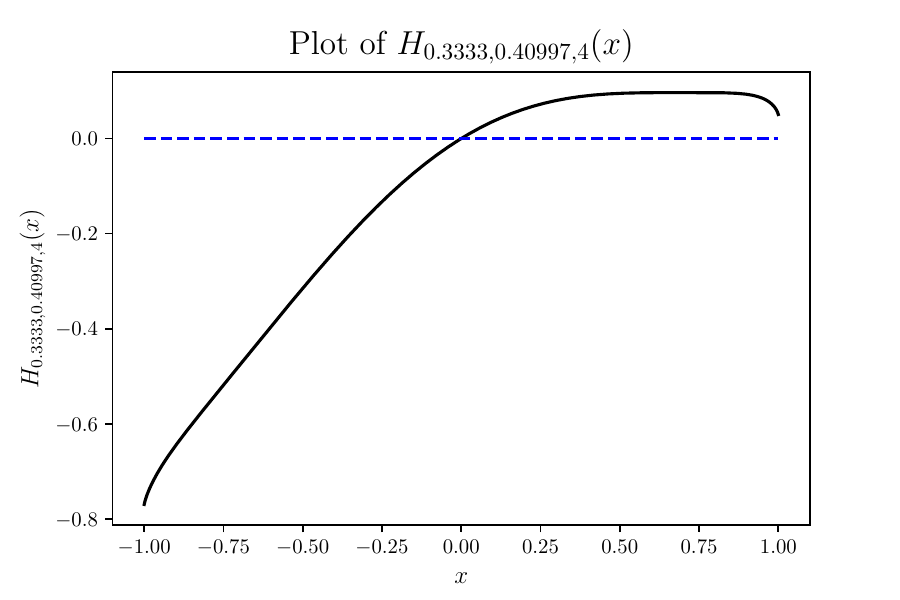}\\
		\small{(b)} 
	\end{minipage}
	\caption{\small{(a) Plot of the function $H_{\beta, h, p}$  at the $4$-regular point $(\beta, h)=(0.2, 0.1)$, where the function $H_{\beta, h, p}$ has a single global maximizer and the second derivative is negative at the maximizer; (b) plot of the function $H_{\beta, h, p}$  at the $4$-special point $(\beta, h)=(0.3333,0.40997)$, where the function $H_{\beta, h, p}$ has a single global maximizer, but the second derivative is zero at the maximizer.}}
	\label{fig:example_I}
\end{figure}
%%%%%%%%%%%%%%%%%%%%%%%%%%%%%%%%%%%%%%%%%%%%%%%%%%%%%%%%%%%%%%%%%%%%%%%%%%%%%%%%%%%%%%%%%%%%%%%%%%%%%%%%%%  
It remains to describe the structure of $p$-special points. To this end, fix $p \geq 3$ and define the following quantities: 
\begin{align}\label{eq:beta_h_special}
\check{\beta}_p := \frac{1}{2(p-1)} \left(\frac{p}{p-2}\right)^{\frac{p-2}{2}}\quad\textrm{and}\quad\check{h}_p := \tanh^{-1}\left(\sqrt{\frac{p-2}{p}}\right) - \check{\beta}_p p \left(\frac{p-2}{p}\right)^{\frac{p-1}{2}}. 
\end{align}
Again, depending on whether $p$ is even or odd there are two cases: 

\begin{itemize}
	
	\item $p \geq 3$ odd: In this case, Lemma \ref{derh22} shows that there is only one $p$-special point $\tau_p:=(\check{\beta}_p,\check{h}_p)$.
	
	\item $p \geq 4$ even: Here, again from Lemma \ref{derh22} and the symmetry of the model about $h=0$, there are two $p$-special points $\tau_p^+:=(\check{\beta}_p, \check{h}_p)$ and $\tau_p^-:=(\check{\beta}_p, -\check{h}_p)$. 
	
\end{itemize} 
These points are especially interesting, because, as we will see in a moment, here the sample mean has fluctuations of order $N^{\frac{1}{4}}$ and a non-Gaussian limiting distribution.

%%%%%%%%%%%%%%%%%%%%%%%%%%%%%%%%%%%%%%%%%%%%%%%%%%%%%%%%%%%%%%%%%%%%%%%%%%%%%%%%
\begin{figure}[h]
	\centering
	\begin{minipage}[c]{0.48\textwidth}
		\centering
		\includegraphics[width=3.25in]
		{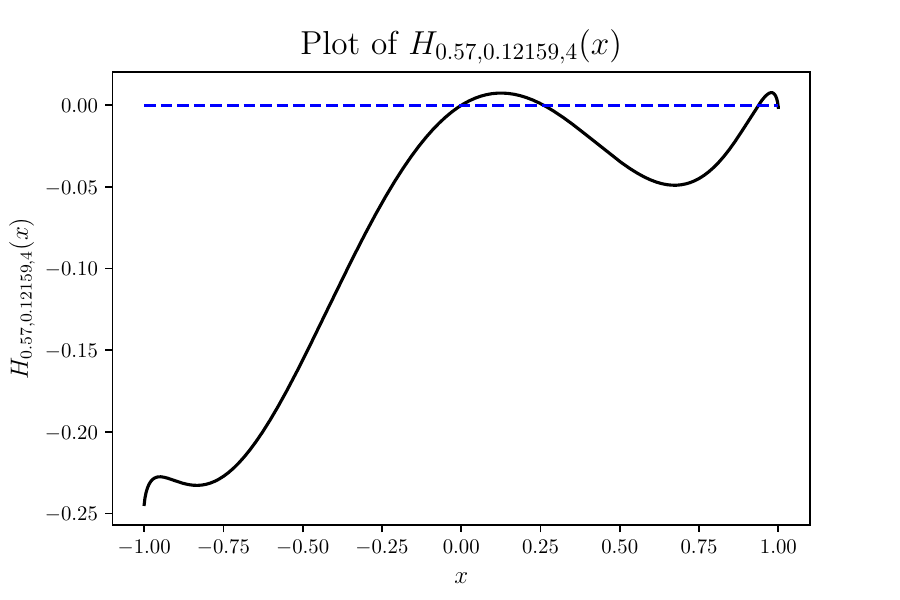}\\
		\small{(a)}
	\end{minipage}
	\begin{minipage}[c]{0.48\textwidth}
		\centering
		\includegraphics[width=3.25in]
		{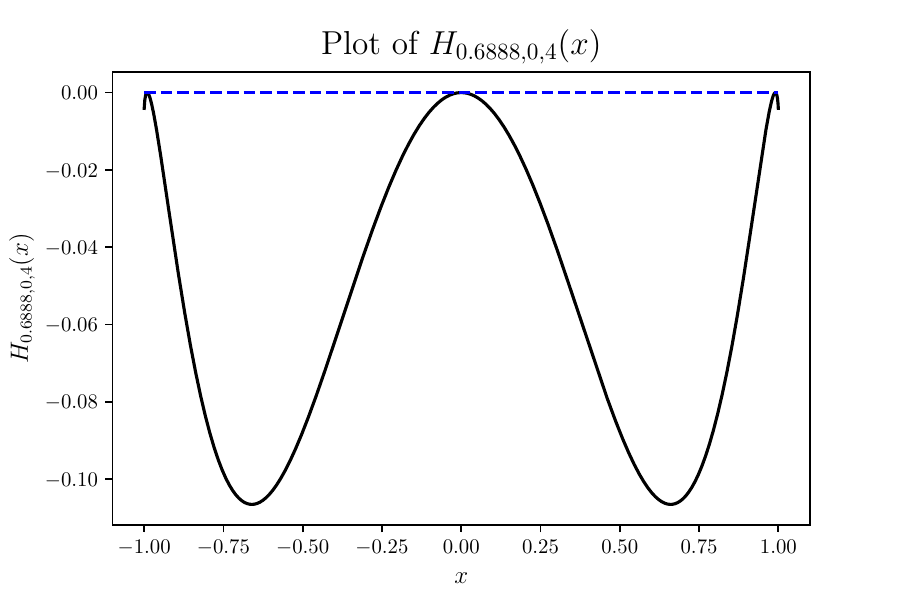}\\
		\small{(b)} 
	\end{minipage}
	\caption{\small{Plots of the function $H_{\beta, h, p}$  at $p$-critical points. For the plot in (a) $p=4$ and $(\beta, h)=(0.57,0.12159)$ and the function $H_{\beta, h, p}$ has two global  maximizers; and for (b) $p=4$ and $(\beta, h)=(0.688, 0)$ and the function $H_{\beta, h, p}$ has three global  maximizers, that is, the point $(0.688, 0)$ is $4$-strongly critical.  }}
	\label{fig:example_II}
\end{figure}
%%%%%%%%%%%%%%%%%%%%%%%%%%%%%%%%%%%%%%%%%%%%%%%%%%%%%%%%%%%%%%%%%%%%%%%%%%%%%%%%%%%%%%%%%%%%%%%%%%%%%%%%%%

The plots in Figure \ref{fig:example_I} and Figure \ref{fig:example_II} show instances of the different cases described above: Figure \ref{fig:example_I}(a) shows the plot of the function $H_{\beta, h, p}$ at the $4$-regular point $(\beta, h)=(0.2, 0.1)$, and Figure \ref{fig:example_I}(b) shows the plot of the function $H_{\beta, h, p}$ at the $4$-special point $(\beta, h)=(0.3333,0.40997)$.  On the other hand, Figure \ref{fig:example_II}(a) shows the plot of the function $H_{\beta, h, p}$ at the $3$-critical point  $(\beta, h)=(0.57, 0.12159)$,  which has two global maximizers, and  Figure \ref{fig:example_II}(b)  shows the plot of the function at the $4$-strongly critical point $(\beta, h)=(0.688, 0)$, where the function $H_{\beta, h, p}$ has three global  maximizers. In fact, recalling that $\cR_p$ denotes the set of all $p$-regular points and $\sC_p^+$ the set of points $(\beta, p)$ where $H_{\beta, h, p}$ has exactly two maximizers, the discussion above can be summarized as follows: 
\begin{align}\label{eq:parameter_space}
\Theta= 
\left\{
\begin{array}{cc}
\cR_p \bigcup \sC_p^+ \bigcup \{\tau_p\}  &  \text{ for } p \geq 3  \text{ odd},    \\
\cR_p \bigcup \sC_p^+ \bigcup  \{\lambda_p, \tau_p^+, \tau_p^-\}  &  \text{ for } p \geq 4 \text{ even}. 
\end{array}
\right. 
\end{align}
Figure \ref{figure:ordering1} and Figure \ref{figure:ordering2} illustrates this decomposition of the parameter space for $p=4$ and $p=5$, respectively.

%\iffalse

\begin{rem}\label{remark:theta} Note that \eqref{eq:parameter_space} provides a complete characterization of the parameter space for $p\geq 3$. As mentioned before, in the well-studied case of $p=2$, the situation is relatively simpler \cite{dm_gibbs_notes,ellis_book}. In this case, $H_{\beta, h, p}$ can have at most two global maximizers, that is, it has no strongly critical points, hence, $\cC_2=\cC_2^+$. In fact, it follows from \cite{ellis_book} that the set of points $(\beta, h)$ with exactly two global maximizers $\cC_2^+$ is the open half-line $(0.5, \infty)\times \{0\}$. Moreover, there is a single 2-special point $(0.5, 0)$ (where there the function $H$ has a unique maximum, but the double derivative is zero), and all the remaining points $\Theta \backslash [0.5, \infty)$ are 2-regular. This shows that for $p=2$ there is no point in $\Theta$ with $h \ne 0$ that is critical. In contrast, for $p \geq 3$ odd, the set of critical points is a continuous curve in $\Theta$ which intersects the line $h=0$ at a single point, and for $p \geq 4$ even, the set of critical points is a continuous curve in $\Theta$ which has two arms that intersect the line $h=0$ in the half-line $[\tilde \beta_p, \infty)$ (see Lemma \ref{derh33} for the precise statement and Figures \ref{figure:ordering1} and  \ref{figure:ordering2} for an illustration.) Moreover, this curve has exactly one limit point (if $p \geq 3$ is odd) and exactly two limit points (if $p \geq 4$ is even) outside it, which is (are) precisely the $p$-special point(s).  
\end{rem}

Having described the behavior of the function $H_{\beta, h, p}$, we can now state the limiting 
distribution of $\os$, which depends on whether the point $(\beta, h)$ is regular, critical, or special.

\begin{thm}[Asymptotic distribution of the sample mean]\label{cltintr1} Fix $p \geq 3$ and $(\beta, h) \in \Theta$, and suppose $\bs \sim \mathbb{P}_{\beta,h,p} $. Then with $H=H_{\beta,p,h}$ as defined in \eqref{eq:H}, the following hold: 
	
	\begin{itemize}
		
		\item[$(1)$] Suppose $(\beta, h)$ is $p$-regular and denote the unique maximizer of $H$ by $m_*= m_*(\beta,h,p)$. Then, as $N \rightarrow \infty$,
		\begin{align}\label{eq:meanclt_I}
		N^{\frac{1}{2}}\left(\os - m_*\right)\xrightarrow{D} N\left(0,-\frac{1}{H''(m_*)}\right).
		\end{align}

		\item[$(2)$] Suppose $(\beta, h)$ is $p$-critical and denote the $K \in \{2, 3\}$ maximizers  of $H$ by $m_1:=m_1(\beta,h,p)< \ldots < m_K:=m_K(\beta,h,p)$. Then, as $N \rightarrow \infty$, 
		\begin{align}\label{eq:meanclt_II_mixture}
		\os \xrightarrow{D} \sum_{k=1}^K p_k \delta_{m_k}, 
		\end{align}
		where for each $1\leq k\leq K$,\footnote{Note that all the global maximizers of the function $H$ belong to the open interval $(-1, 1)$, and if $(\beta, p)$ is $p$-critical and $m_1, \ldots, m_K$ are the global maximizers of $H$, for some $K \in \{2, 3\}$, then $H''_{\beta, h, p}(m_i) < 0$, for all $1 \leq i \leq K$. These statements are proved in Lemma \ref{derh11} and Lemma \ref{derh22}, respectively. This implies that the probabilities $p_1, \ldots, p_K$ in \eqref{eq:p1} are well-defined. Moreover, when $(\beta, h)$ is $p$-strongly critical, that is, $H_{\beta, h, p}$ has three global maximizers, the symmetry of the model about $h=0$ (recall that $p \geq 4$ is even and $h=0$ for a strongly critical point), implies that  the three maximizers are $m_1, 0, -m_1$, for some $m_1=m_1(\beta, h, p) < 0 $.}
		\begin{equation}\label{eq:p1}
		p_k := \frac{\left[(m_k^2-1)H''(m_k)\right]^{-1/2}}{\sum_{i=1}^K \left[(m_i^2-1)H''(m_i)\right]^{-1/2}}.
		\end{equation} 
		Moreover, if $A \subseteq [-1,1]$ is an interval containing $m_k$ in its interior for some $1\leq k \leq K$, such that $H(m_k) > H(x)$ for all $x\in A\setminus \{m_k\}$, then 
		\begin{align}\label{eq:meanclt_II}
		N^{\frac{1}{2}}\left(\os - m_k\right)\Big\vert \{\os \in A\}  \xrightarrow{D} N\left(0,-\frac{1}{H''(m_k)}\right).\end{align}

		\item[$(3)$] Suppose $(\beta,h)$ is $p$-special and denote the unique maximizer of $H$ by $m_*= m_*(\beta,h,p)$. Then, as $N \rightarrow \infty$, 
		$$N^\frac{1}{4}(\os - m_*) \xrightarrow{D} F,$$ 
		where the density of $F$ with respect to the Lebesgue measure is given by 
		\begin{align}\label{eq:meanclt_III}
		\mathrm dF(x) = \frac{2}{\Gamma(\tfrac{1}{4})}\left(-\frac{H^{(4)}(m_*)}{24}\right)^{\frac{1}{4}}\exp\left(\frac{H^{(4)}(m_*)}{24} x^4\right)\mathrm d x,
		\end{align}
		with $H^{(4)}$ denoting the fourth derivative of the function $H$. 
	\end{itemize}
\end{thm}

The result in Theorem \ref{cltintr1} follows from a slightly more general version (see Theorem \ref{cltun} in Section \ref{cltmagnetization}), where, instead of deriving the limiting distribution of $\os$ at a fixed point $(\beta, h)$, we compute the limits at appropriately perturbed parameter values $(\beta_N, h_N)$, with $\beta_N \rightarrow \beta$ and $h_N \rightarrow h$. This generalization will be required for deriving the asymptotic distribution of the ML estimates of $\beta$ and $h$, described in the following section. Deferring the technical details for later, we describe below the key ideas involved in the proof of Theorem \ref{cltintr1}: 

\begin{itemize}
	
	\item In the $p$-regular case, the proof has three main steps: The first step is to prove a concentration inequality of $\os$ in an asymptotically vanishing neighborhood $m_*$ (Lemma \ref{conc}). This not only shows that $m_*$ is the typical value of $\os$, but also implies that  the partition function $Z_N$ (which is the sum over all $\bm X \in \cC_N$ as in \eqref{ptnepn}), can be restricted over those $\bs$ for which $\os$ lies within this concentration interval around $m_*$. The second step is to find an accurate asymptotic expansion of $Z_N$ by first approximating this restricted sum by an integral over the concentration interval, and then applying saddle point techniques to get a further approximation to this integral (Lemma \ref{ex}). The third and final step is to use this approximation of $Z_N$ to compute the limit of the moment generating function of $N^{\frac{1}{2}}(\os- m_*)$, and show that the limit converges to that of the Gaussian distribution appearing in \eqref{eq:meanclt_I}. Details are given in Section \ref{sec:pfcltun_I}.

	\item The proof in the $p$-special case follows the same strategy as the $p$-regular case, with appropriate modifications to deal with the vanishing second derivative at the maximizer. As before, the first step is to prove the concentration of $\os$ within a vanishing neighborhood of $m_*$ which, in this case, requires a higher-order Taylor expansion, since $H_{\beta,h,p}''(m_*) = 0$ (Lemma \ref{irrconch}). The second step, as before, is the approximation of the partition function (Lemma \ref{ex2}).  The proof is completed by calculating the limit of the moment generating function of $N^{\frac{1}{4}}(\os- m_*)$ using this approximation to the partition function. Details are given in Section \ref{sec:pfcltun_II} and Section \ref{irregproof}.

	\item For the $p$-critical case, the basic proof strategy remains the same as above. However, to deal with the presence of multiple maximizers, we need to prove a conditional concentration result for the sample mean,  that is, $\os$ concentrates at one of the maximizers, given that $\os$ lies in a small  neighborhood of that maximizer (Lemma \ref{lem:multiple_max}). Similarly, for the second step, we need to approximate a restricted partition function, where instead of taking a sum over all configurations $\bm X \in \cC_N$ as in \eqref{cwwithmg}, we sum over configurations $\bm X \in \cC_N$ such that $\os$  lies in the neighborhood of one of the maximizers (Lemma \ref{lm:condpart}). Details are given in Section \ref{nonpuniqueclt} and Section \ref{sec:proofofnonun}.

\end{itemize}
%\iffalse

\begin{figure}[h]
	\begin{minipage}[c]{0.48\textwidth}
		\centering
		\includegraphics[width=3.2in,height=2.2in]
		{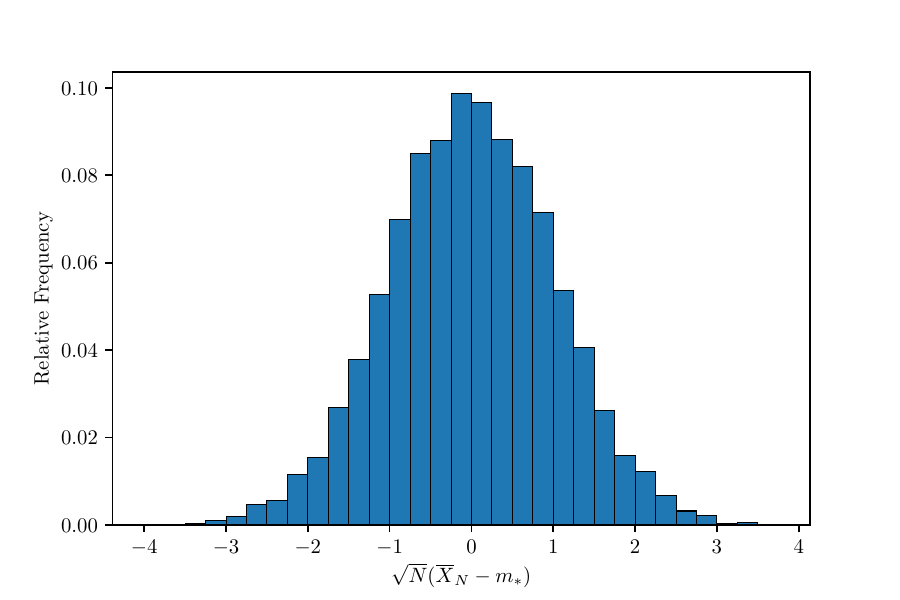}\\
		\small{(a)}
	\end{minipage}
	\begin{minipage}[l]{0.48\textwidth} 
		\centering
		\includegraphics[width=4.15in]
		{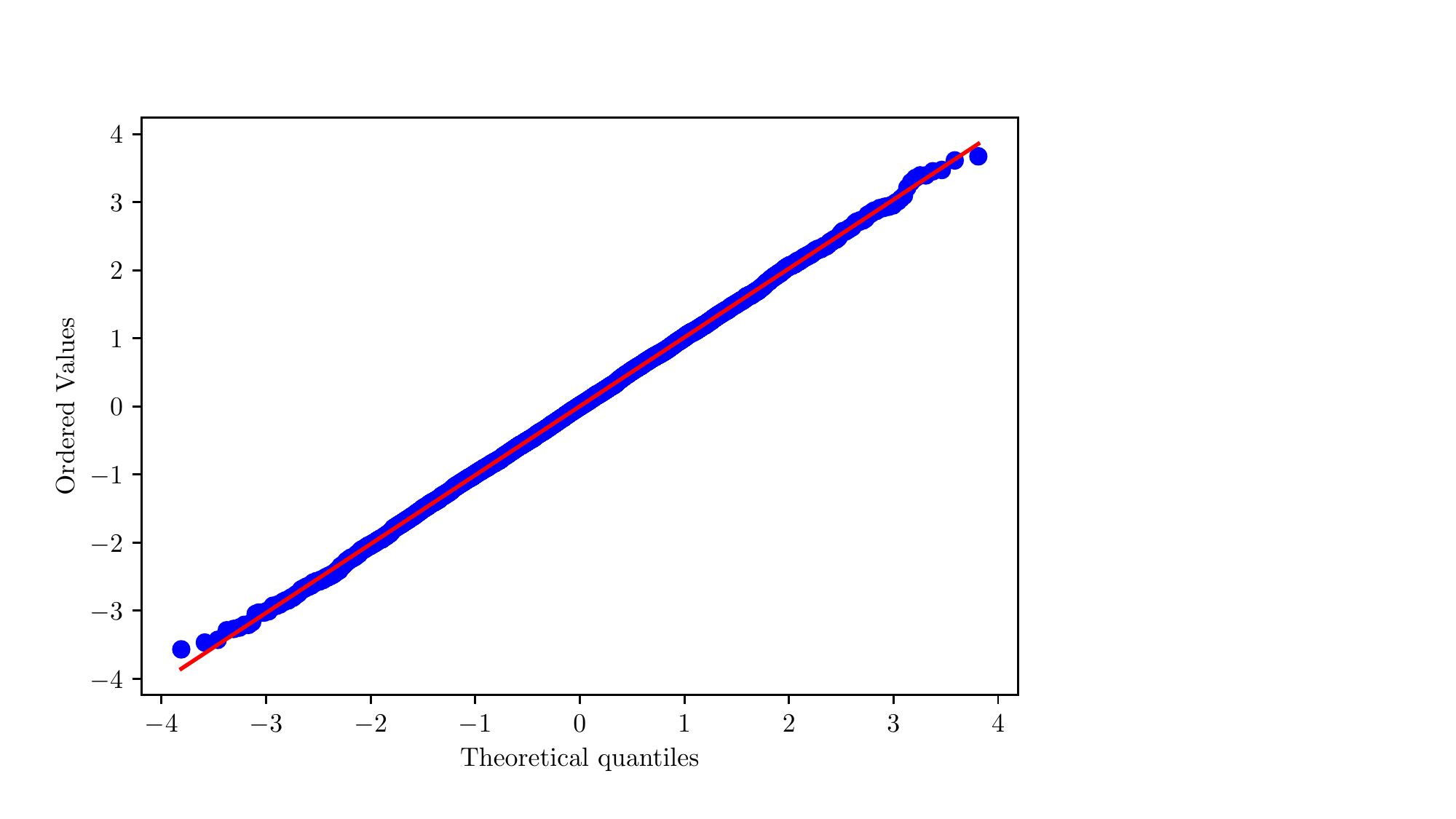}\\
		\small{(b)} 
	\end{minipage}
	\caption{\small{(a) The histogram of $N^{\frac{1}{2}} (\overline{X}_N - m_*)$ at the $4$-regular point $(\beta, h)=(0.2, 0.1)$ and (b) the corresponding quantile-quantile (QQ) plot confirming the asymptotic normality.}}
	\label{fig:sigma_histogram_I_regular}
\end{figure}

\begin{figure}[h]
	\centering
	\begin{minipage}[c]{0.48\textwidth}
		\centering 
		%\vspace{0.05in} 
		\includegraphics[width=3.2in,height=2.2in]
		{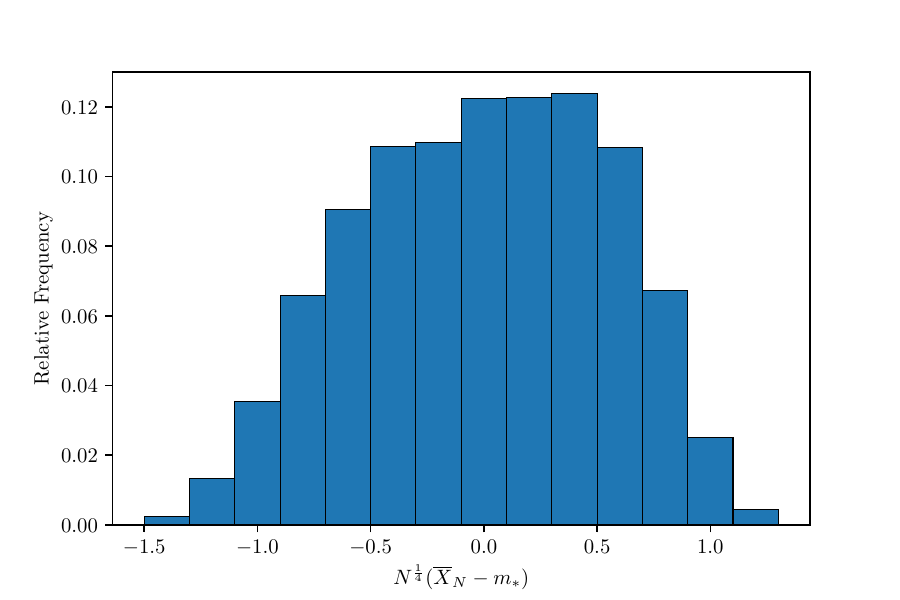}\\
		\small{(a)}
	\end{minipage}
	\begin{minipage}[c]{0.48\textwidth}
		\centering
		%\vspace{-0.1in}
		\includegraphics[width=4.35in,height=2.7in]
		{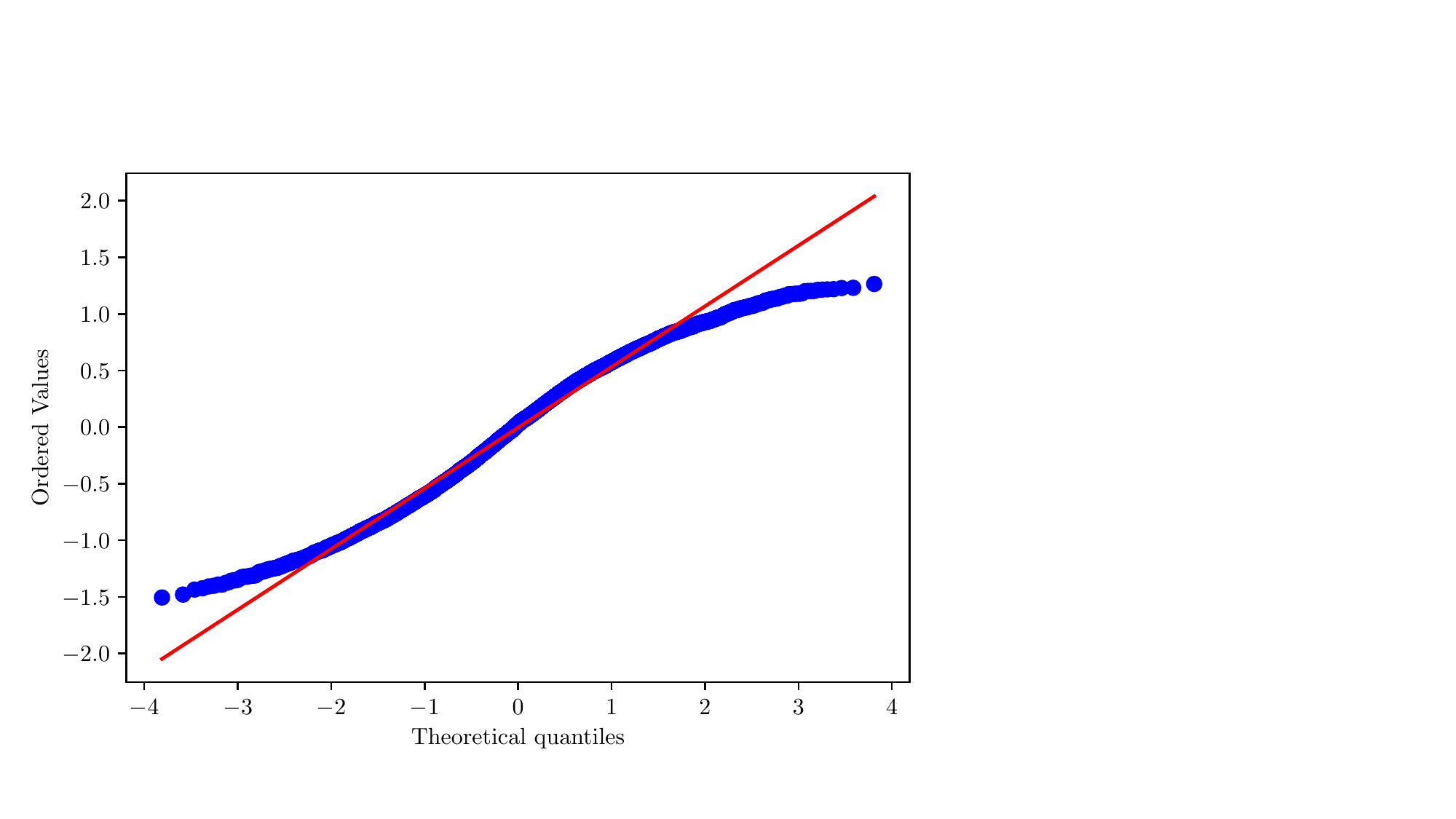}\\
		%\vspace{-0.1in}
		\small{(b)} 
	\end{minipage}
	\caption{\small{(a) The histogram of $N^{\frac{1}{4}}(\overline{X}_N - m_*)$ at the $4$-special point $(\beta, h)=(0.3333,0.40997)$ and (b) the corresponding QQ plot indicating a non-normal distribution.}}
	\label{fig:sigma_histogram_I_special}
\end{figure}

\begin{figure}[h]
	\centering
	\begin{minipage}[c]{1.0\textwidth}
		\centering
		\includegraphics[width=5.65in]
		{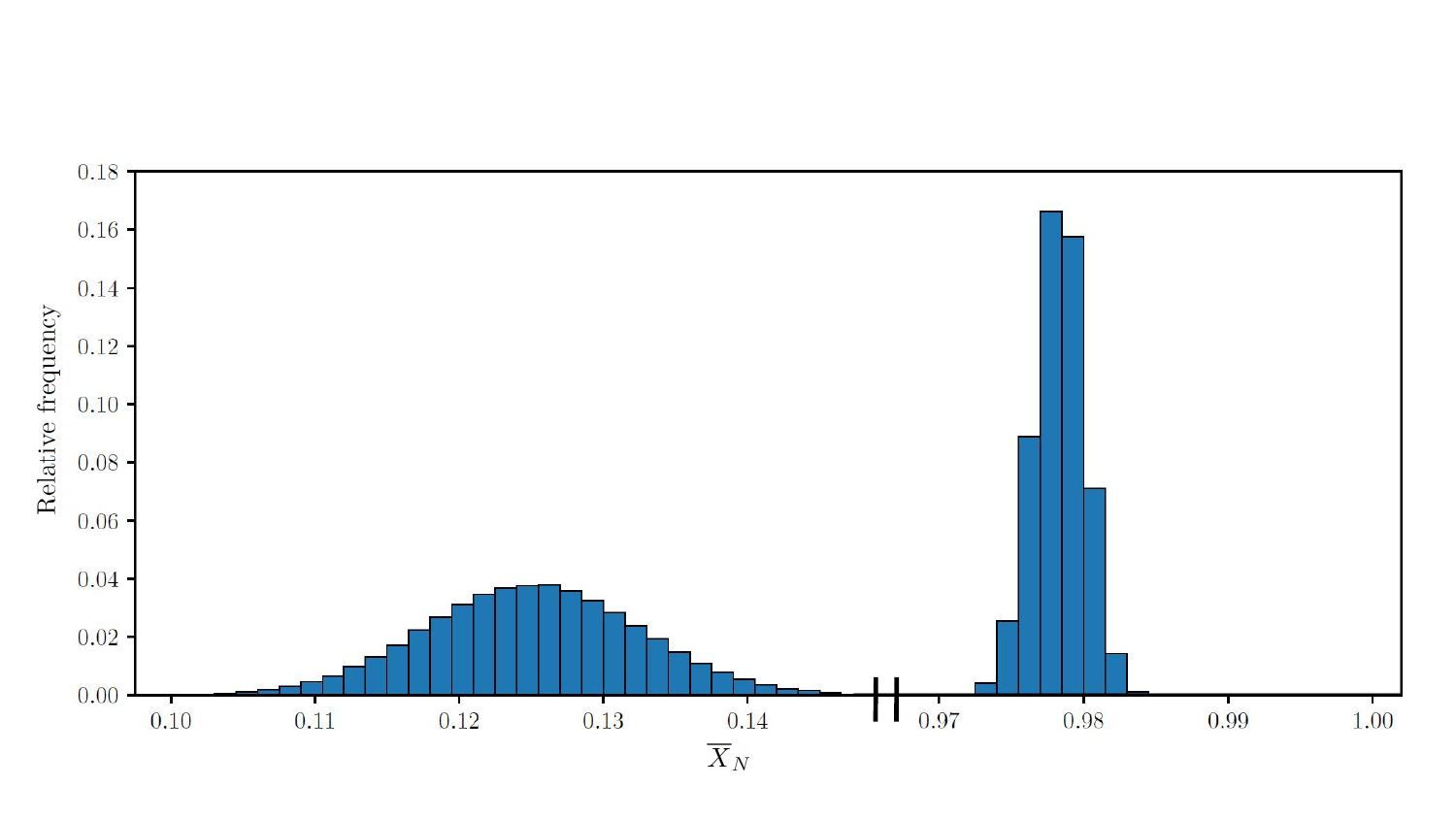}\\
	\end{minipage}
	\caption{\small{Histogram of $\os$ at the $4$-critical point $(0.57,0.12159)$, where the function  $H_{0.57,0.12159, 4}$ has two global maximizers, around which $\os$ concentrates. }}
	\label{fig:sigma_histogram_II}
\end{figure}

\begin{figure}[h]
	\centering
	\begin{minipage}[c]{1.0\textwidth}
		\centering
		\includegraphics[width=5.65in]
		{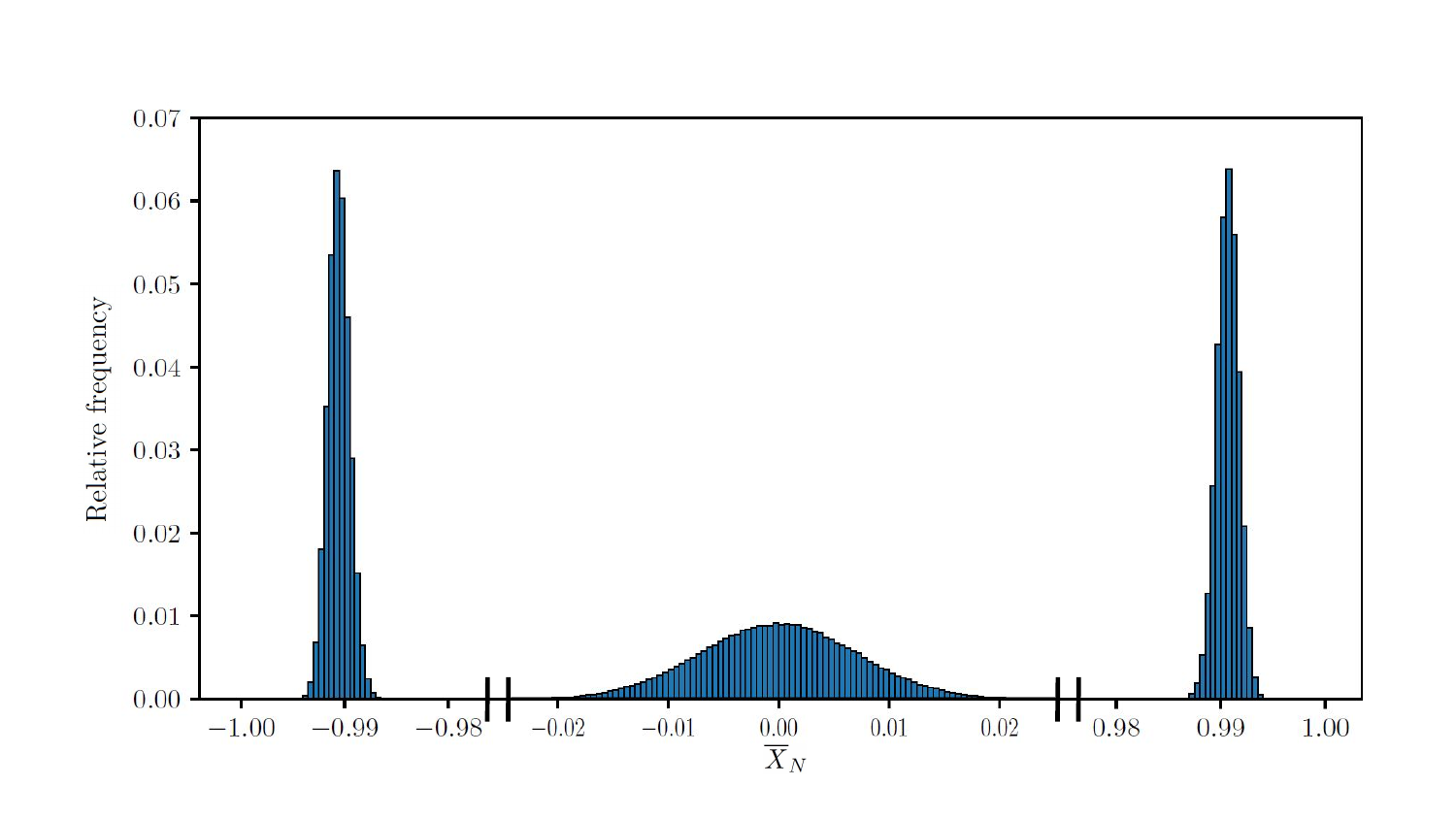}\\
	\end{minipage}
	\caption{\small{Histogram of $\os$ at the non $4$-strongly critical point $(0.6888,0)$, where the function  $H_{0.6888,0, 4}$ has three global maximizers, around which $\os$ concentrates. }}
	\label{fig:sigma_histogram_III}
\end{figure}

To empirically validate the different results in Theorem \ref{cltintr1}, we fix  $p \geq 3$,  some $(\beta, h) \in \Theta$, and $N = 20,000$. Then we  generate $10^6$ replications from $\p_{\beta, h, p}$ and plot the histograms of the sample means.  Figure \ref{fig:sigma_histogram_I_regular}(a) shows the histogram of $N^{\frac{1}{2}} (\overline{X}_N - m_*)$ at the $4$-regular point $(\beta, h)=(0.2, 0.1)$ where, as expected from (\ref{eq:meanclt_I}), we see a limiting normal distribution. This is also confirmed from the corresponding quantile-quantile (QQ) plot in Figure \ref{fig:sigma_histogram_I_regular}(b). Next, Figure \ref{fig:sigma_histogram_I_special}(a) shows the histogram of $N^{\frac{1}{4}} (\overline{X}_N - m_*)$ at the $4$-special point $(\beta, h)=(0.3333,0.40997)$, where a non-normal shape emerges, as predicted by \eqref{eq:meanclt_III}. The non-normality is also confirmed from the QQ plot in Figure \ref{fig:sigma_histogram_I_special}(b).  Figure \ref{fig:sigma_histogram_II} shows the histogram of $\os$ at the $4$-critical point  $(\beta, h)=(0.57,0.12159)$, where the function $H_{0.57,0.12159, 4}$ has two global maximizers (see plot in Figure \ref{fig:example_II}(a)). Hence, the histogram of $\os$ has two peaks located at two maximizers (as shown in \eqref{eq:meanclt_II_mixture}). Finally, in Figure \ref{fig:sigma_histogram_III} we show the histogram of $\os$ at a $4$-strongly critical point $(\beta, h)=(0.688, 0)$. Here, the histogram has three peaks, since the function $H_{\beta, h, p}$ has three global  maximizers (see plot in Figure \ref{fig:example_II}(b)). Note that the histograms of $\os$ both in Figure \ref{fig:sigma_histogram_II} and \ref{fig:sigma_histogram_III}, look like a Gaussian distribution in a neighborhood of each of the maximizers, as predicted by \eqref{eq:meanclt_II} in the theorem above. 

%%%%%%%%%%%%%%%%%%%%%%%%%%%%%%%%%%%%%%%%%%%%%%%%%%%%%%%%%%%%%%%%%%%%%%%%%%%%%%%%
%%%%%%%%%%%%%%%%%%%%%%%%%%%%%%%%%%%%%%%%%%%%%%%%%%%%%%%%%%%%%%%%%%%%%%%%%%%%%%%%%%%%%%%%%%%%%%%%%%%%%%%%%%

%%%%%%%%%%%%%%%%%%%%%%%%%%%%%%%%%%%%%%%%%%%%%%%%%%%%%%%%%%%%%%%%%%%%%%%%%%%%%%%%
%%%%%%%%%%%%%%%%%%%%%%%%%%%%%%%%%%%%%%%%%%%%%%%%%%%%%%%%%%%%%%%%%%%%%%%%%%%%%%%%%%%%%%%%%%%%%%%%%%%%%%%%%%

%%%%%%%%%%%%%%%%%%%%%%%%%%%%%%%%%%%%%%%%%%%%%%%%%%%%%%%%%%%%%%%%%%%%%%%%%%%%%%%%
%%%%%%%%%%%%%%%%%%%%%%%%%%%%%%%%%%%%%%%%%%%%%%%%%%%%%%%%%%%%%%%%%%%%%%%%%%%%%%%%%%%%%%%%%%%%%%%%%%%%%%%%%%

%%%%%%%%%%%%%%%%%%%%%%%%%%%%%%%%%%%%%%%%%%%%%%%%%%%%%%%%%%%%%%%%%%%%%%%%%%%%%%%%
%%%%%%%%%%%%%%%%%%%%%%%%%%%%%%%%%%%%%%%%%%%%%%%%%%%%%%%%%%%%%%%%%%%%%%%%%%%%%%%%%%%%%%%%%%%%%%%%%%%%%%%%%%

%\iffalse
\subsection{Asymptotics of the ML Estimates}
\label{sec:mle_beta_h_I}

In this section we consider the problem of estimating the parameters $\beta$ and $h$ given a single sample $\bs \sim \p_{\beta,h,p} $ using the method of maximum likelihood. Note that the distribution of  the $p$-tensor Curie-model \eqref{cwwithmg} has a single sufficient statistic $\os$. This suggests, as mentioned before, that the parameters $(\beta, h)$ cannot be estimated simultaneously. In fact, one can show that the joint ML estimates for $(\beta, h)$ might not exist with probability 1 (see Lemma \ref{mleexist} for details). However, it is possible to marginally estimate one of the parameters assuming that the other is known. Hereafter, given $\bs \sim \p_{\beta,h,p} $, we denote by $\hat{\beta}_N$ and $\hat{h}_N$ the maximum likelihood estimators of $\beta$ and $h$, respectively. Note that, for fixed $h \in \R$, $\hat{\beta}_N$ is a solution of the equation (in $\beta$)
\begin{equation}\label{eqmle}
\e_{\beta,h,p} \left(\overline{X}^p_N \right) = \overline{X}^p_N. 
\end{equation}
Similarly, for fixed $\beta$, $\hat{h}_N $ is a solution of the equation (in $h$)
\begin{equation}\label{eqmleh}
\e_{\beta,h,p} \left(\os \right) = \os,
\end{equation}
The limiting properties of the ML estimates of $h$ and $\beta$ are presented in Section \ref{sec:mle_h} and Section \ref{sec:mle_beta}, respectively. The full phase diagrams summarizing the results are given in Section \ref{sec:mle_beta_h_II}.

\subsubsection{ML Estimate of $h$} 
\label{sec:mle_h}

In order to describe the asymptotic distribution of the ML estimate of $h$, we need the following definition: 

\begin{definition}\label{defn:halfnormal}
	For $\sigma > 0$, the {\it positive half-normal distribution} $N^+(0,\sigma^2)$ is defined as the distribution of $|Z|$, where $Z \sim N(0,\sigma^2)$. The {\it negative half-normal distribution} $N^-(0,\sigma^2)$ is defined as the distribution of $-|Z|$, where $Z \sim N(0,\sigma^2)$.
\end{definition}

The asymptotic distribution of the ML estimate of $h$ is summarized in the theorem below. As expected, the results depend on whether $(\beta, h)$ is regular, critical, or special, which we state separately in the theorems below. In this regard, denote by $\delta_x$ the point mass at $x$.  We begin with the case when $(\beta, h)$ is regular. Throughout, $H=H_{\beta,p,h}$ will be as defined in \eqref{eq:H}.

\begin{thm}[Asymptotic distribution of $\hat{h}_N$ at $p$-regular points]\label{cltintr3_1}  Fix $p \geq 3$ and suppose $(\beta, h) \in \Theta$ is $p$-regular. Assume $\beta$ is known and $\bs \sim \mathbb{P}_{\beta,h,p} $. Then denoting the unique maximizer of $H$ by $m_*= m_*(\beta,h,p)$, as $N \rightarrow \infty$,
	\begin{align}\label{eq:h_1}
	N^{\frac{1}{2}}(\hat{h}_N  - h) \xrightarrow{D} N\left(0, -H''(m_*)\right).
	\end{align} 
\end{thm}

%The proof the above theorem is given in \ref{}. 
This result shows that $\hat{h}_N $ is $N^{\frac{1}{2}}$-consistent and asymptotically normal at the regular points. Before discussing more about the implications of this theorem, we state the result for the asymptotic distribution of $\hat{h}_N$ when $(\beta, h)$ is $p$-special.

\begin{thm}[Asymptotic distributions of $\hat{h}_N$ at $p$-special points]\label{cltintr3_2}	
	Fix $p \geq 3$ and suppose $(\beta, h) \in \Theta$ is $p$-special. Assume $\beta$ is known and $\bs \sim \mathbb{P}_{\beta,h,p} $. Then denoting the unique maximizer of $H$ by $m_*= m_*(\beta,h,p)$, as $N \rightarrow \infty$, 
	\begin{align}\label{eq:h_5} 
	N^\frac{3}{4}(\hat{h}_N - h) \xrightarrow{D} G_1, 
	\end{align} 
	where the distribution function of $G_1$ is given by  
	$$G_1(t) = F_{0,0}\left(\int_{-\infty}^{\infty} u ~\mathrm d F_{0,t}(u)\right),$$ with $F_{0, t}$ as defined in \eqref{eq:beta_h_distribution} below. 
\end{thm}

Finally, we consider the case $(\beta, h)$ is $p$-critical. Here, it is convenient to consider the cases $p$ is odd or even separately. 

\begin{thm}[Asymptotic distribution of $\hat{h}_N$ at $p$-critical points]\label{cltintr3_III}
	Fix $p \geq 3$ and suppose $(\beta, h) \in \Theta$ is $p$-critical. Assume $\beta$ is known and $\bs \sim \mathbb{P}_{\beta,h,p} $. Denote the $K \in \{2,3\}$ maximizers of $H$ by $m_1 := m_1(\beta, h, p) < \ldots < m_K := m_K(\beta, h, p)$, and let $p_1, \ldots, p_K$ be as in \eqref{eq:p1}. 
	
	\begin{itemize}

		\item[$(1)$] Suppose $p \geq 3$ is odd. In this case, the function $H$ has exactly two (asymmetric) maximizers $m_1 < m_2$ and, as $N \rightarrow \infty$,
		\begin{align}\label{eq:h_31}
		N^{\frac{1}{2}}(\hat{h}_N - h) \xrightarrow{D} \tfrac{p_1}{2} N^{-}\left(0, -H''(m_1)\right) + \tfrac{1-p_1}{2} N^{+}\left(0, -H''(m_2)\right) + \tfrac{1}{2} \delta_0, 
		\end{align} 
		where $N^{\pm}$ are the half-normal distributions as in Definition $\ref{defn:halfnormal}$.

		\item[$(2)$] Suppose $p \geq 4$ is even.  Then the following hold:  
		
		\begin{enumerate}
			
			\item[$\bullet$] If $h \ne 0$,  then the function $H$ has exactly two (asymmetric) maximizers $m_1 < m_2$ and, as $N \rightarrow \infty$,
			\begin{align}\label{eq:h_3}
			N^{\frac{1}{2}}(\hat{h}_N - h) \xrightarrow{D} \tfrac{p_1}{2} N^{-}\left(0, -H''(m_1)\right) + \tfrac{1-p_1}{2} N^{+}\left(0, -H''(m_2)\right) + \tfrac{1}{2} \delta_0.\nonumber\\  
			\end{align}

			\item[$\bullet$] If $h=0$ and $\beta> \tilde{\beta}_p$, then the function $H$ has exactly two symmetric maximizers $m_1, m_2$, where $m_2 = -m_1 = m_*$, for some $m_* = m_*(\beta, h, p)> 0$. Then, as $N \rightarrow \infty$, 
			\begin{align}\label{eq:h_4}
			N^{\frac{1}{2}} \hat{h}_N  \xrightarrow{D} \tfrac{1}{2} N\left(0, -H''(m_*)\right)  + \tfrac{1}{2} \delta_0.
			\end{align}
			
			\item[$\bullet$] If $h=0$ and $\beta = \tilde{\beta}_p$, the function $H$ has three maximizers $m_1=-m_*$, $m_2=0$, and $m_3=m_*$, where $m_*=m_*(\beta, h, p) > 0$.  Then, as $N \rightarrow \infty$,
			\begin{align}\label{eq:h_2}
			N^{\frac{1}{2}}\hat{h}_N \xrightarrow{D} p_1 N\left(0, -H''(m_1)\right) + (1-p_1)\delta_0,
			\end{align} 
			where $p_1$ is as defined in \eqref{eq:p1}.  
		\end{enumerate} 
	\end{itemize}
	
\end{thm}

%\iffalse
The proofs of these results are given in Section \ref{sec:pf_beta_h} (a short roadmap of the proof is given in Section \ref{sec:pfsketch_mle}). The results above show that for all points in the parameter space, the ML estimate $\hat{h}_N$ is a consistent estimate of $h$, that is, $\hat{h}_N \pto h$. Moreover, the rate of convergence is $N^{\frac{1}{2}}$, except at the $p$-special points. However, at the $p$-special point(s), 
%(recall that there is only one such point when $p \geq 3$ is odd and two such points for $p \geq 4$ is even), 
the rate improves to $N^\frac{3}{4}$, that is, the ML estimate of $h$ at these point(s) is {\it superefficient}, converging to the true value of $h$ faster than the usual $N^{\frac{1}{2}}$ rate at the neighboring points. Another interesting feature is that, while at the regular points $\hat h_N$ has a simple Gaussian limit, at the critical points it has a mixture distribution, consisting of (half) normals and a point mass at $0$. The reason the limiting distribution has a point mass at $0$ is because the sample mean $\os$ is  ``discontinuous"  under the perturbed measure $\p_{\beta,h+ t/\sqrt N,p}$, as $t$ transitions from negative to positive. In fact, Lemma \ref{redunsup} (in Section \ref{sec:pf_concentration}) shows that under the measure $\p_{\beta,h+ t/\sqrt{N},p}$, the point where $\os$ concentrates depends on the sign of the perturbation factor $t$. Therefore, since the distribution function of $N^{\frac{1}{2}}(\hat{h}_N - h)$ evaluated at $t$ depends on the law of $\os$ under the perturbed measure $\p_{\beta,h+ t/\sqrt N, p}$ (see the calculations in Section \ref{sec:pf_beta_h_2} for details), it has a discontinuity at the point $t=0$, and, hence, a point mass at $0$ appears in the limit. 

Another interesting revelation are the results in \eqref{eq:h_31} and \eqref{eq:h_3}, where the $H$ function has two (asymmetric) maximizers. In this case, the ML estimate $\hat h_N$ converges to a three component mixture, which has a point mass at zero with probability $\frac{1}{2}$ and is a mixture of two half normal distributions, with probabilities $\frac{p_1}{2}$ and $\frac{1-p_1}{2}$, respectively. This corresponds to the region of the critical curve where $h \ne 0$ (and also the point $(\tilde \beta_p, 0)$, for $p \geq 3$ odd), a striking new phenomena that emerges only when $p \geq 3$. Note that, this does not happen for $p=2$, because, in this case, $\sC_p^+=(0.5, \infty)\times \{0\}$, hence, the two maximizers at any 2-critical point are symmetric about zero, and the two half normal mixing components combine to form a single Gaussian, and the resulting limit is the mixture of a single normal and a point mass at zero, as is the case in \eqref{eq:h_4} above.

%\fi

%%
%%
%%
%%
%%
%%
%%
%%
%%
%%
%%
%%
%%
%%
%%
%%
%%
%%
%%
%%
%%
%%
%%
%%

%\iffalse
\subsubsection{ML Estimate of $\beta$}
\label{sec:mle_beta} 

Here, we consider the ML estimate  $\hat{\beta}_N$ of $\beta$. As before, the results depend on whether $(\beta, h)$ is regular, critical, or special. However, the analysis here is more involved, and each of these cases breaks down into further cases, depending on the value of the maximizers, the parity of $p$, and the sign of $h$. We begin with the case when $(\beta, h)$ is regular. As always, $H=H_{\beta,p,h}$ will be as defined in \eqref{eq:H}.

\begin{thm}[Asymptotic distributions of $\hat{\beta}_N$ at $p$-regular points]\label{thmmle1} 
	Fix $p \geq 3$ and suppose $(\beta, h) \in \Theta$ is $p$-regular. Assume $h$ is known and $\bs \sim \mathbb{P}_{\beta,h,p} $. Then denoting the unique maximizer of $H$ by $m_*= m_*(\beta,h,p)$, the following hold, 
	\begin{itemize}
		
		\item[$\bullet$]  If $m_* \neq 0$, then, as $N \rightarrow \infty$, 
		\begin{align}\label{eq:bmle_m_1}
		N^{\frac{1}{2}}(\hat{\beta}_N-\beta) \xrightarrow{D} N\left(0, -\frac{H''(m_*)}{p^{2}m_*^{2p-2}} \right). 
		\end{align}
		
		\item[$\bullet$] If $m_* = 0$, (equivalently, $h = 0$ and $\beta < \tilde{\beta}_p$), then, as $N \rightarrow \infty$,
		\begin{align}\label{eq:bmle_m_2} 
		\hat{\beta}_N \xrightarrow{D} 
		\begin{cases}
		\frac{1}{2}\delta_{\tilde{\beta}_p} + \frac{1}{2} \delta_{-\tilde{\beta}_p} &\quad\text{if}~p~\textrm{is odd},\\
		\gamma_p \delta_{-\infty} + (1-\gamma_p)\delta_{\tilde{\beta}_p} &\quad\text{if}~p~\textrm{is even},\\
		\end{cases}
		\end{align}
		where $\gamma_p := \p(Z^p \leq \e Z^p)$ with $Z \sim N(0,1)$. 
	\end{itemize}
\end{thm}

We will discuss the various implications of the above theorem later in this section. Now, we state the result for the asymptotic distribution of $\hat{\beta}_N$ when $(\beta, h)$ is $p$-special.

\begin{thm}[Asymptotic distributions of $\hat{\beta}_N$ at $p$-special points]\label{thmmle2}	
	Fix $p \geq 3$ and suppose $(\beta, h) \in \Theta$ is $p$-special. Assume $h$ is known and $\bs \sim \mathbb{P}_{\beta,h,p} $. Then denoting the unique maximizer of $H$ by $m_*= m_*(\beta,h,p)$, as $N \rightarrow \infty$,
	\begin{align}\label{eq:bmle_m_3} 
	N^\frac{3}{4}(\hat{\beta}_N - \beta) \xrightarrow{D} G_2, 
	\end{align} 
	where the distribution function of $G_2$ is given by 
	$$G_2(t) = F_{0,0}\left(\int_{-\infty}^{\infty} u ~\mathrm d F_{t,0}(u)\right),$$ with $F_{t,0}$ as defined in  \eqref{eq:beta_h_distribution} below. 
\end{thm}

Finally, we consider the case $(\beta, h)$ is $p$-critical. The situation here is quite delicate, depending on various things like weak and strong criticality, parity of $p$, and the sign of the field $h$. 

\begin{thm}[Asymptotic distribution of $\hat{\beta}_N$ at $p$-critical points]\label{thmmle_III}
	Fix $p \geq 3$ and suppose $(\beta, h) \in \Theta$ is $p$-critical. Assume $h$ is known and $\bs \sim \mathbb{P}_{\beta,h,p} $. Denote the $K \in \{2,3\}$ maximizers of $H$ by $m_1 := m_1(\beta, h, p) < \ldots < m_K := m_K(\beta, h, p)$, and let $p_1, \ldots, p_K$ be as in \eqref{eq:p1}. 
	
	\begin{itemize}
		
		\item[$(1)$] Suppose $p \geq 3$ is odd. In this case, the function has exactly two maximizers  $m_1 < m_2$. Then, as $N \rightarrow \infty$, the following hold: 
		
		\begin{enumerate}
			\item[$\bullet$] If $(\beta,h) \neq (\tilde{\beta}_p,0)$, where $\tilde{\beta}_p$ is defined in \eqref{eq:betatilde}, then 
			\begin{align}\label{eq:bmle_multimax_1} 
			N^{\frac{1}{2}}(\hat{\beta}_N - \beta)\xrightarrow{D} \tfrac{p_1}{2}N^-\left(0,-\frac{H''(m_1)}{p^2m_1^{2p-2}}\right)+ \tfrac{1-p_1}{2}N^+\left(0,-\frac{H''(m_2)}{p^2m_2^{2p-2}}\right)+\tfrac{1}{2}\delta_0.\nonumber\\
			\end{align} 
			
			\item[$\bullet$] If $(\beta,h) = (\tilde{\beta}_p,0)$, then 
			\begin{align}\label{eq:bmle_multimax_2} 
			N^{\frac{1}{2}}(\hat{\beta}_N - \beta)\xrightarrow{D} \tfrac{p_1}{2}\delta_{-\infty} + \tfrac{1-p_1}{2}N^+\left(0,-\frac{H''(m_2)}{p^2m_2^{2p-2}}\right)+\tfrac{1}{2}\delta_0. 
			\end{align} 
		\end{enumerate}
		
		\item[$(2)$] Suppose $p \geq 4$ is even. Then the following hold, as $N \rightarrow \infty$:   
		\begin{enumerate}
			\item[$\bullet$] If $h > 0$, then 
			\begin{align}\label{eq:bmle_multimax_3}
			N^{\frac{1}{2}}(\hat{\beta}_N - \beta)\xrightarrow{D} \tfrac{p_1}{2}N^-\left(0,-\frac{H''(m_1)}{p^2m_1^{2p-2}}\right)+ \tfrac{1-p_1}{2}N^+\left(0,-\frac{H''(m_2)}{p^2m_2^{2p-2}}\right)+\tfrac{1}{2}\delta_0.\nonumber\\
			\end{align}
			
			\item[$\bullet$]  If $h < 0$, then 
			\begin{align}\label{eq:bmle_multimax_4} 
			N^{\frac{1}{2}}(\hat{\beta}_N - \beta)\xrightarrow{D} \tfrac{p_1}{2}N^+\left(0,-\frac{H''(m_1)}{p^2m_1^{2p-2}}\right)+ \tfrac{1-p_1}{2}N^-\left(0,-\frac{H''(m_2)}{p^2m_2^{2p-2}}\right)+\tfrac{1}{2}\delta_0.\nonumber\\ 
			\end{align}
			
			\item[$\bullet$] If $h=0$ and $\beta> \tilde{\beta}_p$, there are exactly two maximizers $m_1=-m_*$ and $m_2=m_*$ of $H$, where $m_* = m_*(\beta, h, p)> 0$. In this case, 
			\begin{align}\label{eq:bmle_multimax_5} 
			N^{\frac{1}{2}}(\hat{\beta}_N - \beta)\xrightarrow{D} N\left(0,-\frac{H''(m_*)}{p^2m_*^{2p-2}}\right).
			\end{align}
			
			\item[$\bullet$] If $h=0$ and $\beta = \tilde{\beta}_p$, there are exactly three maximizers $m_1 = -m_*$, $m_2=0$, and $m_3 = m_*$ of $H$, where $m_*=m_*(\beta, h, p) > 0$. In this case, \begin{align}\label{eq:bmle_multimax_6} 
			N^{\frac{1}{2}}(\hat{\beta}_N - \beta)\xrightarrow{D} p_2 \gamma_p \delta_{-\infty} + p_1N^+\left(0,-\frac{H''(m_*)}{p^2m_*^{2p-2}}\right) + (1-p_1-p_2\gamma_p) \delta_0,\nonumber\\
			\end{align} 
			where $\gamma_p := \p(Z^p \leq \e Z^p)$ and $Z$ is a standard normal random variable.
		\end{enumerate}
	\end{itemize}	
\end{thm}

The proofs of the above results are given in Section \ref{sec:pf_beta_h}. One of the main technical ingredients is a strengthening of Theorem \ref{cltintr1}, which requires obtaining the asymptotic distribution of $\os$ when the parameters $(\beta,h)$ are perturbed by an $o(N)$ term. A more detailed overview of the proof technique is given in Section \ref{sec:pfsketch_mle}. Here, we summarize the main consequences of the above results and highlight the various new phenomena that emerge as one moves from the matrix $(p=2)$ to the tensor $(p \geq 3)$ case.

\begin{itemize}

	\item  For $p$-regular points, Theorem \ref{thmmle1} shows that when the unique maximizer $m_* \ne 0$, then $\hat{\beta}_N$ is consistent at rate $N^{\frac{1}{2}}$ with a limiting normal distribution. On the other hand, when $m_* = 0$, which happens in the interval  $[0, \tilde \beta_p)$, the ML estimate $\hat{\beta}_N$ is inconsistent. In this regime, when $p \geq 3$ is odd, then $\hat{\beta}_N$ concentrates at $\pm \bm \tilde \beta_p$ with probability $\frac{1}{2}$, irrespective of the value of true value of $\beta \in  [0, \tilde \beta_p)$. The situation is even more strange when $p \geq 4$ is even. Here,  $\hat{\beta}_N$ concentrates at either $\bm \tilde \beta_p$ or escapes to negative infinity, that is, with positive probability $\hat{\beta}_N$ is unbounded, when $p \geq 4$ and $\beta \in  [0, \tilde \beta_p)$. The corresponding results for $p=2$ are similar in the sense that, for $\beta \in [0, 0.5)$ (recall that $\tilde{\beta_2}=0.5$), the ML estimate $\hat{\beta}_N$ is inconsistent. However, unlike in the case for $p\geq 4$ even, the ML estimate $\hat{\beta}_N$, when $p=2$, is always finite and converges to a (properly centered and rescaled) chi-squared distribution \cite[Theorem 1.4]{comets}.

	\item For $p$-special points Theorem \ref{thmmle2} shows that $\hat{\beta}_N$  converges to $\beta$ at rate $N^{-\frac{3}{4}}$, that is, it is superefficient. Recall that the same thing happens for $\hat h_N$ at $p$-special points (Theorem \ref{cltintr3_2}). In comparison, for $p=2$ at the only 2-special $(0.5, 0)$,  $\hat{h}_N$ is superefficient with rate $N^{-\frac{3}{4}}$ \cite[Theorem 1.3]{comets}, but $\hat{\beta}_N$ remains $N^{\frac{1}{2}}$-consistent \cite[Theorem 1.4]{comets}. This is because when $(\beta,h)$ is $p$-special, the unique maximizer $m_*$ of $H_{\beta,h,p}$ is $0$ when $p=2$, but non-zero, for $p \geq 3$. This creates a difference in the rate of convergence of the maxima of $H_{\beta_N,h_N,p}$ towards the maximum of $H_{\beta,h,p}$ for some suitably chosen perturbation $(\beta_N,h_N)$ of $(\beta,h)$, which is an important step in deriving the asymptotic rate of convergence of the ML estimates. Another interesting difference is that for $p=2$, the only 2-special point $(0.5, 0)$ coincides with the thermodynamic threshold of the 2-tensor Curie-Weiss model. However, for $p\geq 3$, the $p$-special points (the point $(\check{\beta}_p,\check{h}_p)$, for $p \geq 3$ odd, and the points 
	$(\check{\beta}_p, \pm \check{h}_p)$, for $p \geq 4$ even), where we get the non-Gaussian limits of $\os$, ${\hat h}_N$, and ${\hat \beta}_N$, have nothing to do with the thermodynamic threshold of the $p$-tensor Curie-Weiss model, but rather depends on the vanishing property of the second derivative of $H$ at its maximizer. On the contrary,  quite remarkably, the thermodynamic threshold $(\tilde{\beta}_p, 0)$ of the $p$-tensor Curie-Weiss model (recall definition in \eqref{eq:betatilde}) turns out to be a $p$-weakly critical point for $p \geq 3$ odd, and the only $p$-strongly critical point for $p \geq 4$ even, another unexpected phenomenon unearthed by our results.

	\item  The landscape is much more delicate for $p$-critical points, as can be seen from Theorem \ref{thmmle_III}. In this case, the limiting distribution of $\hat{\beta}_N$ converges to various mixture distributions, depending on, among other things, the sign of $h$ and the parity of $p$. As in the case of $\hat{h}_N$, a particularly interesting new phenomena is the three component mixture that arises in the limiting distribution of $\hat{\beta}_N$ when the critical curve $\sC_p^+$ intersects the region $h \ne 0$. This corresponds to the result \eqref{eq:bmle_multimax_1} for $p \geq 3$ odd, and results in \eqref{eq:bmle_multimax_3} and \eqref{eq:bmle_multimax_4}  for $p \geq 4$ even. Recall, from the discussion following Theorem \ref{cltintr3_III}, that this does not happen for $p=2$, because, in this case, $\sC_p^+=(0.5, \infty)\times \{0\}$, hence, the two maximizers at any 2-critical point are symmetric about zero, and the two half normal mixing components combine to form a single Gaussian. As a result, the limit is the mixture of a single normal and a point mass at zero. Interestingly, this also happens for $p \geq 4$ even, when the critical curve intersects the line $h=0$ and is strictly above the threshold $\tilde{\beta}_p$, as seen in  \eqref{eq:bmle_multimax_5} above.

	\item The final bit in the puzzle is the point of thermodynamic phase transition $(\tilde \beta_p, 0)$. Here, the ML estimate $\hat{\beta}_N$  is not $N^{\frac{1}{2}}$-consistent. More precisely, in the limit, $N^{\frac{1}{2}}(\hat \beta_N-\beta)$ has a point mass at negative infinity with positive probability,  and is a mixture of a folded normal and a point mass with the remaining probability (as described in \eqref{eq:bmle_multimax_2} and \eqref{eq:bmle_multimax_6}). In contrast, as explained in the second case above, when $p=2$, then at the point of thermodynamic phase transition ($\tilde{\beta}_2=0.5$) the ML estimate ${\hat \beta}_N$ is $N^{\frac{1}{2}}$-consistent. 
	
\end{itemize}

%\iffalse

\subsection{Summarizing the Phase Diagram}  
\label{sec:mle_beta_h_II}

\begin{figure}
	\centering
	\begin{minipage}[c]{1.0\textwidth}
		\centering
		\includegraphics[width=6.25in,height=3.05in]
		{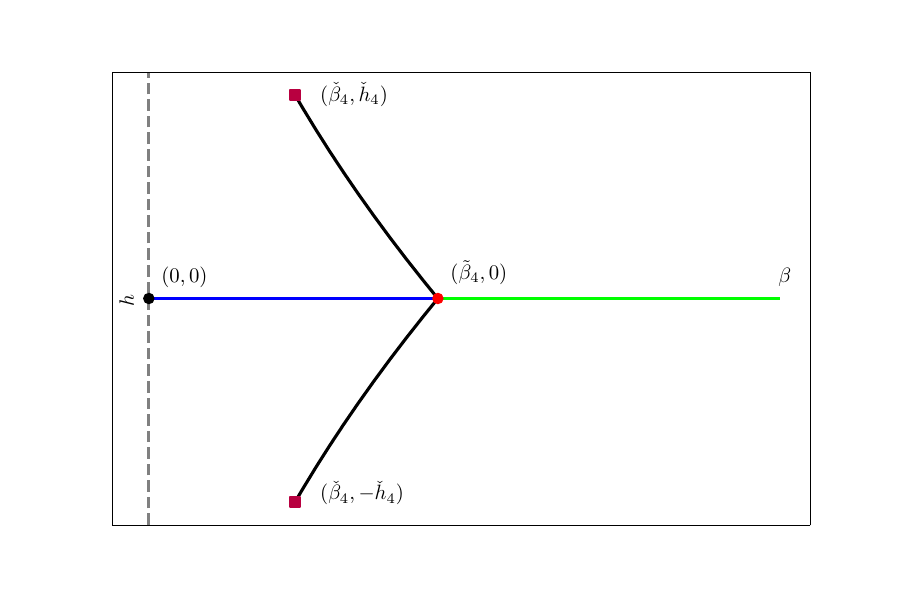}\\
		\caption{\small{The phase diagram for $p=4$: The properties of the ML estimates in the  different regions of the parameter space $\Theta= [0, \infty) \times \R$ are as follows: }}
		\label{figure:ordering1}
		\small{
			\begin{itemize} 
				\item[] 
				\begin{itemize}
					
					\item The $\mysquare[white]$ (white) region: These are the $p$-regular points where $H$ has a unique global maximizer $m_* \ne 0$  and $H''(m_*) < 0$.  Hence, $\hat \beta_N$ and $\hat h_N$  are both $N^{\frac{1}{2}}$-consistent and asymptotically normal, by \eqref{eq:bmle_m_1} and \eqref{eq:h_1}, respectively. 
					
					\item The \textcolor{blue}{\rule{0.75cm}{0.75mm}} line: These are the $p$-regular points where $H$ has a unique global maximizer $m_* = 0$  and $H''(0) < 0$. Hence, $\hat \beta_N$ is inconsistent by \eqref{eq:bmle_m_2},  but $\hat h_N$ is $N^{\frac{1}{2}}$-consistent and asymptotically normal by \eqref{eq:h_1}. 
					
					\item The $\mysquare[purple]$ points: These are the $p$-special points. Here, $H$ has a unique maximizer $m_*$, but $H''(m_*)=0$. Hence,  $\hat \beta_N$ and $\hat h_N$  are both superefficient, converging at rate $N^{\frac{3}{4}}$ to non-Gaussian distributions, by \eqref{eq:bmle_m_3} and \eqref{eq:h_5}, respectively.

					\item The \textcolor{black}{\rule{0.75cm}{0.75mm}} curve: These are $p$-weakly critical points where $h \ne 0$. Here, $H$ has two global (non-symmetric) maximizers. Both $\hat \beta_N$ and $\hat h_N$  are $N^{\frac{1}{2}}$-consistent and asymptotically a three component mixture (comprising of two half normal normal distributions and a point mass at zero), by  \eqref{eq:bmle_multimax_3}, \eqref{eq:bmle_multimax_4},  and \eqref{eq:h_3}, respectively.

					\item The \textcolor{green}{\rule{0.75cm}{0.75mm}} line: These are $p$-weakly critical points where $h = 0$. Here, $H$ has two global symmetric maximizers. Hence, $\hat \beta_N$ is $N^{\frac{1}{2}}$-consistent and asymptotically normal by \eqref{eq:bmle_multimax_5}, and $\hat h_N$ is $N^{\frac{1}{2}}$-consistent and asymptotically a mixture of a normal distribution and a point mass at zero, by \eqref{eq:h_4}. 
					
					\item The \tikz\draw[red,fill=red] (0,0) circle (.8ex); point: This is the $p$-strongly critical point. Here, $H$ has three global maximizers. Hence, $\hat{\beta}_N$ is not $N^{\frac{1}{2}}$-consistent, by \eqref{eq:bmle_multimax_6}, but $\hat{h}_N$ is $N^{\frac{1}{2}}$-consistent and asymptotically a mixture of normal distribution and point mass at 0, by \eqref{eq:h_2}. 	
				\end{itemize}
			\end{itemize}
		}
	\end{minipage}
\end{figure}

\begin{figure}
	\centering
	\begin{minipage}[c]{1.0\textwidth}
		\centering
		\includegraphics[width=6.25in,height=3.05in]
		{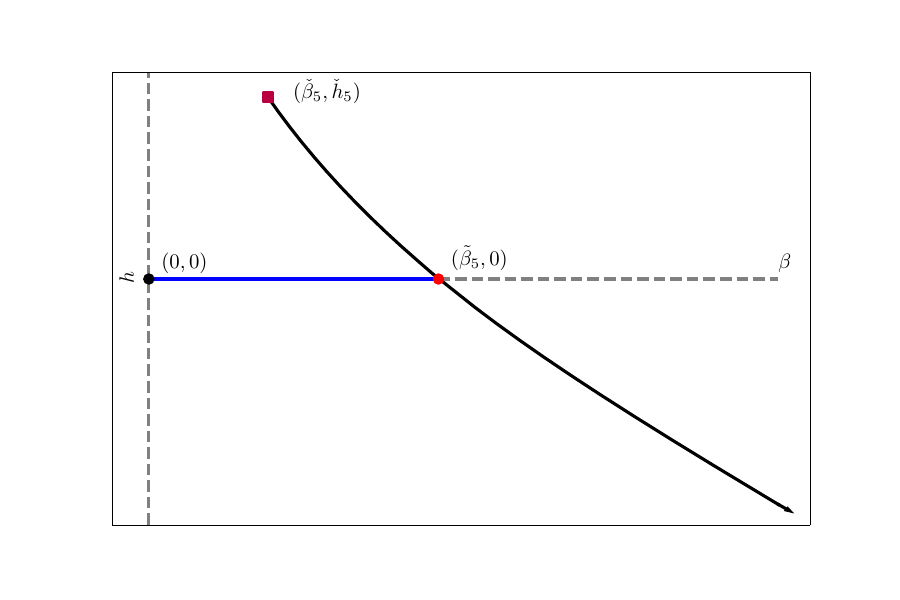}\\
		\caption{\small{The phase diagram for $p=5$: The properties of the ML estimates in the  different regions of the parameter space $\Theta= [0, \infty) \times \R$ are as follows: }}
		\label{figure:ordering2}
		\small{
			\begin{itemize} 
				\item[] 
				\begin{itemize}
					
					\item The $\mysquare[white]$ (white) region: These are the $p$-regular points where $H$ has a unique global maximizer $m_* \ne 0$  and $H''(m_*) < 0$.  Hence, $\hat \beta_N$ and $\hat h_N$  are both $N^{\frac{1}{2}}$-consistent and asymptotically normal, by \eqref{eq:bmle_m_1} and \eqref{eq:h_1}, respectively. 
					
					\item The \textcolor{blue}{\rule{0.75cm}{0.75mm}} line: These are the $p$-regular points where $H$ has a unique global maximizer $m_* = 0$  and $H''(0) < 0$. Hence, $\hat \beta_N$ is inconsistent by \eqref{eq:bmle_m_2},  but $\hat h_N$   $N^{\frac{1}{2}}$-consistent and asymptotically normal by \eqref{eq:h_1}. 
					
					\item The $\mysquare[purple]$ point: This is the only $p$-special point. Here, $H$ has a unique maximizer $m_*$, but $H''(m_*)=0$. Hence,  $\hat \beta_N$ and $\hat h_N$  are both superefficient, converging at rate $N^{\frac{3}{4}}$ to non-Gaussian distributions, by \eqref{eq:bmle_m_3} and \eqref{eq:h_5}, respectively.

					\item The \textcolor{black}{\rule{0.75cm}{0.75mm}} curve: These are $p$-weakly critical points where $h \ne 0$. Here, $H$ has two global (non-symmetric) maximizers. Both, $\hat \beta_N$ and $\hat h_N$  are $N^{\frac{1}{2}}$-consistent and asymptotically a three component mixture (comprising two half normal normal distributions and a point mass at zero), by  \eqref{eq:bmle_multimax_1} and \eqref{eq:h_31}, respectively.

					\item The \tikz\draw[red,fill=red] (0,0) circle (.8ex); point: This is the $p$-weakly critical point with $h=0$. Here, $H$ has two (non-symmetric) global maximizers. Hence, $\hat h_N$  is asymptotically a three component mixture by  \eqref{eq:h_31}, but $\hat{\beta}_N$ is not $N^{\frac{1}{2}}$-consistent by \eqref{eq:bmle_multimax_2}. 
				\end{itemize}
			\end{itemize}
		}
	\end{minipage}
\end{figure}

The results above can be compactly summarized and better visualized in a phase diagram, which shows the  partition of the parameter space described in \eqref{eq:parameter_space}. The phase diagrams for $p=4$ and $p=5$, obtained by numerical optimization of the function $H$ over a fine grid of parameter values, are shown in Figure \ref{figure:ordering1} and Figure \ref{figure:ordering2}, respectively.   The limiting distributions that arise in the different regions of the phase diagram are described in the figure legends.

%\iffalse

\section{Asymptotic Distribution of the Sample Mean: Proof of Theorem \ref{cltintr1}} 
\label{cltmagnetization}

In this section, we will prove  Theorem \ref{cltintr1} in the $p$-regular case and present proof roadmaps for the remaining two cases. To this end, note that the model \eqref{cwwithmg} can be written more compactly as  
\begin{align*}%\label{cwwithmg}
\mathbb{P}_{\beta,h,p}(\bs) & = \frac{1}{2^{N} Z_N(\beta,h,p)} \exp\Big\{N\left(\beta \overline{X}^p_N + h\os  \right) \Big\}, 
\end{align*}
where $\os := \frac{1}{N} \sum_{i=1}^N X_i$ is the sample mean. Therefore, the sample mean has the probability mass function, 
\begin{align*}%\label{eq:m}
\p_{\beta,h,p}(\overline{X}_N=m)= \frac{1}{2^{N}Z_N(\beta,h,p)} \binom{N}{\frac{N(1+m)}{2}}e^{N(\beta m^p + hm)}
\end{align*} 
$\text{for} ~m \in \left\{-1,-1+\dfrac{2}{N},\ldots, 1-\dfrac{2}{N},1\right\}.$

Observe that the expression for the probability mass function of $\os$ involves the partition function $Z_N(\beta,h,p)$, which does not have a closed form.\footnote{Note that even though the partition function has no closed form, for a given value of $(\beta, h)$ it can be easily computed in $O(N)$ time. This is one of the major advantages of the Curie-Weiss model, which 
	allows, among other things, efficient computation of the ML estimates. } Therefore, obtaining limiting properties of $\os$ requires accurate estimation of $Z_N(\beta,h,p)$. Moreover, as mentioned before, with the goal of deriving the limiting distribution of the ML estimates of $\beta$ and $h$, we will need to prove the limiting distribution of $\os$ at slightly perturbed parameter values $(\beta_N, h_N)$, for some sequences $\beta_N \rightarrow \beta$ and $h_N \rightarrow h$ to be chosen later. Hereafter, we will denote $\p_{\beta_N, h_N,p}, ~Z_N(\beta_N,h_N,p)$, and $F_N(\beta_N,h_N,p)$,  by $\bar{\p}, \bar{Z}_N$,  and $\bar{F}_N$, respectively. The asymptotic distribution of $\os$ in the different cases at the appropriately perturbed parameter values is summarized below:

\begin{thm}[Asymptotic distribution of $\os$ under perturbed parameters]\label{cltun} Fix $p \geq 3$, $(\beta, h) \in \Theta$, and $\bar{\beta},\bar{h}\in \mathbb{R}$. Then with $H=H_{\beta,p,h}$ as defined in \eqref{eq:H} the following hold: 
	
	\begin{itemize}
		
		\item[$(1)$] Suppose $(\beta, h)$ is $p$-regular and denote the unique maximizer of $H$ by $m_*= m_*(\beta,h,p)$. Then, for $\bs \sim \mathbb{P}_{\beta+N^{-\frac{1}{2}}\bar{\beta},~h + N^{-\frac{1}{2}}\bar{h},~p}$, as $N \rightarrow \infty$,
		\begin{align}\label{eq:cltun_I} 
		N^{\frac{1}{2}}\left(\os - m_*(\beta,h,p)\right)\xrightarrow{D} N\left(-\frac{\bar{h}+\bar{\beta}pm_*(\beta,h,p)^{p-1}}{H''(m_*)},~-\frac{1}{H''(m_*)}\right). 
		\end{align} 
		
		\item[$(2)$] Suppose $(\beta, h)$ is $p$-critical and denote the $K \in \{2, 3\}$ maximizers  of $H$ denoted by $m_1:=m_1(\beta,h,p),$$\ldots,m_K:=m_K(\beta,h,p)$. Then, for $\bs \sim \p_{\beta, h, p}$, as $N \rightarrow \infty$, 
		\begin{align}\label{eq:cltun_p1}
		\os \xrightarrow{D} \sum_{k=1}^K p_k \delta_{m_k}, 
		\end{align}
		where $p_1, \ldots, p_K$ are as defined in \eqref{eq:p1}. Moreover, if $m$ is any local maximizer of $H$ contained in the interior of an interval $A \subseteq [-1,1]$, such that $H(m) > H(x)$ for all $x\in A\setminus \{m\}$, then for $\bs \sim \mathbb{P}_{\beta+N^{-\frac{1}{2}}\bar{\beta},~h + N^{-\frac{1}{2}}\bar{h},~p}$, as $N \rightarrow \infty$, 
		\begin{equation}\label{eq:cltun_II}
		N^{\frac{1}{2}}\left(\os - m\right)\Big\vert \{ \os \in A\} \xrightarrow{D} N\left(-\frac{\bar{h}+\bar{\beta}pm^{p-1}}{H''(m)},~-\frac{1}{H''(m)}\right).
		\end{equation}

		\item[$(3)$] Suppose $(\beta, h)$ is $p$-special and denote the unique maximizer of $H$ by $m_*= m_*(\beta,h,p)$. Then, for $\bs \sim \mathbb{P}_{\beta+ N^{-\frac{3}{4}}\bar{\beta},~h + N^{-\frac{3}{4}}\bar{h},~p}$, as $N \rightarrow \infty$, 
		\begin{align}\label{eq:cltun_III} 
		N^\frac{1}{4}(\os - m_*(\beta,h,p)) \xrightarrow{D} F_{\bar{\beta},\bar{h}},
		\end{align} 
		where density of $F$ with respect to the Lebesgue measure is given by 
		\begin{align}\label{eq:beta_h_distribution}
		\frac{\mathrm dF_{\bar{\beta},\bar{h}} (x) }{\mathrm d x} ~\propto~ \exp\left(\frac{H^{(4)}(m_*)}{24} x^4 + (\bar{\beta}pm_*^{p-1} + \bar{h})x\right). 
		\end{align}
	\end{itemize} 
\end{thm}

Note that Theorem \ref{cltintr1} follows directly from the above by taking $\bar{\beta}=0$ and $\bar{h}=0$. We prove Theorem \ref{cltun} in the $p$-regular case in Section \ref{sec:pfcltun_I} below. The roadmaps of the remaining two cases, which follow a similar strategy but requires more a delicate analysis, are described in Section \ref{sec:pfcltun_II} and  Section \ref{nonpuniqueclt}. The complete proofs of \eqref{eq:cltun_II} and \eqref{eq:cltun_III} are given in Section \ref{sec:proofofnonun} and \ref{irregproof}, respectively. 

\subsection{Proof of Theorem \ref{cltun} when $(\beta, h)$ is $p$-regular} 
\label{sec:pfcltun_I}

Fix a $p$-regular point $(\beta, h) \in \Theta$ and consider a sequence $(\beta_N, h_N) \in \Theta$ (to be specified later) such that $\beta_N \rightarrow \beta$ and $h_N \rightarrow h$. It has been shown in  Lemma \ref{derh44} that the function $H_N(x) := H_{\beta_N,h_N,p}(x)$ will have a unique global maximizer $m_*(N)$, for all large $N$, and $m_*(N)\rightarrow m_*$ as $N \rightarrow \infty$. Choose this maximizer $m_*(N)$ and define, for $\alpha \in (0, 1)$, 
\begin{align}\label{eq:an_maximizer}
A_{N,\alpha} := \left(m_*(N)-N^{-\frac{1}{2} + \alpha},m_*(N) +N^{-\frac{1}{2} + \alpha}\right). 
\end{align}
The first step in the  proof of Theorem \ref{cltun} when $(\beta, h)$ is $p$-regular,  is to show that under $\bar{\p}$, the sample mean $\os$ concentrates around $m_*(N)$ at rate $N^{-\frac{1}{2} + \alpha}$, for any $\alpha > 0$. 

\begin{lem}\label{conc}
	Suppose  $(\beta, h) \in \Theta$ is $p$-regular. Then for $\alpha \in \left(0,\frac{1}{6}\right]$ and  $A_{N,\alpha}$ as defined above in \eqref{eq:an_maximizer},\footnote{For any set $A$, $A^c$ denotes the complement of the set $A$.}
	\begin{equation*}
	\bar{\mathbb P}\left(\os \in A_{N,\alpha}^c\right) = \exp\left\{\frac{1}{3}N^{2\alpha} H''(m_*) \right\}O(N^{\frac{3}{2}}).
	\end{equation*}
\end{lem}

\begin{proof} Note that the support of the magnetization $\os$ is the set  $$\mathcal{M}_N:=\left\{-1,-1+\dfrac{2}{N},\ldots, 1-\dfrac{2}{N},1\right\}.$$	It follows from \cite{talagrand}, Equation (5.4), that for any $m \in \mathcal{M}_N$, the cardinality of the set
	$$A_m := \left\{\bs \in \sa_N: \os = m\right\}$$ can be bounded by
	\begin{equation}\label{boundtal}
	\frac{2^N}{LN^{\frac{1}{2}}} \exp\left\{-N I(m) \right\}\leq |A_m| \leq 2^N\exp\left\{-N I(m) \right\}
	\end{equation}
	for some universal constant $L$ (recall that $I(\cdot)$ is the binary entropy function). Hence, we have from \eqref{boundtal},
	\begin{align}\label{sb}
	\bar{\p}(\os \in A_{N,\alpha}^c) & =  \frac{\sum_{m \in \mathcal{M}_N \bigcap A_{N,\alpha}^c} |A_m|\exp\left\{N(\beta_N m^p + h_Nm) \right\}}{\sum_{m \in \mathcal{M}_N} |A_m|\exp\left\{N(\beta_N m^p + h_Nm) \right\}}\nonumber\\& \leq  \frac{LN^{\frac{1}{2}}(N+1) \sup_{x \in A_{N,\alpha}^c}e^{NH_N(x)}}{\sup_{x\in [-1,1]}e^{NH_N(x)} }\nonumber\\& =  \exp\left\{N\left(\sup_{x\in A_{N,\alpha}^c} H_N(x) - H_N\left(m_*(N)\right)\right)\right\} O(N^{\frac{3}{2}}).
	\end{align}
	By Lemma \ref{mest}, we know that for all large $N$, $\sup_{x\in A_{N,\alpha}^c} H_N(x)$ is either $H_N(m_*(N)-N^{-\frac{1}{2} + \alpha})$ or $H_N(m_*(N)+N^{-\frac{1}{2} + \alpha})$. Since $H_N'\left(m_*(N)\right) = 0$ and the functions $H_N^{(3)}$ are uniformly bounded on any closed interval contained in $(-1,1)$, Taylor's theorem gives us:
	\begin{align}
	H_N\left(m_*(N)\pm N^{-\frac{1}{2} + \alpha}\right) - H_N(m_*(N)) & =  \frac{1}{2}N^{-1+2\alpha} H_N''(m_*(N)) + O\left(N^{-\frac{3}{2} + 3\alpha}\right)\label{taylorH1}\\& \leq  \frac{1}{3}N^{-1+2\alpha} H''(m_*) +  O\left(N^{-\frac{3}{2} + 3\alpha}\right)\label{taylorH2}.
	\end{align}
	Note that \eqref{taylorH2} follows from \eqref{taylorH1} since $H_N''(m_*(N)) \rightarrow H''(m_*) <0$.  The proof of Lemma \ref{conc} is now complete, in view of \eqref{sb}. 
\end{proof}

Lemma \ref{conc} shows that almost all contribution to $\bar{Z}_N$ comes from configurations whose average lies in a vanishing neighborhood of the maximizer $m_*(N)$ of $H_N$. This enables us to accurately approximate  the partition function $\bar{Z}_N$. This involves a Riemann approximation of the sum of the mass function $\p_{\beta_N,h_N,p}(\bs)$ over all $\bm X$ whose mean lies in a vanishing neighborhood of $m_*$, followed by a further saddle-point approximation of the resulting integral. 

\begin{lem}\label{ex} Suppose  $(\beta, h) \in \Theta$ is $p$-regular. Then for $\alpha > 0$ and $N$ large enough, the partition function can be expanded as,
	\begin{equation}\label{partex}
	\bar{Z}_N =  \frac{e^{NH_N(m_*(N))}}{\sqrt{(m_*(N)^2-1)H_N''(m_*(N))}}\left(1+O\left(N^{-\frac{1}{2} +\alpha}\right)\right),
	\end{equation}
	where $m_*(N)$ is the unique maximizer of the function $H_N$. Moreover, for $N$ large enough, the log-partition function can be expanded as, 
	\begin{equation}\label{logpartex}
	\bar{F}_N =  N H_N(m_*(N)) - \tfrac{1}{2}\log \left[(m_*(N)^2-1)H_N''(m_*(N))\right] + O\left(N^{-\frac{1}{2} +\alpha}\right).
	\end{equation}
\end{lem}

\begin{proof}
	Without loss of generality, let $\alpha \in \left(0,\frac{1}{6}\right]$ and note that
	\begin{equation}\label{fststep}
	\bar{\mathbb P}(\os \in A_{N,\alpha}) = \bar{Z}_N^{-1}\sum_{m \in \mathcal{M}_N\bigcap A_{N,\alpha}} \binom{N}{N(1+m)/2} \exp\left\{N(\beta_N m^p + h_Nm-\log 2) \right\}.
	\end{equation}
	By Lemma \ref{conc}, $\bar{\mathbb P}(\os \in A_{N,\alpha}) = 1-O\left(e^{-N^\alpha}\right)$ and hence \eqref{fststep} gives us
	\begin{align}\label{secstep}
	\bar{Z}_N & =  \left(1+ O\left(e^{-N^\alpha}\right) \right)\sum_{m \in \mathcal{M}_N\bigcap A_{N,\alpha}} \binom{N}{N(1+m)/2} \exp\left\{N(\beta_N m^p + h_Nm-\log 2) \right\}\nonumber\\& =  \left(1+ O\left(e^{-N^\alpha}\right) \right) \sum_{m \in \mathcal{M}_N\bigcap A_{N,\alpha}} \zeta(m)
	\end{align}
	where $\zeta:[-1,1] \rightarrow \mathbb{R}$ is defined as
	\begin{equation}\label{xidef}
	\zeta(x) := \binom{N}{N(1+x)/2} \exp\left\{N(\beta_N x^p + h_Nx - \log 2)\right\},
	\end{equation}
	where $\binom{N}{N(1+x)/2}$ is interpreted as a continuous binomial coefficient (refer to Section \ref{mathfunc} for the definition of continuous binomial coefficients). The next step is to approximate the sum in \eqref{secstep} by an integral, using Lemma \ref{Riemann}. Note that Lemma \ref{Riemann} can be applied with $n = \Theta(N^{\frac{1}{2}+\alpha})$ to obtain (using Lemma \ref{imp1}),
	\begin{align}\label{rax}
	\left|\int_{A_{N,\alpha}} \zeta(x)\mathrm d x - \frac{2}{N} \sum_{m \in \mathcal{M}_N\bigcap A_{N,\alpha}} \zeta(m)\right| &\leq \Theta(N^{-\frac{1}{2}+\alpha})N^{-1}\sup_{x\in A_{N,\alpha}}|\zeta'(x)|\nonumber\\
	&= O\left(N^{-\frac{1}{2}+\alpha}\cdot N^{-1}\cdot N^{\frac{1}{2}+\alpha}\right)\zeta(m_*(N))\nonumber\\
	&= O\left(N^{-1 + 2\alpha}\right)\zeta(m_*(N)).
	\end{align}
	It now follows from \eqref{rax}, Lemma \ref{stir}, Lemma \ref{Laplace} and Lemma \ref{xiact}, that
	\begin{align}\label{fin1}
	&\sum_{m \in \mathcal{M}_N\bigcap A_{N,\alpha}} \zeta(m)\nonumber\\ 
	&= \frac{N}{2} \int_{A_{N,\alpha}} \zeta(x)\mathrm d x + O(N^{2\alpha}) \zeta(m_*(N))\nonumber\\
	&= \frac{N^{\frac{1}{2}}}{2}\left(1+O(N^{-1})\right)\int_{A_{N,\alpha}} e^{NH_N(x)}\sqrt{\frac{2}{\pi (1-x^2)}}\mathrm d x + O(N^{2\alpha}) \zeta(m_*(N))\nonumber\\
	&= \frac{N^{\frac{1}{2}}}{2} \sqrt{\frac{2\pi}{N|H_N''(m_*(N))|}}\sqrt{\frac{2}{\pi(1-m_*(N)^2)}} e^{N H_N(m_*(N))}\left(1+O\left(N^{-\frac{1}{2}+3\alpha}\right)\right)\nonumber\\ &+ O(N^{2\alpha}) \zeta(m_*(N))\nonumber\\
	&= \frac{e^{N H_N(m_*(N))}}{\sqrt{(m_*(N)^2-1)H_N''(m_*(N))}}\left(1+O\left(N^{-\frac{1}{2}+3\alpha}\right)\right)\nonumber\\ 
	&+ \sqrt{\frac{2}{\pi N (1-m_*(N)^2)}} e^{NH_N(m_*(N))}\left(1+O(N^{-1})\right)O(N^{2\alpha})\nonumber\\ 
	&= \frac{e^{N H_N(m_*(N))}}{\sqrt{(m_*(N)^2-1)H_N''(m_*(N))}}\left(1+O\left(N^{- \frac{1}{2}+3\alpha}\right)\right).
	\end{align}
	Combining \eqref{secstep} and \eqref{fin1}, we have:
	\begin{align}\label{fin2}
	\bar{Z}_N & =  \left(1+ O\left(e^{-N^\alpha}\right) \right)\left(1+O\left(N^{- \frac{1}{2}+3\alpha}\right)\right)\frac{e^{N H_N(m_*(N))}}{\sqrt{(m_*(N)^2-1)H_N''(m_*(N))}}\nonumber\\& =  \left(1+O\left(N^{- \frac{1}{2}+3\alpha}\right)\right)\frac{e^{N H_N(m_*(N))}}{\sqrt{(m_*(N)^2-1)H_N''(m_*(N))}}.
	\end{align}
	This completes the proof of \eqref{partex}. If we take logarithm on all sides in \eqref{fin2} and use the fact that $\log\left(1+ O(a_n)\right) = O(a_n)$ for any sequence $a_n = o(1)$, then we get \eqref{logpartex}, completing the proof. 
\end{proof}

\noindent\textbf{\textit{Completing the Proof of}} \eqref{eq:cltun_I}: We now have all the necessary ingredients in order to derive the CLT for $\os$ when $(\beta, h)$ is $p$-regular. Throughout this subsection, we take 
$$\beta_N = \beta+N^{-\frac{1}{2}}\bar{\beta} \quad \text{and} \quad  h_N = h + N^{-\frac{1}{2}}\bar{h},$$ for some fixed $\beta\geq 0$ and $\bar{\beta}, h, \bar{h} \in \mathbb{R}$.  Now, recall that $H_N := H_{\beta_N,h_N,p}$ and $m_*=m_*(\beta,h,p)$ is the unique maximizer of $H$. To complete the proof we will show that the moment generating function of $N^{\frac{1}{2}}\left(\os - m_*\right)$ under $\p_{\beta_N,h_N,p}$ converges pointwise to the moment generating function of the Gaussian distribution with mean $-\bar{h}/H''(m_*)$ and variance $-1/H''(m_*)$. Towards this, fix $t \in \mathbb{R}$ and note that the moment generating function of $N^{\frac{1}{2}}\left(\os - m_*\right)$ at $t$ can be expressed as 
\begin{equation}\label{cltst1}
\mathbb{E}_{\beta_N,h_N,p}e^{tN^{\frac{1}{2}}\left(\os - m_*\right)} = e^{-tN^{\frac{1}{2}}m_*}\frac{Z_N\left(\beta_N,h_N+N^{-\frac{1}{2}}t,p\right)}{Z_N(\beta_N,h_N,p)}.
\end{equation}
Using Lemma \ref{ex} and the fact that $m_*(N) \rightarrow m_*$, the right side of \eqref{cltst1} simplifies to 
\begin{equation}\label{cltst2}
(1+o(1)) e^{-tN^{\frac{1}{2}}m_* +  N\left\{H_{\beta_N, h_N+N^{-\frac{1}{2}}t,p} \left(m_*\left(\beta_N, h_N+N^{-\frac{1}{2}}t,p\right) \right) - H_{\beta_N, h_N,p} \left(m_*\left(\beta_N, h_N,p\right) \right) \right\} } .
\end{equation}
Now, Lemma \ref{derh2} and a simple Taylor expansion gives us
\begin{align}\label{cltst3}
m_*\left(\beta_N, h_N+N^{-\frac{1}{2}}t,p\right) - m_*\left(\beta_N, h_N,p\right) &= N^{-\frac{1}{2}}t~ \frac{\partial}{\partial{\underline{h}}} m_*(\beta_N,\underline{h},p)\Big|_{\underline{h}=h_N} + O(N^{-1})\nonumber\\
&= -\frac{t}{N^{\frac{1}{2}} H_N''(m_*(\beta_N,h_N,p))} + O(N^{-1}).\nonumber\\ 
\end{align}
Using \eqref{cltst3} and a further Taylor expansion, we have
\begin{align}\label{cltst4}
&N\left\{H_{\beta_N, h_N,p} \left(m_*\left(\beta_N, h_N+N^{-\frac{1}{2}}t,p\right) \right) - H_{\beta_N, h_N,p} \left(m_*\left(\beta_N, h_N,p\right) \right) \right\}\nonumber\\
&\quad \quad \quad \quad = \frac{N}{2}\left\{m_*\left(\beta_N, h_N+N^{-\frac{1}{2}}t,p\right) - m_*\left(\beta_N, h_N,p\right)\right\}^2 H_N''\left(m_*\left(\beta_N, h_N,p\right)\right) + o(1)\nonumber\\ 
&\quad \quad \quad \quad = \frac{t^2}{2H_N''\left(m_*\left(\beta_N, h_N,p\right)\right)} + o(1) \nonumber\\ 
&\quad \quad \quad \quad = \frac{t^2}{2H''(m_*)} + o(1).
\end{align}
Next, we have by Lemma \ref{derh2} and a Taylor expansion,
\begin{align}\label{cltst5}
tN^{\frac{1}{2}} m_*\left(\beta_N, h_N+N^{-\frac{1}{2}}t,p\right) &= tN^{\frac{1}{2}} m_*(\beta_N,h,p) + t(t+\bar{h}) \frac{\partial}{\partial{\underline{h}}} m_*(\beta_N,\underline{h},p)\Big|_{\underline{h} = h} + o(1)\nonumber\\
&= tN^{\frac{1}{2}} m_* + t\bar{\beta}\frac{\partial}{\partial \underline{\beta}} m_*(\underline{\beta},h,p)\Big|_{\underline{\beta} = \beta}  - \frac{t(t+\bar{h})}{H''(m_*)} + o(1)\nonumber\\ 
&= tN^{\frac{1}{2}}m_* -\frac{t\bar{\beta}pm_*^{p-1}}{H''(m_*)} -\frac{t(t+\bar{h})}{H''(m_*)} + o(1).
\end{align}
Adding \eqref{cltst4} and \eqref{cltst5}, we have:
\begin{align}\label{cltst6}
& N\left\{H_{\beta_N, h_N+N^{-\frac{1}{2}}t,p} \left(m_*\left(\beta_N, h_N+N^{-\frac{1}{2}}t,p\right) \right) - H_{\beta_N, h_N,p} \left(m_*\left(\beta_N, h_N,p\right) \right) \right\}\nonumber\\ 
& \quad \quad \quad \quad = tN^{\frac{1}{2}}m_* -\frac{t(\bar{h}+\bar{\beta}pm_*^{p-1})}{H''(m_*)} -\frac{t^2}{2H''(m_*)} + o(1).
\end{align}
Using \eqref{cltst6}, the expression in \eqref{cltst2} becomes
\begin{equation}\label{cltst7}
\exp\left\{ -\frac{t(\bar{h}+\bar{\beta}pm_*^{p-1})}{H''(m_*)} -\frac{t^2}{2H''(m_*)}\right\} + o(1).
\end{equation}
The constant in expression \eqref{cltst7} is easily recognizable as the moment generating function of $N(-\frac{\bar{h}+\bar{\beta}pm_*^{p-1}}{H''(m_*)},~-\frac{1}{H''(m_*)})$ evaluated at $t$. This completes the proof of Theorem \ref{cltun}. \qed

\subsection{Proof Roadmap of Theorem \ref{cltun} when $(\beta, h)$ is $p$-special}\label{sec:pfcltun_II} 

When $(\beta, h)$ is $p$-special, we consider local perturbations of the parameters 
$$(\beta_N,h_N) :=\left(\beta + \bar{\beta} N^{-\frac{3}{4}}, h + \bar{h} N^{-\frac{3}{4}}\right)$$
as in the statement of Theorem \ref{cltun} (3).  
Note that in this case the function $H_{\beta, h, p}$ still has a unique maximizer $m_*=m_*(\beta, h, p)$, but $H_{\beta,h,p}''(m_*) = 0$. The proof strategy here follows essentially the same roadmap as in the $p$-regular case, with relevant modifications while taking Taylor expansions, since $H_{\beta,h,p}''(m_*) = 0$. As before, the first step is to prove the concentration of $\os$ within a vanishing neighborhood of $m_*$ (Lemma \ref{irrconch}). Here, the concentration window turns out to be a little more inflated, that is, its length is of order $N^{-\frac{1}{4} + \alpha}$, for $\alpha > 0$. Next, we approximate the partition function $\bar Z_N$, where, since the second derivative of $H$ is zero at the maximizer, we need to consider derivatives up to order four to accurately approximate $\bar Z_N$ (Lemma \ref{ex2}). The details of the proof are given in Section \ref{irregproof}.

\subsection{Proof Roadmap of Theorem \ref{cltun} when $(\beta, h)$ is $p$-critical}\label{nonpuniqueclt}

Throughout this section we assume that $(\beta, h) \in \Theta$ is $p$-critical. This means, by definition and Lemma \ref{derh11}, that the function $H = H_{\beta,h,p}$ has $K \in \{2, 3\}$ global maximizers, which we denote by  $m_1< \ldots<m_K$. It also follows from Lemma \ref{derh44}, that for sequences $\beta_N \rightarrow \beta$ and $h_N \rightarrow h$, the function $H_N := H_{\beta_N,h_N,p}$, for all large $N$, have local maximizers $m_1(N), \ldots,m_K(N)$ such that $m_k(N) \rightarrow m_k$, as $N \rightarrow \infty$, for all $1\leq k \leq K$. As before, $\bar{\p}$ and $\bar{Z}_N$ will denote $\p_{\beta_N,h_N,p} $ and $Z_N(\beta_N,h_N,p)$, respectively.

In presence of multiple global maximizers, the magnetization $\os$ will concentrate around the set of all global maximizers.  In fact, we can prove the following stronger result: Consider an open interval $A$ around a local maximizer $m$ such that $m$ is the unique global maximizer of $H$ over $A$. Then conditional on the event  $\os \in A$ (which is a rare event if $m$ is not a global maximizer), $\os$ concentrates around $m$. This is the first step in the proof of Theorem \ref{cltun} when $(\beta, h)$ is $p$-critical. To state the result precisely, assume that $m$ is a local maximizer of $H$ and let $m(N)$ be local maximizers of $H_N$ converging to $m$, which exist by Lemma \ref{derh44}. Define 
\begin{align}\label{eq:An}
A_{N,\alpha}(m(N)) = \left(m(N)-N^{-\frac{1}{2} + \alpha},m(N) +N^{-\frac{1}{2} + \alpha}\right). 
\end{align} 
The following lemma gives the conditional and, hence, the unconditional, concentration result of $\os$ around local maximizers.\footnote{The unconditional concentration derived in  \eqref{eq:multi_max_local_II} is not required in the proof of Theorem \ref{cltun}. Nevertheless, we include it for the sake of completeness.} The proof is given in Section \ref{sec:proofofnonun_I}. 

\begin{lem}\label{lem:multiple_max} Suppose $(\beta, h) \in \Theta$ is $p$-critical. Then for $\alpha \in \left(0,\frac{1}{6}\right]$ fixed and $A_{N,\alpha}(m(N))$ as defined in \eqref{eq:An}, 
	\begin{equation}\label{eq:multi_max_local_I}
	\bar{\mathbb P}\left(\os \in A_{N,\alpha}(m(N))^c\big|
	\os \in A\right) = \exp\left\{\frac{1}{3}N^{2\alpha} H''(m) \right\}O(N^{\frac{3}{2}}),
	\end{equation} 
	for any interval $A \subseteq [-1,1]$ such that $m \in \textrm{int}(A)$ and $H(m) > H(x)$, for all $x \in \mathrm{cl}(A) \setminus \{m\}$\footnote{For any set $A \subseteq \R$, $\textrm{int}(A)$ and $\mathrm{cl}(A)$ denote the topological interior and closure of $A$, respectively.} As a consequence, for $A_{N,\alpha,K} := \bigcup_{k=1}^K A_{N,\alpha}(m_k(N))$, 	
	\begin{equation}\label{eq:multi_max_local_II}
	\bar{\mathbb P}\left(\os \in A_{N,\alpha,K}^c\right) = \exp\left\{\frac{1}{3}N^{2\alpha} \max_{1\leq k\leq K} H''(m_k) \right\}O(N^{\frac{3}{2}}).
	\end{equation}
\end{lem} 

In order to derive a conditional CLT of $\os$ around the local maximizer $m$, given that $m$ is in $A$ (where $A$ is as in Lemma \ref{lem:multiple_max} above), we need precise estimates of the \textit{restricted partition functions} defined as
$$\bar{Z}_N\big|_A := \frac{1}{2^N}\sum_{\bs\in \sa_N: \os \in A} \exp\left\{N(\beta_N \overline{X}^p_N + h_N\os)  \right\}.$$ Note that $\bar{Z}_N\big|_A$ is the partition function of the conditional measure $\bar{\mathbb P}\left(\bs \in \cdot\big|\os \in A \right)$, in the sense that for any $\bt = (\tau_1, \tau_2, \ldots, \tau_N) \in \sa_N$ such that $\bar{\bt} \in A$, we have
$$\bar{\mathbb P}\left(\bs = \bt\big|\os \in A \right) = \frac{1}{2^{N} \bar{Z}_N\big|_{A}} \exp\left\{N(\beta_N \overline{X}^p_N + h_N\os) \right\}.$$ The following lemma gives an approximation of the restricted and, hence, the unrestricted partition functions. To this end, recall that $m(N)$ is a local maximizer of $H_N$ converging to $m$.

\begin{lem}\label{lm:condpart} Suppose  $(\beta, h) \in \Theta$ is $p$-critical. Then for $\alpha > 0$ and $N$ large enough, the restricted partition function can be expanded as
	\begin{equation}\label{eq:multi_max_I}
	\bar{Z}_N\big|_A =  \frac{e^{NH_N(m(N))}}{\sqrt{(m(N)^2-1)H_N''(m(N))}}\left(1+O\left(N^{-\frac{1}{2} +\alpha}\right)\right),
	\end{equation}
	where the set $A$ is as in Lemma \ref{lem:multiple_max}. This implies, for every $\alpha > 0$ and $N$ large enough, the (unrestricted) partition function can be expanded as 
	\begin{equation}\label{eq:multi_max_II}
	\bar{Z}_N =  \sum_{k=1}^K\frac{e^{NH_N(m_k(N))}}{\sqrt{(m_k(N)^2-1)H_N''(m_k(N))}}\left(1+O\left(N^{-\frac{1}{2} +\alpha}\right)\right). 
	\end{equation}
\end{lem}

The proof of this result is given in Section \ref{sec:proofofnonun_II}. We can now use the results above to complete the proof of Theorem \ref{cltun} (2). \\

\noindent\textbf{\textit{Completing the Proof of Theorem $\ref{cltun}$ when $(\beta, h)$ is $p$-critical}}: For each $\varepsilon > 0$ and $1\leq s \leq K$, define $B_{s, \varepsilon} = (m_s -\varepsilon, m_s+\varepsilon)$. Then for all $\varepsilon>0$ small enough, $H(m_s) > H(x)$, for all $x \in B_{s,\varepsilon}\setminus \{m_s\}$. Now, for each $1\leq s \leq K$, we have
\begin{equation}\label{smallprob}
\p_{\beta,h,p} (\os \in B_{s,\varepsilon}) = \frac{Z_N(\beta,h,p)\big|_{B_{s,\varepsilon}}}{Z_N(\beta,h,p)}.
\end{equation}
By Lemma \ref{lm:condpart} we have 
\begin{equation}\label{smallprob1}
Z_N(\beta,h,p)\big|_{B_{s,\varepsilon}} =  \frac{e^{N\sup_{x\in[-1,1]} H(x)}}{\sqrt{(m_s^2-1)H''(m_s)}}\left(1+o(1)\right)\quad\textrm{for all} ~1\leq s \leq K, 
\end{equation} 
and 
\begin{equation}\label{smallprob2}
Z_N(\beta,h,p) = e^{N\sup_{x\in[-1,1]} H(x)} \sum_{s=1}^K \frac{1}{\sqrt{(m_s^2-1)H''(m_s)}}\left(1+o(1)\right).
\end{equation}
The result in \eqref{eq:cltun_p1} now follows from \eqref{smallprob}, \eqref{smallprob1} and \eqref{smallprob2}.

Now, we proceed we prove \eqref{eq:cltun_II}. Hereafter, let $\beta_N := \beta+ N^{-\frac{1}{2}} \bar{\beta}$ and $h_N := h + N^{-\frac{1}{2}} \bar{h}$. A direct calculation reveals that 
\begin{equation}\label{condexp}
\e_{\beta_N,h_N,p}  \left[e^{tN^{\frac{1}{2}}(\os-m)}\Big| \os \in A  \right] =  e^{-tN^{\frac{1}{2}}m} \frac{Z_N(\beta_N,h_N + N^{-\frac{1}{2}} t,p)\big|_A}{Z_N(\beta_N,h_N,p)\big|_A}.
\end{equation}
Using Lemma \ref{lm:condpart}, the right side of \eqref{condexp} simplifies to
\begin{equation*}
(1+o(1)) e^{-tN^{\frac{1}{2}}m +  N\left\{H_{\beta_N, h_N+N^{-\frac{1}{2}}t,p} \left(m\left(\beta_N, h_N+N^{-\frac{1}{2}}t,p\right) \right) - H_{\beta_N, h_N,p} \left(m\left(\beta_N, h_N,p\right) \right)  \right\}},
\end{equation*}
where $m(\beta_N, h_N,p )$ and $m(\beta_N, h_N+N^{-\frac{1}{2}}t,p)$ are the local maximizers of the functions $H_{\beta_N, h_N,p}$ and $H_{\beta_N, h_N+N^{-\frac{1}{2}}t,p}$ respectively, converging to $m$. We can mimic the proof of Theorem \ref{cltun} verbatim from this point onward, to conclude that as $N \rightarrow \infty$,   
\begin{equation}\label{moment generating functionnonun}
\e_{\beta_N,h_N,p}  \left[e^{tN^{\frac{1}{2}}(\os-m)}\Big| \os \in A  \right] \rightarrow 	\exp\left\{ -\frac{t(\bar{h}+\bar{\beta}pm^{p-1})}{H''(m)} -\frac{t^2}{2H''(m)}\right\}.
\end{equation}
The result in \eqref{eq:cltun_II} now follows from \eqref{moment generating functionnonun}. \qed

\section{Asymptotic Distribution of the ML Estimates: Proof Overview} 
\label{sec:pfsketch_mle}

Here, we provide an overview on how the limiting distributions of $\hat{\beta}_N$ and $\hat{h}_N$, described in Section \ref{sec:mle_beta_h_I} above, can be obtained from the distribution of $\os$ presented in Theorem \ref{cltun}. To this end, recall the ML equations \eqref{eqmle} and \eqref{eqmleh}, and for notational convenience, we introduce the following definition , for $m \geq 1$: $$u_{N,m}(\beta, h, p) := \e_{\beta,h,p}  \overline{X}^m_N. $$

\begin{itemize}
	
	\item  The first step is to express the distribution functions of $\hat{\beta}_N$ and $\hat{h}_N$ in terms of the sample mean $\os$. This follows very easily from the ML equations \eqref{eqmle} and \eqref{eqmleh} and the monotonicity of the function $u_{N,m}$ (proved in Lemma \ref{increasing}). To this end, define $a_N = N^{\frac{1}{2}}$, if $(\beta,h)$ not $p$-special and $a_N=N^{\frac{3}{4}}$ if $(\beta,h)$ is $p$-special. Now, note that, fixing $t \in \R$,  
	\begin{align}\label{intui1}
	\left\{a_N(\hat{h}_N -h) \leq t\right\} = \left\{\hat{h}_N  \leq h + \frac{t}{a_N}\right\} 
	&= \left\{\os \leq \e_{\beta, h+ \frac{t}{a_N},p}  (\os)\right\},
	\end{align}
	by the monotonicity of the function $u_{N, 1}(\beta, \cdot, p)$ (using Lemma \ref{increasing}) and the ML equation \eqref{eqmleh}. Similarly, 
	\begin{align}\label{intui2}
	\left\{a_N(\hat{\beta}_N-\beta) \leq t\right\} = \left\{\hat{\beta}_N \leq \beta + \frac{t}{a_N}\right\} 
	&= \left\{\overline{X}^p_N \leq \e_{\beta+ \frac{t}{a_N},h,p} (\overline{X}^p_N)\right\}.
	\end{align}
	by the monotonicity of the function $u_{N, p}(\cdot, h, p)$ (using Lemma \ref{increasing}) and the ML equation \eqref{eqmle}. 
	
	\item The next step is to write the event in \eqref{intui1}  as $$\left\{\frac{N}{a_N}(\os-c) \leq \e_{\beta, h+ \frac{t}{a_N},p}  \left[ \frac{N}{a_N} (\os-c)\right]\right\},$$ 
	for some appropriately chosen centering $c$, and and similarly, for the event \eqref{intui2}. Now, if the point $(\beta,h)$ is $p$-regular or $p$-special, the centering $c$ will be the unique global maximizer of $H_{\beta,h,p}$, around which $\os$ concentrates. However, if $(\beta,h)$ is $p$-critical, then the situation is more tricky. In that case, one needs to look at the sign of $t$, and choose the centering $c$ to be that global maximizer of $H_{\beta,h,p}$ around which $\os$ concentrates, under the measure $\p_{\beta, h+ t/a_N, p}$ (for $\hat{h}_N$), and the measure $\p_{\beta+t/a_N, h,p}$ (for $\hat{\beta}_N$). The proofs are now completed by computing the asymptotic probabilities of the events in the centered and scaled forms, written above, by applying the results in Section \ref{cltmagnetization}.  
	
\end{itemize}

The details of the proof are given in Section \ref{sec:pf_beta_h}. The proofs in the $p$-regular and $p$-special cases (which includes Theorems \ref{cltintr3_1}, \ref{cltintr3_2}, \ref{thmmle1}, and \ref{thmmle2}) are given in Section \ref{sec:pf_beta_h_1}. The asymptotic distributions of 
$\hat{h}_N$ in the $p$-critical case (Theorem \ref{cltintr3_III}) are given in Section \ref{sec:pf_beta_h_2}. The results for $\hat{\beta}_N$ in the $p$-critical case (Theorem \ref{thmmle_III}) are proved in Section \ref{sec:pf_beta_III}.

\section{Constructing Confidence Intervals}
\label{sec:applications}

In this section, we discuss how the limiting distributions for the ML estimates $\hat{\beta}_N$ and $\hat{h}_N$ obtained above can be used to construct asymptotically valid confidence intervals for the respective parameters. One complication towards using the above results directly is that the limiting distributions $\hat{\beta}_N$ and $\hat{h}_N $ depend on the actual position of the true parameter $(\beta,h) \in \Theta$.  However, if there were an oracle that told us that the unknown parameter $(\beta,h)$ is $p$-regular, then using the results in \eqref{eq:h_1} and \eqref{eq:bmle_m_1} we would be able to easily construct  confidence intervals for the parameters with asymptotic coverage probability $1-\alpha$, as follows:

\begin{itemize}

	\item {\it Confidence interval for $h$ at the regular points}: Suppose $\beta \geq 0$ is known and $(\beta, h)$ is $p$-regular. Denote the unique maximizer of the function $H$ by $m_*= m_*(\beta, h, p)$. Note that, by  \eqref{eq:meanclt_I}, the sample mean $\os \pto m_*$, under $\p_{\beta, h, p}$. Therefore, by \eqref{eq:h_1}, 
	\begin{equation}\label{eq:cih1}
	I_{\mathrm{reg}} := \left(\hat{h}_N - \sqrt{\frac{-H''(\os)}{N}} z_{1-\frac{\alpha}{2}},~\hat{h}_N + \sqrt{\frac{-H''(\os)}{N}} z_{1-\frac{\alpha}{2}}\right), 
	\end{equation} 
	is an interval which contains $h$ with asymptotic coverage probability $1-\alpha$, whenever $(\beta, h)$- is $p$-regular.\footnote{Note that $z_\alpha$ is the $\alpha$-th quantile of the standard normal distribution, that is, $\mathbb{P}(N(0, 1) \leq z_\alpha) = \alpha$.} More precisely, $\p_{\beta, h, p}(h \in I_{\mathrm{reg}}) \rightarrow 1-\alpha$, for  $(\beta, h) \in \Theta$ which is regular.

	\item {\it Confidence interval for $\beta$ at the regular points}: Suppose $h \ne 0$ is known and $(\beta, h)$ is $p$-regular. As before, denote the unique maximizer of the function $H$ by $m_*= m_*(\beta, h, p)$. Therefore, by \eqref{eq:bmle_m_1}, 
	\begin{equation}\label{eq:cibeta1}
	J_{\mathrm{reg}}: = \left(\hat{\beta}_N - \frac{\os^{1-p}}{p}\sqrt{\frac{-H''(\os)}{N}} z_{1-\frac{\alpha}{2}},~\hat{\beta}_N + \frac{\os^{1-p}}{p}\sqrt{\frac{-H''(\os)}{N}} z_{1-\frac{\alpha}{2}}\right),
	\end{equation}   
	is an interval which contains $\beta$ with asymptotic coverage probability $1-\alpha$, whenever $(\beta, h)$ is $p$-regular. Note that the assumption $h \ne 0$ is essential, since \ref{eq:bmle_m_2} shows that the ML estimate $\hat{\beta}_N$ may be inconsistent otherwise. 
\end{itemize}
Note that the length of $I_{\mathrm{reg}}$ does not depend upon the true value of $\beta$, and 
the length of $J_{\mathrm{reg}}$ does not depend on the true $h$. 

Now, we discuss how the intervals in \eqref{eq:cih1} and \eqref{eq:cibeta1} can be modified so that they are valid for all parameter points. To this end, let $\overline{\cp}$ denote the closure of the curve $\cp$ with respect to the Euclidean topology on $\Theta$, that is, the union of $\cp$ with the $p$-special point(s) (recall \eqref{eq:beta_h_special}).

\begin{itemize}
	
	\item {\it Confidence interval for $h$ for all points}: Suppose $\beta \geq 0$ is known, and define $S_p(\beta) :=  \{\underline{h} \in \mathbb{R}: (\beta,\underline{h}) \in  \overline{\cp}\}$. Note that $S_p(\beta)$ is either empty, a singleton or a doubleton (recall Figures \ref{figure:ordering1} and  \ref{figure:ordering2}), and is free of the unknown parameter $h$. Then 
	\begin{align}\label{eq:cih2}
	I := I_{\mathrm{reg}} \bigcup S_p(\beta), 
	\end{align}  
	is an interval with the same length (Lebesgue measure) as the regular interval $I_{\mathrm{reg}}$, and contains $h$ with asymptotic probability at least $1-\alpha$, {\it for all $(\beta, h) \in \Theta$}.  This is because the asymptotic coverage probability is guaranteed to be $1-\alpha$ when $(\beta, h)$ is $p$-regular by the discussion following \eqref{eq:cih1} above. On the other hand, 
	if $(\beta, h)$ is not $p$-regular, by definition $h \in S_p(\beta)$, and hence, by \eqref{eq:cih2}, $\p_{\beta,h, p}(I \ni h) = 1$.

	\item {\it Confidence interval for $\beta$ for all points with $h \ne 0$}: 
	Fix $h \neq 0$, and define $T_p(\beta) :=  \{\underline{\beta} \geq 0: (\underline{\beta},h) \in  \overline{\cp}\}$. Note that $T_p(h)$ is either empty or a singleton, and is free of the unknown parameter $\beta$. Then, as above,  
	\begin{align}\label{eq:cibeta2}
	J:= J_{\mathrm{reg}} \bigcup T_p(h), 
	\end{align}   
	is an interval with the same length (Lebesgue measure) as the regular interval $J_{\mathrm{reg}}$, and contains $\beta$ with asymptotic probability at least $1-\alpha$, {\it for all $(\beta, h) \in \Theta$.} 
	
\end{itemize}

\begin{figure}
	\centering
	\begin{minipage}[c]{1.0\textwidth}
		\centering
		\includegraphics[width=5.25in,height=2.75in]
		{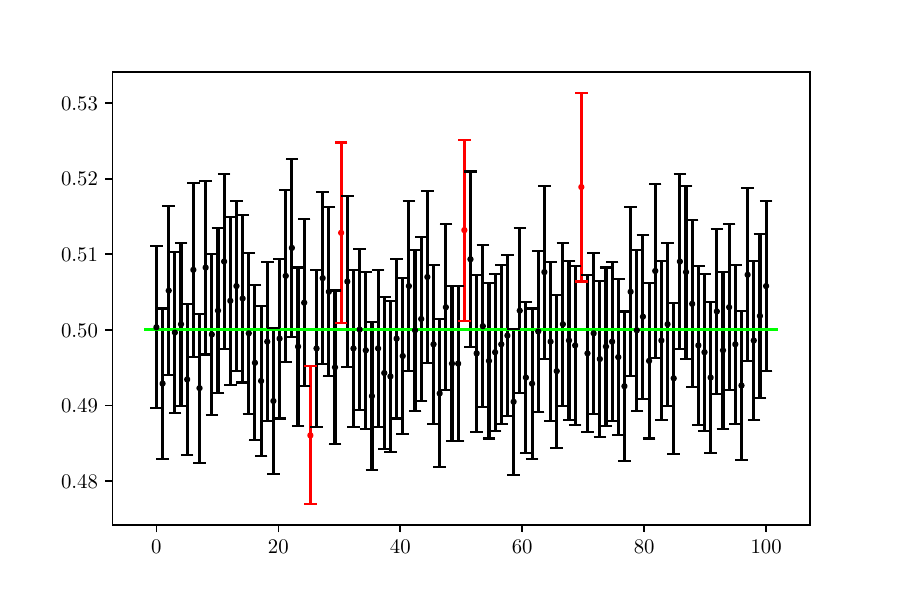}\\
		%\small{(a)}
	\end{minipage}
	\caption{\small{100 realizations of the 95\% confidence interval for $\beta$ at the $3$-regular point $(\beta, h)=(0.5, 0.2)$, with $N=10,000$. The intervals not containing the true parameter $\beta=0.5$ are shown in red. }}
	\label{fig:sigma_interval_beta}
\end{figure}

Figure \ref{fig:sigma_interval_beta} shows 100 realizations of the 95\% confidence interval for $\beta$ at the $3$-regular point $(\beta, h)=(0.5, 0.2)$, with $N=10,000$. The green horizontal line represents the true parameter $\beta=0.5$ and the intervals not containing the true parameter are shown in red. 

%%%%%%%%%%%%%%%%%%%%%%%%%%%%%%%%%%%%%%%%%%%%%%%%%%%%%%%%%%%%%%%%%%%%%%%%%%%%%%%%
%%%%%%%%%%%%%%%%%%%%%%%%%%%%%%%%%%%%%%%%%%%%%%%%%%%%%%%%%%%%%%%%%%%%%%%%%%%%%%%%%%%%%%%%%%%%%%%%%%%%%%%%%%

%\fi

%%
%
%
%%
%
%
%
%
%
%
%
%
%
%
%
%
%
%
%
%
%
%
%
%
%
%
%
%
%
%
%
%
%
%

\section{Proof of Theorem \ref{cltun} when $(\beta, h)$ is $p$-special}\label{irregproof}

Throughout this section, as usual, we will denote $H_{\beta,h,p}$ by $H$, $H_{\beta_N, h_N,p}$ by $H_N$, the unique global maximizer of $H_{\beta_N, h_N,p}$ (for large $N$) by $m_*(N)$, $\p_{\beta_N,h_N,p} $ by $\bar{\p}$, $Z_N(\beta_N,h_N,p)$ by $\bar{Z}_N$ and $F_N(\beta_N,h_N,p)$ by $\bar{F}_N$. As outlined in Section \ref{sec:pfcltun_II}, the first step in the proof of Theorem \ref{cltun} when $(\beta, h)$ is $p$-special, is to show the concentration of $\os$ within a vanishing neighborhood of $m_*(N)$. In the $p$-special case, this is more delicate, because it requires Taylor expansions up to the fourth order term. Here, the concentration window turns out to be a bit more inflated as well, and it is given by:
\begin{equation}\label{concwind}
\ma_{N,\alpha} := (m_*(N)-N^{-\frac{1}{4} + \alpha},m_*(N) +N^{-\frac{1}{4} + \alpha} ). 
\end{equation}

\begin{lem}\label{irrconch} Suppose  $(\beta, h) \in \Theta$ is $p$-special. Fix $\alpha \in \left(0,\frac{1}{20}\right]$ and let $\ma_{N,\alpha}$ be as in \eqref{concwind}. Then,
	\begin{equation*}
	\bar{\mathbb P}\left(\os \in \ma_{N,\alpha}^c\right) = \exp\left\{\frac{1}{24}N^{4\alpha} H^{(4)}(m_*)(1+o(1)) \right\}O(N^{\frac{3}{2}}).
	\end{equation*}	
\end{lem}

\begin{proof}
	It follows from the proof of Lemma \ref{conc} and the fact $H_N''(m_*(N)) \leq 0$, that
	\begin{align}\label{simtoconc}
	\bar{\mathbb P} & (\os \in \ma_{N,\alpha}^c)   \nonumber \\ 
	&\leq  \exp\left\{N\left( H_N\left(m_*(N) \pm N^{-\frac{1}{4}+\alpha}\right) - H_N\left(m_*(N)\right)\right)\right\} O(N^{\frac{3}{2}})\nonumber\\
	&\leq  \exp\left\{\frac{1}{6} N^{\frac{1}{4}+3\alpha}H_N^{^{(3)}} (m_*(N)) + \frac{1}{24}N^{4\alpha} H_N^{(4)}(m_*(N)) + O\left(N^{-\frac{1}{4} + 5\alpha}\right)\right\} O(N^{\frac{3}{2}})\nonumber\\
	\end{align}	
	Now, it follows from Lemma \ref{mnminusm}, that $|H_N^{(3)}(m_*(N))| =  O(N^{-1/4})$. Hence, $$N^{(1/4) + 3\alpha} H_N^{(3)}(m_*(N)) + N^{4\alpha}H_N^{(4)}(m_*(N)) = N^{4\alpha}H^{(4)}(m_*)(1+o(1)),$$ and Lemma \ref{irrconch} follows from \eqref{simtoconc}. 
\end{proof}

The next step in the proof of Theorem \ref{cltun} when $(\beta, h)$ is $p$-special is the approximation of the partition function.

\begin{lem}\label{ex2} Suppose  $(\beta, h) \in \Theta$ is $p$-special, and let $(\beta_N,h_N) = (\beta+ N^{-\frac{3}{4}}\bar{\beta},~h + N^{-\frac{3}{4}}\bar{h})$. Then for $N$ large enough, the partition function $\bar{Z}_N$ can be expanded as
	\begin{equation*}%\label{partex2}
	\bar{Z}_N =  \frac{N^\frac{1}{4}e^{N H_N(m_*(N))}}{\sqrt{2\pi(1-m_*(N)^2)}}\int_{-\infty}^\infty e^{\eta_{\bar{\beta},\bar{h},p}(y)} \mathrm d y \left(1+o(1)\right),
	\end{equation*}
	where $\eta_{\bar{\beta},\bar{h},p}(y)= a y^2 +  b y^3 + c y^4$, with 	
	\begin{equation*}
	a:=\frac{(6(\bar{\beta}pm_*^{p-1} + \bar{h}))^{\frac{2}{3}}\left(H^{(4)}(m_*)\right)^{\frac{1}{3}}}{4}, ~ b:= - \frac{(6(\bar{\beta}pm_*^{p-1} + \bar{h}))^{\frac{1}{3}}\left(H^{(4)}(m_*)\right)^{\frac{2}{3}}}{6},
	\end{equation*}
	\begin{equation*}
 ~c:= \frac{H^{(4)}(m_*)}{24}.
	\end{equation*}

\end{lem}

\begin{proof}
	Once again, as in the proof of Lemma \ref{ex}, it follows from Lemma \ref{irrconch}, that for $\alpha \in \left(0,\frac{1}{20}\right]$,
	\begin{equation}\label{secstep22}
	\bar{Z}_N = \left(1+ O\left(e^{-N^\alpha}\right) \right) \sum_{m \in \mathcal{M}_N\bigcap \mathcal{A}_{N,\alpha}} \zeta(m),
	\end{equation}
	where $\zeta:[-1,1] \rightarrow \mathbb{R}$ is defined in \eqref{xidef} and $\mathcal{A}_{N,\alpha}$ is defined in \eqref{concwind}.
	It also follows from Lemma \ref{Riemann} and Lemma \ref{xiprimeirreg}, exactly as in the proof of Lemma \ref{ex}, that
	\begin{equation}\label{rax2}
	\left|\int_{\ma_{N,\alpha}} \zeta(x)\mathrm d x - \frac{2}{N} \sum_{m \in \mathcal{M}_N\bigcap \ma_{N,\alpha}} \zeta(m)\right| = O\left(N^{-1 + 4\alpha}\right)\zeta(m_*(N)).
	\end{equation}
	Hence, we have from \eqref{rax2}, Lemma \ref{xiact}, Lemma \ref{Laplaceirreg} and Lemma \ref{mnminusm},
	\begin{align}\label{fin3} 
	& \sum_{m \in \mathcal{M}_N\bigcap \ma_{N,\alpha}} \zeta(m)\nonumber\\
	& =  \frac{N}{2} \int_{\ma_{N,\alpha}} \zeta(x)\mathrm d x + O(N^{4\alpha}) \zeta(m_*(N))\nonumber\\
	& =  \frac{N^{\frac{1}{2}}}{2}\left(1+O(N^{-1})\right)\int_{\ma_{N,\alpha}} e^{NH_N(x)}\sqrt{\frac{2}{\pi (1-x^2)}}\mathrm d x + O(N^{4\alpha}) \zeta(m_*(N))\nonumber\\
	& =   \frac{N^\frac{1}{4}}{\sqrt{2\pi(1-m_*(N)^2)}} e^{N H_N(m_*(N))}  \int_{-N^\alpha}^{N^\alpha}e^{\eta_{\bar{\beta},\bar{h},p}(y)}\mathrm d y \left(1+O\left(N^{-\frac{1}{4}+5\alpha}\right)\right) \nonumber\\
	& +O(N^{4\alpha}) \zeta(m_*(N))\nonumber\\
	& =  \frac{N^\frac{1}{4}e^{N H_N(m_*(N))}}{\sqrt{2\pi (1-m_*(N)^2)}}\int_{-\infty}^\infty e^{\eta_{\bar{\beta},\bar{h},p}(y)}\mathrm d y  (1+o(1))\left(1+O\left(N^{-\frac{1}{4}+5\alpha}\right)\right)\nonumber\\
	&+  \sqrt{\frac{2}{\pi N (1-m_*(N)^2)}} e^{NH_N(m_*(N))}\left(1+O(N^{-1})\right)O(N^{4\alpha})\nonumber\\
	& =  \frac{N^\frac{1}{4}e^{N H_N(m_*(N))}}{\sqrt{2\pi(1-m_*(N)^2)}}\int_{-\infty}^\infty e^{\eta_{\bar{\beta},\bar{h},p}(y)}\mathrm d y \left(1+o(1)\right).
	\end{align}
	Combining \eqref{secstep22} and \eqref{fin3}, we have:
	\begin{align}\label{fin4}
	\bar{Z}_N & =  \left(1+ O\left(e^{-N^\alpha}\right) \right)\left(1+o(1)\right)\frac{N^\frac{1}{4}e^{N H_N(m_*(N))}}{\sqrt{2\pi(1-m_*(N)^2)}}\int_{-\infty}^\infty e^{\eta_{\bar{\beta},\bar{h},p}(y)}\mathrm d y \nonumber\\& =  \left(1+o(1)\right)\frac{N^\frac{1}{4}e^{N H_N(m_*(N))}}{\sqrt{2\pi(1-m_*(N)^2)}}\int_{-\infty}^\infty e^{\eta_{\bar{\beta},\bar{h},p}(y)}\mathrm d y .
	\end{align}
	This completes the proof of Lemma \ref{ex2}. 
\end{proof}

\noindent\textbf{\textit{Completing the Proof of}} \eqref{eq:cltun_III}: As before, we start by computing the limiting moment generating function of $$N^{\frac{1}{4}}\left(\os - m_*(\beta,h,p)\right),$$ in the following lemma. 

\begin{lem}\label{cltun133pr}
	For every $p$-special point $(\beta,h)\in \Theta$ and $\bar{\beta}, \bar{h} \in \mathbb{R}$, if $\bs \sim \mathbb{P}_{\beta+ N^{-\frac{3}{4}}\bar{\beta},~h + N^{-\frac{3}{4}}\bar{h},~p} $, then
	\begin{align}\label{complexp}
	\lim_{N \rightarrow \infty} & \e_{\beta+ N^{-\frac{3}{4}}\bar{\beta},~h + N^{-\frac{3}{4}}\bar{h},~p}  \left[e^{t N^{\frac{1}{4}}\left(\os - m_*(\beta,h,p)\right)} \right] \nonumber \\ 
	& \quad \quad = C_p(\bar \beta, \bar h, t) \exp\Big\{- t R_p(\bar \beta, \bar h, t)  + \eta_{\bar{\beta},\bar{h},p}\left( R_p(\bar \beta, \bar h, 0) -  R_p(\bar \beta, \bar h, t) \right) \Big\},
	\end{align}
	where $\eta_{\bar{\beta},\bar{h},p}$ is defined in the statement of Lemma \ref{ex2}, 
	$$C_p(\bar \beta, \bar h, t):=\frac{\int_{-\infty}^\infty e^{ \eta_{\bar{\beta},\bar{h}+t,p}(y) } \mathrm d y }{\int_{-\infty}^\infty e^{\eta_{\bar{\beta},\bar{h},p}(y) } \mathrm d y },$$
	and $R_p(\bar \beta, \bar h, t):=\left(\frac{6(\bar{\beta}pm_*^{p-1} +\bar{h}+t)}{H^{(4)}(m_*)}\right)^\frac{1}{3}$.
\end{lem}
\begin{proof} Once again, throughout this proof, we will denote $m_*(\beta,h,p)$ by $m_*$, $\beta + N^{-\frac{3}{4}} \bar{\beta}$ by $\beta_N$, and $h + N^{-\frac{3}{4}} \bar{h}$ by $h_N$. Fix $t \in \mathbb{R}$ and note that the moment generating function of $N^{\frac{1}{4}}\left(\os - m_*\right)$ at $t$ can be expressed as 
	\begin{equation}\label{cltstr1}
	\mathbb{E}_{\beta_N,h_N,p}e^{tN^{\frac{1}{4}}\left(\os - m_*\right)} = e^{-tN^{\frac{1}{4}}m_*}\frac{Z_N\left(\beta_N,h_N+N^{-\frac{3}{4}}t,p\right)}{Z_N(\beta_N,h_N,p)}.
	\end{equation}
	Using Lemma \ref{ex2} and the facts that $m_*(\beta_N,h_N,p) \rightarrow m_*$ and $m_*(\beta_N,h_N+N^{-\frac{3}{4}}t,p) \rightarrow m_*$, the right side of \eqref{cltstr1} simplifies to 
	\begin{align}\label{cltstr2}
	C_p(\bar \beta, \bar h, t)e^{-tN^{\frac{1}{4}}m_* + N\left\{H_{\beta_N, h_N+N^{-\frac{3}{4}}t,p} \left(m_*\left(\beta_N, h_N+N^{-\frac{3}{4}}t,p\right) \right) - H_{\beta_N, h_N,p} \left(m_*\left(\beta_N, h_N,p\right) \right) \right\} } (1+o(1)).
	\end{align}
	By Lemma \ref{mnminusm}, we have:
	\begin{equation}\label{cltstr3}
	N^\frac{1}{4} \left(m_*\left(\beta_N, h_N+N^{-\frac{3}{4}}t,p\right) - m_*\right) =  -R_p(\bar{\beta},\bar{h},t) + o(1).
	\end{equation}
	By a further Taylor expansion and using \eqref{mnminusm}, we have (denoting $H_N=H_{\beta_N, h_N,p} $), 
	\begin{align}\label{cltstr4}
	N\left\{H_{N} \left(m_*\left(\beta_N, h_N+N^{-\frac{3}{4}}t,p\right) \right) - H_{N} \left(m_*\left(\beta_N, h_N,p\right) \right) \right\}  = T_1 + T_2 + T_3 + T_4,\nonumber\\
	\end{align} 
	where 
	\begin{align*}
	T_1&:=\tfrac{N}{2} \left\{ m_*\left(\beta_N, h_N+N^{-\frac{3}{4}}t,p\right) - m_*\left(\beta_N, h_N,p\right) \right \} ^2 H_{\beta_N,h_N,p}''\left(m_*\left(\beta_N, h_N,p\right)\right)\nonumber\\
	& = \tfrac{1}{2}\left\{ R_p(\bar \beta, \bar h, 0) - R_p(\bar \beta, \bar h, t) \right\} ^2\cdot \frac{1}{2}(6(\bar{\beta}pm_*^{p-1} + \bar{h}))^\frac{2}{3}\left(H^{(4)}(m_*)\right)^\frac{1}{3} + o(1), 
	\end{align*} 
	\begin{align*} 
	T_2 & := \tfrac{N}{6}\left\{ m_*\left(\beta_N, h_N+N^{-\frac{3}{4}}t,p\right) - m_*\left(\beta_N, h_N,p\right) \right \} ^3 H_{\beta_N,h_N,p}^{(3)}\left(m_*\left(\beta_N, h_N,p\right)\right)\nonumber\\
	& = - \tfrac{1}{6} \left\{ R_p(\bar \beta, \bar h, 0) - R_p(\bar \beta, \bar h, t) \right\} ^3(6(\bar{\beta}pm_*^{p-1} + \bar{h}))^\frac{1}{3}\left(H^{(4)}(m_*)\right)^\frac{2}{3} + o(1) , 
	\end{align*} 
	\begin{align*} 
	T_3 & := \tfrac{N}{24}\left\{ m_*\left(\beta_N, h_N+N^{-\frac{3}{4}}t,p\right) - m_*\left(\beta_N, h_N,p\right)\right\} ^4 H_{\beta_N,h_N,p}^{(4)}\left(m_*\left(\beta_N, h_N,p\right)\right)\nonumber\\ 
	& = \tfrac{1}{24} \left\{ R_p(\bar \beta, \bar h, 0) - R_p(\bar{\beta},\bar{h},t) \right\} ^4 H^{(4)}(m_*) + o(1), 
	\end{align*}
	and  
	$$T_4  := O(N\{ m_*(\beta_N, h_N+N^{-\frac{3}{4}}t,p) - m_*(\beta_N, h_N,p)\} ^5) = o(1).$$
	
	Now, using both \eqref{cltstr3} and \eqref{cltstr4}, we have 
	\begin{align}%\label{cltstr5}
	& N\left[H_{\beta_N, h_N+N^{-\frac{3}{4}}t,p} \left(m_*\left(\beta_N, h_N+N^{-\frac{3}{4}}t,p\right) \right) - H_{\beta_N, h_N,p} \left(m_*\left(\beta_N, h_N,p\right) \right) \right]\nonumber\\ 
	&= tN^{\frac{1}{4}}m_* - tR_p(\bar \beta, \bar h, t)  + \eta_{\bar{\beta},\bar{h},p}\left(R_p(\bar \beta, \bar h, 0) - R_p(\bar \beta, \bar h, t)\right) + o(1). \nonumber 
	\end{align}
	Using the above with  \eqref{cltstr1} and \eqref{cltstr2} Lemma \ref{cltun133pr} follows. 
\end{proof}

Although \eqref{complexp} is not readily recognizable as the moment generating function of any probability distribution, we will show below that it is indeed the moment generating function of the distribution $F_{\bar{\beta},\bar{h}}$ defined in \eqref{eq:beta_h_distribution}.

\begin{lem}\label{wint}
	Let $F_{\bar{\beta},\bar{h}}$ be the distribution defined in \eqref{eq:beta_h_distribution}. Then, 
	\begin{align} 
	\label{cltstr6pr}
	&\int e^{tx} \mathrm d F_{\bar{\beta},\bar{h}}(x) =  C_p(\bar \beta, \bar h, t) \exp\Big\{- t R_p(\bar \beta, \bar h, t)  + \eta_{\bar{\beta},\bar{h},p}\left( R_p(\bar \beta, \bar h, 0) -  R_p(\bar \beta, \bar h, t) \right) \Big\},\nonumber\\ 
	\end{align}
	with notations as in Lemma \ref{cltun133pr}.  
\end{lem}
\begin{proof} 
	Let us denote the right side of \eqref{cltstr6pr} by $M(t)$. Define
	\begin{align*}
	\Delta(t,y) :=  - t R_p(\bar \beta, \bar h, t)  + \eta_{\bar{\beta},\bar{h},p}\left( R_p(\bar \beta, \bar h, 0) -  R_p(\bar \beta, \bar h, t) \right)  + \eta_{\bar{\beta},\bar{h}+t,p}(y),
	\end{align*}
	Note that
	\begin{equation}\label{num}
	M(t) =\frac{\int_{-\infty}^\infty e^{\Delta(t,y)}\enskip \!\! \mathrm d y }{\int_{-\infty}^\infty e^{\eta_{\bar{\beta},\bar{h},p}(y) } \mathrm d y }.
	\end{equation}
	Using the change of variables $u = y- R_p(\bar \beta, \bar h, t)$ and a straightforward algebra, we have 
	\begin{equation}\label{integrand}
	\Delta(t,y) = \frac{H^{(4)}(m_*)}{24} u^4 + (\bar{\beta} p m_*^{p-1}+\bar{h})u + tu + \frac{(6(\bar{\beta} p m_*^{p-1}+\bar{h}))^{\frac{4}{3}}}{8 (H^{(4)}(m_*))^{\frac{1}{3}}} 
	\end{equation}
	and 
	\begin{equation}\label{integrand2}
	\eta_{\bar{\beta},\bar{h},p}(y)=\frac{H^{(4)}(m_*)}{24}u^4 + (\bar{\beta} p m_*^{p-1}+\bar{h}) u + \frac{\left(6(\bar{\beta}pm_*^{p-1}+\bar{h})\right)^{\frac{4}{3}}}{8(H^{(4)}(m_*))^{\frac{1}{3}}}.
	\end{equation}
	
	Lemma \ref{wint} now follows on substituting \eqref{integrand} and \eqref{integrand2} in \eqref{num}.
\end{proof}

The proof of \eqref{eq:cltun_III} now follows from Lemmas \ref{cltun133pr} and \ref{wint}. This completes the proof of Theorem \ref{cltun} when $(\beta, h)$ is $p$-special. 

\section{Missing Details in the Proof of Theorem \ref{cltun} when $(\beta, h)$ is $p$-critical}\label{sec:proofofnonun}

In this section we prove Lemma \ref{lem:multiple_max} and Lemma \ref{lm:condpart}, in Section \ref{sec:proofofnonun_I} and Section \ref{sec:proofofnonun_II}, respectively. These lemmas where used in Section \ref{nonpuniqueclt} in the proof of Theorem \ref{cltun} when $(\beta, h)$ is $p$-critical. 

\subsection{Proof of Lemma \ref{lem:multiple_max}}\label{sec:proofofnonun_I}

It follows from Lemma \ref{derh44}, that for all $N$ sufficiently large, $H_N(m(N)) > H_N(x)$ for all $x \in \mathrm{cl}(A) \setminus \{m(N)\}$, whence we can apply Lemma \ref{mest} to conclude that $$\sup_{x\in A\setminus A_{N,\alpha}(m(N))} H_N(x) = H_N\left(m(N) \pm N^{-\frac{1}{2}+\alpha}\right),$$ for all large $N$ such that $A_{N,\alpha}(m(N))\subset A$, as well. Following the proof of Lemma \ref{conc}, we have for all large N ,
\begin{align}\label{coinnonun1}
&\bar{\mathbb P}\left(\os \in A_{N,\alpha}(m(N))^c\big|
\os \in A\right)\nonumber\\ & \leq   \exp\left\{N\left(\sup_{x\in A\setminus A_{N,\alpha}(m(N))} H_N(x) - \sup_{x\in A} H_N(x)\right)\right\} O(N^{\frac{3}{2}})\nonumber\\& =  \exp\left\{N\left( H_N\left(m(N) \pm N^{-\frac{1}{2}+\alpha}\right) - H_N\left(m(N)\right)\right)\right\} O(N^{\frac{3}{2}})\nonumber\\& \leq  \exp\left\{\frac{N}{3}\left(N^{-1+2\alpha} H''(m) +  O\left(N^{-\frac{3}{2} + 3\alpha}\right) \right)\right\} O(N^{\frac{3}{2}}).
\end{align} 
The result \eqref{eq:multi_max_local_I} now follows from \eqref{coinnonun1}.

Next, we proceed to prove \eqref{eq:multi_max_local_II}. Let $A_1 := [-1,(m_1+m_2)/2)$, $A_K := [(m_{K-1}+m_K)/2,1]$ and for $1< k< K$, $A_k := [(m_{k-1} + m_k)/2,(m_k+m_{k+1})/2)$. Then, $A_1,A_2,\ldots, A_K$ are disjoint intervals uniting to $[-1,1]$, $m_k \in \textrm{int}(A_k)$, and $H(m_k) > H(x)$ for all $x\in \mathrm{cl}(A_k) \setminus \{m_k\}$ and all $1\leq k\leq K$. Hence, by Lemma \ref{lem:multiple_max}, $$\bar{\mathbb P}\left(\os \in A_{N,\alpha}(m_k(N))^c\big|
\os \in A_k\right) = \exp\left\{\frac{1}{3}N^{2\alpha} H''(m_k) \right\}O(N^{\frac{3}{2}})\quad\textrm{for all}~ 1\leq k\leq K.$$ Since $A_{N,\alpha}(m_k(N)) \subset A_k$ for all $1\leq k \leq K$, for all large $N$, we have $A_{N,\alpha}(m_k(N))^c \bigcap A_k = A_{N,\alpha,K}^c \bigcap A_k$ for all $1\leq k \leq K$, for all large $N$ (recall the definition of $ A_{N,\alpha,K}$ from the statement of Lemma \ref{lem:multiple_max}). Hence, $$\bar{\mathbb P}\left(\os \in A_{N,\alpha}(m_k(N))^c\big|
\os \in A_k\right) = \bar{\mathbb P}\left(\os \in A_{N,\alpha,K}^c\big|
\os \in A_k\right)\quad\textrm{for all}~1\le k \le K$$ for all large $N$. Hence, for all large $N$, we have
\begin{equation}\label{condtouncond}
\bar{\mathbb P}\left(\os \in A_{N,\alpha,K}^c\big|
\os \in A_k\right) = \exp\left\{\frac{1}{3}N^{2\alpha} H''(m_k) \right\}O(N^{\frac{3}{2}})\quad\textrm{for all}~ 1\leq k\leq K.
\end{equation}
It follows from \eqref{condtouncond} that for all large $N$,
\begin{align}\label{lstp1}
\bar{\p}(\os \in A_{N,\alpha,K}^c) & =  \sum_{k=1}^K \bar{\mathbb P}\left(\os \in A_{N,\alpha,K}^c\big|
\os \in A_k\right) \bar{\p} (\os \in A_k)\nonumber\\& \leq  \exp\left\{\frac{1}{3}N^{2\alpha} \max_{1\leq k\leq K}H''(m_k) \right\}O(N^{\frac{3}{2}})\sum_{k=1}^K \bar{\p} (\os \in A_k)\nonumber\\& =  \exp\left\{\frac{1}{3}N^{2\alpha} \max_{1\leq k\leq K}H''(m_k) \right\}O(N^{\frac{3}{2}}).
\end{align}
The result in \eqref{eq:multi_max_local_II} now follows from \eqref{lstp1}, completing the proof of Lemma \ref{lem:multiple_max}. \qed

%\subsection{Proof of Lemma \ref{lm:condpart}, Corollary \ref{exnonun} and Corollary \ref{limitpmf}}

\subsection{Proof of Lemma \ref{lm:condpart}}
\label{sec:proofofnonun_II}

The arguments below are meant for all sufficiently large $N$. Without loss of generality, let $\alpha \in \left(0,\frac{1}{6}\right]$ and note that
\begin{align}\label{exp1}
& \bar{\p}\left(\os \in A_{N,\alpha}(m(N))\Big|\os \in A\right)\nonumber\\ & =  \bar{Z}_N\big|_A^{-1} \sum_{m \in \mathcal{M}_N \bigcap A_{N,\alpha}(m(N))} \binom{N}{N(1+m)/2} \exp\left\{N(\beta_N m^p + h_N m - \log 2) \right\}.
\end{align}
By Lemma \ref{lem:multiple_max}, $\bar{\p}\left(\os \in A_{N,\alpha}(m(N))\Big|\os \in A\right) = 1-O(e^{-N^\alpha})$ and hence \eqref{exp1} gives us
\begin{equation}\label{secc}
\bar{Z}_N\big|_A = \left(1+O(e^{-N^\alpha})\right) \sum_{m \in \mathcal{M}_N \bigcap A_{N,\alpha}(m(N))} \binom{N}{N(1+m)/2} \exp\left\{N(\beta_N m^p + h_N m - \log 2) \right\}.
\end{equation}
Since $m(N)$ is the unique global maximizer of $H_N$ over the interval $A_{N,\alpha}(m(N))$,  by mimicking the proof of Lemma \ref{ex} on the interval $A_{N,\alpha}(m(N))$, it follows that 
\begin{align}\label{nonunmim1}
&\sum_{m \in \mathcal{M}_N\bigcap A_{N,\alpha}(m(N))} \binom{N}{N(1+m)/2} \exp\left\{N(\beta_N m^p + h_Nm-\log 2) \right\}\nonumber\\&= \frac{e^{N H_N(m(N))}}{\sqrt{(m(N)^2-1)H_N''(m(N))}}\left(1+O\left(N^{- \frac{1}{2}+3\alpha}\right)\right).
\end{align} 
The result in \eqref{eq:multi_max_I} now follows from \eqref{secc} and \eqref{nonunmim1}. 

For each $1\leq k\leq K$, (\ref{eq:multi_max_I}) immediately gives us 
\begin{equation}\label{tobeadded}
\bar{Z}_N\big|_{A_k} =  \frac{e^{NH_N(m_k(N))}}{\sqrt{(m_k(N)^2-1)H_N''(m_k(N))}}\left(1+O\left(N^{-\frac{1}{2} +\alpha}\right)\right),
\end{equation}
where the sets $A_1, \ldots,A_K$ are as defined in the proof of \eqref{eq:multi_max_local_II}. The result in \eqref{eq:multi_max_II} now follows from \eqref{tobeadded} on observing that $\bar{Z}_N = \sum_{k=1}^K \bar{Z}_N\big|_{A_k}$. \qed

\section{Perturbative Concentration Lemmas at $p$-critical Points}
\label{sec:pf_concentration}

It was shown in \eqref{eq:cltun_p1} that for $(\beta,h) \in \Theta$ which is $p$-critical, the limiting distribution of $\os$ assigns positive mass to each of the global maximizers $m_1,m_2,\ldots,m_K$. However, to use this result to obtain the limiting distribution of the ML estimates, we need to derive a similar concentration for $\os$ under $\p_{\beta_N,h_N,p}$. In particular, is it the case that $\os$ assigns positive mass to each of $m_1,m_2,\ldots,m_K$, or is the asymptotic support of $\os$ in this case a proper subset of $\{m_1,m_2,\ldots,m_K\}$ (we already know from \eqref{eq:multi_max_local_II} that the asymptotic support of $\os$ is a subset of  $\{m_1,m_2,\ldots,m_K\}$)? The answer to this question depends upon the rate of convergence of $(\beta_N,h_N)$ to $(\beta,h)$. This section is devoted to deriving these concentration results, which will be essential in proving the asymptotic distributions of $\hat{\beta}_N$ and $\hat{h}_N$ at the critical points, presented in Section \ref{sec:pf_beta_h} below. 

In what follows, assume  $(\beta,h) \in \Theta$ which is $p$-critical and let $m_1<m_2<\ldots<m_K$ be the global maximizers of $H_{\beta,h,p}$, and let $A_1,A_2,\ldots,A_K$ be the sets defined in the proof of \eqref{eq:multi_max_local_II} (in Section \ref{sec:proofofnonun_I}).  The following lemma shows that keeping $\beta$ fixed, if $h$ is perturbed at a rate slower than $1/N$, then under the perturbed sequence of measures, $\os$ concentrates around the largest/smallest global maximizer according as the perturbation is in the positive/negative direction, respectively.

\begin{lem}\label{redunsup}
	For any positive sequence $y_N$ satisfying $N^{-1} \ll  y_N \ll 1$, there exist positive constants $C_1$ and $C_2$ not depending on $N$, such that
	$$\p_{\beta,h+\bar{h}y_N,p}\left(\os \in A_{\bm{1}(\bar{h}<0) + K \bm{1}(\bar{h}>0)}^c\right) \leq C_1e^{-C_2 Ny_N}.$$
\end{lem} 

\begin{proof} Let $H_N:= H_{\beta,h+\bar{h}y_N,p}$ and $m_k(N)$ be the local maximizers of $H_N$ converging to $m_k$. In what follows, for two positive sequences $\phi_N$ and $\psi_N$, we will use the notation $\phi_N \lesssim \psi_N$ to denote that $\phi_N \leq C \psi_N$ for all $N$ and some constant $C$ not depending on $N$. Let $t := \bm{1}(\bar{h}<0) + K \bm{1}(\bar{h}>0)$. Then for any $s \neq t$, we have by Lemma \ref{lm:condpart} and Lemma \ref{bdhbest},
	\begin{align}\label{compact}
	\p_{\beta,h+\bar{h}y_N,p}(\os \in A_s) & =  \frac{Z_N(\beta,h+\bar{h}y_N,p)\big|_{A_s}}{Z_N(\beta,h+\bar{h}y_N,p)}\nonumber\\& \leq   \frac{Z_N(\beta,h+\bar{h}y_N,p)\big|_{A_s}}{Z_N(\beta,h+\bar{h}y_N,p)\big|_{A_t}}\nonumber \\&\lesssim  \sqrt{\frac{(m_t(N)^2-1)H_N''(m_t(N))}{(m_s(N)^2-1)H_N''(m_s(N))}} e^{N\left[H_N(m_s(N)) - H_N(m_t(N)) \right]}\nonumber\\&\lesssim  e^{N\left[\bar{h}y_N(m_s-m_t) + O(y_N^2)\right]} = e^{Ny_N\left[\bar{h}(m_s-m_t) + o(1)\right]}.
	\end{align} 
	Lemma \ref{redunsup} now follows from \eqref{compact}, since $\bar{h}(m_s-m_t) <0$ for every $s \neq t$, by definition. 
\end{proof}

The situation becomes a bit trickier when $h$ is fixed and $\beta$ is perturbed, as two cases arise depending upon the parity of $p$. The case $p \geq 3$ is odd, is the easier one, and is exactly similar to the previous setting. Note that in this case, $K = 2$.
\begin{lem}\label{redunsup2}
	Suppose that $p\geq 3$ is odd. Then, for any positive sequence $x_N$ satisfying $N^{-1} \ll  x_N \ll 1$, there exist positive constants $C_1$ and $C_2$ not depending on $N$, such that
	$$\p_{\beta+\bar{\beta}x_N,h,p}\left(\os \in A_{\bm{1}\{\bar{\beta}<0\} + 2 \cdot \bm{1}\{\bar{\beta}>0\}}^c\right) \leq C_1e^{-C_2 Nx_N}.$$
\end{lem}

\begin{proof} 
	Let $H_N:= H_{\beta+\bar{\beta}x_N,h,p}$ and $m_k(N)$ be the local maximizers of $H_N$ converging to $m_k$. Then for any $s \neq t := \bm{1}\{\bar{\beta}<0\} + 2 \cdot \bm{1}\{\bar{\beta}>0\}$, by exactly following the proof of Lemma \ref{redunsup}, one gets
	\begin{equation}\label{compact2}
	\p_{\beta+\bar{\beta}x_N,h,p}(\os \in A_s) \lesssim e^{Nx_N\left[\bar{\beta}(m_s^p-m_t^p) + o(1)\right]}.
	\end{equation}
	Lemma \ref{redunsup2} now follows from \eqref{compact2}, since $\bar{\beta}(m_s^p-m_t^p) <0$ for every $s \neq t$, by definition.  
\end{proof}

In the following lemma, we deal with the case $p\geq 4$ even. The result is presented in two cases, depending upon whether $h=0$ or not. Note that, if $h \neq 0$, then $K=2$. On the other hand, if $h = 0$, then we may assume that $\beta \geq \tilde{\beta}_p$, since otherwise, $(\beta,h)$ is $p$-regular.  In this case, $K=2$ if $\beta > \tilde{\beta}_p$ and $K=3$ if $\beta = \tilde{\beta}_p$.

\begin{lem}\label{redunsup3}
	The following hold when $p \geq 4$ is even.
	\begin{enumerate}
		\item[$(1)$]  Suppose that $h \neq 0$. Then, for any positive sequence $x_N$ satisfying $N^{-1} \ll  x_N \ll 1$, there exist positive constants $C_1$ and $C_2$ not depending on $N$, such that the following hold.
		\begin{itemize}
			\item[$\bullet$] If $h>0$, then
			$$\p_{\beta+\bar{\beta}x_N,h,p}\left(\os \in A_{\bm{1}\{\bar{\beta}<0\} + 2 \cdot \bm{1}\{\bar{\beta}>0\}}^c\right) \leq C_1e^{-C_2 Nx_N}. $$
			\item[$\bullet$] If $h < 0$, then
			$$\p_{\beta+\bar{\beta}x_N,h,p}\left(\os \in A_{\bm{1}\{\bar{\beta}>0\} + 2 \cdot \bm{1}\{\bar{\beta}<0\}}^c\right) \leq C_1e^{-C_2 Nx_N}. $$
		\end{itemize}
		\item[$(2)$] Suppose that $h = 0$. 
		\begin{itemize}
			\item[$\bullet$] If $\beta > \tilde{\beta}_p$, then for any sequence $(\beta_N,h_N) \rightarrow (\beta,h)$, there exists a positive constant $C$ not depending on $N$, such that
			\begin{equation}\label{cc1}
			\max\left\{\left|\p_{\beta_N,h_N,p}\left(\os \in A_1\right) -\frac{1}{2}\right| ,~ \left|\p_{\beta_N,h_N,p}\left(\os \in A_2\right) -\frac{1}{2}\right|\right\} \leq Ce^{-N^{\frac{1}{6}}}.
			\end{equation}
			\item[$\bullet$] If $\beta = \tilde{\beta}_p$ and $\bar{\beta} > 0$, then for any positive sequence $x_N$ satisfying $N^{-1} \ll  x_N \ll 1$, there exist positive constants $C_1$ and $C_2$ not depending on $N$, such that
			\begin{equation}\label{cc2}
			\max\left\{\left|\p_{\beta_N,h,p}\left(\os \in A_1\right)-\frac{1}{2}\right|,~\left|\p_{\beta_N, h,p}\left(\os \in A_3\right) - \frac{1}{2}\right|\right\} \leq C_1e^{-C_2 Nx_N},
			\end{equation}
			where $\beta_N=\beta + \bar{\beta} x_N$. 
			\item[$\bullet$] If $\beta = \tilde{\beta}_p$ and $\bar{\beta} < 0$, then for any positive sequence $x_N$ satisfying $N^{-1} \ll  x_N \ll 1$, there exist a positive constants $C_1$ and $C_2$ not depending on $N$, such that
			\begin{equation}\label{cc3}
			\p_{\beta+\bar{\beta}x_N,h,p}\left(\os \in A_2^c\right) \leq C_1 e^{-C_2 N x_N}.
			\end{equation}
		\end{itemize}
	\end{enumerate}
	
\end{lem} 

\begin{proof}
	The proof of (1) is exactly similar to that of Lemma \ref{redunsup2}, and hence we ignore it. One only needs to observe that $m_1<m_2<0$ if $h <0$, and $0<m_1<m_2$ if $h> 0$. Hence, $m_1^p < m_2^p$ if $h >0$, and $m_1^p > m_2^p$ if $h<0$.

	Next, we prove (2). Note that \eqref{cc1} follows directly from \eqref{eq:multi_max_local_II} in Lemma \ref{lem:multiple_max} (taking $\alpha = \frac{1}{6}$) and using the fact that for even $p$ and $h=0$, $\os \stackrel{D}{=} -\os$. Next, note that if $\beta = \tilde{\beta}_p$ and $\bar{\beta} > 0$, then from \eqref{compact2},
	$\p_{\beta + \bar{\beta} x_N,h,p}\left(\os \in A_2\right) \lesssim e^{-C_2 Nx_N}$ for some positive constant $C_2$ not depending on $N$, and \eqref{cc2} follows from the symmetry of the distribution of $\os$. Finally, if $\beta = \tilde{\beta}_p$ and $\bar{\beta} < 0$, then once again from \eqref{compact2}, $\p_{\beta + \bar{\beta} x_N,h,p}\left(\os \in A_k\right) \lesssim e^{-C_2 N x_N}$ for some positive constant $C_2$ not depending on $N$ and $k \in \{1,3\}$. This gives \eqref{cc3} and completes the proof of Lemma \ref{redunsup3}. 
\end{proof}

%\iffalse

\section{Proofs from Section \ref{sec:mle_beta_h_I}} 
\label{sec:pf_beta_h}

In this section we derive the limiting distribution of the ML estimates as presented in Section \ref{sec:mle_beta_h_I}. The proofs of Theorems \ref{cltintr3_1}, \ref{cltintr3_2}, \ref{thmmle1} and \ref{thmmle2} are given Section \ref{sec:pf_beta_h_1}. The proof of Theorem \ref{cltintr3_III} is given in Section \ref{sec:pf_beta_h_2}, and the proof of the Theorem \ref{thmmle_III} is Section \ref{sec:pf_beta_III}.

\subsection{Proofs of Theorems \ref{cltintr3_1}, \ref{cltintr3_2}, \ref{thmmle1} and \ref{thmmle2}} 
\label{sec:pf_beta_h_1} 

We will only prove the case $(\beta,h)$ is $p$-regular, which includes  Theorems \ref{cltintr3_1} and \ref{thmmle1}. The proofs for the $p$-special case, that is,  Theorems\ref{cltintr3_2} and \ref{thmmle2}, follow similarly from part (3) of Theorem \ref{cltun}.

\subsubsection{Proof of Theorem \ref{cltintr3_1}}

For any $t \in \mathbb{R}$, we have by \eqref{eqmleh}, Lemma \ref{increasing}, \eqref{eq:meanclt_I} and \eqref{eq:cltun_I}, together with the fact that pointwise convergence of moment generating functions on $\mathbb{R}$ imply convergence of moments, 
\begin{align}\label{trclt2}
\p_{\beta,h,p}  \left(N^{\frac{1}{2}}(\hat{h}_N -h) \leq t\right) &= \p_{\beta,h,p} \left(\hat{h}_N  \leq h + \frac{t}{N^{\frac{1}{2}}}\right)\nonumber\\
&= \p_{\beta,h,p} \left(u_{N,1}(\beta,\hat{h}_N ,p) \leq u_{N,1}\left(\beta, h + \frac{t}{N^{\frac{1}{2}}},p\right) \right)\nonumber\\
&= \p_{\beta,h,p} \left(\os \leq \e_{\beta, h+ N^{-\frac{1}{2}} t,p}  (\os) \right)\nonumber\\
&= \p_{\beta,h,p} \left(N^{\frac{1}{2}}(\os-m_*) \leq \e_{\beta, h+ N^{-\frac{1}{2}} t,p}  (N^{\frac{1}{2}}(\os-m_*)) \right)\nonumber\\
& \rightarrow \p_{\beta,h,p} \left(N\left(0, -\frac{1}{H''(m_*)}\right)\leq -\frac{t}{H''(m_*)} \right)\nonumber\\
&= \p_{\beta,h,p} \left(N\left(0, -H''(m_*)\right)\leq t \right).
\end{align}
Now, the proof of Theorem \ref{cltintr3_1} follows from \eqref{trclt2}.

\subsubsection{Proof of Theorem \ref{thmmle1}}

We begin with the case $m_* \ne 0$. By Theorem \ref{cltun}, $\left(\os - m_*\right)^s =  O_{P}(N^{-\frac{s}{2}}) = O_{P}(N^{-1})$, for every $s \geq 2$ under  $\p = \p_{\beta+ \bar{\beta}/\sqrt N, h, p}$. Further, since pointwise convergence of moment generating functions on $\mathbb{R}$ imply convergence of moments, we also have $\e_{\beta+ \bar{\beta}/\sqrt N, h, p}(\os - m_*)^s  = O(N^{-1})$, for every $s \geq 2$. Now,
\begin{equation}\label{binoexp}
N^{\frac{1}{2}}(\overline{X}^p_N - m_*^p) =N^{\frac{1}{2}}pm_*^{p-1} (\os-m_*) + N^{\frac{1}{2}}\sum_{s=2}^p \binom{p}{s} m_*^{p-s}(\os-m_*)^s.
\end{equation}
It follows from Theorem \ref{cltun} and \eqref{binoexp} that under $\p_{\beta+{N^{-\frac{1}{2}}\bar{\beta}},h,p} $,
\begin{equation}\label{mpl1}
N^{\frac{1}{2}}(\overline{X}^p_N - m_*^p) \xrightarrow{D} N\left(-\frac{\bar{\beta}p^2m_*^{2p-2}}{H''(m_*)},-\frac{p^2m_*^{2p-2}}{H''(m_*)}\right) ,
\end{equation} 
and 
\begin{equation}\label{mpl1_N}
\e\left[N^{\frac{1}{2}}(\overline{X}^p_N - m_*^p)\right] \rightarrow -\frac{\bar{\beta}p^2m_*^{2p-2}}{H''(m_*)}.
\end{equation} 
Now, note that for any $t \in \mathbb{R}$, we have by \eqref{eqmle} and the monotonicity of the function $u_{N, p}(\cdot, h, p)$ (Lemma \ref{increasing}), we have
\begin{align}\label{trclt}
\p_{\beta,h,p}  \left(N^{\frac{1}{2}}(\hat{\beta}_N-\beta) \leq t\right) &= \p_{\beta,h,p} \left(\hat{\beta}_N \leq \beta + \frac{t}{N^{\frac{1}{2}}}\right)\nonumber\\
&= \p_{\beta,h,p} \left(u_{N,p}(\hat{\beta}_N,h,p) \leq u_{N,p}\left(\beta + \frac{t}{N^{\frac{1}{2}}}, h,p\right) \right)\nonumber\\
&= \p_{\beta,h,p} \left(\overline{X}^p_N \leq \e_{\beta + N^{-\frac{1}{2}} t, h,p}  (\overline{X}^p_N) \right)\nonumber\\ 
&= \p_{\beta,h,p} \left(N^{\frac{1}{2}}(\overline{X}^p_N-m_*^p) \leq \e_{\beta + N^{-\frac{1}{2}} t, h,p}  (N^{\frac{1}{2}}(\overline{X}^p_N-m_*^p)) \right).\nonumber\\
\end{align} 
Now, weak convergence to a continuous distribution implies uniform convergence of the distribution functions, by \eqref{mpl1}, \eqref{mpl1_N}, and \eqref{trclt}, it follows that under $\p_{\beta,h,p} $, 
\begin{align*}%\label{finmleconc}
\p_{\beta,h,p}  \left(N^{\frac{1}{2}}(\hat{\beta}_N-\beta) \leq t\right) &\rightarrow \p_{\beta,h,p} \left(N\left(0,-\frac{p^2m_*^{2p-2}}{H''(m_*)}\right) \leq -\frac{tp^2m_*^{2p-2}}{H''(m_*)} \right)  \nonumber \\ 
& = \p_{\beta,h,p}\left(N\left(0, -\frac{H''(m_*)}{p^2m_*^{2p-2}}\right)\leq t \right) .
\end{align*}
This completes the proof of \eqref{eq:bmle_m_1}. 

Next, we consider the case $m_* = 0$. This implies that $\sup_{x \in [-1,1]} H_{\beta,h,p}(x) = 0$. Hence, by part (1) of Lemma \ref{derh11}, $h=0$, and then, \eqref{eq:betatilde} implies $\beta \leq \tilde{\beta}_p$. However, the point $(\tilde{\beta}_p,0)$ is $p$-critical, and hence, we must have $\beta < \tilde{\beta}_p$. Now, for every $t \in \mathbb{R}$,we have by \eqref{eqmle} and Lemma \ref{increasing}, 
\begin{equation}\label{fct1}
\p_{\beta,0,p}\left(\hat{\beta}_N > t\right) = \p_{\beta,0,p}\left({\left(N^{\frac{1}{2}}\os\right)}^p > N^{\frac{p}{2}}u_{N,p}(t,0,p)\right).
\end{equation}
First, fix $t \in (\tilde{\beta}_p,\infty)$ and note that:
\begin{equation}\label{adnlst1}
u_{N,p}(t,0,p) = \frac{1}{N} \frac{\partial}{\partial \underline{\beta}}F_N(\underline{\beta},0,p)\Big|_{\underline{\beta} = t}. 
\end{equation}
Now, by the mean value theorem and the fact that $F_N(0,0,p) = 0$, we have:
\begin{equation}\label{adnlstp2}
F_N(t,0,p) = t\frac{\partial}{\partial \underline{\beta}}F_N(\underline{\beta},0,p)\Big|_{\underline{\beta} = \xi}
\end{equation}
for some $\xi \in (0,t)$. By Lemma \ref{increasing}, we have:
\begin{equation}\label{adnlstp3}
\frac{\partial}{\partial \underline{\beta}}F_N(\underline{\beta},0,p)\Big|_{\underline{\beta} = \xi} \leq \frac{\partial}{\partial \underline{\beta}}F_N(\underline{\beta},0,p)\Big|_{\underline{\beta} = t}. 
\end{equation}
Combining \eqref{adnlst1}, \eqref{adnlstp2} and \eqref{adnlstp3}, we have:
\begin{equation}\label{adnlstp4}
u_{N,p}(t,0,p) \geq t^{-1} N^{-1} F_N(t,0,p). 
\end{equation}
Now, \eqref{logpartex} in  Lemma \ref{ex} (for odd $p$) and \eqref{eq:multi_max_II} in Lemma \ref{lm:condpart} (for even $p$) implies that\footnote{For two positive sequences $\{a_n\}_{n \geq 1}$ and $\{b_n\}_{n \geq 1}$, $a_n = \Omega(b_n)$, if there exists a positive constant $C$, such that $a_n \geq C b_n$, for all large $n$.} $$N^{-1} F_N(t,0,p) = \Omega(1).$$
This, together with \eqref{adnlstp4} implies that:
\begin{equation}\label{udifc}
u_{N,p}(t,0,p) = \Omega(1).
\end{equation}
Since, by \eqref{eq:meanclt_I}, $N^{1/2}\os \xrightarrow{D} N(0,1)$ under $\p_{\beta,0,p}$,   \eqref{fct1} and \eqref{udifc} implies, as $N \rightarrow \infty$,
\begin{equation}\label{tgeqbet}
\p_{\beta,0,p}\left(\hat{\beta}_N > t\right) \rightarrow 0.
\end{equation}
Next, fix $t \in [0,\tilde{\beta}_p)$. Since we have pointwise convergence of moment generating functions in part (1) of Theorem \ref{cltintr1}, we get:
\begin{equation}\label{adl}
N^{\frac{p}{2}}u_{N,p}(t,0,p) = \e_{t,0,p}[(N^{\frac{1}{2}}\os)^p] \rightarrow \e Z^p. 
\end{equation}
Hence by \eqref{fct1},	
\begin{equation}\label{tgeqbet1}
\p_{\beta,0,p}\left(\hat{\beta}_N \leq t\right) \rightarrow \gamma_p.
\end{equation}
Finally, fix $t\in (-\infty,0)$. If $p$ is odd, the function $\beta \mapsto F_N(\beta,0,p)$ becomes an even function (recall \eqref{ptnepn}), and hence, its partial derivative with respect to $\beta$ becomes an odd function. Consequently, $$u_{N,p}(t,0,p) = -u_{N,p}(-t,0,p). $$ Now, if $t <  -\tilde{\beta}_p$, then $-t \in ( \tilde{\beta}_p,\infty)$, so by \eqref{udifc}, $N^{\frac{p}{2}}u_{N,p}(-t,0,p)$ converges to $\infty$ , i.e. $$\lim_{N\rightarrow\infty} N^{\frac{p}{2}}u_{N,p}(t,0,p) = -\infty. $$ If $t > -\tilde{\beta}_p$, then $-t \in (0,\tilde{\beta}_p)$, and hence, by \eqref{adl} (note that $\mathbb{E} Z^p = 0$ when $p$ is odd)
$$\lim_{N\rightarrow\infty} N^{\frac{p}{2}}u_{N,p}(t,0,p) = -\lim_{N\rightarrow\infty} N^{\frac{p}{2}}u_{N,p}(-t,0,p) = 0. $$                                
Hence, we have from \eqref{fct1}, as $N \rightarrow \infty$,

\[   
\p_{\beta,0,p}\left(\hat{\beta}_N \leq t\right) \rightarrow 
\begin{cases}
0 &\quad\text{if}~t < -\tilde{\beta}_p,\\
\frac{1}{2} &\quad\text{if}~t > -\tilde{\beta}_p\\
\end{cases}
\]
This, combined with \eqref{tgeqbet} and \eqref{tgeqbet1}, shows that $\hat{\beta}_N \xrightarrow{D} \frac{1}{2}\delta_{\tilde{\beta}_p} + \frac{1}{2} \delta_{-\tilde{\beta}_p}$ if $p$ is odd.

Now, assume that $p \geq 4$ is even. Then, for $t\in (-\infty,0)$, 
\begin{align}\label{myt1}
\bigg|N^{\frac{p}{2}-1} \frac{\partial}{\partial \underline \beta} F_N(\underline \beta,0,p)\Big|_{\underline \beta = t} - & \mathbb{E}_{0,0,p} \left[N^{\frac{p}{2}}\overline{X}^p_N\right]\bigg|   \nonumber\\ 
&= N^{\frac{p}{2}-1}\left\{\frac{\partial}{\partial \underline \beta} F_N(\underline \beta,0,p)\Big|_{\underline \beta = 0} - \frac{\partial}{\partial \underline \beta} F_N(\underline \beta,0,p)\Big|_{\underline \beta = t}\right \} \nonumber\\
&   \leq -tN^{\frac{p}{2}-1} \sup_{\zeta \in [t,0]} \frac{\partial^2}{\partial \underline \beta^2} F_N(\underline \beta, 0, p)\Big|_{\underline \beta=\zeta}  \nonumber\\
&=-tN^{\frac{p}{2}-1} \sup_{\zeta \in [t,0]} \mathrm{Var}_{\zeta,0,p} \left(N\overline{X}^p_N\right)\nonumber\\ 
&= -tN^{1-\frac{p}{2}}\sup_{\zeta \in [t,0]} \mathrm{Var}_{\zeta,0,p} \left(N^{\frac{p}{2}}\overline{X}^p_N\right)\nonumber\\
& \leq -tN^{1-\frac{p}{2}}\sup_{\zeta \in [t,0]} \mathbb{E}_{\zeta,0,p}  (N^{p}\overline{X}^{2p}_N) \nonumber\\
&= -t N^{1-\frac{p}{2}}\sup_{\zeta \in [t,0]} \mathbb{E}_{0,0,p}  \left[N^{p}\overline{X}^{2p}_N e^{\zeta N\overline{X}^p_N-F_N(\zeta,0,p)}\right].\nonumber\\
\end{align}
Next, for every $\zeta \in [t,0]$, since the map $\beta \mapsto \frac{\partial}{\partial \beta} F_N(\beta,0,p)$ is increasing,
$$-F_N(\zeta,0,p) \leq -\zeta \frac{\partial}{\partial \underline \beta} F_N(\underline \beta,0,p)\Big|_{\underline \beta=0} \leq -t N\mathbb{E}_{0,0,p} \overline{X}^p_N = o(1)$$ $$\implies \sup_{N \geq 1}\sup_{\zeta \in [t,0]} e^{-F_N(\zeta,0,p)} := B < \infty.$$ We thus have from \eqref{myt1},
\begin{equation*}
\left|N^{\frac{p}{2}-1} \frac{\partial}{\partial \underline \beta} F_N(\underline \beta,0,p)\Big|_{\underline \beta = t} - \mathbb{E}_{0,0,p} \left[N^{\frac{p}{2}}\overline{X}^p_N\right]\right| \leq -tBN^{1-\frac{p}{2}} \mathbb{E}_{0,0,p}\left[N^p\overline{X}^{2p}_N\right] = o(1).
\end{equation*}
Hence, $N^{\frac{p}{2}} u_{N,p}(t,0,p) = N^{\frac{p}{2}-1} \frac{\partial}{\partial \tilde{\beta}} F_N(\tilde{\beta},0,p)\Big|_{\tilde{\beta} = t} \rightarrow \e Z^p$ as $N \rightarrow \infty$. Consequently, as $N \rightarrow \infty$,
\begin{equation}\label{lstpeven}
\p_{\beta,0,p}\left(\hat{\beta}_N \leq t\right) \rightarrow \gamma_p.
\end{equation}
We conclude from \eqref{tgeqbet}, \eqref{tgeqbet1} and \eqref{lstpeven}, that $\hat{\beta}_N \xrightarrow{D} \gamma_p \delta_{-\infty} + (1-\gamma_p) \delta_{\tilde{\beta}_p}$ if $p$ is even. This completes the proof of \eqref{eq:bmle_m_2}. \qed

\begin{rem}[Efficiency of the ML estimates at $p$-regular points]\label{remark1} 
	An interesting consequence of the results proved above is that, at the $p$-regular points, the  limiting variance of the ML estimates equals the limiting inverse Fisher information, that is, the ML estimates are asymptotically efficient. To see this, note that the Fisher information of $\beta$ and $h$ (scaled by $N$) in the model \eqref{cwwithmg} are given by 
	$$I_N(\beta) = \frac{1}{N} \e_{\beta,h,p}\left[\left(\frac{\partial}{\partial \beta} \log \p_{\beta,h,p}(\bs)\right)^2\right] = \mathrm{Var}_{\beta,h,p}(N^{\frac{1}{2}} \os^p)$$ 
	and 
	$$I_N(h) = \frac{1}{N} \e_{\beta,h,p}\left[\left(\frac{\partial}{\partial h} \log \p_{\beta,h,p}(\bs)\right)^2\right] = \mathrm{Var}_{\beta,h,p}(N^{\frac{1}{2}} \os),$$ 
	respectively. It follows from the proof of Theorem \ref{cltintr1}, that for a $p$-regular point $(\beta,h)$, the moment generating of $\sqrt{N}\left(\os- m_*\right)$ converges pointwise to that of the centered Gaussian distribution with variance $-\left[H''(m_*)\right]^{-1}$. Hence, 
	\begin{equation}\label{hinfconv}
	\lim_{N \rightarrow \infty}I_N(h) =  -\left[H''(m_*)\right]^{-1}~.
	\end{equation} 
	Also, it follows from \eqref{binoexp} and \eqref{mpl1} and the fact $\e_{\beta,h,p}\left[(\os- m_*)^s\right] = O(N^{-s/2})$, for each $s \geq 1$, that
	\begin{equation}\label{betainfconv}
	\lim_{N \rightarrow \infty}I_N(\beta) =  -\frac{p^2 m_*^{2p-2}}{H''(m_*)}~.
	\end{equation} 
	Therefore, by Theorem \ref{cltintr3_1}, at a $p$-regular point $(\beta,h)$, $\hat{h}_N$ is an efficient estimator of $h$, and by Theorem \ref{thmmle1}, if $(\beta,h)$ is a $p$-regular point with $m_*\neq 0$, then $\hat{\beta}_N$ is an efficient estimator of $\beta$.
\end{rem}

\subsection{Proof of Theorem \ref{cltintr3_III}}
\label{sec:pf_beta_h_2} 

Recall the definitions of the sets $A_k~(1\leq k\leq K)$ from the proof of  Lemma \ref{lem:multiple_max} (in Section \ref{sec:proofofnonun_I}). Now, fixing $t <0$, we have similar to the proof of \eqref{eq:h_1},
\begin{align*}
\p_{\beta,h,p}  \left(N^{\frac{1}{2}}(\hat{h}_N -h) \leq t\right)  &= \p_{\beta,h,p} \left(N^{\frac{1}{2}}(\os-m_1) \leq \e_{\beta, h+ N^{-\frac{1}{2}} t,p}  (N^{\frac{1}{2}}(\os-m_1)) \right) \nonumber\\
& = T_1 + T_2, %\label{m2}
\end{align*}
where 
\begin{align*}
T_1&= \p_{\beta,h,p} \left(N^{\frac{1}{2}}(\os-m_1) \leq \e_{\beta, h+ N^{-\frac{1}{2}} t,p}  (N^{\frac{1}{2}}(\os-m_1))\Big| \os \in A_1\right) \p_{\beta,h,p} (\os \in A_1) , \\%\label{m1} \\
T_2 & = \p_{\beta,h,p} \left(N^{\frac{1}{2}}(\os-m_1) \leq \e_{\beta, h+ N^{-\frac{1}{2}} t,p}  (N^{\frac{1}{2}}(\os-m_1))\Big| \os \in A_1^c\right) \p_{\beta,h,p} (\os \in A_1^c). %\label{m2} 
\end{align*} 
Now, by the law of iterated expectations, we have
\begin{align}\label{lie1}
\e_{\beta, h+ N^{-\frac{1}{2}} t,p}  (N^{\frac{1}{2}}(\os-m_1))  & =  S_1 + S_2 , 
\end{align}
where 
\begin{align}\label{lie2}
S_1:= \e_{\beta, h+ N^{-\frac{1}{2}} t,p}  \left(N^{\frac{1}{2}}(\os-m_1)\Big|\os \in A_1\right)  \p_{\beta, h+ N^{-\frac{1}{2}} t,p} (\os \in A_1) 
\end{align}
and
\begin{align}\label{lie3}
S_2:= \e_{\beta, h+ N^{-\frac{1}{2}} t,p}  \left(N^{\frac{1}{2}}(\os-m_1)\Big|\os \in A_1^c\right) \p_{\beta, h+ N^{-\frac{1}{2}} t,p} (\os \in A_1^c). 
\end{align}  
Note that by \eqref{moment generating functionnonun}, $$ \e_{\beta, h+ N^{-\frac{1}{2}} t,p} \left(N^{\frac{1}{2}}(\os-m_1)\Big|\os \in A_1\right) \rightarrow - \frac{t}{H''(m_1)},$$ 
as $N \rightarrow \infty$. Also, by Lemma \ref{redunsup}, $\p_{\beta, h+ t/\sqrt N,p} (\os \in A_1^c) \leq C_1 e^{-C_2 N^{\frac{1}{2}}}$ for positive constants $C_1,C_2$ not depending on $N$. Hence, \eqref{lie2} converges to $-t/H''(m_1)$ and \eqref{lie3} converges to $0$. Consequently, \eqref{lie1} converges to $-t/H''(m_1)$. 

Next, under $\p_{\beta, h,p}(~\cdot~|\os \in A_1^c)$, $N^{\frac{1}{2}}(\os - m_1) \xrightarrow{P} \infty,$ by Lemma \ref{lem:multiple_max}. Hence, $T_2 \rightarrow 0$. Also, by Theorem \ref{cltintr1}, $N^{\frac{1}{2}}(\os - m_1) \xrightarrow{D} N\left(0,-1/H''(m_1)\right)$ under $\p_{\beta,h,p} (~\cdot~|\os \in A_1)$. Hence, $T_1$ converges to
$$p_1\p(N\left(0,-\frac{1}{H''(m_1)}) \leq -\frac{t}{H''(m_1)}\right) = p_1\p\left(N(0,-H''(m_1)) \leq t\right).$$ Hence,  
\begin{equation}\label{htl0}
\p_{\beta,h,p}  \left(N^{\frac{1}{2}}(\hat{h}_N -h) \leq t\right)  \rightarrow p_1\p\left(N(0,-H''(m_1)) \leq t\right)\quad\quad \textrm{for all}\quad t<0.
\end{equation}
Next, fix $t >0$, whence we have
\begin{equation*}
\p_{\beta,h,p}  \left(N^{\frac{1}{2}}(\hat{h}_N -h) > t\right) = T_3+ T_4,
\end{equation*}
where 
\begin{align*}
T_3&= \p_{\beta,h,p} \left(N^{\frac{1}{2}}(\os-m_K) > \e_{\beta, h+ N^{-\frac{1}{2}} t,p}  (N^{\frac{1}{2}}(\os-m_K))\Big| \os \in A_K\right) \p_{\beta,h,p} (\os \in A_K), \\%\label{n1} \\
T_4 & = \p_{\beta,h,p} \left(N^{\frac{1}{2}}(\os-m_K) > \e_{\beta, h+ N^{-\frac{1}{2}} t,p}  (N^{\frac{1}{2}}(\os-m_K))\Big| \os \in A_K^c\right) \p_{\beta,h,p} (\os \in A_K^c).%\label{n2}
\end{align*} 
By the same arguments as before, it follows that $$\e_{\beta, h+  N^{-\frac{1}{2}} t, p}  (N^{\frac{1}{2}}(\os-m_K)) \rightarrow -\frac{t}{H''(m_K)}.$$ Next, under $\p_{\beta, h,p} (~\cdot~|\os \in A_K^c )$, $N^{\frac{1}{2}}(\os - m_K) \xrightarrow{P} -\infty$ by Lemma \ref{lem:multiple_max}. Hence, $T_4\rightarrow 0$. Also, by Theorem \ref{cltintr1}, $N^{\frac{1}{2}}(\os - m_K) \xrightarrow{D} N(0,-1/H''(m_K))$ under $\p_{\beta,h,p} (~\cdot~|\os \in A_K)$. Hence, $T_3$ converges to
$$p_K\p\left(N\left(0,-\frac{1}{H''(m_K)}\right) > -\frac{t}{H''(m_K)}\right) = p_K\p\left(N(0,-H''(m_K)) > t\right).$$ Hence,
\begin{equation}\label{htl1}
\p_{\beta,h,p}  \left(N^{\frac{1}{2}}(\hat{h}_N -h) > t\right) \rightarrow p_K\p\left(N(0,-H''(m_K)) > t\right)\quad\quad \textrm{for all}\quad t>0.
\end{equation}
Combining \eqref{htl0} and \eqref{htl1}, we conclude that for all $p$-critical points $(\beta,h)$, under $\p_{\beta,h,p} $,
\begin{equation}\label{htl3}
N^{\frac{1}{2}}(\hat{h}_N -h) \xrightarrow{D} \tfrac{p_1}{2} N^{-}(0,-H''(m_1)) + \tfrac{p_K}{2} N^+(0,-H''(m_K)) + \left(1-\frac{p_1+p_K}{2}\right)\delta_0.	
\end{equation}
Theorem \ref{cltintr3_III} follows from \eqref{htl3} on observing that if $p\geq 4$ is even and $(\beta,h) = (\tilde{\beta}_p,0)$, then $K=3$, $m_3 = -m_1$ and $p_1=p_3$, and otherwise, $K=2$.

\subsection{Proof of Theorem \ref{thmmle_III}} 
\label{sec:pf_beta_III} 

We first deal with the case $p \geq 3$ is odd.

\noindent\textit{Proof of \eqref{eq:bmle_multimax_1}:} In this case, $0$ is not a global maximizer of $H_{\beta,p,h}$. Fixing $t < 0$, we have
\begin{align*}
\p_{\beta,h,p}  \left(N^{\frac{1}{2}}(\hat{\beta}_N-\beta) \leq t\right) &= \p_{\beta,h,p} \left(N^{\frac{1}{2}}(\overline{X}^p_N-m_1^p) \leq \e_{\beta+N^{-\frac{1}{2}} t, h ,p}  (N^{\frac{1}{2}}(\overline{X}^p_N-m_1^p)) \right)\nonumber\\&= T_5+T_6,
\end{align*}
where
\begin{align*}
T_5&= \p_{\beta,h,p} \left(N^{\frac{1}{2}}(\overline{X}^p_N-m_1^p) \leq \e_{\beta+N^{-\frac{1}{2}} t, h ,p}  (N^{\frac{1}{2}}(\overline{X}^p_N-m_1^p))\Big| \os \in A_1\right)\p_{\beta,h,p} (\os \in A_1), \\%\label{m1} \\
T_6 & = \p_{\beta,h,p} \left(N^{\frac{1}{2}}(\overline{X}^p_N-m_1^p) \leq \e_{\beta+N^{-\frac{1}{2}} t, h ,p}  (N^{\frac{1}{2}}(\overline{X}^p_N-m_1^p))\Big| \os \in A_1^c\right)\p_{\beta,h,p} (\os \in A_1^c). %\label{m2} 
\end{align*} 
Now, by the law of iterated expectations, we have
\begin{align}
\e_{\beta+N^{-\frac{1}{2}} t, h ,p}  (N^{\frac{1}{2}}(\overline{X}^p_N-m_1^p)) = S_3 + S_4, \label{lie11} 
\end{align} 
where 
\begin{align} 
S_3:=\e_{\beta+N^{-\frac{1}{2}} t, h,p}  \left(N^{\frac{1}{2}}(\overline{X}^p_N-m_1^p)\Big|\os \in A_1\right) \p_{\beta+ N^{-\frac{1}{2}} t, h,p} (\os \in A_1)\label{lie21}
\end{align} 
and 
\begin{align} 
S_4:=\e_{\beta+N^{-\frac{1}{2}} t, h,p}  \left(N^{\frac{1}{2}}(\overline{X}^p_N-m_1^p)\Big|\os \in A_1^c\right) \p_{\beta+ N^{-\frac{1}{2}} t, h,p} (\os \in A_1^c).\label{lie31}
\end{align} 

\noindent From Theorem \ref{cltun} and a simple binomial expansion (see \eqref{binoexp}), it follows that
\begin{equation*}%\label{eqe11}
\e_{\beta+N^{-\frac{1}{2}} t, h,p}  \left(N^{\frac{1}{2}}(\overline{X}^p_N-m_1^p)\Big|\os \in A_1\right) \rightarrow -\frac{tp^2m_1^{2p-2}}{H''(m_1)},
\end{equation*}
and under $\p_{\beta,h,p} (~\cdot~|\os \in A_1)$,
\begin{equation}\label{eqel2}
N^{\frac{1}{2}}(\overline{X}^p_N - m_1^p) \xrightarrow{D} N\left(0,-\frac{p^2m_1^{2p-2}}{H''(m_1)}\right). 
\end{equation}

\noindent By Lemma \ref{redunsup2}, $\p_{\beta+ t/\sqrt N, h,p} (\os \in A_1^c) \leq C_1 e^{-C_2 \sqrt N}$ for positive constants $C_1,C_2$ not depending on $N$. Hence, \eqref{lie21} converges to $-tp^2m_1^{2p-2}/H''(m_1)$ and \eqref{lie31} converges to $0$. Consequently, \eqref{lie11} converges to $-tp^2m_1^{2p-2}/H''(m_1)$.  

Next, under $\p_{\beta, h,p} (~\cdot~|\os \in A_1^c )$, $N^{\frac{1}{2}}(\overline{X}^p_N - m_1^p) \xrightarrow{P} \infty$ by Lemma \ref{lem:multiple_max}. Hence, $T_6 \rightarrow 0$. Then, by \eqref{eqel2}, $T_5$ converges to
$$p_1\p\left(N\left(0,-\frac{p^2m_1^{2p-2}}{H''(m_1)}\right) \leq -\frac{tp^2m_1^{2p-2}}{H''(m_1)}\right) = p_1\p\left(N\left(0,-\frac{H''(m_1)}{p^2 m_1^{2p-2}}\right) \leq t\right).$$ Hence, for all $t<0$, we have:
\begin{equation}\label{htl03}
\p_{\beta,h,p}  \left(N^{\frac{1}{2}}(\hat{\beta}_N-\beta) \leq t\right) \rightarrow p_1\p\left(N\left(0,-\frac{H''(m_1)}{p^2 m_1^{2p-2}}\right)\leq t\right).
\end{equation}
Next, fix $t >0$, whence we have
\begin{equation*}
\p_{\beta,h,p}  \left(N^{\frac{1}{2}}(\hat{\beta}_N-\beta) > t\right) = T_7 + T_8,
\end{equation*}
where
\begin{align*}
T_7&= \p_{\beta,h,p} \left(N^{\frac{1}{2}}(\overline{X}^p_N-m_2^p) > \e_{\beta+ N^{-\frac{1}{2}} t, h,p}  (N^{\frac{1}{2}}(\overline{X}^p_N-m_2^p))\Big| \os \in A_2\right)\p_{\beta,h,p} (\os \in A_2),\\%\label{m1} \\
T_8 & = \p_{\beta,h,p} \left(N^{\frac{1}{2}}(\overline{X}^p_N-m_2^p) > \e_{\beta+ N^{-\frac{1}{2}} t, h,p}  (N^{\frac{1}{2}}(\overline{X}^p_N-m_2^p))\Big| \os \in A_2^c\right)\p_{\beta,h,p} (\os \in A_2^c). %\label{m2} 
\end{align*}
By the same arguments as before, it follows that $$\e_{\beta+ N^{-\frac{1}{2}} t, h,p}  (N^{\frac{1}{2}}(\overline{X}^p_N-m_2^p)) \rightarrow -\frac{tp^2 m_2^{2p-2}}{H''(m_2)}.$$ Next, under $\p_{\beta, h,p} (~\cdot~\big|\os \in A_2^c )$, $N^{\frac{1}{2}}(\overline{X}^p_N - m_2^p) \xrightarrow{P} -\infty$ by Lemma \ref{lem:multiple_max}. Hence, $T_8 \rightarrow 0$. Also, we know that  $N^{\frac{1}{2}}(\overline{X}^p_N - m_2^p) \xrightarrow{D} N\left(0,-p^2m_2^{2p-2}/H''(m_2)\right)$ under $\p_{\beta,h,p} \left(~\cdot~\big|\os \in A_2\right)$. Hence, $T_7$ converges to
$$p_2\p\left(N\left(0,-\frac{p^2m_2^{2p-2}}{H''(m_2)}\right) > -\frac{tp^2m_2^{2p-2}}{H''(m_2)}\right) = p_2\p\left(N\left(0,-\frac{H''(m_2)}{p^2 m_2^{2p-2}}\right) > t\right).$$ Hence,
\begin{equation}\label{htl13}
\p_{\beta,h,p}  \left(N^{\frac{1}{2}}(\hat{\beta}_N-\beta) > t\right) \rightarrow p_2\p\left(N\left(0,-\frac{H''(m_2)}{p^2 m_2^{2p-2}}\right) > t\right), \quad \textrm{for all } t>0.
\end{equation}
Combining \eqref{htl03} and \eqref{htl13}, we conclude that if $p \geq 3$ is odd, then for all $p$-critical points $(\beta,h)$, under $\p_{\beta,h,p} $,
\begin{equation}\label{htl33}
N^{\frac{1}{2}}(\hat{\beta}_N-\beta) \xrightarrow{D} \frac{p_1}{2} N^{-}\left(0,-\frac{H''(m_1)}{p^2 m_1^{2p-2}}\right) + \frac{p_2}{2} N^+\left(0,-\frac{H''(m_2)}{p^2 m_2^{2p-2}}\right) + \left(1-\frac{p_1+p_2}{2}\right)\delta_0.	
\end{equation}
\eqref{eq:bmle_multimax_1} now follows from \eqref{htl33} on observing that $p_2=1-p_1$. \\ 

\medskip

\noindent \textit{Proof of \eqref{eq:bmle_multimax_2}:} In this case, $m_1 = 0$. We can write for any $t < 0$,
\begin{align}
& \p_{\beta,h,p}  \left(N^{\frac{1}{2}}(\hat{\beta}_N-\beta) \leq t\right)\nonumber\\& =  \p_{\tilde{\beta}_p,0,p} \left(N^\frac{p}{2}\overline{X}^p_N \leq \e_{\tilde{\beta}_p+N^{-\frac{1}{2}} t, 0 ,p}  (N^\frac{p}{2}\overline{X}^p_N)\Big| \os \in A_1\right)\p_{\tilde{\beta}_p,0,p} (\os \in A_1)\label{m111}\\&+  \p_{\tilde{\beta}_p,0,p} \left(N^\frac{p}{2}\overline{X}^p_N \leq \e_{\tilde{\beta}_p+N^{-\frac{1}{2}} t, 0 ,p}  (N^\frac{p}{2}\overline{X}^p_N)\Big| \os \in A_1^c\right)\p_{\tilde{\beta}_p,0,p} (\os \in A_1^c)\label{m211}.
\end{align}
By Theorem \ref{cltun} under both $\p_{\tilde{\beta}_p,0,p} (~\cdot~\big|\os \in A_1)$ and $\p_{\tilde{\beta}_p + t/\sqrt N, 0,p} (~\cdot~\big|\os \in A_1)$, $N^\frac{p}{2} \overline{X}^p_N$ converges to $Z^p$ in distribution and in moments, where $Z \sim N(0,1)$. Consequently, $\e_{\tilde{\beta}_p+ t/\sqrt N, 0 ,p}  (N^\frac{p}{2}\overline{X}^p_N) \rightarrow 0$ by arguments similar to before, since $\p_{\tilde{\beta}_p+  t/\sqrt N, 0 ,p} (\os \in A_1^c)$ decays to $0$ exponentially fast. Hence, \eqref{m111} converges to $p_1/2$. Also, under $\p_{\tilde{\beta}_p,0,p} (~\cdot~\big| \os \in A_1^c )$, $N^\frac{p}{2}\overline{X}^p_N \xrightarrow{P} \infty$ and hence, \eqref{m211} converges to $0$. This shows that for all $t < 0$,
\begin{equation}\label{consprob}
\p_{\beta,h,p}  \left(N^{\frac{1}{2}}(\hat{\beta}_N-\beta) \leq t\right) \rightarrow \frac{p_1}{2}.
\end{equation}
Of course, \eqref{htl13} still remains valid. \eqref{eq:bmle_multimax_2} now follows from \eqref{htl13} and \eqref{consprob}.

Now, assume that $p \geq 4$ is even. If $h \neq 0$, then $K=2$. Also, $m_1<m_2<0$ if $h <0$ and $0<m_1<m_2$ if $h> 0$. Hence, $m_1^p < m_2^p$ if $h >0$ and $m_1^p > m_2^p$ if $h<0$. We can now use Lemma \ref{redunsup3} to derive \eqref{eq:bmle_multimax_3} and \eqref{eq:bmle_multimax_4}, and the proof is so similar to that for the $p \geq 3$ odd case, that we skip it. We now prove \eqref{eq:bmle_multimax_5} and \eqref{eq:bmle_multimax_6}. \\ 

\noindent\textit{Proof of \eqref{eq:bmle_multimax_5}:}~By Theorem \ref{cltun} and a standard binomial expansion (see \eqref{binoexp}), it follows that for any $\bar{\beta} \in \mathbb{R}$ and $i \in \{1,2\}$, under the conditional measure $\p_{\beta+N^{-\frac{1}{2}}\bar{\beta}, 0,p}\left(~\cdot~\big|\os \in A_i\right)$,  
\begin{equation}\label{iin12}
N^{\frac{1}{2}}(\overline{X}^p_N - m_*^p) \xrightarrow{D} N\left(-\frac{\bar{\beta}p^2 m_*^{2p-2}}{H''(m_*)}, -\frac{p^2 m_*^{2p-2}}{H''(m_*)}\right) \textrm{ and }  \e\left(N^{\frac{1}{2}}(\overline{X}^p_N - m_*^p)\right) \rightarrow -\frac{\bar{\beta}p^2m_*^{2p-2}}{H''(m_*)}.
\end{equation}
Since $A_1 \sqcup A_2 = [-1,1]$, \eqref{iin12} also holds under the unconditional measure $\p_{\beta+\bar{\beta}/\sqrt N, 0,p}$. The result in \eqref{eq:bmle_multimax_5} now follows easily, since for every $t \in \mathbb{R}$, we have
\begin{align*}
\p_{\beta,0,p}\left(N^{\frac{1}{2}}(\hat{\beta}_N - \beta) \leq t\right) &= \p_{\beta,0,p}\left(N^{\frac{1}{2}}(\overline{X}^p_N - m_*^p) \leq \e_{\beta+ N^{-\frac{1}{2}} t, 0,p} \left(N^{\frac{1}{2}}(\overline{X}^p_N - m_*^p) \right)\right)\nonumber\\&\rightarrow \p\left(N\left(0, -\frac{p^2 m_*^{2p-2}}{H''(m_*)}\right) \leq -\frac{tp^2m_*^{2p-2}}{H''(m_*)}\right)\nonumber\\&= \p\left( N\left(0, -\frac{H''(m_*)}{p^2 m_*^{2p-2}}\right) \leq t\right).
\end{align*}
\noindent \textit{Proof of \eqref{eq:bmle_multimax_6}:}~Fix $t < 0$. By \eqref{cc3} in Lemma \ref{redunsup3} and Theorem \ref{cltun}, we have
\begin{equation*}
\e_{\beta + N^{-\frac{1}{2}}t,0,p} \left(N^\frac{p}{2}\overline{X}^p_N\right) = \e_{\beta + N^{-\frac{1}{2}}t,0,p} \left(N^\frac{p}{2}\overline{X}^p_N\Big| \os \in A_2\right)(1-o(1)) + o(1) = \e Z^p + o(1).
\end{equation*}
This, together with the fact that $N^\frac{p}{2} \overline{X}^p_N \xrightarrow{P} \infty$ under $\p_{\beta,0,p}(~\cdot~\big|\os \in A_2^c)$, implies that
\begin{align}\label{neg}
\p_{\beta,0,p}\left(N^{\frac{1}{2}}(\hat{\beta}_N - \beta) \leq t\right) &= \p_{\beta,0,p}\left(N^\frac{p}{2} \overline{X}^p_N \leq \e Z^p + o(1)\Big|\os \in A_2\right)\p_{\beta,0,p}(\os \in A_2) + o(1)\nonumber\\&\rightarrow p_2\gamma_p.
\end{align}
Next, fix $t > 0$. Note that for any $\bar{\beta} \in \mathbb{R}$ and $i \in \{1,3\}$, we have the following under $\p_{\beta+\bar{\beta}/\sqrt N, 0,p}\left(~\cdot~\big|\os \in A_i\right)$,  
\begin{equation}\label{iin121}
N^{\frac{1}{2}}(\overline{X}^p_N - m_*^p) \xrightarrow{D} N\left(-\frac{\bar{\beta}p^2 m_*^{2p-2}}{H''(m_*)}, -\frac{p^2 m_*^{2p-2}}{H''(m_*)}\right) \textrm{ and } \e\left(N^{\frac{1}{2}}(\overline{X}^p_N - m_*^p)\right) \rightarrow -\frac{\bar{\beta}p^2m_*^{2p-2}}{H''(m_*)}.
\end{equation} 
By\eqref{cc2} in Lemma \ref{redunsup3}, $\p_{\beta+ t/\sqrt N, 0, p}(\os \in A_2) \leq Ce^{-DN^{\frac{1}{2}}}$ for some positive constants $C$ and $D$. It thus follows from the second convergence in \eqref{iin121}, that
\begin{equation}\label{expconv}
\e_{\beta+N^{-\frac{1}{2}}t, 0, p} \left(N^{\frac{1}{2}}(\overline{X}^p_N - m_*^p)\right) \rightarrow -\frac{tp^2m_*^{2p-2}}{H''(m_*)}.
\end{equation} 
Next, observe that $N^{\frac{1}{2}}(\overline{X}^p_N - m_*^p) \xrightarrow{P} - \infty$ under $\p_{\beta, 0,p}\left(~\cdot~\big|\os \in A_2\right)$. Combining this with \eqref{expconv} and using the fact that $p_1=p_3$, we have by the first convergence in \eqref{iin121},
\begin{align}\label{final}
&\p_{\beta,0,p}\left(N^{\frac{1}{2}}(\hat{\beta}_N - \beta)> t\right)\nonumber\\ 
& =  \p_{\beta,0,p} \left(N^{\frac{1}{2}}(\overline{X}^p_N - m_*^p) > -\frac{tp^2m_*^{2p-2}}{H''(m_*)} + o(1)\Bigg|\os \in A_1\right)\p_{\beta,0,p}(\os \in A_1)\nonumber\\
&+  \p_{\beta,0,p} \left(N^{\frac{1}{2}}(\overline{X}^p_N - m_*^p) > -\frac{tp^2m_*^{2p-2}}{H''(m_*)} + o(1)\Bigg|\os \in A_3\right)\p_{\beta,0,p}(\os \in A_3) + o(1) \nonumber\\
&\rightarrow 2p_1\p_{\beta,0,p}\left(N\left(0, -\frac{p^2 m_*^{2p-2}}{H''(m_*)}\right) > -\frac{tp^2m_*^{2p-2}}{H''(m_*)}\right)\nonumber\\& =  p_1 \p_{\beta,0,p}\left(N^+\left(0,-\frac{H''(m_*)}{p^2m_*^{2p-2}}\right) > t\right).
\end{align}
The result in \eqref{eq:bmle_multimax_6} now follows from \eqref{neg} and \eqref{final}.

%\iffalse
%\fi

\chapter{Inference in General Ising Models: the Maximum Pseudo-likelihood Method}\footnotetext{This chapter is a joint work with Jaesung Son and Bhaswar B. Bhattacharya}\label{ch:generalmple}
%\import{Chapters/}{intro}

 The increasing popularity of the Ising model as a foundational tool for understanding nearest-neighbor interactions in network data, has made it imperative to develop computationally tractable algorithms for learning the model parameters and understanding their rates of convergence (statistically efficiencies). In this chapter, we are interested in estimating the parameters of a general $p$-tensor model given a single sample of binary outcomes from an underlying network. This problem was classically studied in the $p=2$ case, when the  underlying network was a spatial lattice, where consistency and optimality of the maximum likelihood (ML) estimates were derived~\cite{comets_exp,gidas,guyon,discrete_mrf_pickard}. In Chapter \ref{curiech}, we discussed about ML estimation in the $p$-tensor case, when the underlying network was a complete hypergraph, i.e. all $p$-tuples of nodes interact with equal strength. However, as mentioned before, for general networks, parameter estimation using the  ML method turns out to be notoriously hard due to the appearance of an intractable normalizing constant in the likelihood. To circumvent this issue, Chatterjee \cite{chatterjee} proposed using the maximum pseudolikelihood (MPL) estimator \cite{besag_lattice,besag_nl}, which is a computationally efficient algorithm for estimating the parameters of a  Markov random field, that maximizes an approximation to the likelihood function (a `pseudo-likelihood') based on conditional distributions. This method and results in \cite{chatterjee} were later generalized in \cite{BM16} and \cite{pg_sm} to obtain rates of estimation for Ising models on general weighted graphs and joint estimation of parameters, respectively. These techniques were recently  used in Daskalakis et al. \cite{cd_ising_II,cd_ising_I} to obtain rates of convergence of the MPLE in general logistic regression models with dependent observations. Very recently, Dagan et al. \cite{cd_ising_estimation} considered the problem of parameter estimation in a more general model where the binary outcomes can be influenced by various underlying networks, and, as a consequence, improved some of the results in \cite{BM16}. Related problems in hypothesis testing given a single sample from the Ising model are considered in \cite{bresler_II,ising_testing,rm_sm}.

While the results above are promising, both from a practical and a mathematical standpoint, much is still left desired. For instance, in most real-world scenarios, the dependencies between nodes in a network are not consequences of just pairwise interactions, but arise due to peer-group effects. Higher-order relational data, which arise naturally in a variety of applications \cite{hypergraph_learning,hypergraph_applications,hypergraph_complex_data,hypergraph_image,hypergraph_multimedia,hypergraph_gene}, are generally  modeled using hypergraphs/tensors. In order to understand dependencies of binary variables in such datasets, it is natural to consider tensor Ising models, where the interaction matrix is replaced by a tensor (hypergraph) which encodes the strength of the interactions between, not just pairwise, but groups of individuals. 
%found several applications such as image retrieval \cite{hypergraph_image},  recommender systems \cite{hypergraph_multimedia} folksonomy \cite{hypergraph_applications}, bioinformatics \cite{hypergraph_network,hypergraph_gene}
%The $p$-tensor Ising model, which is a natural generalization of the matrix Ising model, provides a convenient mathematical framework for modeling dependencies in such datasets. 
%While the thermodynamic properties of these models have been studied in statistical physics under the purview of higher-order spin systems \cite{ferromagnetic_mean_field,ising_general,turban} for several years, parameter estimation for these distributions, arising from their applications in modeling  higher-order network data, has only been considered recently \cite{cd_ising_II}.  
The $p$-tensor Ising model provides a useful primitive for modeling such dependencies, where given a vector of binary outcomes $\bm X := (X_1,\ldots,X_N) \in \cC_N:=\{-1, 1\}^N$ and a $p$-tensor $\bm J_N = ((J_{i_1 \ldots i_p}))_{1 \leq i_1 \ldots i_p \leq N}$, encoding the strength of interactions between $p$-tuples of individuals, the joint distribution of $\bm X$ takes the following form: 
%the $p$-tensor (spin) Ising model, which is a discrete exponential family on $\mathcal{C}_N := \{-1,1\}^N$ with the following probability mass function: For $\bm X := (X_1,\ldots,X_N) \in \mathcal{C}_N$, 
\begin{equation}\label{model}
\p_{\beta,p}(\bm X) = \frac{1}{2^N Z_N(\beta,p)} e^{\beta H_N(\bm X)}  ,
\end{equation}
where the sufficient statistic (Hamiltonian)
\begin{equation}\label{eq:HN}
H_N(\bm X) := \sum_{1\leq i_1,\ldots,i_p \leq N} J_{i_1\ldots i_p} X_{i_1}\ldots X_{i_p}, 
\end{equation} 
and the parameter $\beta \geq 0$ (referred to as the {\it inverse temperature} in statistical physics)  measures the overall magnitude of dependency in the model across the tensor network (referred to as the `peer-group' effect in \cite{cd_ising_II}). The  {\it normalizing constant} $Z_N(\beta,p)$ (also referred to as the {\it partition function}) is determined by the condition $\sum_{\bm X \in \mathcal{C}_N} \p_{\beta,p}(\bm X) = 1$, that is,  
\begin{equation*}%\label{partition}
Z_N(\beta,p) = \frac{1}{2^N} \sum_{\bm X \in \mathcal{C}_N} \exp\left\{ \beta \sum_{1\leq i_1,\ldots,i_p \leq N} J_{i_1\ldots i_p} X_{i_1}\ldots X_{i_p} \right\}
\end{equation*} 
We will denote by $F_N(\beta,p) := \log Z_N(\beta,p)$ the {\it log-partition} function of the model. 
%The sufficient statistic $$H_N(\bm X) := \frac{1}{p!} \sum_{1\leq i_1,\ldots,i_p \leq N} J_{i_1\ldots i_p} X_{i_1}\ldots X_{i_p}$$ is called the \textit{Hamiltonian} of the model. 
%The $(i_1,\ldots,i_p)^{\textrm{th}}$ entry of the tensor $\bm J_N$ can be thought of as a measure of interaction between the observations $X_{i_1}, \ldots, X_{i_p}$, which is why $\bm J_N$ is also called the \textit{interaction tensor}. 
%We will always assume that the tensor $\bm J_N$ is non-zero, i.e. at least one of its entries is non-zero. 
Moreover, unless mentioned otherwise, we will assume that the tensor $\bm J_N$ satisfies the following two properties: 
\begin{itemize} 
	
	\item[(1)] The tensor $\bm J_N$ is {\it symmetric}, that is, $J_{i_1\ldots i_p} = J_{i_{\sigma(1)}\ldots i_{\sigma(p)}}$ for every $1\leq i_1<\cdots<i_p\leq N$ and every permutation $\sigma$ of $\{1,\ldots,p\}$,  and 
	
	\item[(2)] The tensor $\bm J_N$ has zeros on the `diagonals', that is, $J_{i_1\ldots i_p} = 0$, if $i_s = i_t$ for some $1\leq s<t\leq p$.
	
\end{itemize}

In this paper, we consider the problem of estimating the parameter $\beta$ given a single sample $\bm X = (X_1, X_2, \ldots, X_n)$ from the $p$-tensor Ising model \eqref{model}. Extending the results of Chatterjee \cite{chatterjee} on MPL estimation in matrix ($p=2$) Ising models, we obtain a general theorem which gives conditions under which the MPL estimate is $\sqrt N$-consistent in the $p$-tensor Ising model, for any $p \geq 3$.\footnote{A sequence of estimators $\{\hat{\beta}_N\}_{N\geq 1}$ is said to be {\it consistent} at $\beta$, if $\hat{\beta}_N \pto \beta$ under $\p_{\beta}$, that is, for every $M > 0$, $\p_{\beta}(|\hat \beta_N(\bm X) - \beta| \leq M ) \rightarrow 1$ as $N \rightarrow \infty$. Moreover, a sequence of estimators $\{\hat{\beta}_N\}_{N\geq 1}$ is said to be $\sqrt N$-{\it consistent} at $\beta$, if for every $\delta > 0$, there exists $M:=M(\delta, \beta) > 0$ such that $\p_{\beta}(\sqrt N |\hat \beta_N(\bm X) - \beta| \leq M ) > 1-\delta$, for all $N$.} The main bottleneck in extending the results from the matrix to the tensor case, is the lack of a natural spectral condition that is strong enough to control the fluctuations of the MPL function, but still verifiable in natural examples. To this end, we introduce the notion of a {\it local interaction matrix} which, given a configuration $\bm x \in \{-1, 1\}^n$,  measures the strength of the interaction between pairs of vertices (Definition \ref{defn:JN_interaction}). Our result shows that the MPL estimate is $\sqrt N$-consistent, whenever we have an appropriate moment bound on the local interaction matrix, and if the normalized log-partition function stays bounded away from zero (Theorem \ref{chextension}). We illustrate the robustness and generality of our result by verifying the conditions of the theorem in various commonly studied tensor Ising models. This includes the $\sqrt N$-consistency of the MPL estimate in the well-known $p$-spin Sherrington-Kirkpatrick (SK) model \cite{bovier,panchenko_book} (Corollary \ref{skthreshold}), and in Ising models on $p$-uniform hypergraphs under appropriate conditions on the adjacency tensors (Corollary \ref{boundeddeg}). The latter is also related to the recent work of Daskalakis et al. \cite{cd_ising_II}, where, as alluded to earlier, a general model  for logistic regression with dependent observations using higher-order Ising models was proposed, which includes as a special case the model in \eqref{model}. However, the conditions in \cite{cd_ising_II} are based directly on the interaction tensor, hence, cannot handle models where the rate of  estimation undergoes a phase transition. This is understandable because  \cite{cd_ising_II} considered the problem of jointly estimating multiple parameters in a more general model, hence, stronger assumptions were necessary for ensuring consistency. Our goal, on the other hand, is to pin down the precise conditions necessary for estimating the single parameter $\beta$ and develop methods for verifying those conditions in natural examples. To this end, our general theorem recovers as a corollary, the results in \cite{cd_ising_II} when specialized to the model \eqref{model}. 
%For instance, the conditions in \cite{cd_ising_II} are based directly on the interaction tensor, hence, cannot handle models where the rate of  estimation undergoes a phase transition, which happens whenever the underlying hypergraph becomes dense. 
More importantly, our results can handle models where the rate of estimation has phase transitions, which happens whenever the underlying hypergraph becomes dense. 
%, which to the best of our knowledge has not been done before for tensor Ising models. 
To illustrate this phenomenon we consider the Ising model on a hypergraph stochastic block model (HSBM), a natural generalization of the widely studied (graph) stochastic block model, that serves as a natural model for capturing higher-order relational data \cite{hypergraph_learning,hypergraph_image,hypergraph_multimedia,hypergraph_gene}.  In this case, we show there is a critical value $\beta_{\mathrm{HSBM}}^*$, such that if $\beta > \beta_{\mathrm{HSBM}}^*$ then the MPL estimate is $\sqrt N$ consistent, while if $\beta < \beta_{\mathrm{HSBM}}^*$ there is no consistent estimator for $\beta$ (Theorem \ref{sbmthr}). While it is relatively straightforward to show the $\sqrt N$-consistency of the MPL estimate above the threshold using our general theorem, proving that estimation is impossible below the threshold is more challenging. This is one of the technical highlights of the paper, which requires careful combinatorial estimates that go beyond the standard mean-field approximation techniques. Finally, we consider the special case of the $p$-tensor Curie-Weiss model, which is the Ising model on the complete $p$-uniform hypergraph. Here, using the special structure of the interaction tensor we are able to obtain the exact limiting distribution of the MPL estimate for all points above the estimation threshold (Theorem \ref{thm:cwmplclt}). In fact, it turns out that the  asymptotic variance of the MPL estimate saturates the Cramer-Rao lower bound, that is, the MPL estimate attains the best asymptotic variance among the class of consistent estimates. The formal statements of the results and their various consequences are given below in Section \ref{sec:statements}. 

\begin{rem}
	A related area of active research is the problem of structure learning in Ising models and Markov Random Fields. Here, one is given access to {\it multiple} i.i.d. samples from an Ising model, or a more general graphical model, and the goal is to estimate the underlying graph structure. Efficient algorithms and statistical lower bounds for this problem has been developed over the years under various structural assumptions on the underlying graph (cf.~\cite{structure_learning,discrete_tree,highdim_ising,graphical_models_binary,ising_nonconcave} and the references therein). Bresler \cite{bresler} made the first breakthrough for general bounded degree graphs, giving an efficient algorithm for structure learning, which required only logarithmic samples in the number of nodes of the graph.  This result has been subsequently generalized to Markov-random fields with higher-order interactions and alphabets with more than two elements (cf. \cite{graphical_models_algorithmic,multiplicative} and the references therein). The related problems of goodness-of-fit and independence testing given multiple samples from an Ising model has been studied in Daskalakis et al. \cite{cd_testing}. Recently, Neykov and Liu \cite{neykovliu_property} and Cao et al. \cite{high_tempferro} considered the problem of testing graph properties, such as connectivity and presence of cycles or cliques, using multiple samples from the Ising model on the underlying graph.

	All these results, however, are in contrast with the present thesis, where the underlying graph structure is assumed to be known and the goal is to  estimate the natural parameters given a {\it single} sample from the model. This is motivated by the applications described earlier, where it is more common to have access to only a single sample of node activities across the whole network, such as in disease modeling or social network interactions, where it is unrealistic, if not impossible, to generate many independent samples from the underlying model within a reasonable amount of time. 
\end{rem}

\section{Main Results}\label{sec:statements}
%\iffalse 
In this section we state our main results related to maximum pseudolikelihood estimation in general Ising models. The general result about the $\sqrt N$-consistency of the MPL estimate in tensor Ising models is discussed in Section \ref{sec:2.1}. Applications of this result to the  $p$-spin   SK model and Ising models on various hypergraphs are discussed in Section \ref{sec:2.2}.  Finally, in Section \ref{cltsec} we obtain the limiting distribution of the MPL estimate in the $p$-spin   Curie-Weiss model. Hereafter, we will often omit the dependence on $p$ and abbreviate $\p_{\beta,p}$, $Z_N(\beta,p)$, and $F_N(\beta,p)$ by $\p_\beta$, $Z_N(\beta)$, and $F_N(\beta)$, respectively, when there is no scope of confusion.

\subsection{Rate of Consistency of the MPL Estimator}\label{sec:2.1}

The maximum pseudo-likelihood (MPL) method, introduced by Besag \cite{besag_lattice,besag_nl}, provides a way to conveniently approximate the joint distribution of $\bm X \sim \p_{\beta, p}$ that avoids calculations with the normalizing constant. 

\begin{definition}{\em\cite{besag_lattice,besag_nl}} {\em Given a discrete random vector $\bm X= (X_1, X_2, \ldots, X_N)$ whose joint distribution is parameterized by a parameter $\beta \in \R$,  the MPL estimate of $\beta$ is defined as
		\begin{equation*}%\label{eq:mple_defn}
		\hat {\beta}_N(\bm X):=\arg\max_{\beta \in \R}\prod_{i=1}^N f_i(\beta, \bm X),
		\end{equation*}
		where $f_i(\beta, \bm X)$ is the conditional probability mass function of $X_i$ given $(X_j)_{j \ne i}$. }
\end{definition}

To compute the MPL estimate in the $p$-tensor Ising model \eqref{model}, fix $\beta > 0$ and consider $\bm X \sim \p_{\beta}$. Then from  \eqref{model}, the conditional distribution of $X_i$ given $(X_j)_{j\neq i}$ can be easily computed as: 
\begin{equation}\label{eq:conditional}
\p_{\beta}\left(X_i\big|(X_j)_{j\neq i}\right) = \frac{e^{p \beta  X_i m_i(\bm X) }}{ e^{p \beta  m_i(\bm X)  } + e^{-p \beta  m_i(\bm X)  } },
\end{equation}
where $m_i(\bm X) := \sum_{1\leq i_2,\ldots,i_p \leq N} J_{i i_2 \ldots i_p} X_{i_2}\cdots X_{i_p}$, is the local effect at the node $1 \leq i \leq N$ (often referred to as the local magnetization of the vertex $i$ in the statistical physics literature).  Then the pseudolikelihood estimate of $\beta$ (as defined in \eqref{model}) in the $p$-tensor Ising model \eqref{eq:conditional} is obtained by maximizing the function below, with respect to $b$,
\begin{align*}
L(b|\bm X) := \prod_{i=1}^N \p_{b}\left(X_i\big|(X_j)_{j\neq i}\right) & = \frac{1}{2^N} \exp\left\{\sum_{i=1}^N \left\{ p b   X_i m_i(\bm X) -  \log \cosh\left(p b  m_i(\bm X) \right) \right\} \right\} .
\end{align*}
Now, since $\log L(b|\bm X)$ is concave in $b$, the MPL estimator $\hat{\beta}_N(\bm X)$ can be obtained by solving the gradient equation $\frac{\partial \log L(b|\bm X)}{\partial b} = 0$, which simplifies to 
\begin{equation}\label{solution}
H_N(\bm X) -  \sum_{i=1}^N m_i(\bm X) \tanh\left(p b m_i(\bs)\right) = 0.
\end{equation} %\iffalse
To ensure well-definedness, in case \eqref{solution} does not have a solution or has more than one solution, the MPL estimate $\hat{\beta}_N(\bm X)$ is more formally defined as: 
\begin{equation}\label{mple}
\hat{\beta}_N(\bm X):= \inf\left\{b \geq 0 : H_N(\bm X) =  \sum_{i=1}^N m_i(\bm X) \tanh\left( p b m_i(\bs)\right)\right\},
\end{equation}
where the infimum of an empty-set is defined to be $+\infty$. Note that the expression in the RHS of the equality in \eqref{mple} is an increasing function of $t$, hence $\hat{\beta}_N(\bm X)$ can be very easily computed by the Newton-Raphson method or even a simple grid search.

Our first result is about the rate of consistency of the MPL estimate in general tensor Ising models. In particular,  we show in the proposition below that the MPL estimate $\hat \beta_N(\bm X)$, based on a single sample $\bm X\sim \p_{\beta}$ converges to the true parameter $\beta$ at rate $1/\sqrt N$, whenever the interaction tensor $\bm J_N$ satisfies a certain spectral-type condition and the log-partition function is $\Omega(N)$\footnote{For positive sequences $\{a_n\}_{n\geq 1}$ and $\{b_n\}_{n\geq 1}$, $a_n = O(b_n)$ means $a_n \leq C_1 b_n$, $a_n =\Omega(b_n)$ means $ a_n \geq C_2 b_n$, and  $a_n = \Theta(b_n)$ means $C_1 b_n \leq a_n \leq C_2 b_n$, for all $n$ large enough and positive constants $C_1, C_2$. Moreover, subscripts in the above notation,  for example $O_\square$, denote that the hidden constants may depend on the subscripted parameters.} at the true parameter value. To state our result formally, we need the following definition:

\begin{definition}\label{defn:JN_interaction} {\em Given a $p$-tensor $\bm J_N=((J_{i_1 i_2\ldots i_p}))_{1 \leq i_1, i_2, \ldots, i_p \leq N}$ and $\bm x = (x_1, x_2, \ldots, x_N) \in \cC_N$, define the {\it local interaction matrix} of $\bm J_N$ at the point $\bm x$ as the $N \times N$ matrix 
		$\bm J_N(\bt) := ((J_{i_1 i_2}(\bt)))_{1\leq i_1, i_2 \leq N}$, where the entries are given by: 
		\begin{align}\label{eq:Jx}
		J_{i_1 i_2}(\bt) := \sum_{1\leq i_3,\ldots,i_p \leq N}J_{i_1 i_2 i_3\ldots i_p}x_{i_3}\cdots x_{i_p}. 
		\end{align} 
		(Note that in the case $p=2$, $J_{i_1 i_2}(\bt) =J_{i_1 i_2}$, that is, the  local interaction matrix $\bm J_N(\bm x)$  is same as the interaction matrix $\bm J_N$, for all $\bm x \in \cC_N$.) }
\end{definition}

We are now ready to state our result on the convergence rate of the MPL estimate in a tensor Ising model.\footnote{For a vector $\bm v \in \R^N$, $\|\bm v\|$ will denote the Euclidean norm of $\bm v$. Moreover, for a $N \times N$ matrix $A$, $\|A\| := \sup_{\|\bm x\|=1} \|A\bm x\|$ denotes the operator norm of $A$.} 

\begin{thm}\label{chextension}
	Fix $p \geq 2$, $\beta > 0$ and a sequence of $p$-tensors $\{\bm J_N\}_{N \geq 1}$ such that the following two conditions hold: 
	\begin{enumerate}
		\item[$(1)$] $\sup_{N\geq 1}\mathbb{E}_{\beta}[ \|\bm J_N(\bm Z)\|^4] <\infty$, where the expectation is taken with respect to $\bm Z \sim \p_{\beta}$, 
		\item[$(2)$] $\liminf_{N\rightarrow \infty} \frac{1}{N}F_N(\beta) > 0$.
	\end{enumerate}
	Then given a single sample $\bm X$ from the model \eqref{model} with interaction tensor $\bm J_N$, the MPL estimate $\hat{\beta}_N(\bm X)$, as defined in \eqref{mple}, is $\sqrt N$-consistent for $\beta$, that is, for every $\delta > 0$, there exists $M:=M(\delta, \beta) > 0$ such that $$\p_{\beta}(\sqrt N |\hat \beta_N(\bm X) - \beta| \leq M ) > 1-\delta,$$ for all $N$ large enough.  
\end{thm}

The proof of this theorem is given in Section \ref{sec:3}. The proof has two main steps: In the first step we use the method of exchangeable pairs to show that the derivative of the log-pseudolikelihood (the LHS of \eqref{solution}) is concentrated around zero at the true model parameter (see Lemma \ref{lm:boundsecondmoment} for details). The proof adapts the method of exchangeable pairs introduced in \cite{chatterjee} where a similar result was proved for matrix (2-spin) Ising models. The main technical challenge as one goes from the matrix to the tensor case, is the absence of a natural spectral condition in tensor models. To this end, we introduce condition (1), which requires that the fourth-moment of the spectral norm of the local interaction matrix is uniformly bounded. This condition allows us to prove the desired concentration of the  log-pseudolikelihood, and, as we will see below, can be easily verified for a large class of natural tensor models. The second step in the proof of Theorem \ref{chextension} is to show that the log-pseudolikelihood is strongly concave, that is, its second derivative is strictly negative with high probability. Here, we use condition (2) to first show that the Hamiltonian is $\Omega(N)$ with high-probability, which then implies the strong concavity of the log-pseudolikelihood by a truncated second-moment argument.\footnote{Recalling the discussion in Definition \ref{defn:JN_interaction}, note that when $p=2$, condition $(1)$ simplifies to $\sup_{N \geq 1}||\bm J_N|| < \infty$, hence Theorem \ref{chextension} recovers Chatterjee's result on $\sqrt N$-consistency of MPL estimates in 2-spin Ising models \cite[Theorem 1.1]{chatterjee}.}

\begin{rem}\label{remark:condition} {\em The $L_4$-condition (condition (1)) on the local interaction matrix in Theorem \ref{chextension} can be replaced by the following stronger $L_\infty$-condition, which is often easier to verify in examples: 
		\begin{equation}\label{eq:JNcondition}
		\sup_{N\geq 1}\sup_{\bt \in \sa_N} \|\bm J_N(\bt)\| <\infty.
		\end{equation} 
		Condition \eqref{eq:JNcondition}, hence condition $(1)$ in Theorem \ref{chextension}, is also weaker than the `bounded-degree condition': 
		\begin{align}\label{eq:JNcondition_II}
		\sup_{1 \leq i_1 \leq N} \sum_{1 \leq i_2, i_3, \ldots, i_p \leq N} |J_{i_1 i_2 i_3 \ldots i_p} | =O(1). 
		\end{align} 
		In particular, condition \eqref{eq:JNcondition} allows us to handle the $p$-spin Sherrington-Kirkpatrick model, an example where the bounded-degree condition \eqref{eq:JNcondition_II} fails to hold. } 
\end{rem}

%\iffalse

\subsection{Applications}\label{sec:2.2}

In this section we discuss the consequences of Theorem \ref{chextension} to the $p$-spin SK model (Section \ref{sec:sk}), spin systems of on general hypergraphs (Section \ref{nonstoch}), and the hypergraph stochastic block model (Section \ref{stoch}).

\subsubsection{The $p$-Spin Sherrington-Kirkpatrick Model}
\label{sec:sk}

In the $p$-spin Sherrington-Kirkpatrick (SK) model \cite{bovier}, the interaction tensor is of the form  
\begin{align}\label{eq:JN_sk}
J_{i_1\ldots i_p} = N^{\frac{1-p}{2}} g_{i_1\ldots i_p}, 
\end{align}
where  $(g_{i_1\ldots i_p})_{1\leq i_1<\ldots<i_p< \infty}$ is a fixed realization of a collection of independent standard Gaussian random variables, and $g_{i_1\ldots i_p} = g_{\sigma(i_1)\ldots \sigma(i_p)}$, for any permutation $\sigma$ of $\{1, 2, \ldots, p\}$. This is a canonical example of a spin glass model which has remarkable thermodynamic properties \cite{spinglass_book}. A whole new discipline has emerged from the study of this object, with many beautiful theorems that have unearthed deep connections between diverse areas in mathematics and statistical physics  (cf.~\cite{bovier,sk_optimization,panchenko_book,talagrand_sk} and the references therein). The problem of parameter estimation in the SK model was initiated by Chatterjee \cite{chatterjee}, where $\sqrt N$-consistent of the MPL estimate for all $\beta > 0$ was proved for the 2-spin SK model. The following corollary extends this to all $p \geq 3$.

\begin{cor}\label{skthreshold} In the $p$-spin SK model, the MPL estimate $\hat{\beta}_N(\bm X)$  is $\sqrt{N}$-consistent for all $\beta> 0$. 
\end{cor}

The proof of this result is given in Section \ref{skproof}. In this case, condition (2) in Theorem \ref{chextension} can be easily verified using monotonicity and the well-known asymptotics of $F_N(\beta)$ in the high-temperature (small $\beta$) regime: In particular, we know from \cite[Theorem 1.1]{bovier} that, almost surely, $\lim_{N \rightarrow \infty}\frac{1}{N}F_N(\beta) = \tfrac{\beta^2}{2}$, for $\beta >0$ small enough. Hence, by the monotonicity of $F_N(\beta)$, we have $\lim_{N \rightarrow \infty}\frac{1}{N}F_N(\beta) > 0 $ for all $\beta > 0$, which establishes (2).  However, unlike when $p=2$, verifying condition $(1)$ in Theorem \ref{chextension} when $p \geq 3$ requires more work.\footnote{Note that when $p=2$, $\bm J_N$ is a Wigner matrix, and hence, by \cite[Theorem 2.12]{spectral_norm}  $\sup_{N \geq 1} ||\bm J_N|| < \infty$, thus verifying condition $(1)$ of Theorem \ref{chextension}.} To this end, note that for $p \geq 3$ and every fixed $\bm x \in \cC_N$, the local interaction matrix $\bm J_N(\bm x)$ is a Gaussian random matrix, but the elements are now dependent because of the symmetry of the tensor $\bm J_N$. This dependence, however, is relatively weak and using standard Gaussian process machinery we can show the validity of \eqref{eq:JNcondition}, and, hence, that of condition $(1)$ in Theorem \ref{chextension}.

\subsubsection{Ising Models on Hypergraphs}\label{nonstoch}

The $p$-tensor model \eqref{model} can be interpreted as a spin system on a weighted $p$-uniform hypergraph, where the entries of the tensor correspond to the weights of the hyperedges. More precisely, given a symmetric tensor $\bm J_N = ((J_{i_1 i_2 \ldots i_p}))_{1 \leq i_1, i_2, \ldots, i_p \leq N}$, construct a weighted $p$-uniform hypergraph $H_N$ with vertex set $[N]:=\{1, 2, \ldots, N\}$ and edge weights $w(\bm e)= J_{i_1 i_2 \ldots i_p}$, for $\bm e =(i_1, i_2, \ldots, i_p) \in [N]_p$.\footnote{For the set $[N]=\{1, 2, \ldots, N\}$, $[N]^p$ denotes the $p$-fold Cartesian product $[N]\times [N] \times \cdots \times [N]$, and $[N]_p$ is the collection of $p$-tuples in $[N]^p$ with distinct entries.} The model \eqref{model} is then a spin system on $H_N$ where the Hamiltonian \eqref{eq:HN} can be rewritten as 
$$H_N(\bm X) = \sum_{\bm e \in {[N]_p}} w(\bm e) X_{\bm e},$$
where $\bm X= (X_1, X_2, \ldots, X_N) \in \cC_N$ and $X_{\bm e} =X_{i_1} X_{i_2} \ldots X_{i_p}$, for $\bm e = (i_1, i_2, \ldots, i_p)$. For a tensor $\bm J_N=((J_{i_1 i_2 \ldots i_N}))$, define the (weighted) degree of the vertex $i_1$ as $$d_{\bm J_N}({i_1}):= \frac{1}{(p-1)!}\sum_{1\leq i_2, i_3,\ldots,i_p \leq N}|J_{i_1 i_2 i_3\ldots i_p}|,$$ 
which is the sum of the absolute values weights of the hyperedges passing through the vertex $i_1$. Similarly, define the weighted {\it co-degree} of the vertices $i_1, i_2$ as 
\begin{align}\label{eq:dJN}
d_{\bm J_N}(i_1, i_2):=\frac{1}{(p-2)!}\sum_{1\leq i_3,\ldots,i_p \leq N} |J_{i_1 i_2 i_3\ldots i_p}|, 
\end{align}
which is the sum of the absolute values of weights of the hyperedges incident on both $i_1$ and $i_2$. Denote by $\bm D_{\bm J_N} = ((d_{\bm J_N}(i_1, i_2)))_{1 \leq i_1, i_2 \leq N}$, the {\it co-degree matrix} corresponding to the tensor $\bm J_N$. The following corollary provides useful sufficient conditions under which the MPL estimate is $\sqrt N$-consistent at all   temperatures. The proof is given in Section \ref{sec:boundeddegpf}.

\begin{cor}\label{boundeddeg} Suppose $\{\bm J_N\}_{N \geq 1}$ is a sequence of $p$-tensors such that the following two conditions hold: 
	\begin{enumerate}
		\item[$(1)$] $\sup_{N\geq 1} \|\bm D_{\bm J_N} \| <\infty$,
		\item[$(2)$] $\liminf_{N\rightarrow \infty} \frac{1}{N} \sum_{1 \leq i_1 < i_2 < \ldots < i_p \leq N} J_{i_1 i_2 \ldots i_p}^2> 0$.
	\end{enumerate}
	Then the MPL estimate $\hat{\beta}_N(\bm X)$ is $\sqrt N$-consistent for all $\beta > 0$. 
\end{cor}

%\iffalse
\begin{rem}\label{remark:condition_II} 
	{\em Note that, since the $L_2$-operator norm of a symmetric matrix is bounded by its $L_\infty$-operator norm,\footnote{For two sequences $a_n$ and $b_n$, $a_n \lesssim_{\Box} b_n$ means that there exists a positive constant $C(\Box)$ depending only on the subscripted parameters $\Box$, such that $a_n \leq C(\Box) b_n$ for all $n$ large enough.}     
		\begin{align}\label{eq:DN_bound}
		\|\bm D_{\bm J_N} \| \leq \max_{1 \leq i_1 \leq N} \sum_{i_2=1}^N d_{\bm J_N}(i_1, i_2) & = \frac{1}{(p-2)!} \max_{1 \leq i_1 \leq N} \sum_{1 \leq i_2, i_3, \ldots, i_p \leq N} |J_{i_1 i_2 \ldots i_p}| \nonumber \\   
		& \lesssim_p \max_{1 \leq i_1 \leq N} d_{\bm J_N}(i_1), 
		\end{align}
		that is, if a tensor has bounded maximum degree, then condition $(1)$ of Theorem \ref{chextension} holds. This shows that Corollary \ref{boundeddeg} recovers the general theorem of \cite{cd_ising_II}, where $\sqrt N$-consistency of the MPL was proved, albeit for a more general model, under condition (2) and condition $(1)$ replaced by the bounded degree assumption $\max_{1 \leq i_1 \leq N} d_{i_1} = O(1)$. }   
\end{rem}

As mentioned earlier, the conditions in Corollary \ref{boundeddeg}, neither of which depend on the true parameter $\beta$, cannot hold for hypergraphs where the rate of estimation undergoes a phase transition. In fact, as explained in Remark \ref{remark:hypergraph}, the scope of this corollary is really only restricted to Ising models on hypergraphs which are sparse. The importance of the second condition in  Theorem \ref{chextension} becomes evident when the  hypergraph becomes dense, where $F_N(\beta)$ ceases to be $\Omega(N)$ for all $\beta$, and the rate of estimation changes as $\beta$ varies. This is illustrated in Section \ref{stoch} below, where the exact location of the phase transition is derived for Ising models on block hypergraphs.

\begin{rem}\label{remark:hypergraph} {\em Suppose $H_N=(V(H_N), E(H_N))$ is a sequence of unweighted $p$-uniform hypergraphs with vertex set $V(H_N)=[N] = \{1, 2, \ldots, N\}$ and edge set $E(H_N)$, with no isolated vertex. Denote by $\bm A_{H_N} = ((a_{i_1 i_2 \ldots i_p}))_{1 \leq i_1, i_2, \ldots, i_p \leq N}$ the adjacency tensor of $H_N$, that is, $a_{i_1 i_2 \ldots i_p} =1$ if $(i_1, i_2, \ldots, i_p) \in E(H_N)$ and zero otherwise. Then in order to ensure that a $p$-spin-system on $H_N$, as in \eqref{model}, has a non-trivial scaling limit, one needs to consider the scaled tensor,  
		$$\bm J_{H_N}= \frac{N}{|E(H_N)|} \bm A_{H_N}.$$
		In this case, the Frobenius norm condition in Corollary \ref{boundeddeg} simplifies to,   
		\begin{align}\label{eq:JHN_condition}
		\frac{1}{N} ||\bm J_{H_N}||_F^2 = \frac{N}{|E(H_N)|^2} \sum_{1 \leq i_1, i_2, \ldots, i_p \leq N} a_{i_1 i_2 
			\ldots i_p} = \Theta\left(\frac{N}{|E(H_N)|} \right) = \Omega(1).
		\end{align}
		This implies, $|E(H_N)| = \Theta(N)$, since $H_N$ has no isolated vertex. Moreover, condition $(1)$ can be written as,  
		\begin{align}\label{eq:DHN}
		||\bm D_{\bm A_{H_N} }|| = O\left(\frac{|E(H_N)|}{N} \right). 
		\end{align}
		Therefore, combining \eqref{eq:JHN_condition}, \eqref{eq:DHN}, and Corollary \ref{boundeddeg}, shows that for any sequence of (unweighted) $p$-uniform hypergraphs $H_N= (V(H_N), E(H_N))$, such that $||\bm D_{\bm A_{H_N} }|| = O(1)$ and $|E(H_N)| = O(N)$, the MPL estimate $\hat{\beta}_N(\bm X)$ in the Ising model \eqref{model} with interaction tensor $\bm J_{H_N}$,  is $\sqrt N$-consistent for all $\beta > 0$.  In particular, by the bound in \eqref{eq:DN_bound} applied to the adjacency tensor $\bm A_{H_N}$, the MPL estimate $\hat{\beta}_N(\bm X)$ is $\sqrt N$-consistent for all $\beta > 0$, whenever $H_N$ has bounded maximum degree and $O(N)$ edges. } 
\end{rem}

\subsubsection{Hypergraph Stochastic Block Models}\label{stoch}

The hypergraph stochastic block model (HSBM) is a random hypergraph model where each hyperedge is present independently with probability depending on the membership of the vertices to various blocks (see \cite{hypergraph_block,hypergraph_clustering_II,hypergraph_phase_transition} and the references therein for more on the HSBM and its applications in higher-order community detection).

\begin{definition}\label{defn:block} {\em (Hypergraph Stochastic Block Model) Fix $p \geq 2$, $K\geq 1$,  a vector of community proportions $\bm \lambda := (\lambda_1,\ldots,\lambda_K)$, such that $\sum_{j=1}^K \lambda_j = 1$, and a symmetric probability tensor $\bm \Theta := ((\theta_{j_1\ldots j_p}))_{1\leq j_1,\ldots,j_p\leq K}$,  where $\theta_{j_1\ldots j_p} \in [0, 1]$, for $1\leq i_1,\ldots, i_p\leq K$. The hypergraph stochastic block model with proportion vector $\bm \lambda$ and probability tensor $\bm \Theta$ is a $p$-uniform hypergraph $H_N$ on $[N]=\{1, 2, \ldots, N\}$ vertices with adjacency tensor $\bm A_{H_N}=((a_{i_1 i_2 \ldots i_p}))_{1 \leq i_1, i_2, \ldots, i_p \leq N}$, where 
		$$a_{i_1\ldots i_p} \sim \mathrm{Ber}\left(\theta_{j_1\ldots j_p}\right)\quad\textrm{for } i_1<\ldots<i_p \textrm{ and }  (i_1,\ldots,i_p) \in \cB_{j_1}\times\cdots\times \cB_{j_p},$$ where $\cB_j := (N \sum_{i=1}^{j-1}\lambda_i, N \sum_{i=1}^j \lambda_i] \bigcap [N]$, for $j \in \{1,\ldots,K\}$, and $\{a_{i_1\ldots i_p}\}_{1\leq i_1<\ldots<i_p \leq |V|}$ are independent. We denote this model by $\cH_{p, K, N} (\bm \lambda, \bm \Theta)$ and a realization from this model as $H_N \sim \cH_{p, K, N} (\bm \lambda, \bm \Theta)$. }
\end{definition}

In this section, we consider the problem of parameter estimation given a sample from an Ising model on a HSBM.  The following theorem shows that for the $p$-tensor Ising models on a HSBM, there is a critical value of $\beta$, below which estimation is impossible, and above which the MPL estimate is $\sqrt N$-consistent. The location of the phase transition is determined by the first time the maximum of a certain variational problem, which arises from the mean-field approximation of the partition function,  becomes non-zero. More formally, this is defined as,   
\begin{equation}\label{eq:beta_threshold}
\beta_{\mathrm{HSBM}}^* :=  \sup\left\{\beta \geq 0: \sup_{(t_1,\ldots,t_K)\in [0,1]^K} \phi_\beta(t_1,\ldots,t_K) = 0\right\},
\end{equation}
where the function 
$\phi_\beta: [-1,1]^K \mapsto \mathbb{R}$ is: 
\begin{align}\label{eq:threshold_function}
\phi_{\beta}(t_1, t_2, \ldots, t_K) := \beta \left(\sum_{1\leq j_1,\ldots,j_p\leq K} \theta_{j_1\ldots j_p} \prod_{\ell=1}^p \lambda_{j_\ell} t_{j_\ell} \right) - \sum_{j=1}^K \lambda_j I(t_j),
\end{align}
and $I(t) := \frac{1}{2}\left\{(1+t)\log(1+t) + (1-t)\log(1-t)\right\}$ is the binary entropy function. 

%\iffalse

\begin{thm}\label{sbmthr}  Fix $p \geq 2$ and a realization of a HSBM $H_N \sim \cH_{p, K, N}(\bm \lambda, \bm \Theta)$ on $N$ vertices, where $\bm \lambda$ is a proportion vector and $\bm \Theta$ is a symmetric probability tensor as in Definition \ref{defn:block}.  Then given a sample $\bm X \sim \p_{\beta}$ from the model \eqref{model}, with adjacency tensor $\bm J_N= \frac{1}{N^{p-1}} \bm A_{H_N}$, the following hold: 
	\begin{itemize}
		
		\item[(1)] The MPL estimate $\hat{\beta}_N(\bm X)$ is $\sqrt{N}$-consistent for $\beta > \beta_{\mathrm{HSBM}}^*$. 
		
		\item[(2)] There does not exist any consistent sequence of estimators for any $\beta < \beta_{\mathrm{HSBM}}^*$.
		
	\end{itemize}
\end{thm}

The proof of the above result is given in Section \ref{proof6}. To show the result in $(1)$ we verify the conditions of Theorem \ref{chextension}. Here, we invoke the standard mean-field lower bound to the Gibbs variational representation of the partition function \cite{CD16}, from which it can be easily verified that $F_N(\beta) = \Omega(N)$, whenever $\beta > \beta_{\mathrm{HSBM}}^*$. Perhaps the more interesting consequence of Theorem \ref{sbmthr} is the result in (2), which shows that not only is the MPL estimate not $\sqrt N$-consistent below the threshold, no estimator is consistent in this regime, let alone $\sqrt N$-consistent. The main argument in this proof is to show that 
\begin{align}\label{eq:Fbeta_bounded_I}
F_N(\beta) = O(1), \quad \text{ for } \beta < \beta_{\mathrm{HSBM}}^*. 
\end{align} 
Once this is proved, then it can be easily verified that the Kullback-Leibler (KL) divergence between the measures $\p_{\beta_1, p}$ and $\p_{\beta_2, p}$, for any two $0 < \beta_1 < \beta_2 < \beta_{\mathrm{HSBM}}^*$ remains bounded, which in turn implies that the measures $\p_{\beta_1, p}$ and $\p_{\beta_2, p}$ are untestable, and hence inestimable. The main technical difficulty in proving an estimate like \eqref{eq:Fbeta_bounded_I} in tensor models, is the absence of `Gaussian' techniques \cite{BM16,comets}, which allows one to compare the partition function of Ising models with quadratic Hamiltonians with an appropriately chosen Gaussian model. This method, unfortunately, does not apply when $p \geq 3$, hence, to estimate the partition function we take the following more direct approach: We first consider the {\it averaged model} where the interaction tensor is replaced by the expected interaction tensor $\tilde{\bm J}_N := \e {\bm J_N}$. Using the block structure of the tensor $\tilde{\bm J}_N$ the Hamiltonian in the averaged model can be written in terms of the average of the spins in the different blocks, and hence, the partition function in the averaged model can be accurately estimated using bare-hands combinatorics (Lemma \ref{lm:logpartition_I}). We then move from the averaged model to the actual model using standard concentration arguments (Lemma \ref{st2}).

\begin{rem} {\em Using the machinery of non-linear large deviations developed in \cite{CD16}, we can in fact show that for the HSBM,  
		\begin{align}\label{Fbeta_variational_problem}
		\lim_{N \rightarrow \infty}\frac{1}{N} F_N(\beta) = \sup_{(t_1,\ldots,t_K)\in [0,1]^K} \phi_\beta(t_1,\ldots,t_K) , 
		\end{align}
		with probability 1.  Although the proof of this result has not been included in the paper, because for proving Theorem \ref{sbmthr} $(1)$ we only need to establish a lower bound on $\frac{1}{N} F_N(\beta)$, this is worth mentioning as it motivates the definition of the threshold $\beta_{\mathrm{HSBM}}^*$ and corroborates the result in Theorem \ref{sbmthr} (1).  The result in \eqref{Fbeta_variational_problem} is, however, not strong enough to show that estimation is impossible below the threshold $\beta_{\mathrm{HSBM}}^*$. Here, we need to understand the asymptotic behavior of $F_N(\beta)$ itself (without scaling by $N$), which is a more delicate matter that require arguments beyond the purview of non-linear large deviations and mean-field approximations, as discussed above. In this case, the proof of Theorem \ref{sbmthr} (2) shows that whenever the log-partition function is $o(N)$, which happens when $\beta < \beta_{\mathrm{HSBM}}^*$, it is actually $O(1)$, and hence,  there is a sharp transition from inestimability to $\sqrt N$-consistency. } 
\end{rem}

An important special case of the HSBM is the Erd\H{o}s-R\'enyi random hypergraph model, where every hyperedge is present independently with the same fixed probability. 
%\iffalse
\begin{example}\label{example:random_hypergraph} {\em (Erd\H{o}s-R\'enyi random hypergraphs) The HSBM reduces to the classical Erd\H{o}s-R\'enyi random $p$-hypergraph model when the number of blocks $K=1$. In this case, each hyperedge is present independently with probability $\theta \in (0, 1]$, and the variational problem \eqref{eq:beta_threshold} for the threshold simplifies to 
		\begin{equation}\label{hypergraph_random}
		\beta_{\mathrm{ER}}^*(p, \theta) :=  \sup\left\{\beta \geq 0: \sup_{t \in [0,1]} \left\{ \beta \theta t^p - I(t) \right\} = 0\right\}.
		\end{equation} 
		We will denote this hypergraph model by $\mathscr{G}_p(N, \theta)$. In this case, Theorem \ref{sbmthr} gives the following:
		\begin{itemize}
			
			\item In the  Erd\H{o}s-R\'enyi random $p$-hypergraph model $\mathscr{G}_p(N, \theta)$, the MPL estimate $\hat{\beta}_N(\bm X)$ is $\sqrt{N}$-consistent for all $\beta > \beta_{\mathrm{ER}}^*(p, \theta)$. 
			
			\item On the other hand, there does not exist any consistent sequence of estimators for any $\beta < \beta_{\mathrm{ER}}^*(p, \theta)$.  
		\end{itemize} 
		Note that by the change of variable $\kappa = \beta \theta$, it follows that $\beta_{\mathrm{ER}}^*(p, \theta) = \beta_{\mathrm{ER}}^*(p, 1)/\theta$. A simple analysis shows $\beta_{\mathrm{ER}}^*(2, 1)=0.5$, and hence, $\beta_{\mathrm{ER}}^*(2, \theta)  = 0.5/\theta$. For higher values of $p$, $\beta_{\mathrm{ER}}^*(p, 1)$ can be easily computed numerically. In particular, we have $\beta_{\mathrm{ER}}^*(3, 1) \approx 0.672$ and $\beta_{\mathrm{ER}}^*(4, 1) \approx 0.689$. In fact, $\beta_{\mathrm{ER}}^*(p, 1)$ is strictly increasing in $p$ and $\lim_{p \rightarrow \infty} \beta_{\mathrm{ER}}^*(p, 1) = \log 2$ (see Appendix \ref{cwtp} for a proof). } 
\end{example}

Another example is that of random $p$-partite $p$-uniform hypergraphs, which are natural extensions of random bipartite graphs.

\begin{example}[Random $p$-partite $p$-uniform hypergraphs] {\em A $p$-uniform hypergraph is said to be  $p$-partite if the vertex set of the hypergraph can be partitioned into $p$-nonempty sets in such a way that every edge intersects every set of the partition in exactly one vertex.  A random $p$-partite $p$-uniform hypergraph, is a  $p$-partite $p$-uniform hypergraph where each edge is present independently with some fixed probability $\theta \in (0, 1]$ \cite{multipartite_random_hypergraph}. More formally, given a vector $\bm N = (N_1, N_2, \ldots, N_p)$ of positive integers, such that $\sum_{j=1}^p N_j = N$ and $\theta \in (0, 1]$, in the random $p$-partite $p$-uniform hypergraph $\cH_p(\bm N, \theta)$, the vertex set  $[N]=\{1, 2, \ldots, N\}$ is partitioned into $p$ disjoint sets $S_1,\ldots, S_p$, such that $|S_j| = N_j$ for $1 \leq j \leq p$, and each edge $\bm e \in V_1 \times V_2 \times \cdots \times V_p$ is present independently with probability $\theta$.  If $\bm N$ is such that $\frac{1}{N} \bm N \rightarrow \bm \lambda = (\lambda_1, \lambda_2, \ldots, \lambda_p)$, as $N \rightarrow \infty$, then this is a special case of the hypergraph stochastic block model and the threshold \eqref{eq:beta_threshold} simplifies to,  
		\begin{equation}\label{rurpthr}
		\beta_{\mathrm{partite}}^*(p, \bm \lambda, \theta) :=  \sup\left\{\beta \geq 0: \sup_{(t_1,\ldots,t_p)\in [0,1]^p}   \left\{ \beta \theta \prod_{j=1}^p \lambda_j t_j - \sum_{j=1}^p \lambda_j I(t_j) \right\}  = 0\right\}.
		\end{equation}
		Theorem \ref{sbmthr} then implies that the MPL estimate is $\sqrt{N}$-consistent for all $\beta > \beta_{\mathrm{partite}}^*(p, \bm \lambda, \theta) $, and consistent estimation is impossible for $\beta < \beta_{\mathrm{partite}}^*(p, \bm \lambda, \theta) $. In case the $p$ partitioning sets have asymptotically equal size, that is, $\lambda_j = \frac{1}{p}$ for all $1 \leq j \leq p$, the threshold in \eqref{rurpthr} simplifies further to: 
		\begin{align}\label{rurpthr_equal}
		\beta_{\mathrm{equipartite}}^*(p,  \theta) :=  \sup\left\{\beta \geq 0: \sup_{(t_1,\ldots,t_p)\in [0,1]^p}   \left\{ \beta \theta p^{-p} \prod_{j=1}^p  t_j - \frac{1}{p} \sum_{j=1}^p I(t_j) \right\}  = 0\right\}.
		\end{align}
		Now, a simple analysis shows that $\beta_{\mathrm{equipartite}}^*(p,  \theta)=p^p\beta_{\mathrm{ER}}^*(p,  \theta)$. The upper bound $$\beta_{\mathrm{equipartite}}^*(p,  \theta) \leq p^p\beta_{\mathrm{ER}}^*(p,  \theta)$$ follows by substituting $t_1=t_2  \cdots = t_p = t \in [0, 1]$ in \eqref{rurpthr_equal} and relating it to \eqref{hypergraph_random}. For the lower bound, note by the convexity of the function $I(x)$ and the AM-GM inequality, that  
		$$\beta \theta p^{-p} \prod_{j=1}^p  t_j - \frac{1}{p}\sum_{j=1}^p I(t_j) \leq \beta \theta p^{-p}\left(\frac{1}{p}\sum_{j=1}^p t_j \right)^p - I\left(\frac{1}{p}\sum_{j=1}^p t_j\right).$$
		Then, by the change of variable $\kappa=\beta p^{-p}$, it follows that $\beta_{\mathrm{equipartite}}^*(p,  \theta) \geq p^p\beta_{\mathrm{ER}}^*(p,  \theta)$. } 
\end{example}

\subsection{Precise Fluctuations in the Curie-Weiss Model}\label{cltsec}

The $p$-tensor Curie-Weiss model is the Ising model on the complete $p$-uniform hypergraph,\footnote{In the complete $p$-uniform hypergraph with vertex set $[N]=\{1, 2, \ldots, N\}$ the set of hyperedges is the collection of all the $p$-element subsets of $[N]$.} where all the $p$-tuples of interactions are present \cite{ferromagnetic_mean_field}. In other words, this is the Ising model on the Erd\H os-R\'enyi $p$-hypergraph with $\theta=1$. Denoting $\beta_{\mathrm{CW}}^*(p) := \beta_{\mathrm{ER}}^*(p, 1)$, we know from the discussion in Example \ref{example:random_hypergraph}, that for  $\beta< \beta_{\mathrm{CW}}^*(p)$ consistent estimation is impossible, while for  $\beta > \beta_{\mathrm{CW}}^*(p)$ the MPL estimate $\hat{\beta}_N(\bm X)$  is $\sqrt{N}$-consistent. Given that we know the rate of consistency, the next natural question is to wonder whether anything can be said about the limiting distribution of the MPL estimate above the threshold. While tackling this question appears to be extremely difficult, if not impossible, for general models, the special structure of the Curie-Weiss model allows us to say much more. This begins with the observation that in the  Curie-Weiss model  the MPL estimate can be written as a function of the sample mean $\bar X_N = \frac{1}{N} \sum_{i=1}^N X_i$. Then combining the recent results on the   asymptotic distribution of $\bar X_N$ \cite{mlepaper} and the delta theorem, we can get the precise fluctuations of the MPL estimate at all points above the estimation threshold $\beta_{\mathrm{CW}}^*(p)$. This is formalized in the theorem below:

\begin{thm}\label{thm:cwmplclt} 
	Fix $p \geq 2$ and consider the $p$-spin Curie-Weiss model with interaction tensor $\bm J_N = ((J_{i_1 \ldots i_p}))_{1 \leq i_1, \ldots, i_p \leq N}$, where $J_{i_1  \ldots i_p} = \frac{1}{N^{p-1}}$, for all $1 \leq i_1, \ldots, i_p \leq N$. Then for every $\beta > \beta_{\mathrm{CW}}^*(p)$, as $N \rightarrow \infty$,
	\begin{equation}\label{statement1}
	\sqrt{N}(\hat{\beta}_N(\bm X) -\beta) \xrightarrow{D} N\left(0, -\frac{g''(m_*)}{p^2m_*^{2p-2}}\right),
	\end{equation}
	where $g(t) := \beta t^p - I(t)$, for $t \in [-1, 1]$, and $m_*=m_*(\beta, p)$ is the unique positive global maximizer of $g$.
\end{thm}

\begin{figure*}\vspace{-0.15in}
	\centering
	\begin{minipage}[l]{1.0\textwidth}
		\centering
		\includegraphics[height=2.65in,width=4.5in]
		{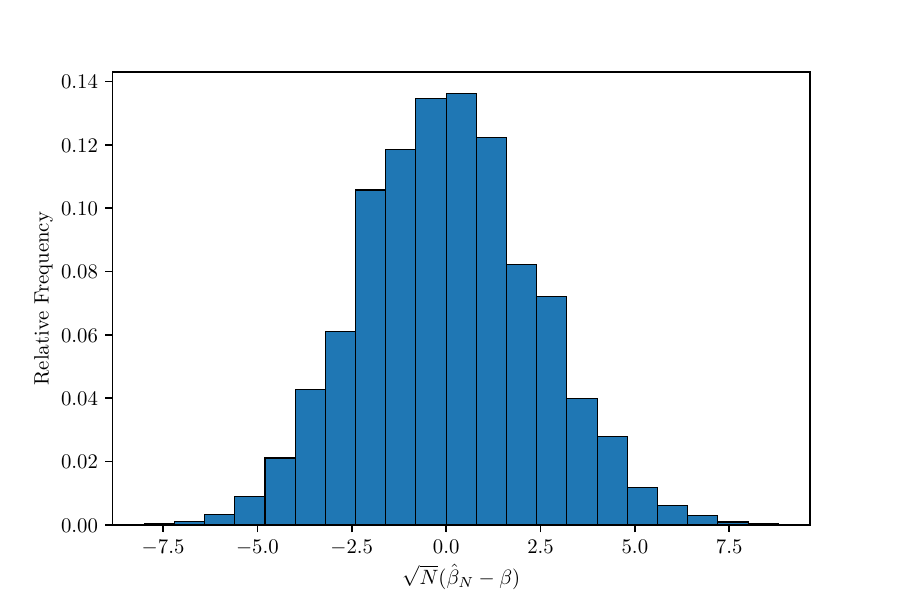}\\
		%\scriptsize{(a)}
	\end{minipage}
	\caption{\small{The histogram $\sqrt N (\hat{\beta}_N(\bm X)-\beta)$ in the 4-tensor Curie-Weiss model at $\beta= 0.75 > \beta_{\mathrm{CW}}^*(4) \approx 0.689$ (above the estimation threshold).}}
	\label{fig:histogram_I}
\end{figure*}

The proof of this result is given in Section \ref{sec:pf_cwmplclt}. Figure \ref{fig:histogram_I} shows the histogram (over $10^6$ replications) of $\sqrt N (\hat{\beta}_N(\bm X)-\beta)$ with $p=4$, $\beta= 0.75$, and $N=20000$. As predicted by the result above, we see a limiting Gaussian distribution, since $\beta= 0.75 >\beta_{\mathrm{CW}}^*(4) \approx 0.689$ is above the estimation threshold.

%%%%%%%%%%%%%%%%%%%%%%%%%%%%%%%%%%%%%%%%%%%%%%%%%%%%%%%%%%%%%%%%%%%%%%%%%%%%%%%
%%%%%%%%%%%%%%%%%%%%%%%%%%%%%%%%%%%%%%%%%%%%%%%%%%%%%%%%%%%%%%%%%%%%%%%%%%%%%%%%%%%%%%%%%%%%%%%%%%%%%%%%%

The result in Theorem \ref{thm:cwmplclt} can be used to construct a confidence interval for the parameter $\beta$ for all points above the estimation threshold. Towards this, note, by \cite[Theorem 2.1]{mlepaper}, that $|\bar X_N| \pto m_*$ under $P_{\beta, p}$, when $\beta > \beta_{\mathrm{CW}}^*(p)$. The result in \eqref{statement1} then implies that 
\begin{equation*}%\label{eq:cibeta1}
\left(\hat{\beta}_N(\bm X) - \frac{|\bar X_N|^{1-p}}{p}\sqrt{\frac{-g''(|\bar X_N|)}{N}} z_{1-\frac{\alpha}{2}},~\hat{\beta}_N(\bm X) + \frac{|\bar X_N|^{1-p}}{p}\sqrt{\frac{-g''(|\bar X_N|)}{N}} z_{1-\frac{\alpha}{2}}\right),
\end{equation*}   
is an interval which contains $\beta$ with asymptotic coverage probability $1-\alpha$, whenever $\beta > \beta_{\mathrm{CW}}^*(p)$.\footnote{For $\alpha \in (0, 1)$, $z_\alpha$ is the $\alpha$-th quantile of the standard normal distribution, that is, $\p_\beta(N(0, 1) \leq z_\alpha) = \alpha$.}
%\iffalse

\begin{rem}\label{remark:information} {\em (Efficiency of the MPL estimate) An interesting consequence of Theorem \ref{thm:cwmplclt} is that the limiting variance in \eqref{statement1} saturates the Cramer-Rao (information) lower bound of the model, when $\beta > \bw$. To see this, note that the (scaled) Fisher information in the model \eqref{model} (recall that the Cramer-Rao lower bound is the inverse of the Fisher information) is given by,  
		$$I_N(\beta) = \frac{1}{N} \e_{\beta}\left[\left(\frac{\mathrm d}{\mathrm d \beta} \log \p_{\beta}(\bs)\right)^2\right] = \mathrm{Var}_\beta (N^\frac{1}{2} \bar{X}_N^p) \rightarrow -\frac{p^2 m_*^{2p-2}}{g''(m_*)},$$ as $N \rightarrow \infty$,  where the last step follows from the asymptotics of $\bar X_N$ derived in \cite{mlepaper}. 
		This implies, for $\beta > \bw$, the MPL estimate $\hat \beta_N(\bm X)$ is {\it asymptotically efficient}, which means that no other consistent estimator can have lower asymptotic mean squared error than $\hat \beta_N(\bm X)$ above the estimation threshold. While this has been shown for the maximum likelihood (ML) estimate \cite{comets,mlepaper}, that the MPL estimate, which only maximizes an approximation to the true likelihood, also has this property, is particularly encouraging, as it showcases the effectiveness of the MPL method, both computationally as well as in terms of statistical efficiency. }   
\end{rem}

The results above show that the MPL estimate is $\sqrt N$-consistent and asymptotic efficient whenever $\beta >  \beta_{\mathrm{CW}}^*(p)$. On the other hand, for $\beta <  \beta_{\mathrm{CW}}^*(p)$, we know from Theorem \ref{sbmthr} that consistent estimation is impossible. In particular, this means that  the MPL estimate is inconsistent for $\beta <  \beta_{\mathrm{CW}}^*(p)$. Therefore, the only case that remains is at the threshold  $\beta = \beta_{\mathrm{CW}}^*(p)$. Here, the situation is much more delicate. We address this case in the theorem below, which shows that the MPL is $\sqrt N$-consistent for $p=2$ (with a non-Gaussian limiting distribution), but inconsistent for $p \geq 3$.

\begin{thm}\label{thm:cw_threshold}(Asymptotics of the MPL estimate at the threshold) 
	Fix $p \geq 2$ and consider the $p$-spin Curie-Weiss model with interaction tensor $\bm J_N = ((J_{i_1 \ldots i_p}))_{1 \leq i_1, \ldots, i_p \leq N}$, where $J_{i_1  \ldots i_p} = \frac{1}{N^{p-1}}$, for all $1 \leq i_1, \ldots, i_p \leq N$. Suppose $\beta = \beta_{\mathrm{CW}}^*(p)$. Denote by $m_*=m_*(\beta, p) \in (0, 1)$ the unique positive maximizer of the function $g:= \beta t^p - I(t)$, for $t \in [-1, 1]$, and define  
	\begin{align}\label{eq:proportion}
	\alpha := 
	\begin{cases}
	\frac{1}{1+2[(m_*^2-1)g''(m_*)]^{-\frac{1}{2}}} &\quad\text{if}~p~\textrm{is even},\\
	\frac{1}{1+[(m_*^2-1)g''(m_*)]^{-\frac{1}{2}}} &\quad\text{if}~p~\textrm{is odd}.\\
	\end{cases}
	\end{align}
	Then, the following hold as $N \rightarrow \infty$, 
	
	\begin{itemize}
		
		\item[(1)]  If $p = 2$ $($recall $\beta_{\mathrm{CW}}^*(2)=\frac{1}{2}$$)$, then for every $t\in \mathbb{R}$, 
		\begin{align}\label{eq:threshold_F}
		\lim_{N\rightarrow \infty} \p_\beta\left (N^\frac{1}{2}\left(\hat{\beta}_N - \frac{1}{2}\right) \leq t\right) = 
		\begin{cases}
		F(\sqrt{6t}) - F(-\sqrt{6t}) &\quad\text{if}~t\geq 0\\
		0 &\quad\text{if}~t < 0\\
		\end{cases}
		\end{align}
		where $F$ is a probability distribution function with density given by $d F(t) \propto \exp\left(-\frac{t^4}{12}\right) \mathrm dt$.

		\item[(2)] If $p \geq 3$, then  
		\begin{align}\label{eq:threshold_mple}
		\sqrt{N}(\hat{\beta}_N(\bm X) -\beta) \xrightarrow{D} (1-\alpha) N\left(0, -\frac{g''(m_*)}{p^2m_*^{2p-2}}\right) + \alpha \delta_\infty, 
		\end{align}
		where $g(\cdot)$ is as defined Theorem \ref{thm:cwmplclt} and $\delta_\infty$ denotes the point mass at $\infty$.

		\item[(3)] Moreover, at a finer scaling, the following hold: 
		\begin{itemize}
			\item[(a)] If $p \geq 4$ is even, then 
			\begin{align}\label{eq:threshold_mple_I}			
			N^{1-\frac{p}{2}} \hat{\beta}_N  \xrightarrow{D} \alpha \left(\frac{1}{pZ^{p-2}}\right) + (1-\alpha)\delta_0,
			\end{align}
			where $Z \sim N(0,1)$. 
			
			\item[(b)] If $p \geq 3$ is odd, then 
			\begin{align}\label{eq:threshold_mple_II}			
			N^{1-\frac{p}{2}} \hat{\beta}_N  \xrightarrow{D} \frac{\alpha}{2} \left(\frac{1}{p |Z|^{p-2}}\right) +\frac{\alpha}{2}\delta_\infty + (1-\alpha)\delta_0. 
			\end{align}
		\end{itemize}
	\end{itemize} 
\end{thm}
%\iffalse

The proof of this result is given in Section \ref{sec:pf_cwmplclt}. As in the proof of Theorem \ref{thm:cwmplclt}, the main ingredient in the proof of the above result is the asymptotic distribution of the sample mean at the threshold derived in \cite{comets,mlepaper}. The reason there is a change in the consistency rates of the MPLE as one moves from the 2-spin model to the $p$-spin model, for $p \geq 3$, is because the rate of convergence of the sample mean $\cs$ in the Curie-Weiss model depends on the value of $p$ at the threshold. More precisely, for $p=2$ and $\beta=\beta_{\mathrm{CW}}^*(2)=\frac{1}{2}$, $N^{\frac{1}{4}} \cs \dto F$, where $F$ is as defined in Theorem \ref{thm:cw_threshold} (1) (see \cite[Proposition 4.1]{comets}).  On the other hand, when  $p \geq 3$ and $\beta=\beta_{\mathrm{CW}}^*(p)$, $N^{\frac{1}{2}} \cs$ converges to a mixture of point masses with two or three components depending on whether $p$ is odd or even, respectively (see \cite[Theorem 1.1]{mlepaper}).

%%%%%%%%%%%%%%%%%%%%%%%%%%%%%%%%%%%%%%%%%%%%%%%%%%%%%%%%%%%%%%%%%%%%%%%%%%%%%%%
\begin{figure*}[h]\vspace{-0.15in}
	\centering
	\begin{minipage}[l]{0.49\textwidth}
		\centering
		\includegraphics[width=3.3in]
		{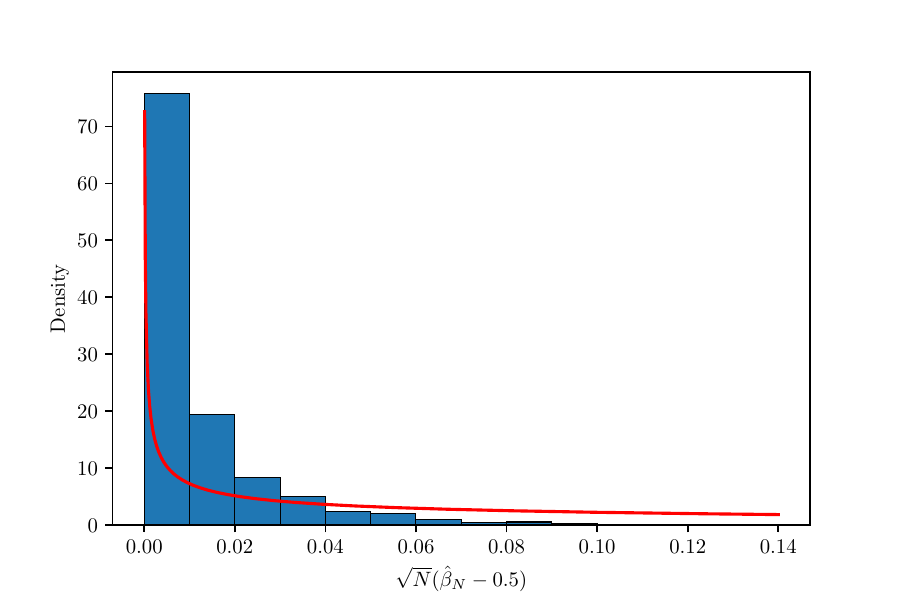}\\
		\small{(a)}
	\end{minipage} 
	\begin{minipage}[l]{0.49\textwidth}
		\centering
		\includegraphics[width=3.3in]
		{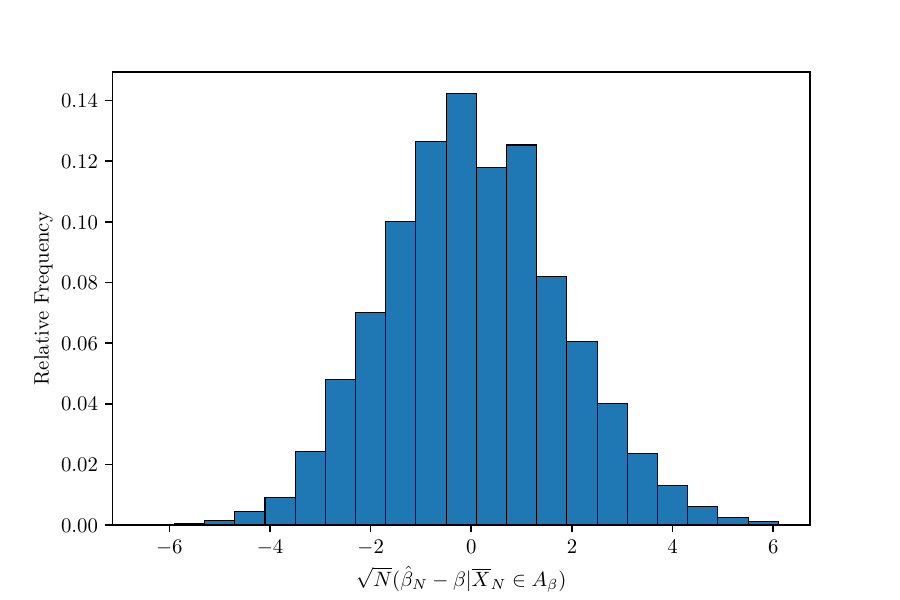}\\ 
		\small{(b)}
	\end{minipage} 
	\caption{\small{(a) The histogram of  $\sqrt N (\hat{\beta}_N(\bm X)-\beta)$  in the 2-tensor Curie-Weiss model at the estimation threshold ($\beta=\frac{1}{2}$) and the limiting density function (in red); and (b) the histogram of the conditional distribution $\sqrt N (\hat{\beta}_N(\bm X)-\beta)|\{\cs \in A_\beta\}$, where $A_\beta$ is the interval $[-1, 1]$ minus a small neighborhood around zero, in the 4-tensor Curie-Weiss model at the estimation threshold, which has a limiting normal distribution.}}
	\label{fig:histogram_II}
\end{figure*}
%%%%%%%%%%%%%%%%%%%%%%%%%%%%%%%%%%%%%%%%%%%%%%%%%%%%%%%%%%%%%%%%%%%%%%%%%%%%%%%%%%%%%%%%%%%%%%%%%%%%%%%%%

Taking derivatives in \eqref{eq:threshold_F} shows that for $p=2$ the MPL estimate has a limiting Gamma distribution with density $g(a) \propto \frac{1}{\sqrt a }e^{-3 a^2} \mathrm da$.  Figure \ref{fig:histogram_II}~(a) shows the histogram of the quantity $\sqrt N (\hat{\beta}_N(\bm X)-\beta)$ for $p=2$ and $\beta= \beta_{\mathrm{CW}}^*(2)= 0.5$, and the limiting density function (plotted in red). On the other hand, for $p \geq 3$, Theorem \ref{thm:cw_threshold} (3)
shows that the MPL estimate is inconsistent at the threshold (in fact, $\hat \beta_N(\bm X) \pto \infty$, for $p \geq 3$ and  $\beta=\beta_{\mathrm{CW}}^*(p)$). However, even though for $p \geq 3$ the MPL estimate is inconsistent when $\beta=\beta_{\mathrm{CW}}^*(p)$, Theorem \ref{thm:cw_threshold} (2) shows $\sqrt N (\hat{\beta}_N(\bm X)-\beta)$ has a Gaussian limit with probability $1-\alpha$,  that is, MPL estimate is $\sqrt N$-consistent at this point with probability $1-\alpha$. In fact, the proof of Theorem \ref{thm:cw_threshold} (2) shows that $\hat{\beta}_N(\bm X)$ is  not $\sqrt N$-consistent at the threshold for $p \geq 3$, only when $\cs$ is close to zero. More precisely, the proof shows that $\sqrt N (\hat{\beta}_N(\bm X)-\beta)|\{\cs \in A_\beta \} \dto N(0, -\frac{g''(m_*)}{p^2m_*^{2p-2}})$, if $A_\beta=[-1, 1] \backslash B_0$, where $B_0$ is a small neighborhood of zero. This is illustrated in Figure \ref{fig:histogram_II}~(b) which plots the histogram of this conditional distribution for $p=4$ and $\beta=0.6888 \approx \beta_{\mathrm{CW}}^*(4)$ .

\subsection{Organization}
The rest of the paper is organized as follows. In Section \ref{millionsong}, we demonstrate through a real data analysis, a scenario where the classical $2$-spin Ising model is not a good fit, and one needs to consider higher order Ising models. In Section \ref{sec:3} we prove Theorem \ref{sec:3}. The proofs of Corollary \ref{skthreshold} and Corollary \ref{boundeddeg} are given in Section \ref{skproof} and Section \ref{sec:boundeddegpf}, respectively. The proof of Theorem \ref{sbmthr} is given in Section \ref{proof6}. The proofs of Theorem \ref{thm:cwmplclt} and Theorem \ref{thm:cw_threshold} are given in Section \ref{sec:pfcwclt}. Additional properties of the Curie-Weiss threshold are given in Appendix \ref{cwtp}.

\section{The Last.fm Dataset}\label{millionsong}
The Last.fm dataset (\url{http://millionsongdataset.com/lastfm/}), a part of the Million Song Dataset (\url{http://millionsongdataset.com/}) contains a list of $1892$ users, their friendship network, and a list of their most favorite artists (see \cite{lastfm,kostisconc}). We wish to investigate if users' preference for a particular artist depends only on pairwise interactions in the user friendship network, or if it is affected by peer group effects. To formulate this precisely, for each artist, we form a vector $\bm X := (X_1,\ldots,X_N)$ where $N$ is the total number of users, and $X_i = +1$ if user $i$ has that artist in his favorite list, and $X_i = -1$ otherwise. We are interested in testing whether the vector $\bm X$ follows a $2$-spin Ising model or not.

We chose four of the most popular artists (and bands) from the dataset, namely Lady Gaga, Britney Spears, Rihanna and the Beatles, and for each of them, implemented the following procedure. Assuming the true model to be a $2$-spin Ising model (without external magnetic field) on the user friendship network, we estimated the parameter $\beta$ from the data $\bm X$, using the MPLE $\hat{\beta}$. We then simulated $100$ observations $\bm X^{(1)},\ldots, \bm X^{(100)}$ from the $2$-spin Ising model on the user friendship network, with parameter $\hat{\beta}$. We decided to accept the null hypothesis of a $2$-spin Ising model if and only if the actual value of the sufficient statistic $H(\bm X) := \sum_{i\sim j} X_i X_j$ \footnote{Here, $i\sim j$ denotes that users $i$ and $j$ are friends.} lies within the $2.5^{\mathrm{th}}$ and $97.5^{\mathrm{th}}$ percentiles of the empirical distribution of $H(\bm X^{(1)}),\ldots,H(\bm X^{(100)})$.

\begin{figure}[h]
	\begin{minipage}[c]{0.5\textwidth}
		\centering
		\includegraphics[width=3.35in]
		{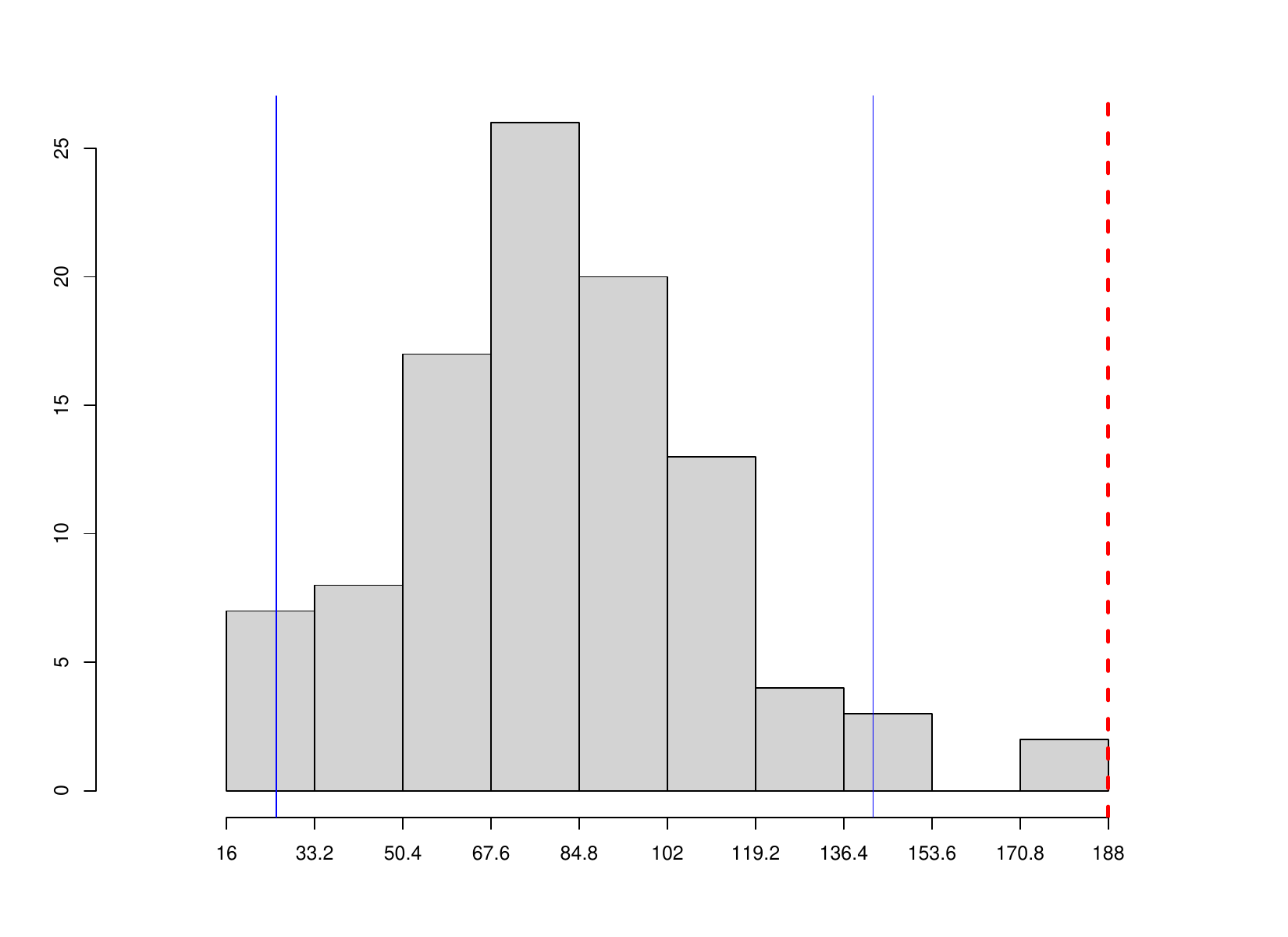}\\
		\small{(a)}
	\end{minipage}
	\begin{minipage}[l]{0.5\textwidth} 
		\centering
		\includegraphics[width=3.35in]
		{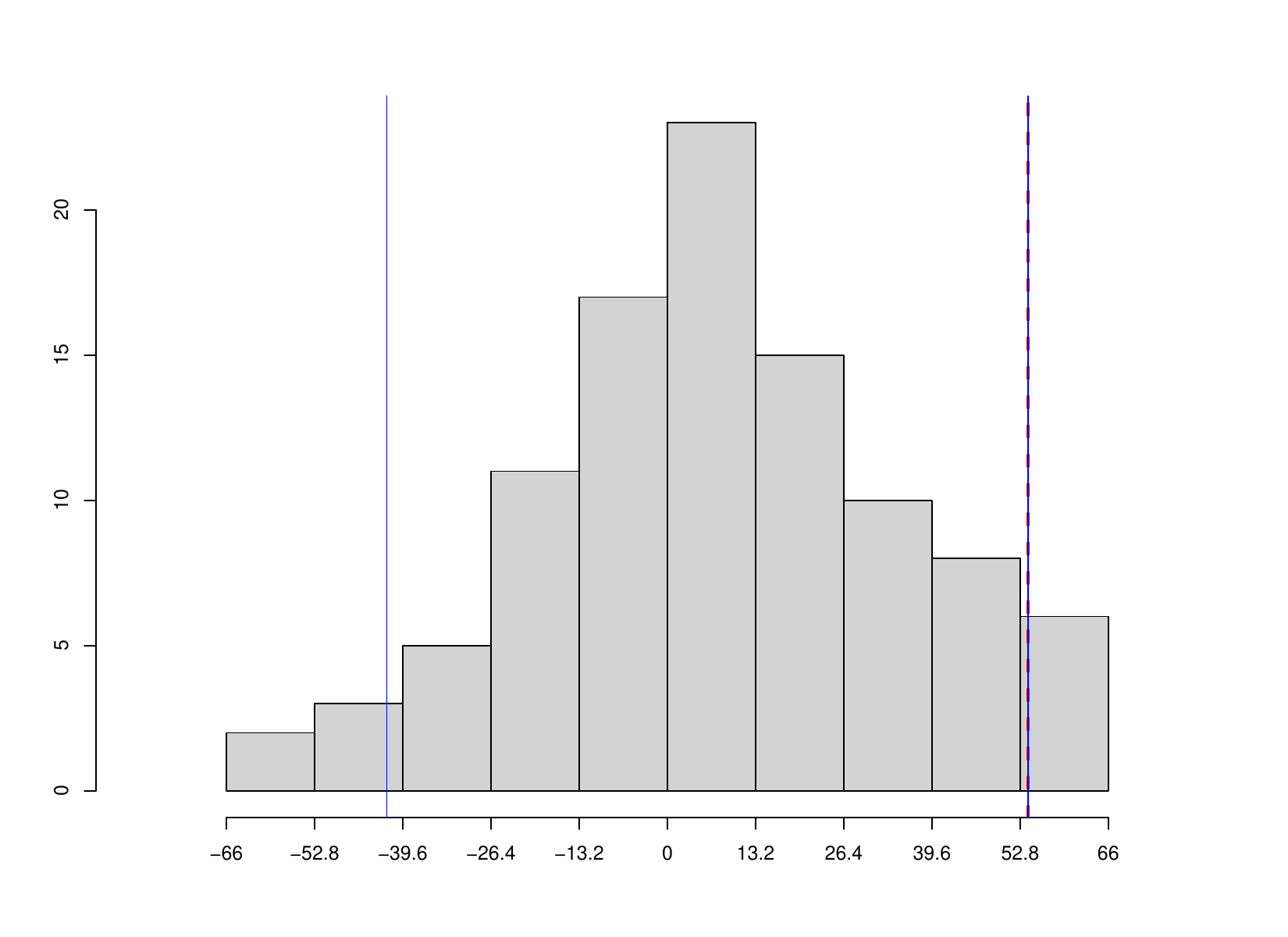}\\
		\small{(b)} 
	\end{minipage}
	\caption{\small{(a) $2$-spin Ising model fit, and (b) $3$-spin Ising model fit on the user preference vector for Lady Gaga.}}
	\label{fig:test1}
\end{figure}

\begin{figure}[h]
	\begin{minipage}[c]{0.5\textwidth}
		\centering
		\includegraphics[width=3.35in]
		{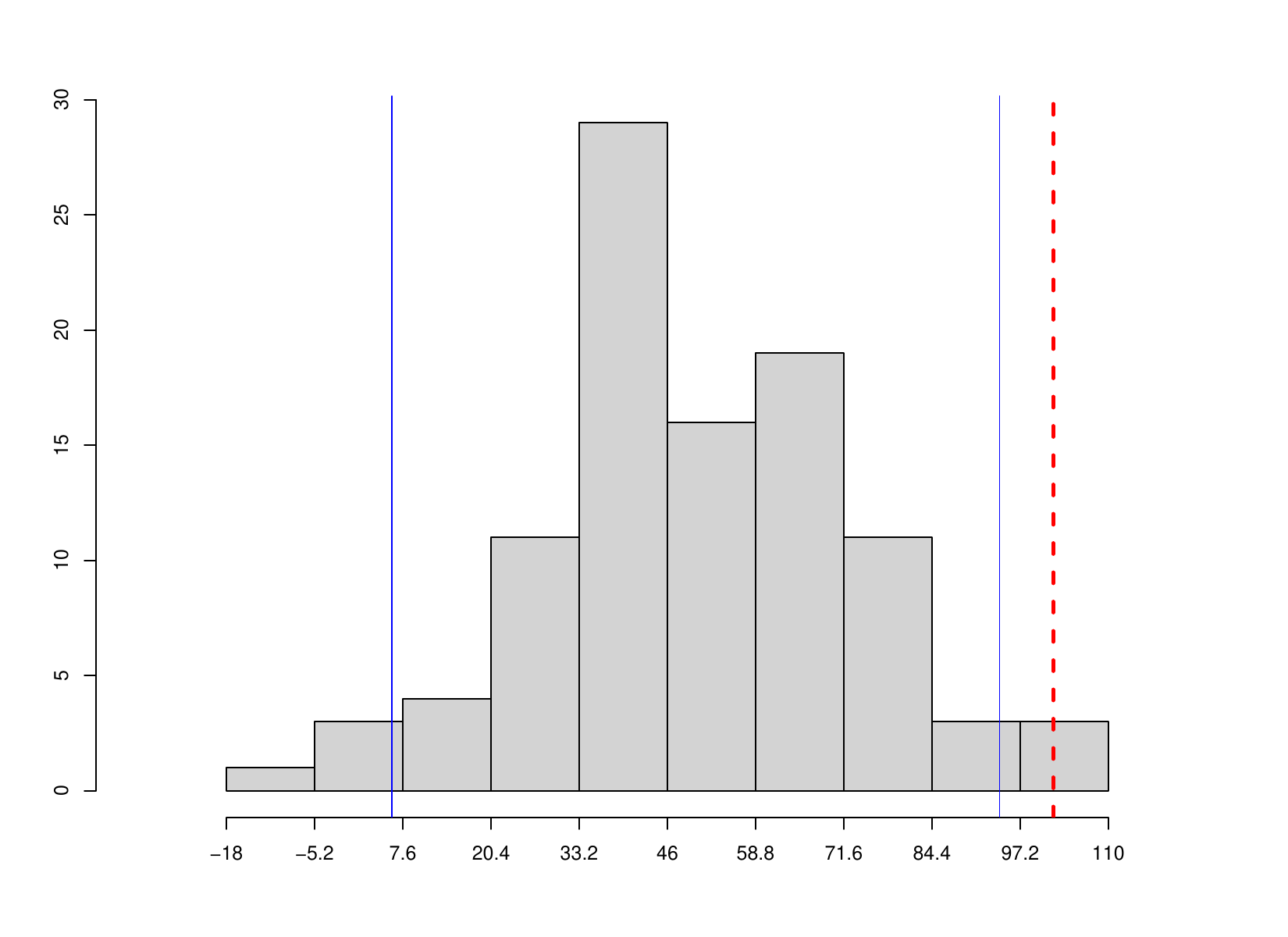}\\
		\small{(a)}
	\end{minipage}
	\begin{minipage}[l]{0.5\textwidth} 
		\centering
		\includegraphics[width=3.35in]
		{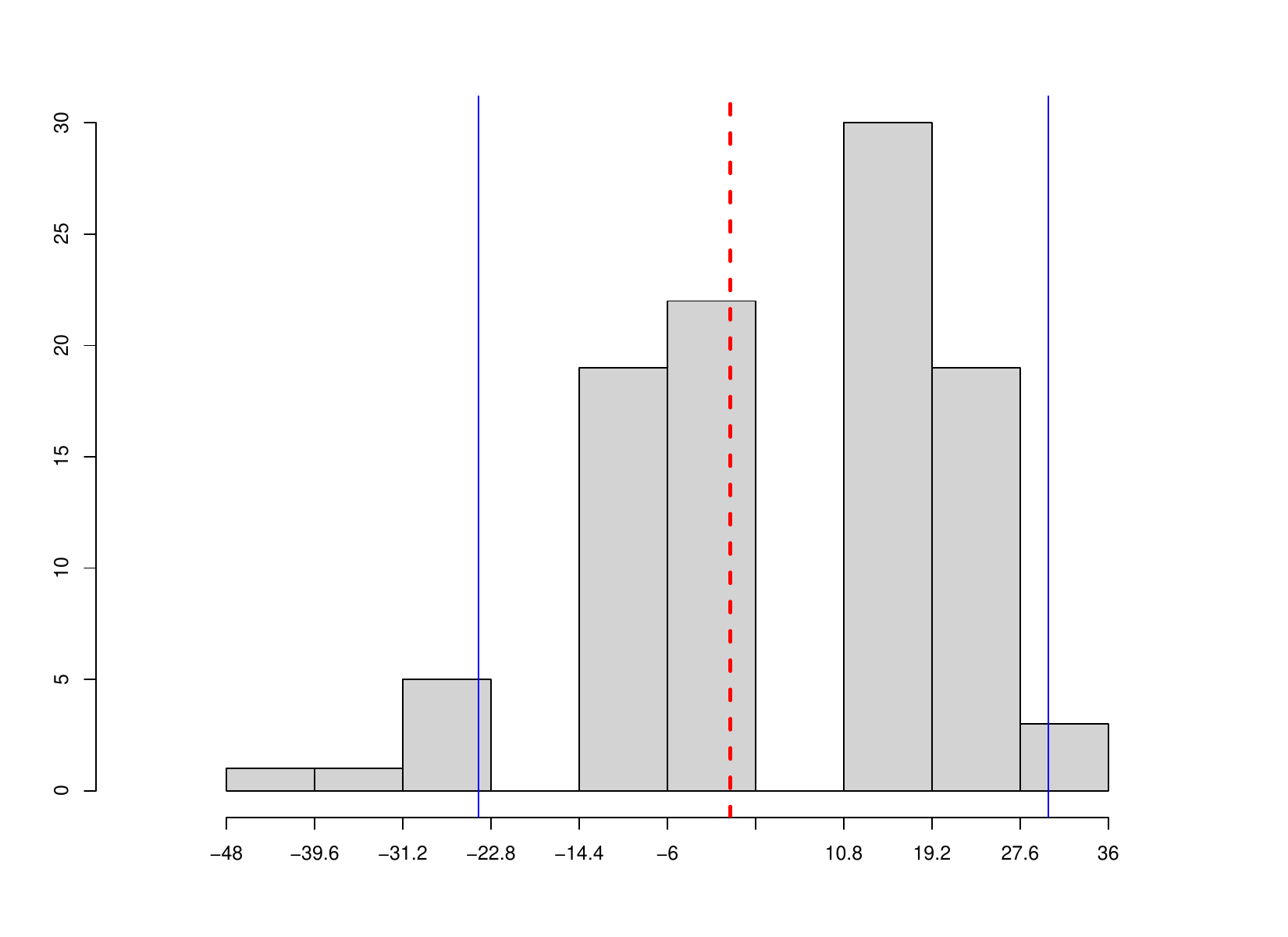}\\
		\small{(b)} 
	\end{minipage}
	\caption{\small{(a) $2$-spin Ising model fit, and (b) $3$-spin Ising model fit on the user preference vector for Britney Spears.}}
	\label{fig:test2}
\end{figure}

\begin{figure}[h]
	\begin{minipage}[c]{0.5\textwidth}
		\centering
		\includegraphics[width=3.35in]
		{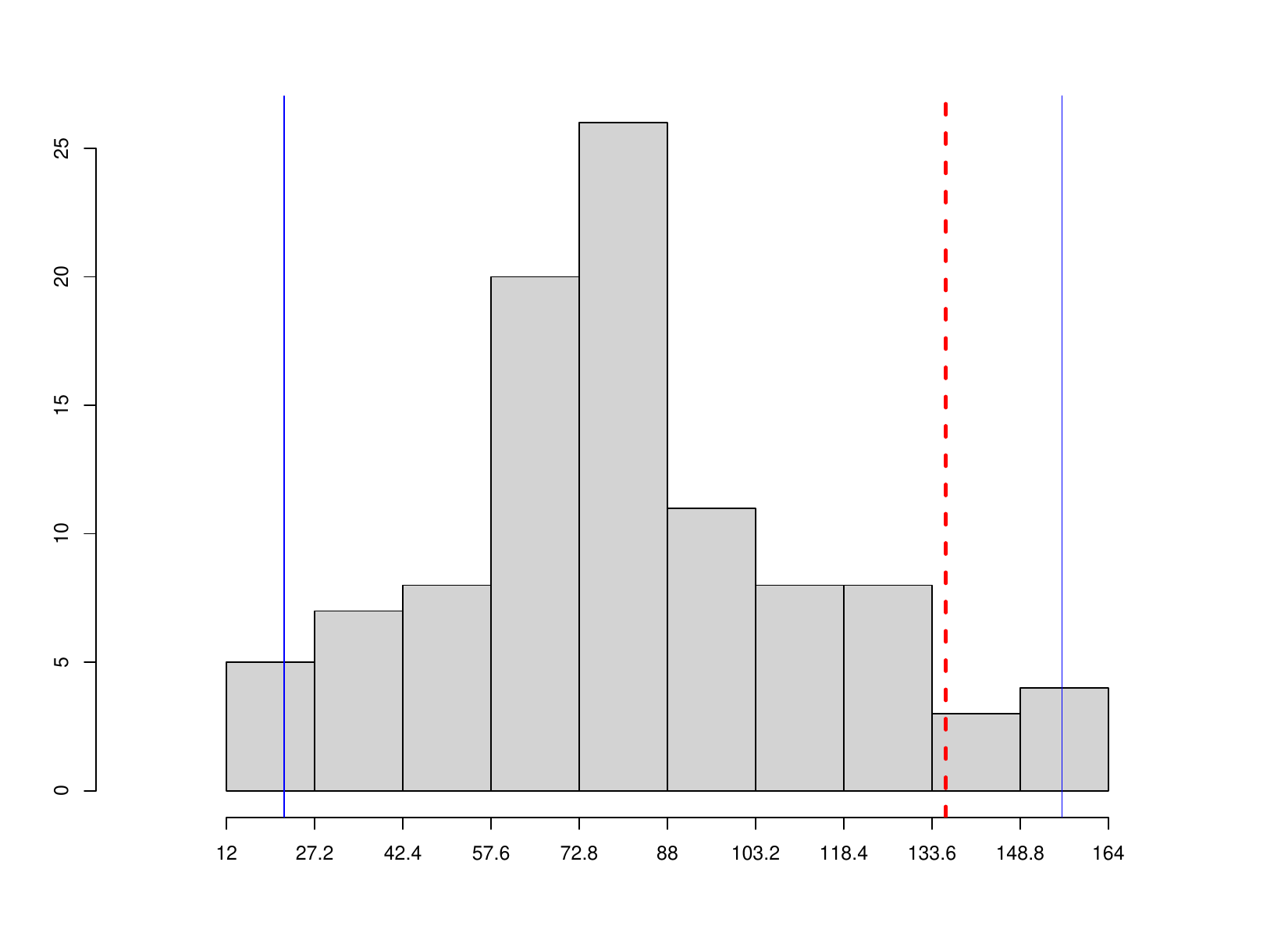}\\
		\small{(a)}
	\end{minipage}
	\begin{minipage}[l]{0.5\textwidth} 
		\centering
		\includegraphics[width=3.35in]
		{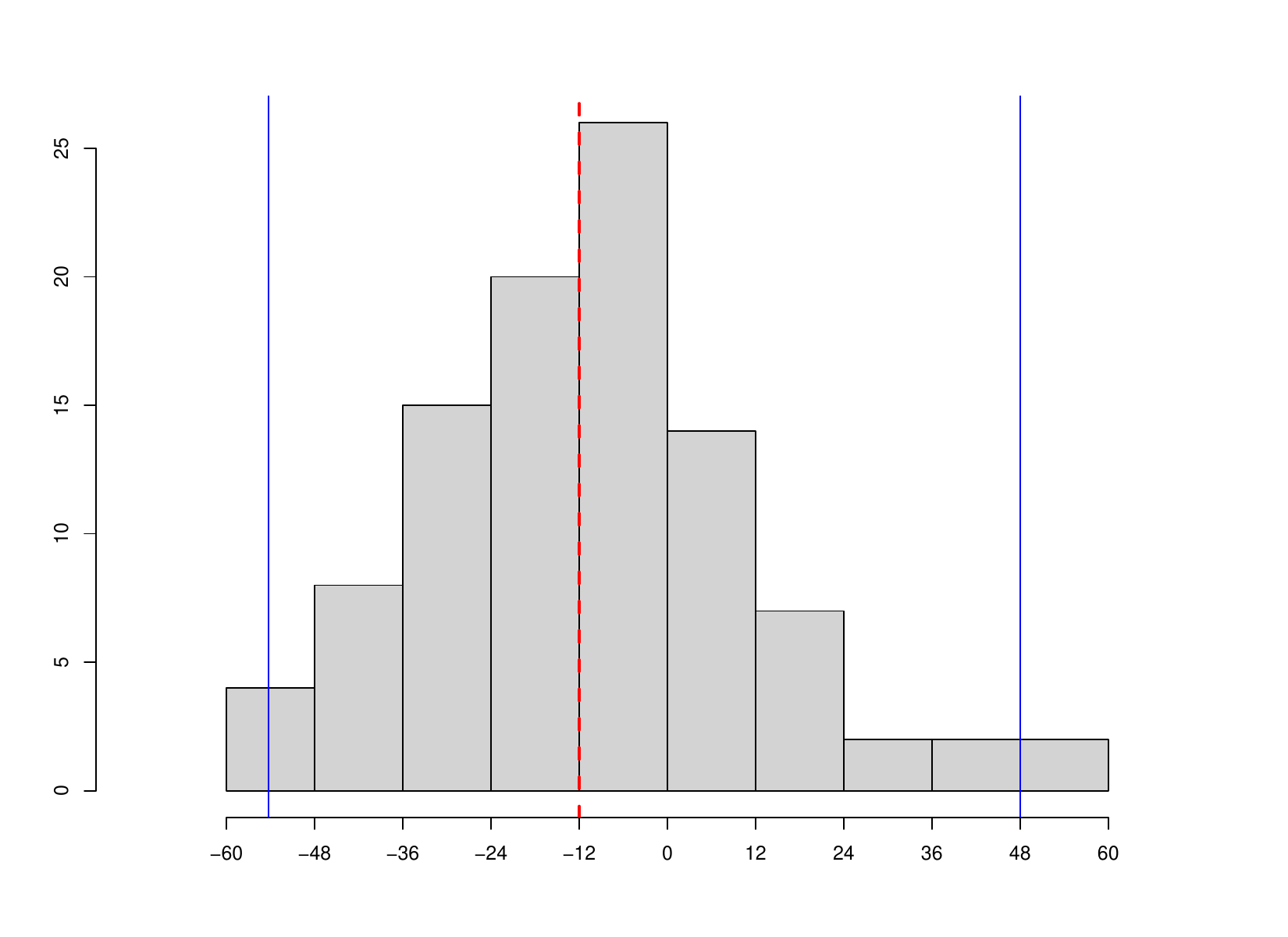}\\
		\small{(b)} 
	\end{minipage}
	\caption{\small{(a) $2$-spin Ising model fit, and (b) $3$-spin Ising model fit on the user preference vector for Rihanna.}}
	\label{fig:test3}
\end{figure}

\begin{figure}[h]
	\begin{minipage}[c]{0.5\textwidth}
		\centering
		\includegraphics[width=3.35in]
		{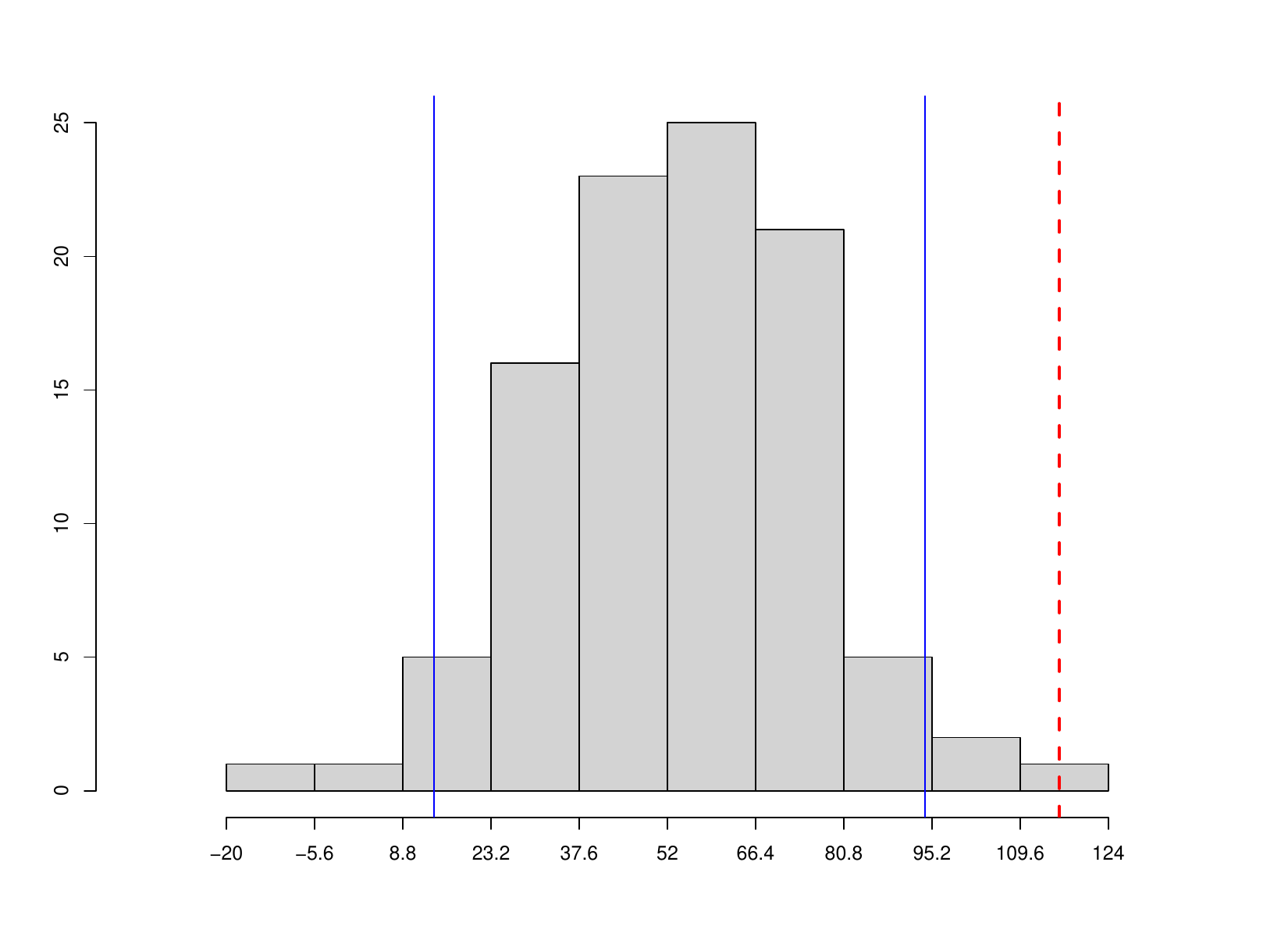}\\
		\small{(a)}
	\end{minipage}
	\begin{minipage}[l]{0.5\textwidth} 
		\centering
		\includegraphics[width=3.35in]
		{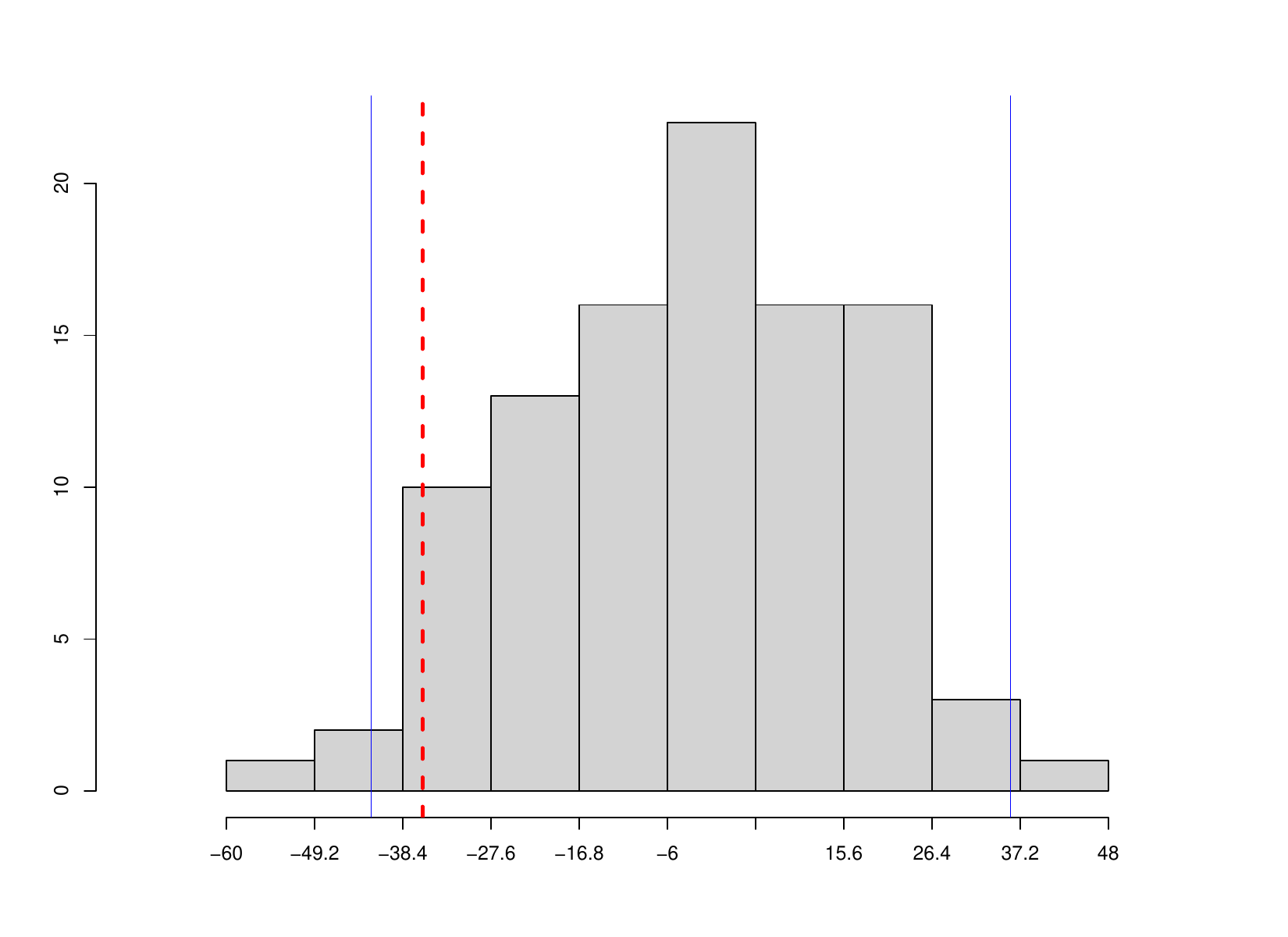}\\
		\small{(b)} 
	\end{minipage}
	\caption{\small{(a) $2$-spin Ising model fit, and (b) $3$-spin Ising model fit on the user preference vector for the Beatles.}}
	\label{fig:test4}
\end{figure}

From Figures \ref{fig:test1} (a), \ref{fig:test2} (a), \ref{fig:test3} (a) and \ref{fig:test4} (a), we observe that the actual value of $H(\bm X)$ lies outside the acceptance thresholds ($2.5$ to $97.5$ percentile of the empirical distribution of the simulated Hamiltonians) for the user preference data corresponding to Lady Gaga, Britney Spears and the Beatles. For Rihanna, the actual value of $H(\bm X)$ lies inside the acceptance thresholds, but only marginally. Hence, there is a strong evidence against our null hypothesis of a $2$-spin Ising model, thereby indicating that pairwise interaction is not enough to explain the dependency in the user friendship network.    

Next, we fit a $3$-tensor Ising model on these data. To be specific, we take $J_{ijk}$ to be the indicator that the users $i,j$ and $k$ form a triangle in the user friendship network (i.e. any two of $i,j$ and $k$ are friends). We then simulate $100$ observations $\bm X^{(1)},\ldots, \bm X^{(100)}$ from the $3$-tensor Ising model with parameter $\hat{\beta}$, where $\hat{\beta}$ is the MPLE of $\beta$ based on the original data $\bm X$ assuming a $3$-tensor Ising model. Similar to the previous setup, we decide to accept the null hypothesis of a $3$-tensor Ising model if and only if the actual value of the sufficient statistic $H(\bm X) := \sum_{i\sim j, j\sim k, k\sim  i} X_i X_j X_k$ lies within the $2.5^{\mathrm{th}}$ and $97.5^{\mathrm{th}}$ percentiles of the empirical distribution of $H(\bm X^{(1)}),\ldots,H(\bm X^{(100)})$.

From Figures \ref{fig:test1} (b), \ref{fig:test2} (b), \ref{fig:test3} (b) and \ref{fig:test4} (b), we observe that the actual value of $H(\bm X)$ lies inside the acceptance thresholds ($2.5$ to $97.5$ percentile of the empirical distribution of the simulated Hamiltonians) for the user preference data corresponding to Lady Gaga (marginally), Rihanna, Britney Spears and the Beatles. Hence, there is no evidence against our null hypothesis of a $3$-spin Ising model, thereby indicating that $3$ or even higher order interactions explain the dependency in the user friendship network much better than pairwise interactions.

\section{Proof of Theorem \ref{chextension}}\label{sec:3}

In this section, we prove of Theorem \ref{chextension}.  We first state the two main technical estimates required in the proof, and show how these results can be used to complete the proof of Theorem \ref{chextension}. As mentioned before, the first step in the proof of Theorem \ref{chextension}  is to show that the (scaled) log-pseudolikelihood concentrates around zero at the true parameter value $\beta > 0$ at the desired rate. This is achieved by proving the following second-moment estimate on the scaled log-pseudolikelihood function. The proof of this lemma is given in Section \ref{sec:lem1}.

\begin{lem}\label{lm:boundsecondmoment} Let $\beta > 0$ be such that assumption $(1)$ of Theorem \ref{chextension} holds. Then 
	$$\e_\beta\left[s_\bs^2(\beta)\right] = O_{\beta, p}\left(\frac{1}{N}\right),$$
	where $s_{\bm X}(b) = \frac{1}{pN} \frac{\partial \log L(b|\bm X)}{\partial b} = \frac{1}{N} (H_N(\bm X) -  \sum_{i=1}^N m_i(\bm X) \tanh (p b m_i(\bs) ) )$. 
\end{lem}

The next step of the proof is to show the strong concavity of the log-pseudolikelihood, that is, $-\frac{\partial}{\partial \beta}s_{\bm X}(\beta ) $ is strictly positive and bounded away from $0$ with high probability. To this end, note that for any $M > 0$, 
\begin{align}\label{eq:sbeta_II}
-\frac{\partial}{\partial \beta}s_{\bm X}(\beta ) & =\frac{p}{N}\sum_{i=1}^N m_i(\bm X)^2\text{sech}^2( p \beta m_i(\bm X)) \nonumber \\ 
& \geq  \frac{p}{N}\text{sech}^2(p \beta M)\sum_{i=1}^N m_i(\bm X)^2 \bm  1\{|m_i(\bm X)|\leq M\}.
\end{align}
Therefore, to show that $-\frac{\partial}{\partial \beta}s_{\bm X}(\beta )$ is strictly positive, it suffices to show that $$\sum_{i=1}^N m_i(\bm X)^2 \bm  1\{|m_i(\bm X)|\leq M\} = \Omega(N)$$ with high probability. This is formalized in the following lemma which is proved in Section \ref{sec:lem2}.

\begin{lem}\label{derivative} Fix $0 < \delta < 1$. Then under the assumptions in Theorem \ref{chextension},  there exists $\varepsilon=\varepsilon(\delta, \beta)>0$ and $M=M(\delta, \beta)<\infty$  such that 
	\begin{align*}
	\p_{\beta}\left(\sum_{i=1}^N m_i(\bm X)^2\bm 1\{|m_i(\bm X)|\leq M\}\geq  \varepsilon N \right)\geq  1-\delta,
	\end{align*}
	for all $N$ large enough. 
\end{lem}

The proof of Theorem \ref{chextension} can now be easily completed using the above lemmas. To this end, note that for any $M_1 > 0$, by Chebyshev's inequality and Lemma \ref{lm:boundsecondmoment} we have, 
\begin{align*}
\p_{\beta}\left(|s_{\bm X}(\beta )|> \frac{M_1}{\sqrt N}\right)\lesssim_{\beta,p} \frac{1}{M_1^2}.
%\label{chebyshev}
\end{align*}
Now, fix $\delta>0$. Therefore, it is possible to choose $M_1=M_1(\delta, \beta)$ such that the RHS above is less than $\delta$. Next, by Lemma \ref{derivative} there exists $\varepsilon=\varepsilon(\delta, \beta)>0$ and $M_2=M_2(\varepsilon,\delta, \beta)<\infty$  such that 
\begin{align*}
\p_{\beta}\left(\sum_{i=1}^N m_i(\bm X)^2\bm 1\{|m_i(\bm X)|\leq M_2\}\geq  \varepsilon N \right)\geq  1-\delta, 
\end{align*}
for $N$ large enough. Thus, defining 
$$T_N:=\left\{\bm X \in \cC_N: |s_{\bm X}(\beta )| \leq \frac{M_1}{\sqrt N},~ \sum_{i=1}^N m_i(\bm X)^2 \bm  1\{|m_i(\bm X)|\leq M_2\}\geq  \varepsilon N \right\},$$ gives $\p_{\beta}(T_N)\geq  1-2\delta$, for $N$ large enough. For $\bm X \in T_N$, recalling \eqref{eq:sbeta_II}, gives 
\begin{align*}
-\frac{\partial}{\partial \beta}s_{\bm X}(\beta ) \geq  \frac{p}{N}\text{sech}^2(p \beta M_2)\sum_{i=1}^N m_i(\bm X)^2 \bm  1\{|m_i(\bm X)|\leq M_2\} &\geq p \varepsilon \text{sech}^2(p \beta M_2). 
\end{align*}			
Therefore, for $\bm X \in T_N$, 
\begin{align*}
\frac{M_1}{\sqrt N}\geq  |s_{\bm X}(\beta)|=|s_{\bm X}(\beta)-s_{\bm X}(\hat \beta_N(\bm X))| &  \geq   - \int_{\beta\wedge \hat \beta_N(\bm X)}^{\beta\vee \hat \beta_N(\bm X)} \frac{\partial}{\partial \beta}s_{\bm X}(\beta ) \mathrm d\beta \nonumber\\
& \geq \frac{\varepsilon }{M_2 }|\tanh(p M_2\hat \beta_N(\bm X))-\tanh(p M_2\beta)| . 
\end{align*}
Then, defining $M=M(\delta, \beta):=\frac{M_2}{M_1\varepsilon}$, shows that 
$$\p_{\beta}\left(\sqrt{N}|\tanh(p M_2\hat \beta_N(\bm X))-\tanh (p M_2\beta)|\leq R \right)\geq 1- 2\delta.$$
The proof of Theorem \ref{chextension} now follows by inverting the $\tanh$ function.

\subsection{Proof of Lemma \ref{lm:boundsecondmoment}}\label{sec:lem1}

For $\bm x, \bm x' \in \cC_N$ define 
\begin{equation*}
F(\bt,\bt')=\frac{1}{2}\sum_{i=1}^N \left(m_i(\bt)+m_i(\bt')\right)(x_i-x_i'),
\end{equation*}
where $m_i$ is as defined in \eqref{eq:conditional}. Note that $F$ is antisymmetric, that is,  $F(\bt, \bt') = -F(\bt',\bt)$.

Now, choose a coordinate $I \in \{1, 2, \ldots, N\}$ uniformly at random and replace the $I$-th coordinate of $\bm X \sim \p_{\beta}$ by a sample drawn from the conditional distribution of $X_u$ given $(X_v)_{v \ne I}$. Denote the resulting vector by $\bm X'$. Note that $F(\bs,\bs')=m_I(\bs)(X_I-X_I')$. Then 
\begin{align}\label{eq:sx_II}
f(\bs):=\mathbb{E}_\beta\left(F(\bs,\bs')|\bs\right) &= \frac{1}{N}\sum_{i=1}^N m_i(\bs) \left\{ X_i- \e_{\beta}\left(X_i|(X_j)_{j\neq i}\right) \right\} \nonumber \\ 
&=\frac{1}{N}\sum_{i=1}^N m_i(\bs)\left(X_i-\tanh\big(p \beta m_i(\bs) \big)\right) \nonumber \\
&=s_\bs(\beta).
\end{align}
Now, since $(\bs,\bs')$ is an exchangeable pair, 
\begin{equation*}
\mathbb{E}_\beta\left(f(\bs)F(\bs,\bs')\right)=\mathbb{E}_\beta\left(f(\bs')F(\bs',\bs)\right).
\end{equation*}
Again, because $F$ is antisymmetric, we have $\mathbb{E}_\beta\left(f(\bs')F(\bs',\bs)\right)=-\mathbb{E}_\beta\left(f(\bs')F(\bs,\bs')\right).$ Hence,  
\begin{align}
\mathbb{E}_\beta\left(f(\bs)^2\right)=\e_{\beta}\left(f(\bs)\e_{\beta}\left[F(\bs,\bs')|\bs\right)\right] & = \mathbb{E}_\beta\left(f(\bs)F(\bs,\bs')\right) \nonumber \\ 	 		
\label{expectation f squared} &=\tfrac{1}{2}\mathbb{E}_\beta\left((f(\bs)-f(\bs'))F(\bs,\bs')\right).
\end{align}

Now, for any $1 \leq t \leq N$ and $\bm x \in \cC_N$, let 
$$\bm x^{(t)}=(x_1, x_2, \ldots, x_{t-1}, 1- x_t, x_{t+1}, \ldots, x_N),$$ 
and 
\begin{equation*}
p_t(\bt):=\mathbb{P}_\beta(X_t'=-x_t|\bs=\bt,I=t)=\frac{e^{-p\beta x_t m_t(\bt)}}{e^{-p\beta  m_t(\bt)} + e^{p\beta m_t(\bt)}}.
\end{equation*}
This implies, 		
\begin{align}\label{cond2step} 
\E_{\beta}((f( \bm X)-f(\bm X'))F( \bm X,  \bm X')| \bm X) \nonumber &=\frac{1}{N}  \sum_{t=1}^N (f( \bm X)-f( \bm X^{(t)})) F( \bm X,  \bm X^{(t)}) p_t( \bm X) \nonumber \\ 	
&=\frac{2}{N}\sum_{t=1}^N m_t(\bs)X_t p_t(\bs) (f(\bs)-f(\bm X^{(t)})).
\end{align} 
For $1 \leq s, t \leq N$, let $a_s(\bt):=x_s -\tanh(p \beta m_s(\bt))$ and $b_{s t}(\bt):=\tanh( p \beta m_s(\bt))-\tanh(p \beta m_s(\bt^{(t)}))$. Then, noting that $f(\bt)= \frac{1}{N} \sum_{s=1}^N m_s(\bt) a_s(\bt)$ gives 
\begin{align} \label{fxx} 
f(\bs)-f(\bm X^{(t)}) &=\frac{1}{N}\sum_{s=1}^N (m_s(\bs)-m_s(\bm X^{(t)}))a_s(\bs)+\frac{1}{N}\sum_{s=1}^N m_s(\bm X^{(t)}) (a_s(\bs)-a_s(\bm X^{(t)}) ) \nonumber  \\ 
&= A_t + B_t + C_t,
\end{align} 
where $$A_t := \frac{2(p-1)X_t}{N}\sum_{s=1}^N  J_{s t}(\bs) a_{s}(\bs), B_t:= \frac{2 m_t(\bs)X_t}{N},\quad\textrm{and}~C_t :=-\frac{1}{N}\sum_{s=1}^Nm_s(\bm X^{(t)}) b_{s t}(\bs).$$ Then using \eqref{expectation f squared} and \eqref{cond2step}, we have 
\begin{equation*} 
\mathbb{E}_\beta\left(f(\bs)^2\right)= \frac{1}{N}\sum_{t=1}^N\e_{\beta}\left[\left(A_t + B_t + C_t \right)m_t(\bs)X_t p_t(\bs)\right]. 
\end{equation*}

Now, define the following three quantities: $$\bm a(\bs) := \left(a_1(\bs),\ldots,a_N(\bs)\right), \quad \bm m (\bs) := \left(m_1(\bs),\ldots,m_N(\bs)\right),$$ and $\bm M(\bs) := \left(m_1(\bs)p_1(\bs),\ldots,m_N(\bs)p_N(\bs)\right)$. Note that $\bm m(\bs) = \bs \bm J_N(\bs)^\top$. Also, observe that each entry of $\bm a(\bs)$ is bounded in absolute value by $2$, hence, $\|\bm a(\bs)\| \leq 2\sqrt{N}$. Moreover, using $p_t(\bm X) \leq 1$, 
\begin{align}\label{eq:norm_bound}
\|\bm M(\bs)\| = \left(\sum_{i=1}^N m_i(\bs)^2 p_i(\bs)^2\right)^\frac{1}{2} \leq \|\bm m(\bs)\| = \|\bm J_N(\bs)\bs^\top\|\leq \sqrt N \|\bm J_N(\bs)\| . 
\end{align}
Hence, recalling the definition of $A_t$ from \eqref{fxx} gives 
\begin{align}\label{firstterm}
\Bigg|\frac{1}{N}\sum_{t=1}^N A_t m_t(\bs) X_t p_t(\bs)\Bigg| & = \frac{2(p-1)}{N^2} \Big|\bm a(\bs) \bm J_N(\bs) \bm M(\bs)^\top\Big|\nonumber\\
&\lesssim_p \frac{1}{N^2}\|\bm J_N(\bs)\|\|\bm a(\bs)\|\|\bm M(\bs)\|\nonumber\\
&\lesssim_p \frac{||\bm J_N(\bm X)||^2}{N} .
\end{align}
Next, we  consider the term corresponding to $B_t$: 
\begin{align}\label{secondterm}
\Bigg|\frac{1}{N}\sum_{t=1}^N B_t m_t(\bs) X_t p_t(\bs)\Bigg| = \frac{2}{N^2}\Bigg|\sum_{t=1}^N m_t^2(\bs)p_t(\bs)\Bigg| \lesssim \frac{1}{N^2} \|\bm m (\bs)\|^2 \leq \frac{||\bm J_N(\bm X)||^2}{N},
\end{align}
where the last step uses \eqref{eq:norm_bound}.

Finally, we consider the term corresponding to $C_t$. Let us define the matrix $\bm J_{N, 2}(\bs) := ((J_{ij}(\bs)^2))_{1\leq i,j\leq N}$. Then, denoting by $\bm e_i$ the vector in $\R^N$ with the $i$-th entry 1 and 0 everywhere else, we get 
\begin{equation*}
\|\bm J_{N, 2}(\bs)\| \leq \max_{1 \leq i \leq N} \sum_{j=1}^N J_{ij}(\bs)^2 \leq \max_{1 \leq i \leq N}
\|\bm e_{i}^\top \bm J_N(\bs)\|^2 \leq  \|\bm J_N(\bs)\|^2.
\end{equation*}
Let $h(x):= \tanh\left(p \beta x \right)$. It is easy to check that $\|h''\|_\infty \leq \beta^2$. Hence, by a Taylor expansion, for $1 \leq s \leq N$, 
\begin{align}\label{taylorh}
\big|h(m_s(\bs))- h(m_s(\bm X^{(t)})) - (m_s(\bs) - m_s(\bm X^{(t)}) ) h'(m_s(\bs))\big| \lesssim_{\beta} \left(m_s(\bs) - m_s(\bm X^{(t)})\right)^2.
\end{align}
Note that $$m_s(\bs)-m_s(\bm X^{(t)})=2(p-1)J_{ s t}(\bs)X_t$$ and $$h(m_s(\bs))- h(m_s(\bm X^{(t)})) = \tanh( p \beta m_s(\bs))-\tanh(p \beta m_s(\bs^{(t)})) = b_{s t}(\bs)$$ Hence, \eqref{taylorh} can be rewritten as:
\begin{equation*}
\left|b_{s t}(\bs)-2(p-1)J_{s t}(\bs)X_t h'(m_s(\bs))\right|\lesssim_{\beta, p} J_{s t}(\bs)^2.
\end{equation*} 
Using the above bounds, we have the following for any two vectors $\bm x=(x_1, x_2, \ldots, x_N)$ and $\bm y = (y_1, y_2, \ldots, y_N) \in\mathbb{R}^N$, 
\begin{align}
& \left|\sum_{1 \leq s, t \leq N} x_s y_t b_{s t}(\bs)\right| \nonumber \\ 
&\leq \left|\sum_{1 \leq s, t \leq N}2(p-1)x_s y_t J_{s t}(\bs)X_t h'(m_s(\bs))\right|\nonumber\\ &+\left|\sum_{1 \leq s, t \leq N}x_s y_t \big(b_{s t}(\bs)-2(p-1)J_{s t}(\bs)X_t h'(m_s(\bs))\big)\right| \nonumber \\
\nonumber & \lesssim_{\beta, p}  \|\bm J_N(\bs)\| \left(\sum_{s=1}^N (x_s h'(m_s(\bs)))^2 \right)^{\frac{1}{2}} \left(\sum_{t=1}^N (y_t X_t)^2 \right)^{\frac{1}{2}}\nonumber\\ &+ \sum_{1 \leq s, t \leq N}|x_s y_t|\ J_{s t}(\bs)^2\\
\nonumber &\lesssim_{\beta, p} \|\bm J_N(\bs)\| \|\bm x\| \|\bm y\| +  \|\bm J_{N, 2}(\bs)\| \|\bm x\|\|\bm y\| \\
\label{ints1} &\lesssim_{\beta, p}  \left(  \|\bm J_N(\bs)\|  +  \|\bm J_N(\bs)\|^2 \right)\|\bm x\|\|\bm y\|.
\end{align}
Again, by a Taylor expansion and using the bound $||h'||_{\infty} \lesssim_{\beta, p} 1$, gives 
$$|b_{s t}(\bs)| \lesssim_{\beta, p} |m_s(\bs)-m_s(\bm X^{(t)}) | \lesssim_{\beta, p} |J_{s t}(\bs)|.$$ Consequently,
\begin{equation}\label{ints2}
\left|\sum_{1 \leq s, t  \leq N}x_s y_t  J_{s t}(\bs)b_{s t}(\bs)\right| \lesssim_{\beta, p} \sum_{1 \leq s, t  \leq N}|x_s y_t | {J}_{s t}(\bs)^2 \leq ||\bm J_N(\bm X)||^2 \|\bm x\|\| \bm y\|.
\end{equation}
Now, recalling the definition of $C_t$ from \eqref{fxx} gives, 
\begin{align}\label{thirdterm}
\nonumber&\left|\frac{1}{N}\sum_{t=1}^N C_t m_t(\bs) X_t p_t(\bm X)\right|\\ \nonumber &=\left|\frac{1}{N^2} \sum_{1 \leq s, t \leq N} m_s( \bm X^{(t)} )b_{s t}(\bs)m_t(\bs) X_t p_t(\bs) \right|\\
\nonumber 
&=\left|\frac{1}{N^2}\sum_{1 \leq s, t \leq N} \big(m_s(\bs)-2(p-1)J_{s t}(\bs) X_t \big) b_{s t}(\bs) m_t(\bs) X_t p_t(\bs)\right|\\
& \lesssim_{\beta, p} \frac{||\bm J_N(\bm X)||^3 + ||\bm J_N(\bm X)||^4}{N}, 
\end{align}
where the last step uses \eqref{ints1}, \eqref{ints2}, and $\|\bm m(\bs)\| = \|\bm J_N(\bs)\bs^\top\|\leq \sqrt N \|\bm J_N(\bs)\|$. 

%\iffalse
Combining \eqref{firstterm}, \eqref{secondterm} and \eqref{thirdterm}, it follows that	$\mathbb{E}_\beta(f(\bm X)^2) \lesssim_{\beta, p} \frac{1}{N}$, since by condition of $(1)$ of Theorem \ref{chextension}, $\E_{\beta}(||\bm J_N(\bm X)||^4 )$ is uniformly bounded. This completes the proof of the lemma, since recalling \eqref{eq:sx_II}, $f(\bs) = s_\bs(\beta)$. \qed

\subsection{Proof of Lemma \ref{derivative}}\label{sec:lem2}

We begin with the following simple observation, which says that if $\liminf_{N\rightarrow\infty}\frac{1}{N}{F_N(\beta)} > 0$, we can find a $\gamma$ small enough such that $\liminf_{N\rightarrow\infty}\frac{1}{N}{F_N'(\beta-\gamma)} > 0$.

\begin{obs}\label{obs:Fderivative} Suppose $\beta > 0$ is such that $\liminf_{N\rightarrow\infty} \frac{1}{N} F_{N}(\beta) > 0$. Then $$\lim_{\delta\rightarrow 0}\liminf_{N\rightarrow\infty}\frac{1}{N} F_{N}(\beta-\delta) > 0.$$
\end{obs}

\begin{proof} Denote $K:=\sup_{N \geq 1}  \E_\beta (||J_N(\bm X)||) < \infty $ Then, $$F_{N}'(\beta): = \frac{\mathrm d}{\mathrm d \beta} F_{N}'(\beta) = \E_\beta (H_N(\bm X)) = \E_\beta (\bm X' \bm J_N(\bm X) \bm X ) \leq N  \E_\beta (|| \bm J_N(\bm X)||) \leq K N.$$ 
	Therefore, by a Taylor expansion
	$$\lim_{\delta\rightarrow 0}\liminf_{N\rightarrow\infty}\frac{1}{N} F_{N}(\beta-\delta)  \geq  \lim_{\delta\rightarrow 0}\liminf_{N\rightarrow\infty}\left(\frac{1}{N} F_{N}(\beta) -   \delta  K \right)>0,$$
	as required. 
\end{proof}

Now, note that for any $\varepsilon, \gamma >0$,  
\begin{align*}
&\p_{\beta}(H_N(\bm X)<\varepsilon N)=\p_{\beta}(e^{-\gamma H_N(\bm X)}>e^{-\gamma \varepsilon N}) \leq e^{\gamma \varepsilon N+F_N(\beta-\gamma)-F_N(\beta)}
\end{align*}
which, on taking logarithms, implies that 
\begin{align*}
\log \p_{\beta}(H_N(\bm X)<\varepsilon N)\leq  \varepsilon\gamma N - \int_{\beta-\gamma}^{\beta}F_N'(t)\mathrm dt\leq  \varepsilon \gamma N -F_N'(\beta-\gamma)\gamma ,
\end{align*} 
by the monotonicity of $F_N'(\cdot)$. Dividing both sides by $N$ and taking limits as $N\rightarrow \infty$ followed by $\varepsilon\rightarrow 0$ we have
$$\lim_{\varepsilon\rightarrow0}\limsup_{N\rightarrow\infty}\frac{1}{N}\log \p_{\beta}(H_N(\bm X)<\varepsilon N)\leq -\liminf_{N\rightarrow\infty}\frac{1}{N}{F_N'(\beta-\gamma)} \leq -\liminf_{N\rightarrow\infty}\frac{F_N(\beta-\gamma)}{N(\beta-\gamma)} <0,$$
by choosing $\gamma$ small enough (by Observation \ref{obs:Fderivative}).  This shows that, for every $0 < \delta < 1$ there exists  $\varepsilon=\varepsilon(\delta)>0$ such that, for $N$ large enough,   
\begin{equation}
\p_{\beta}(H_N(\bm X)< 2\varepsilon N )\leq \delta. 
\label{c00}
\end{equation}

Next, by Lemma \ref{lm:boundsecondmoment} and Chebyshev's inequality, there exists  $M_1=M_1(\delta)< \infty$ such that 
\begin{equation}\label{c1}
\p_{\beta}\left(|s_{\bm X}(\beta )|> \frac{M_1}{\sqrt N} \right)  \leq \delta.
\end{equation}
Moreover, note that for any $M_2 > 0$, 
$$\sum_{i=1}^N|m_i(\bm X)|\bm 1\{|m_i(\bm X)|> M_2\} \leq \frac{1}{M_2} \sum_{i=1}^N m_i(\bm X)^2 =\frac{1}{M_2}  ||J_N(\bm X) \bm X^\top||^2 \leq \frac{N ||J_N(\bm X)||^2}{M_2}.$$
Therefore, using Markov's inequality and condition $(1)$ of Theorem \ref{chextension}, we can choose $M_2=M_2(\delta) < \infty$ such that for all $N$ large enough, 
\begin{align}\label{eq:msigma_epsilon}
\p_{\beta}\left(\sum_{i=1}^N|m_i(\bm X)|\bm 1\{|m_i(\bm X)|> M_2\} > \varepsilon N \right) \leq \frac{\E_\beta(||J_N(\bm X)||^2)}{ \varepsilon M_2} \leq \delta. 
\end{align} 

Then, defining 
\begin{align*}
T_{N}:= \left\{\bm X\in \cC_N: H_N(\bm X)\geq  2\varepsilon N, ~  |s_{\bm X}(\beta )| \leq \frac{M_1}{\sqrt N}, ~ \sum_{i=1}^N|m_i(\bm X)|\bm 1\{|m_i(\bm X)|> M_2\} \leq \varepsilon N   \right \}, 
\end{align*}
and combining \eqref{c00}, \eqref{c1}, and \eqref{eq:msigma_epsilon}, gives $\p_{\beta}(T_N )\geq  1-3\delta$, for $N$ large enough. Now, for $\bm X \in T_N$, 
\begin{align*}
& \sum_{i=1}^N m_i(\bm X)^2 \bm  1\{|m_i(\bm X)|\leq M_2\}+\varepsilon N \\ 
& \geq  \frac{1}{p \beta} \sum_{i=1}^N |m_i(\bm X)| \tanh(p \beta |m_i(\bm X)|) \bm  1\{|m_i(\bm X)|\leq M_2\} + \sum_{i=1}^N |m_i(\bm X)| \bm  1\{|m_i(\bm X)| >  M_2\}  \tag*{(using $\tanh (x) \leq x$)} \\ 
& \gtrsim_{p,\beta} \sum_{i=1}^N |m_i(\bm X)| \tanh(p \beta |m_i(\bm X)|) \bm  1\{|m_i(\bm X)|\leq M_2\}\\ &+ \sum_{i=1}^N |m_i(\bm X)|  \tanh(p \beta |m_i(\bm X)| ) \bm  1\{|m_i(\bm X)| > M_2\} \tag*{(using $\tanh (x) \leq 1$)}   \\ 
& \geq \sum_{i=1}^N m_i(\bm X)\tanh(p \beta m_i(\bm X))\\ 
& = H_N(\bm X)-Ns_{\bm X}(\beta ) \geq 2\varepsilon N - M_1\sqrt{N}. 
\end{align*}
Thus, on the set $T_N$,
$$\sum_{i=1}^N m_i(\bm X)^2 \bm  1\{|m_i(\bm X)|\leq M_2\} \gtrsim_{p,\beta} 2\varepsilon N- M_ 1\sqrt{N}> \varepsilon N,$$ for all $N$ large enough. This completes the proof of Lemma \ref{derivative}. \qed

\section{Proof of Corollary \ref{skthreshold}}\label{skproof}

To prove Corollary \ref{skthreshold} we will verify that the conditions in Theorem \ref{chextension} hold with probability 1. As mentioned before, in this case condition (2) is easy to verify. To this end, note that by \cite[Theorem 1.1]{bovier}, $\lim_{N \rightarrow \infty} \frac{1}{N} F_N(\beta) = \beta^2 /2$ almost surely, for $\beta > 0$ small enough. This implies, since $F_N$ on increasing on the positive half-line, $\lim_{N \rightarrow \infty} \frac{1}{N} F_N(\beta) >0$ almost surely, for all $\beta > 0$. This establishes condition (2) in Theorem \ref{chextension}.

We now proceed to verify condition (1). To begin with, fix $\bt \in \sa_N$ and consider the Gaussian process 
\begin{align}\label{eq:sphere_ux}
G_{\bm u}(\bt):= \bm u^\top \bm J_N(\bt) \bm u
\end{align} 
indexed by $\bm u \in S^{N-1} := \{\bm t \in \mathbb{R}^N: \|\bm t\| = 1\}$. Here, $\bm J_N(\bt)$ is the local interaction matrix corresponding to the tensor \eqref{eq:JN_sk} of the $p$-tensor SK model. Note that the maximum eigenvalue of $\bm J_N(\bt)$ can be expressed as $\lambda_{\max}\left(\bm J_N(\bt)\right) = \sup_{\bm u \in S^{N-1}} G_{\bm u}(\bt)$.\footnote{For any $N \times N$  matrix  $\bm A$, $\lambda_{\max}(\bm A)$ and $\lambda_{\min}(\bm A)$ denotes the maximum and the minimum eigenvalue of $\bm A$, respectively.}

\begin{lem}\label{tec2} Fix $\bm x \in \cC_N$ and consider the Gaussian process $\{G_{\bm u}(\bm x) : \bm u \in S^{N-1} \}$ as defined above in \eqref{eq:sphere_ux}. Then, the following hold: 
	\begin{itemize}
		\item[(1)]~For every vector $\bm u \in S^{N-1}$, $\e\left[G_{\bm u}(\bm x)^2\right] \lesssim_p \frac{1}{N}$;
		\item[(2)]~For vectors $\bm u, \bm v \in S^{N-1}$, $\e\left[G_{\bm u}(\bm x) - G_{\bm v}(\bm x)\right]^2 \lesssim_p \frac{1}{N} \sum_{i=1}^N (u_i-v_i)^2$. 
	\end{itemize}
\end{lem}

\noindent\textit{Proof of} $(1)$: Fix $\bm x \in \cC_N$ and $\bm u \in S^{N-1}$. Then  	
\begin{align*}
G_{\bm u}(\bt) &= \sum_{1 \leq i_1,\ldots,i_p \leq N} J_{i_1\ldots i_p} u_{i_1}u_{i_2}x_{i_3}\ldots x_{i_p}\\ &= (p-2)!\sum_{1\leq i_1<\ldots<i_p \leq N} J_{i_1\ldots i_p} \left(\sum_{1\leq s \neq t \leq p} u_{i_s} u_{i_t} \prod_{a\in \{1, 2, \ldots, p\}\setminus\{s,t\}} x_{i_a}\right).
\end{align*} 
Hence, 
\begin{align}
\e \left[G_{\bm u} (\bt)^2\right] &\lesssim_p \frac{1}{N^{p-1}}\sum_{1\leq i_1<\ldots<i_p \leq N} \left(\sum_{1\leq s \neq t \leq p} u_{i_s} u_{i_t} \prod_{a\in \{1, 2, \ldots, p\}\setminus\{s,t\}} x_{i_a}\right)^2\nonumber\\
&\lesssim_p \frac{1}{N^{p-1}} \sum_{1\leq i_1<\ldots<i_p \leq N} \sum_{1\leq s \neq t \leq p} u_{i_s}^2 u_{i_t}^2\label{st1}\\
&\lesssim_p \frac{1}{N^{p-1}} \sum_{1\leq s\neq t \leq p} \sum_{1 \leq i_1,\ldots,i_p \leq N} u_{i_s}^2 u_{i_t}^2\nonumber\\
&= \frac{1}{N^{p-1}} \sum_{1\leq s\neq t\leq p} N^{p-2}\|\bm u\|^4 \lesssim_p \frac{1}{N}. \nonumber
\end{align}
where in \eqref{st1} we used the inequality $\left(\sum_{i=1}^n a_i\right)^2 \leq n \sum_{i=1}^n a_i^2$, for any sequence of real numbers $a_1,\ldots,a_n$.  \\ 

\noindent\textit{Proof of} $(2)$: Fix $\bm x \in \cC_N$ and $\bm u, \bm v \in S^{N-1}$.	Then, 	 		
\begin{align*}
& \e\left[G_{\bm u}(\bt) - G_{\bm v}(\bt)\right]^2\\ 
&= ((p-2)!)^2\e\left[\sum_{1\leq i_1<\ldots<i_p \leq N} J_{i_1\ldots i_p} \left(\sum_{1\leq s \neq t \leq p} (u_{i_s} u_{i_t} - v_{i_s} v_{i_t}) \prod_{a\in \{1, 2, \ldots, p\}\setminus\{s,t\}} \tau_{i_a}\right)\right]^2\\
& 
\lesssim_p \frac{1}{N^{p-1}}\sum_{1\leq i_1<\ldots<i_p \leq N} \left(\sum_{1\leq s \neq t \leq p} (u_{i_s} u_{i_t} - v_{i_s} v_{i_t}) \prod_{a\in \{1, 2, \ldots, p\}\setminus\{s,t\}} \tau_{i_a}\right)^2\\
& \lesssim_p \frac{1}{N^{p-1}}\sum_{1\leq i_1<\ldots<i_p \leq N} \sum_{1\leq s \neq t \leq p} (u_{i_s} u_{i_t} - v_{i_s} v_{i_t})^2\\ 
& \lesssim_p \frac{1}{N} \sum_{1 \leq i,j \leq N} (u_iu_j - v_iv_j)^2 \\ 
& = \frac{2}{N}  \left[1-\left(\sum_{i=1}^N u_i v_i\right)^2\right] \leq \frac{2}{N} \left[2-2\sum_{i=1}^N u_iv_i \right] \lesssim \frac{1}{N} \sum_{i=1}^N (u_i - v_i)^2.    
\end{align*} 
This completes the proof of Lemma \ref{tec2} (2). \qed \\

Using the lemma above we first show that $\e[\sup_{\bm u \in S^{N-1}} G_{\bm u}(\bt)] \lesssim_p 1$. We do this comparing the supremum of the Gaussian process $\{G_{\bm u}(\bt) : \bm u \in S^{N-1}\}$ with the supremum of the Gaussian process $\{H_{\bm u}: \bm u \in S^{N-1}\}$, where $H_{\bm u} = \sum_{i=1}^N g_i u_i$ and $g_1,\ldots,g_N$ are independent standard Gaussians. Now, by Lemma \ref{tec2} (2), there exists a constant $C:=C(p) > 0$, such that for $\bm u, \bm v \in S^{N-1}$, 
$$\e\left[G_{\bm u}(\bt) - G_{\bm v}(\bt)\right]^2 \leq \frac{C}{N} \sum_{i=1}^N (u_i - v_i)^2 =  \frac{C}{N} \e\left[H_{\bm u} - H_{\bm v}\right]^2.$$ Hence, by the Sudakov-Fernique inequality \cite[Theorem 1.1]{error}, 
\begin{align}\label{eq:Gx_expectation}
\e \left[\sup_{\bm u \in S^{N-1}} G_{\bm u}(\bt)\right]  \leq  \left(\frac{C}{N}\right)^{\frac{1}{2}} \e\left[\sup_{\bm u \in S^{N-1}} H_{\bm u}\right] & =   \left(\frac{C}{N}\right)^{\frac{1}{2}} \e\left[\left(\sum_{i=1}^N g_i^2\right)^{\frac{1}{2}}\right] \nonumber \\
&\leq  C^{\frac{1}{2}} : = D. 
\end{align}

Now, by Lemma \ref{tec2} $(1)$ there exists a constant $K : = K(p) > 0$ such that $$\sup_{\bm u \in S^{N-1}} \e \left[G_{\bm u} (\bt)^2\right] \leq \frac{K}{N}~.$$ Hence, by the Borell-TIS inequality \cite[Theorem 2.1.1]{Gaussian}, for any $t > 0$, 
\begin{equation*}
\p\left(\sup_{\bm u \in S^{N-1}} G_{\bm u}(\bt) - \e \left[\sup_{\bm u \in S^{N-1}} G_{\bm u}(\bt)\right] > t\right) \leq e^{-\frac{N t^2}{2K} }.
\end{equation*}
This implies, by \eqref{eq:Gx_expectation},
\begin{equation*}
\p\left(\sup_{\bm u \in S^{N-1}} G_{\bm u}(\bt) > D+t\right) \leq e^{-\frac{N t^2}{2K} }.
\end{equation*} 
Then, taking $t = \sqrt{2K}$ in the inequality above gives, 
$$\p\left(\|\bm J_N(\bt)\|> D + \sqrt{2K}\right) \leq 2e^{-N},$$
since we have $\lambda_{\max}(\bm J_N(\bt)) \stackrel{D}{=} -\lambda_{\min}(\bm J_N(\bt))$, because  $\bm J_N(\bt) \stackrel{D}{=} - \bm J_N(\bt)$. Therefore, by an union bound, 
$$\p\left(\sup_{\bt \in \sa_N}\|\bm J_N(\bt)\|> D + \sqrt{2 K}\right) \leq 2^{N+1}e^{-N} = 2(e/2)^{-N}.$$ Hence, by the Borel-Cantelli lemma, $\limsup_{N\rightarrow \infty} \sup_{\bt \in \sa_N}\|\bm J_N(\bt)\| \leq D + \sqrt{2 K}$ with probability $1$, which establishes \eqref{eq:JNcondition} and hence, condition $(1)$ of Theorem \ref{chextension}.

\section{Proof of Corollary \ref{boundeddeg} } 
\label{sec:boundeddegpf}

We begin by showing that condition $(1)$ of Corollary \ref{boundeddeg} implies $\sup_{N \geq 1} \sup_{\bm x \in \cC_N} ||\bm J_N(\bm x)|| < \infty$ and, hence, condition $(1)$ of Theorem \ref{chextension}. To this end, fix $\bm x \in \cC_N$, and take $\bm u \in S^{N-1} := \{\bm t \in \mathbb{R}^N: \|\bm t\| = 1\}$. Then, 
\begin{align}\label{sh}
|\bm u^\top \bm J_N(\bm x) \bm u| &= \Big|\sum_{1\leq i_1,i_2,\ldots,i_p\leq N} J_{i_1 i_2\ldots i_p} u_{i_1} u_{i_2} x_{i_3} \ldots x_{i_p}\Big| \leq \sum_{1\leq i_1,i_2,\ldots,i_p\leq N} |J_{i_1 i_2\ldots i_p}| |u_{i_1}| |u_{i_2}| \nonumber\\&= (p-2)!~ |\bm u|^\top \bm D_{\bm J_N} |\bm u| \lesssim_p \|\bm D_{\bm J_N}\|,
\end{align} 
since $|\bm u| := (|u_1|,\ldots,|u_N|) \in S^{N-1}$. Taking supremum over all $\bm u \in S^{N-1}$ followed by the supremum over all $\bm x \in \sa_N$ and further followed by the supremum over all $N \geq 1$ throughout \eqref{sh}, we have:
$$\sup_{N \geq 1} \sup_{\bm x \in \cC_N} ||\bm J_N(\bm x)|| \lesssim_p \sup_{N\geq 1}\|\bm D_{\bm J_N}\| < \infty.$$

Next, we verify condition (2) in Theorem \ref{chextension}. To this end, we need the following lemma: 

\begin{lem}\label{oddp} For every $p \geq 2$, under the assumptions of Corollary \ref{boundeddeg}, $ |F_N^{(3)}(0)| = O(N)$, where $F^{(3)}(0)$ denotes the third derivative of $F_N(\beta)$ at $\beta=0$. 
\end{lem}

The proof of the lemma is given below. First we show how it can be used to prove condition (2) in Theorem \ref{chextension}.  To begin with, note that condition (2) of Corollary \ref{boundeddeg} implies that 
\begin{align}\label{eq:FN_II}
F_N''(0) = \mathrm{Var}_0(H_N(\bs)) = \e_0 H_N^2(\bs) = (p!)^2\sum_{1 \leq i_1<\ldots<i_p \leq N} J_{i_1\ldots i_p}^2 = \Omega(N). 
\end{align} 
Hence, $\liminf_{N\rightarrow \infty} \frac{1}{N} F_N''(0) > 0$. Now, since by Lemma \ref{oddp}, $\limsup_{N\rightarrow \infty} \frac{1}{N} |F_N^{(3)}(0)| <\infty$, we can choose $\varepsilon > 0$ small enough, such that for all $N$ large enough, 
\begin{equation*}
\frac{\varepsilon}{6N} |F_N^{(3)}(0)| < \frac{F_N''(0)}{4N}.
\end{equation*} 
Therefore, because the fourth derivative $F_N^{(4)}(b) = \e_b H_N^4(\bs) \geq 0$ for all $b \geq 0$,  a Taylor expansion gives the following for all $\beta \in (0,\varepsilon)$:
\begin{equation*}%\label{largeN}
\frac{F_N(\beta)}{N} \geq \frac{\beta^2}{2N} F_N''(0) +  \frac{\beta^3}{6N} F_N^{(3)} (0) \geq \frac{\beta^2}{2N} F_N''(0) -  \frac{\beta^3}{6N} |F_N^{(3)} (0)| > \frac{\beta^2}{4N}F_N''(0) = \Omega(1),
\end{equation*}
where the last step uses \eqref{eq:FN_II}. This verifies  condition (2) of Theorem \ref{chextension} for all $\beta > 0$, by the monotonicity of $F_N$.

\subsubsection*{Proof of Lemma \ref{oddp}} To begin with observe that $F_N^{(3)}(0)=\e_0 H_N^3(\bs)$. Now, the proof of the lemma for odd $p$ is trivial. This is because, under $\p_0$, $\bs \stackrel{D}{=} - \bs$, and for odd $p$, $H_N(-\bs) = - H_N(\bs)$, which implies  $\e_0 H_N^3(\bs)$. Hence, we will assume that $p = 2 q$, for $q \geq 1$,  throughout the rest of the proof. Now, note that 
\begin{equation}\label{exh3}
\e_0H_N^3(\bs) = \sum_{\substack{1 \leq i_1, \ldots, i_p \leq N \\ \text{distinct}}}  \sum_{\substack{1 \leq j_1, \ldots, j_p \leq N \\ \text{distinct}}} \sum_{\substack{1 \leq k_1, \ldots, k_p \leq N \\ \text{distinct}}}   J_{i_1, \ldots, i_p } J_{j_1, \ldots, j_p }  J_{k_1, \ldots, k_p }   \e_0\left( \prod_{s=1}^p X_{i_s} X_{j_s} X_{k_s} \right). 
\end{equation}
Observe for each term in the sum above, the expectation is non-zero, if only if the multiplicity of each element in the multi-set $\{ i_1, \ldots, i_p \} \bigcup \{ j_1, \ldots, j_p\} \bigcup \{ k_1, \ldots, k_p\}$ is exactly $2$. This implies, the number of distinct elements in $\{ i_1, \ldots, i_p \} \bigcup \{ j_1, \ldots, j_p\} \bigcup \{ k_1, \ldots, k_p\}$ is  $3q$ and every pair of sets among  $\{ i_1, \ldots, i_p \}$, $\{ j_1, \ldots, j_p\}$, $\{ k_1, \ldots, k_p\}$  must have exactly $q$ elements in common. Therefore, from \eqref{exh3} and recalling the definition of the matrix $\bm D_{\bm J_N}$ from \eqref{eq:dJN} we get, 
\begin{align*}%\label{thirdm}
&\big|\e_0 H_N^3(\bs)\big|\\ & \lesssim_p \sum_{\substack{1 \leq i_1, \ldots, i_q, i_{q+1}, \ldots, i_{2 q}, i_{2q+1}, \ldots, i_{3q} \leq N \\ \textrm{distinct}}} \big|J_{ i_1, \ldots, i_q, i_{q+1}, \ldots, i_{2 q}} J_{i_{q+1}, \ldots, i_{2 q} i_{2q+1}, \ldots, i_{3q} } J_{ i_{2q+1}, \ldots, i_{3q} i_1, \ldots i_q}\big| \nonumber \\ 
& \lesssim_p \sum_{\substack{1 \leq i_1,i_{q +1}, i_{2q+1} \leq N \\ \textrm{distinct}}} \bm d_{\bm J_N}(i_1, i_{q+1}) \bm d_{\bm J_N}(i_{q+1}, i_{2q+1}) \bm d_{\bm J_N}(i_1, i_{2q+1})  \nonumber \\  
& = \mathrm{Trace} (\bm D_{\bm J_N}^3)  \leq N \|\bm D_{\bm J_N}\| = O(N),
\end{align*}
where the last step uses the assumption that $\sup_{N \geq 1}\|\bm D_{\bm J_N}\| < \infty$. \qed

\section{Proof of Theorem \ref{sbmthr}}\label{proof6}

We start by proving Theorem \ref{sbmthr} (1). To this end, it suffices to verify Theorem \ref{chextension}. Note that condition $(1)$ is easily satisfied because the bounded maximum degree condition \eqref{eq:DN_bound} holds almost surely: 
$$\max_{1 \leq i \leq N} d_{\bm J_N}(i)= \frac{1}{N^{p-1}} \max_{1 \leq i \leq N} d_{\bm A_{H_N}}(i) = O(1),$$
since $d_{\bm A_{H_N}}(i) \leq N^{p-1}$, for all $1 \leq i \leq N$, in any $p$-uniform hypergraph $H_N$. Next, we verify condition (2) of Theorem \ref{chextension}. To this end, by the lower bound in \cite[Theorem 1.6]{CD16} (which is the mean-field lower bound to the Gibbs variational representation of the log-partition function), we have 
\begin{align}\label{lowbdgibbs}
\e {F}_N(\beta) & \geq \E\left[\sup_{\bm x \in [-1,1]^N}\left\{ \beta  \sum_{1 \leq i_1, \ldots, i_p \leq N} \e(J_{i_1\ldots i_p}) x_{i_1}\ldots x_{i_p} - \sum_{i=1}^N I(x_i) \right\} \right] \nonumber \\ 
& \geq \sup_{\bm x \in [-1,1]^N}\left\{ \frac{\beta}{N^{p-1}} \sum_{1 \leq i_1, \ldots, i_p \leq N} \e({a}_{i_1\ldots i_p}) x_{i_1}\ldots x_{i_p} - \sum_{i=1}^N I(x_i) \right\} . 
\end{align}
Now, take any $\bm t := (t_1,\ldots,t_K) \in [-1,1]^K$, and define $\bm x \in [-1,1]^N$ by taking $x_i := t_j$, if $i \in \cB_j$, where $\cB_1, \cB_2, \ldots, \cB_K$ are as in Definition \ref{defn:block}. Then, the term inside the supremum in the RHS of \eqref{lowbdgibbs} equals $N \phi_\beta(t_1,\ldots,t_K) + O(1)$ (recall the definition of the function $\phi_\beta(t_1,\ldots,t_K)$ from \eqref{eq:threshold_function}).  Hence, \eqref{lowbdgibbs} gives us, 
\begin{equation}\label{blocklb1}
\frac{\e F_N(\beta)}{N} \geq \sup_{(t_1, t_2, \ldots, t_K) \in [0,1]^K} \phi_\beta(t_1,\ldots,t_K) + o(1).
\end{equation}
The bound in \eqref{blocklb1} above combined with the definition of the threshold $\beta_{\mathrm{HSBM}}^*$ in \eqref{eq:beta_threshold} and now implies that for all $\beta >\beta_{\mathrm{HSBM}}^*$, $\liminf_{N\rightarrow \infty} \frac{1}{N} \e F_N(\beta) > 0$. Then by Lemma \ref{mcdiarmid} below, it follows that $\liminf_{N \rightarrow \infty} \frac{1}{N}F_{N}(\beta) > 0$ with probability $1$. This verifies condition (2) of Theorem \ref{chextension}, and shows that the MPL estimate $\hat \beta_N (\bm X)$ is $\sqrt N$-consistent for $\beta >\beta_{\mathrm{HSBM}}^*$.

\begin{lem}\label{mcdiarmid}
	Let $F_N(\beta)$ denote the log-partition function of the $p$-tensor stochastic block model as in Theorem \ref{sbmthr}. Then, for every $\beta > 0$, the sequence $F_N(\beta) - \e F_N(\beta)$ is bounded in probability.  
\end{lem}

\begin{proof}
	To start with, note that $ {F}_N(\beta)$ is a function of the collection of i.i.d. random variables $\mathcal{A} := \{A_{i_1\ldots i_p}\}_{1\leq i_1<\ldots < i_p \leq N}$, and so, it is convenient to denote $ {F}_N(\beta)$ by $ {F}_{N,\beta}(\mathcal{A})$. Let us take $\mathcal{A}' := \{A'_{i_1\ldots i_p}\}_{1\leq i_1<\ldots< i_p\leq N}$, where $A'_{123\ldots p} = 1- A_{123\ldots p}$ and $A'_{i_1\ldots i_p} = A_{i_1\ldots i_p}$ for all $(i_1,\ldots,i_p) \neq (1,2,3,\ldots,p)$.  Note that 
	\begin{equation*}
	\left| \sum_{i_1<\cdots<i_p} A_{i_1\ldots i_p} X_{i_1}\cdots X_{i_p} - \sum_{i_1<\cdots<i_p} A_{i_1\ldots i_p}' X_{i_1}\cdots X_{i_p}\right|= |X_1 X_2\ldots X_p| = 1.
	\end{equation*} 
	Hence,  
	$$\exp\left\{\frac{\beta p! }{N^{p-1}} \sum_{i_1<\cdots<i_p} A_{i_1\ldots i_p} X_{i_1}\cdots X_{i_p}\right\} \leq e^{\frac{\beta p! }{N^{p-1}}}\exp\left\{\frac{\beta p!}{N^{p-1}} \sum_{i_1<\cdots<i_p} A_{i_1\ldots i_p}' X_{i_1}\cdots X_{i_p}\right\}.$$
	The above inequality implies that $ {F}_{N,\beta}(\mathcal{A}) \leq  {F}_{N,\beta}(\mathcal{A}') + \beta p! N^{1-p}$. Similarly, we also have $ {F}_{N,\beta}(\mathcal{A}') \leq  {F}_{N,\beta}(\mathcal{A}) + \beta p!  N^{1-p}$, and hence, $$\left| {F}_{N,\beta}(\mathcal{A}) -  {F}_{N,\beta}(\mathcal{A}')\right| \leq \beta p!  N^{1-p}.$$ Of course, the above arguments hold if $\mathcal{A}'$ is obtained by flipping any arbitrary entry of $\mathcal{A}$ (not necessarily the $(1,2,\ldots,p)$-th entry) and keeping all other entries unchanged. Hence, the assumption of McDiarmid's inequality \cite{mcdiarmid} holds with bounding constants $c_{i_1\ldots i_p} = \beta p! N^{1-p}$. Therefore, for every $t > 0$:
	$$\mathbb{P}\left(\left| {F}_N(\beta) - \mathbb{E}  {F}_N(\beta)\right| \geq t\right) \leq 2\exp\left\{-\frac{2t^2}{\sum_{ i_1<\ldots<i_p}c_{i_1\ldots i_p}^2}\right\} \leq 2\exp\left\{-\frac{2t^2 N^{p-2}}{\beta^2(p!)^2}\right\},$$ 
	which completes the proof of the lemma. 	 
\end{proof}

We will now use Lemma \ref{mcdiarmid} to prove Theorem \ref{sbmthr} (2). To this end, we will show that
\begin{align}\label{eq:FNbeta_H}
\e F_N(\beta) = O(1), \quad \text{ for } \beta < \beta_{\mathrm{HSBM}}^*, 
\end{align} 
the expectation in \eqref{eq:FNbeta_H} being taken with respect to the randomness of the HSBM. To see why this implies Theorem \ref{sbmthr} (2), assume, on the contrary, that there is a sequence of estimates which is consistent for $\beta < \beta_{\mathrm{HSBM}}^*$.  Using this sequence of estimates we can then construct a consistent sequence of tests $\{\phi_N\}_{N \geq 1}$ for the following hypothesis testing problem:\footnote{A sequence of tests $\{\phi_N\}_{N \geq 1}$ is said to be {\it consistent} 
	if both its Type I and Type II errors converge to zero as $N \rightarrow \infty$, that is, $\lim_{N\rightarrow\infty}\E_{H_0}\phi_N=0$, and the power $\lim_{N\rightarrow\infty}\E_{H_1}\phi_N=1$.}  
\begin{align}\label{eq:hypothesis} 
H_0: \beta = \beta_1\quad\quad\textrm{versus}\quad\quad H_1: \beta=\beta_2, 
\end{align}
if $\beta_1 < \beta_2 < \beta_{\mathrm{HSBM}}^*$. To this end, denote by $\mathbb{Q}_{\beta,p}$ the joint distribution of the HSBM and the $p$-tensor Ising model with parameter $\beta$.  Then a simple calculation shows that for any two positive real numbers $\beta_1 < \beta_2$, the Kullback-Leibler (KL) divergence between the joint measures $\mathbb{Q}_{\beta_1,p}$ and $\mathbb{Q}_{\beta_2,p}$ is given by:
\begin{align}\label{eq:DKL}
D_N(\mathbb{Q}_{\beta_1,p}\| \mathbb{Q}_{\beta_2,p}) = \e D_N(\p_{\beta_1,p}\| \p_{\beta_2,p}) =  \e F_N(\beta_2) - \e F_N(\beta_1) -(\beta_2-\beta_1) \e F_N'(\beta_1),
\end{align} 
where, as before, the expectation in \eqref{eq:DKL} is taken with respect to the randomness of the HSBM.  
Now, by the monotonicity of $F_N'(\cdot)$,		
$$0 = (\beta_2-\beta_1)F_N'(0) \leq (\beta_2-\beta_1)F_N'(\beta_1)\leq \int_{\beta_1}^{\beta_2} F_N'(t)\mathrm dt = F_N(\beta_2)-F_N(\beta_1).$$ 
Hence, by  \eqref{eq:FNbeta_H} and \eqref{eq:DKL}, $D_N(\mathbb{Q}_{\beta_1,p}\| \mathbb{Q}_{\beta_2,p}) = O(1)$. Then, by \cite[Proposition 6.1]{BM16}, there cannot exist any sequence of consistent tests for the hypothesis \eqref{eq:hypothesis}, which leads to a contradiction. This completes the proof of Theorem \ref{sbmthr} (2).

\subsection{Proof of (\ref{eq:FNbeta_H})} The proof of \eqref{eq:FNbeta_H} has the following two steps: 
\begin{itemize}
	\item[(I)] Define a new $p$-tensor Ising model on $N$ nodes, with interaction tensor $\tilde{\bm J}_N := \e {\bm J_N}$. We will call this model $\mathcal M_0 $. The first step in the proof of \eqref{eq:FNbeta_H} is to show that the log-partition function $\tilde{F}_N$ of the model $\mathcal M_0 $ is bounded, for every $\beta < \beta_{\mathrm{HSBM}}^*$.
	
	\item[(II)] The second step is to show that the expected log-partition function $\e F_N(\beta)$ of the original model is bounded, for $\beta < \beta_{\mathrm{HSBM}}^*$, by comparing it with the log-partition function $\tilde{F}_N$ of the model $\mathcal M_0 $. The result in \eqref{eq:FNbeta_H} then follows by an application of Lemma \ref{mcdiarmid}. 
\end{itemize}

\subsubsection{Proof of Step (I)}

Throughout this section we fix $\beta < \beta_{\mathrm{HSBM}}^*$ and denote by $\p_{\beta,\mathcal M_0 }$ the probability measure corresponding to the model $\mathcal M_0 $ at the parameter $\beta$ and $\e_{\beta, \mathcal M_0 }$  the expectation with respect to the probability measure $\p_{\beta,\mathcal M_0 }$.

\begin{lem}\label{lm:logpartition_I} Denote by $\tilde{F}_N(\beta)$ the log-partition function of the model $\cM_0$. Then for $\beta < \beta_{\mathrm{HSBM}}^*$, $\limsup_{N \rightarrow \infty} \tilde{F}_N(\beta) < \infty $. 
\end{lem}

\noindent {\it Proof of Lemma} \ref{lm:logpartition_I}:  Denote the Hamiltonian of the model $\mathcal M_0 $ by $\tilde{H}_N(\bs)$, that is, 
$$\tilde{H}_N(\bs) : = \frac{1}{N^{p-1}} \sum_{1 \leq i_1, i_2, \ldots, i_p \leq N} \E(a_{i_1 i_2 \ldots i_p}) X_{i_1} X_{i_2} \ldots X_{i_p}.$$  
For each $\bs \in \sa_N$ and $1 \leq j \leq K$, define $S_j(\bs):= \sum_{i \in \cB_j} X_i$. With these notations, we have $\tilde{H}_N(\bs) = \overline{H}_N(\bs) + O(1)$, where
\begin{align}\label{eq:HN_tensor}
\overline{H}_N(\bs) = \frac{1}{N^{p-1} }  \sum_{1 \leq j_1,\ldots, j_p \leq K} \theta_{j_1\ldots j_p}  \prod_{\ell=1}^{p} S_{j_\ell}(\bs).
\end{align}  
Let us define $\overline{Z}_N(\beta) := \frac{1}{2^N} \sum_{\bs \in \sa_N} e^{\beta \overline{H}_N(\bs)}$. Since $\tilde{F}_N(\beta) = \log \overline{Z}_N(\beta) + O(1)$, it suffices to show that $\overline{Z}_N(\beta) = O(1)$. Towards this, for each $1 \leq j\leq K$, define the sets:
$$I_j:= \left\{-1,~ -1+\frac{2}{|\cB_j|},~ -1+ \frac{4}{|\cB_j|},~\ldots,~ 1-\frac{2}{|\cB_j|},~ 1\right\}$$ and $A_s(j) := \{\bs \in \sa_N: S_j(\bs) = s\}.$
Recall, $\cB_j=(N \sum_{i=1}^{j-1}\lambda_i, N \sum_{i=1}^j \lambda_i] \bigcap [N]$, hence $| |\cB_j| - N \lambda_j |\leq 2$. Now, note that 
\begin{align}
\overline{Z}_N(\beta) & = \frac{1}{2^N} \sum_{(\ell_1,\ldots,\ell_K) \in I_1 \times \cdots \times I_K} e^{\frac{\beta}{N^{p-1}} \sum_{1\leq j_1,\ldots,j_p\leq K} \theta_{j_1 \ldots j_p} \prod_{m=1}^p \ell_{j_m} |\cB_{j_m}|} \left| \bigcap_{j=1}^K A_{\ell_j |\cB_j|}(j)\right| \nonumber \\ 
\label{fe} &\lesssim \frac{1}{2^N} \sum_{(\ell_1,\ldots,\ell_K) \in I_1 \times \cdots \times I_K} e^{N\beta \sum_{1\leq j_1,\ldots,j_p\leq K} \theta_{j_1 \ldots j_p} \prod_{m=1}^p \lambda_{j_m} \ell_{j_m}} \left| \bigcap_{j=1}^K A_{\ell_j |\cB_j|}(j)\right|\\ 
\label{eq:T12} &= T_1 + T_2 ,
\end{align}
where the term $T_1$ is obtained by restricting the sum in the RHS of \eqref{fe} to the set $(I_1\times\cdots \times I_K) \bigcap [-\frac{1}{2} , \frac{1}{2}]^K$ and the term $T_2$ is the sum restricted to the set $(I_1\times\cdots \times I_K) \bigcap ([-\frac{1}{2} , \frac{1}{2}]^K)^c$. 

Let us bound $T_1$ first. Note that  
$$|A_{\ell_j |\cB_j|}(j)|= {|\cB_j| \choose \frac{|\cB_j|(1+\ell_j)}{2}}.$$
Then by the Stirling's approximation of the binomial coefficient (see, for example, \cite[Lemma B.5]{mlepaper}) and using the fact that the sets $\cB_1, \cB_2, \ldots, \cB_K$ are disjoint, we have for all $(\ell_1,\ldots,\ell_K) \in [-\frac{1}{2} , \frac{1}{2}]^K$,
$$\left| \bigcap_{j=1}^K A_{\ell_j |\cB_j|}(j)\right| = 2^N \exp\left\{-\sum_{j=1}^K |\cB_j| I(\ell_j) \right\} O\left(\frac{1}{\sqrt{\prod_{j=1}^K |\cB_j|}}\right).$$
Hence, denoting $\bm \ell := (\ell_1,\ldots,\ell_K)$ and $\cI := (I_1 \times \cdots \times I_K)$ gives, 
\begin{align*}
T_1 &=  \sum_{\bm \ell \in \cI \bigcap \left[-\frac{1}{2} , \frac{1}{2}\right]^K} e^{ N\beta \sum_{1\leq j_1,\ldots,j_p\leq K} \theta_{j_1 \ldots j_p} \prod_{m=1}^p \lambda_{j_m} \ell_{j_m} -\sum_{j=1}^K |\cB_j| I(\ell_j) } O\left(\frac{1}{\sqrt{\prod_{j=1}^K |\cB_j|}}\right)\\ 
&\lesssim \sum_{\bm \ell \in \cI \bigcap \left[-\frac{1}{2} , \frac{1}{2}\right]^K} e^{ N\beta \sum_{1\leq j_1,\ldots,j_p\leq K} \theta_{j_1 \ldots j_p} \prod_{m=1}^p \lambda_{j_m} \ell_{j_m} - N\sum_{j=1}^K \lambda_j I(\ell_j) } O\left(\frac{1}{\sqrt{\prod_{j=1}^K |\cB_j|}}\right)\\&= \sum_{\bm \ell \in \cI\bigcap \left[-\frac{1}{2} , \frac{1}{2}\right]^K} e^{N\phi_\beta(\ell_1,\ldots,\ell_K) } O\left(\frac{1}{\sqrt{\prod_{j=1}^K |\cB_j|}}\right)
\end{align*} 
Now, since $\beta < \beta_{\mathrm{HSBM}}^*$ we can choose $\varepsilon \in (0,1)$, such that $\beta/(1 - \varepsilon) < \beta_{\mathrm{HSBM}}^*$. Hence, recalling \eqref{eq:beta_threshold}, 
$$\phi_{\beta/(1-\varepsilon)}(\ell_1,\ldots,\ell_K) \leq 0,\quad\textrm{that is, }\quad
\phi_{\beta}(\ell_1,\ldots,\ell_K) \leq -\varepsilon \sum_{j=1}^K \lambda_j I(\ell_j) \leq -\frac{\varepsilon}{2} \sum_{j=1}^K \lambda_j \ell_j^2,$$
where the last step uses $I(x) \geq x^2/2$.  Then by a Riemann sum approximation (see \cite[Lemma B.2]{mlepaper}), 
\begin{align} 
\label{eq:T11}T_1 &\lesssim \sum_{\bm \ell \in \cI\bigcap \left[-\frac{1}{2} , \frac{1}{2}\right]^K} e^{-\frac{N\varepsilon}{2} \sum_{j=1}^K \lambda_j \ell_j^2 } O\left(\frac{1}{\sqrt{\prod_{j=1}^K |\cB_j|}}\right)  \\ 
&\leq  \frac{1}{2^{K}} O\left(\sqrt{\prod_{j=1}^K |\cB_j|}\right) \int\limits_{\mathbb{R}^K} \exp\left\{-\frac{N\varepsilon}{2} \sum_{j=1}^K \lambda_j x_j^2 \right\}~\mathrm dx_1\ldots \mathrm dx_K + O(1) \nonumber \\ 
&= \frac{1}{2^{K}} O\left(\sqrt{\prod_{j=1}^K \frac{|\cB_j|}{N}}\right) \int\limits_{\mathbb{R}^K} \exp\left\{-\frac{\varepsilon}{2} \sum_{j=1}^K \lambda_j y_j^2 \right\}~\mathrm dy_1\ldots \mathrm dy_K + O(1) \nonumber \\ 
\label{eq:T1} &\lesssim \int\limits_{\mathbb{R}^K} \exp\left\{-\frac{\varepsilon}{2} \sum_{j=1}^K \lambda_j y_j^2 \right\}~\mathrm dy_1\ldots \mathrm dy_K + O(1) = O(1) .
\end{align} 

Now, we bound $T_2$. For this we need the following combinatorial estimate:
\begin{obs}\cite[Equation (5.4)]{talagrand} For every integer $s$ and positive integer $m$, 
	\begin{equation}\label{new}
	\left|\left\{\bm x \in \{-1, 1\}^m: \sum_{i=1}^m x_i =  s\right\}\right| \leq 2^{m} e^{-m I\left(\frac{s}{m}\right)}. 
	\end{equation}  
\end{obs}	

Using the bound in \eqref{new} and recalling that the sets $\cB_1, \cB_2, \ldots, \cB_K$ are disjoint, we have for every $(\ell_1, \ell_2, \ldots, \ell_K) \in (I_1, I_2, \ldots, I_K)$, 
\begin{equation*}
\left| \bigcap_{j=1}^K A_{\ell_j |\cB_j|}(j)\right|  \leq 2^N   \exp\left\{-\sum_{j=1}^K |\cB_j| I(\ell_j)\right\}.
\end{equation*}
Hence, by following the arguments used to obtain \eqref{eq:T11} we get,  
\begin{align}\label{eq:T2}
T_2 \lesssim \sum_{\bm \ell \in \cI\bigcap \left(\left[-\frac{1}{2} , \frac{1}{2}\right]^K\right)^c} e^{-\frac{N\varepsilon}{2} \sum_{j=1}^K \lambda_j \ell_j^2 } &\leq \left(\prod_{j=1}^K (|\cB_j|+1)\right)\max_{\bm \ell \in \cI\bigcap \left(\left[-\frac{1}{2} , \frac{1}{2}\right]^K\right)^c}  e^{-\frac{N\varepsilon}{2} \sum_{j=1}^K \lambda_j \ell_j^2 } 
\nonumber\\
&\leq \left(\prod_{j=1}^K(N\lambda_j+1)\right)\max_{1\leq j\leq K}  e^{-\frac{N\varepsilon \lambda_j}{8}} = o(1). 
\end{align}
This shows, combining \eqref{eq:T12}, \eqref{eq:T1}, and \eqref{eq:T2},  that $\overline{Z}_N(\beta) = O(1)$ for all $\beta < \beta_{\mathrm{HSBM}}^*$, completing the proof of Lemma \ref{lm:logpartition_I}.

\subsubsection{Proof of Step {(II)}}  We now show that $\e F_N(\beta)$ is bounded, for $\beta < \beta_{\mathrm{HSBM}}^*$, which will allow us to conclude, using Lemma \ref{mcdiarmid}, that  $F_N(\beta) = O(1)$ with probability $1$, for all $\beta < \beta_{\mathrm{HSBM}}^*$, that is, \eqref{eq:FNbeta_H} holds. 

\begin{lem}\label{st2} For every $\beta < \beta_{\mathrm{HSBM}}^*$, $\limsup_{N\rightarrow \infty} \e F_N(\beta) < \infty$. 
\end{lem}

\begin{proof} Fix $\beta < \beta_{\mathrm{HSBM}}^*$. Then the partition function $Z_N(\beta)$ becomes: 
	\begin{align*}
	Z_N(\beta) &= \frac{1}{2^N} \sum_{\bs \in \sa_N} e^{\frac{p! \beta}{N^{p-1}} \sum_{i_1<\cdots<i_p} A_{i_1\ldots i_p} X_{i_1}\cdots X_{i_p}}\\
	&= \frac{1}{2^N} \sum_{\bs \in \sa_N} e^{\frac{p! \beta}{N^{p-1}}  \sum_{i_1<\cdots<i_p} \left(A_{i_1\ldots i_p} - \e A_{i_1\ldots i_p}\right) X_{i_1}\cdots X_{i_p}}e^{\frac{p! \beta}{N^{p-1}}  \sum_{i_1<\cdots<i_p} \e A_{i_1\ldots i_p}  X_{i_1}\cdots X_{i_p}}.
	\end{align*}
	By Hoeffding's inequality, for each $\bs \in \sa_N$ and $i_1<\ldots<i_p$, 
	\begin{align*}
	\mathbb{E} \exp\left\{\frac{p! \beta}{N^{p-1}} \sum_{i_1<\cdots<i_p} \left(A_{i_1\ldots i_p} - \e A_{i_1\ldots i_p}\right) X_{i_1}\cdots X_{i_p}\right\}  
	&\leq \prod_{i_1<\cdots<i_p}\exp\left\{\frac{\beta^2(p!)^2}{8 N^{2p-2}}\right\}\\
	&\leq \exp\left\{\frac{\beta^2 p!}{8 N^{p-2}}\right\}.
	\end{align*}
	This shows that 
	$$\mathbb{E} Z_N(\beta) \leq \exp\left\{\frac{\beta^2 p!}{8 N^{p-2}}\right\} \left[\frac{1}{2^N} \sum_{\bs \in \sa_N}e^{\frac{p! \beta}{N^{p-1}}\sum_{i_1<\cdots<i_p} \e A_{i_1\ldots i_p} X_{i_1}\cdots X_{i_p}}\right] = \exp\left\{\frac{\beta^2 p!}{8 N^{p-2}}\right\}\tilde{Z}_N(\beta),$$
	where $\tilde{Z}_N$ is the partition function of the model $\mathcal M_0 $. Now, taking logarithms and using Lemma \ref{lm:logpartition_I} shows, $\limsup_{N \rightarrow \infty}\log \mathbb{E} Z_N(\beta) < \infty$, for $\beta < \beta_{\mathrm{HSBM}}^*$. Then, by Jensen's inequality, we conclude that $\limsup_{N \rightarrow \infty}\E F_N(\beta) < \infty$, for $\beta < \beta_{\mathrm{HSBM}}^*$, completing the proof of the lemma.  
\end{proof}

\section{Proofs from Section \ref{cltsec}}\label{sec:pfcwclt}

\subsection{Proof of Theorem \ref{thm:cwmplclt}}\label{sec:pf_cwmplclt}

Define the function $\phi_p: [-1,1]\mapsto (-\infty,\infty]$ as:  
\begin{align}\label{eq:2} 
\phi(t) = \phi_p(t) := 
\begin{cases}
p^{-1} t^{1-p}\tanh^{-1} (t) &\quad\text{if}~p~\textrm{is even and}~t \neq 0,\\
p^{-1} t^{1-p}\tanh^{-1} (t) &\quad\text{if}~p~\textrm{is odd and}~t > 0.\\
\infty &\quad\text{if}~p~\textrm{is odd and}~t < 0.\\
0 &\quad\text{if}~t = 0.\\
\end{cases}
\end{align}
Note that for every $t \neq 0$ when $p$ is even, and every $t > 0$ when $p$ is odd, the function $\vp$ is twice differentiable on some compact set containing $t$ in its interior, and $$\vp'(t) = -\frac{1}{pt^{p-1}}\left\{\frac{(p-1)\tanh^{-1}(t)}{t} -\frac{1}{1-t^2}\right\}.$$
Hence,  recalling the definition of the function $g$ from the statement of the theorem we have, 
\begin{equation}\label{derexp}
\vp'(m_*) = -\frac{g''(m_*)}{pm_*^{p-1}}.
\end{equation}
Moreover, from the definition of the MPLE in \eqref{mple} it is easy to see that in the Curie-Weiss model $\hat{\beta}_N(\bm X) = \phi(\cs)$. Note that $\vp(m_*) = \beta$, since $g'(m_*) = 0$. This implies, 
\begin{equation}\label{delstep}
\sqrt{N}(\hat{\beta}_N(\bm X) - \beta) = \sqrt{N}\left(\vp(\cs) - \vp(m_*)\right).
\end{equation}
We now consider the case $p$ is even and $p$ is odd separately. Throughout, we assume $\beta > \beta_{\mathrm{CW}}^*$: 
\begin{itemize}
	
	\item {\it $p \geq 3$ is odd}: In this case, $m_*$ is the unique global maximizer of the function  $g$ on $[-1, 1]$. This implies, by \cite[Theorem 2.1~(1)]{mlepaper}, 
	$$\sqrt{N}(\bar X_N - m_*) \dto N\left(0, -\frac{1}{g''(m_*)}\right). $$
	Hence, by \eqref{delstep} and the delta method \cite[Theorem 1.8.12]{tpe}, 
	$$\sqrt{N}(\hat{\beta}_N(\bm X) -\beta) ) \xrightarrow{D} N\left(0, -\frac{(\phi'(m_*))^2}{g''(m_*)}\right) \stackrel{D}= N\left(0, -\frac{g''(m_*)}{p^2m_*^{2p-2} }\right),$$
	where the last step uses \eqref{derexp}. This completes the proof when $p$ is odd. 
	
	\item {\it $p \geq 2$ is even}: In this case, the function  $g$ has two (non-zero) global maximizers on $[-1, 1]$, which are given by $m_*$ and $-m_*$ (as shown in \cite[Section C]{mlepaper}). This implies, by \cite[Theorem 2.1~(2)]{mlepaper} and the delta method, 
	\begin{align}\label{eq:curie_weiss_mple_I}
	\sqrt{N}(\hat{\beta}_N(\bm X) -\beta ) \Big| \{\cs \in (0, 1] \} \xrightarrow{D} N\left(0, -\frac{g''(m_*)}{p^2m_*^{2p-2} }\right). 
	\end{align}
	Similarly, observing that $\phi(m_*)=-\phi(-m_*)$ and $g''(m_*)=g''(-m_*)$, 
	\begin{align}\label{eq:curie_weiss_mple_II}
	\sqrt{N}(\hat{\beta}_N(\bm X) -\beta ) \Big| \{ \cs \in [-1, 0] \} \xrightarrow{D} N\left(0, -\frac{g''(m_*)}{p^2m_*^{2p-2} }\right). 
	\end{align}
	Combining \eqref{eq:curie_weiss_mple_I} and \eqref{eq:curie_weiss_mple_II}, gives the desired result when $p$ is even.

\end{itemize}

\subsection{Proof of Theorem \ref{thm:cw_threshold}}\label{sec:pf_cw_threshold} Fix $\beta = \beta_{\mathrm{CW}}^*(p)$.  We now consider the three cases in Theorem \ref{thm:cw_threshold} separately. 

\noindent{\it Proof of }(1): In this case, $p=2$ and, hence, by \eqref{eq:2}, $\hat{\beta}_N = \frac{\tanh^{-1} (\cs)}{2\cs}  \bm 1 \{\cs \neq 0\}$. Therefore, on the event $\cE_N:=\{\cs\neq 0, |\cs|< \frac{1}{2}\}$, 
\begin{align}\label{inf1}
N^\frac{1}{2}\left(\hat{\beta}_N -\frac{1}{2}\right) &= N^\frac{1}{2} \left(\frac{\tanh^{-1} (\cs)}{2\cs} - \frac{1}{2}\right)\nonumber\\ &= \tfrac{1}{2} N^\frac{1}{2}\sum_{s=1}^\infty \frac{\cs^{2s}}{2s+1} = \tfrac{1}{6} N^\frac{1}{2} \cs^2 + N^\frac{1}{2} O(\cs^4). 
\end{align}
From \cite[Proposition 4.1]{comets} we know that $N^\frac{1}{4} \cs \xrightarrow{D} F$, where $F$ is as defined in the statement of Theorem \ref{thm:cw_threshold}. The result in \eqref{eq:threshold_F} now follows from \eqref{inf1}, and the observation that $\p_\beta(\cE_N^c) = o(1)$, since $\cs \pto 0$ and $\p_\beta(\cs =0 )=o(1)$ (from the proof of \cite[Lemma C.6]{mlepaper}). \\ 

\noindent{\it Proof of }(2): Assume $p \geq 3$. To begin with, define the the three intervals $A_1 := [-1,~-\frac{m_*}{2}]$ and $A_2 := (-\frac{m_*}{2},~\frac{m_*}{2})$ and $A_3 := [\frac{m_*}{2},~1]$ (recall that $m_*=m_*(p, \beta)$ is the unique positive maximizer of the function $g(t) := \beta t^p - I(t)$. We now consider the following two cases depending on whether $p \geq 3$ is odd or even: 

\begin{itemize}
	\item[$\bullet$] {\it $p \geq 4$ is even:}  In this case, the function  $g$ has three global maximizers on $[-1, 1]$, which are given by $m_* > 0$ and $-m_*$, and $0$. Then, by \cite[Theorem 2.1 (2)]{mlepaper} and arguments as in \eqref{eq:curie_weiss_mple_I}, we have 
	\begin{equation}\label{evp}
	\sqrt{N}(\hat{\beta}_N(\bm X) - \beta)\big|\{\cs \in A_i\} \xrightarrow{D} N\left(0, -\frac{g''(m_*)}{p^2m_*^{2p-2}}\right), 
	\end{equation}
	for $i \in \{1, 3\}$. Next, recalling $\beta_N(\bm X) =\phi(\bar X_N)$, where $\phi(\cdot)$ as in \eqref{eq:2}, note that 
	\begin{equation}\label{evp2}
	\hat{\beta}_N(\bm X)\big|\{\cs \in A_2 \backslash \{0\}\} \xrightarrow{D} \infty,
	\end{equation}
	since $\frac{\tanh^{-1}(t)}{t^{p-1}} \rightarrow \infty$, as $|t| \rightarrow 0$. Now, since $\p (\cs=0) = o(1)$ and 
	$\p_\beta(\cs \in A_2) \rightarrow \alpha$ (by \cite[Theorem 2.1 (2)]{mlepaper}), combining \eqref{evp} and \eqref{evp2} the result in \eqref{eq:threshold_mple} follows, if $p \geq 4$ is even.  
	
	\item[$\bullet$] {\it $p \geq 3$ is odd:}   In this case, the function  $g$ has two global maximizers on $[-1, 1]$ one of which is non-positive and the other is $m_*> 0$. Then, by similar arguments as above, 
	\begin{equation}\label{odp}
	\sqrt{N}(\hat{\beta}_N(\bm X) - \beta)\big|\{\cs \in A_3\} \xrightarrow{D} N\left(0, -\frac{g''(m_*)}{p^2m_*^{2p-2}}\right). 
	\end{equation}
	Moreover, recalling \eqref{eq:2}, 
	\begin{equation}\label{odp1}
	\hat{\beta}_N(\bm X)\big|\{\cs \in \left(A_1\bigcup A_2 \right) \backslash \{0\}\} \xrightarrow{D} \infty.
	\end{equation}
	The result in \eqref{eq:threshold_mple} now follows from \eqref{odp}, \eqref{odp1}, and the fact that $\p_\beta(\cs \in A_1\bigcup A_2) \rightarrow \alpha$, if $p$ is odd (by \cite[Theorem 2.1 (2)]{mlepaper}).  \\ 
\end{itemize} 

\noindent{\it Proof of }(3): Here, we prove the finer asymptotics of $\hat{\beta}_N(\bm X)$. For this, note 
by \cite[Theorem 2.1 (2)]{mlepaper} that $\sqrt{N}\cs \big|\{\cs \in \cB\} \xrightarrow{D} N(0,1)$ for any interval $\cB$ containing $0$, but no other maximizer of $g$. We now consider the following two cases: 
\begin{itemize}
	\item [--] {\it $p\geq 4$ is even:} In this case, $\sqrt{N}\cs \big|\{\cs \in A_2\} \xrightarrow{D} N(0,1)$. Then, \eqref{eq:2} gives, 
	\begin{equation*}%\label{evenmple}
	N^{1-\frac{p}{2}}\hat{\beta}_N(\bm X) = \frac{1}{p(\sqrt{N}\cs)^{p-2}}\cdot \frac{\tanh^{-1}(\cs)}{\cs} \bm 1\{\cs \neq 0\}.
	\end{equation*}
	Now, since $\frac{\tanh^{-1}(\cs)}{\cs} \bm 1 \{\cs \neq 0\} \big|\{\cs \in A_2\setminus\{0\} \} \xrightarrow{P} 1$, we have,   
	\begin{equation}\label{eq:threshold_II}
	N^{1-\frac{p}{2}} \hat{\beta}_N(\bm X) \big|\{\cs \in A_2\setminus \{0\} \} \xrightarrow{D} \frac{1}{p Z^{p-2}}, 
	\end{equation} 
	where $Z \sim N(0, 1)$. The result in now \eqref{eq:threshold_mple_I} follows from \eqref{eq:threshold_II}, by noting $\p (\cs=0) = o(1)$, that $\p_\beta(\cs \in A_2) \rightarrow \alpha$ and $N^{1-\frac{p}{2}} \hat{\beta}_N(\bm X) \big|\{\cs \in A_1\bigcup A_3\} \xrightarrow{P} 0$ (by  \eqref{evp}).

	\item[--] {\it $p\geq 3$ is odd:} In this case, since the function  $g$ has two global maximizers on $[-1, 1]$ one of which is non-positive and the other is $m_*> 0$, $\sqrt{N}\cs |\{\cs \in A_1\bigcup A_2\} \xrightarrow{D} N(0,1)$. Then, \eqref{eq:2} gives, 
	\begin{equation}\label{eq:beta_p_mple}
	N^{1-\frac{p}{2}}\hat{\beta}_N(\bm X) = \frac{1}{p(\sqrt{N}\cs)^{p-2}}\cdot \frac{\tanh^{-1}(\cs)}{\cs} \bm 1\{\cs > 0\} + \infty \bm 1\{\cs < 0\},
	\end{equation}
	where we adopt the convention infinity times zero is zero. It follows from \eqref{eq:beta_p_mple}, that $N^{1-p/2}\hat{\beta}_N(\bm X)$ is a non-negative random variable. Hence, denoting $\cA_{12}:=\{\cs \in A_1\bigcup A_2\}$ and taking $t \in [0,\infty)$ gives 
	\begin{align*}
	\p\left(N^{1-\frac{p}{2}} \hat{\beta}_N(\bm X)\leq t \big|\cA_{12} \right) 
	&= \p\left(N^{1-\frac{p}{2}} \hat{\beta}_N(\bm X)\leq t,~\cs > 0 \big|\cA_{12}  \right) + o(1), %\label{eq:in2}
	\end{align*}
	since $\p_\beta(\cs = 0| \cA_{12}) = o(1)$. Now, note that 
	\begin{align*}
	& \p\left( N^{1-\frac{p}{2}} \hat{\beta}_N(\bm X) \leq t,~\cs > 0 \big|\cA_{12}  \right) \nonumber \\ 
	&= \p\left( \frac{1}{p(\sqrt{N}\cs)^{p-2}}\cdot \frac{\tanh^{-1}(\cs)}{\cs} \leq t,~\cs \neq 0 \big| \cA_{12}  \right) - \p_\beta(\cs < 0 \big| \cA_{12}) ,
	\end{align*}
	since on the the event $\{\cs < 0\}$, $\frac{1}{p(\sqrt{N}\cs)^{p-2}}\cdot \frac{\tanh^{-1}(\cs)}{\cs} < 0 \leq t$. Note that $\p_\beta(\cs < 0 \big| \cA_{12}) \rightarrow \frac{1}{2}$, since $\sqrt{N}\cs | \cA_{12}  \xrightarrow{D} N(0,1)$. Then, by arguments as in \eqref{eq:threshold_II}, it follows that 
	\begin{align*}
	\p\left( N^{1-\frac{p}{2}} \hat{\beta}_N(\bm X) \leq t \big|\cA_{12}  \right)  \rightarrow  \p\left(\frac{1}{pZ^{p-2}} \leq t\right) -\frac{1}{2} 
	& = \frac{1}{2}\p\left(\frac{1}{p|Z|^{p-2}} \leq t\right).
	\end{align*}
	Therefore, by \eqref{eq:beta_p_mple},
	\begin{equation}\label{eq:threshold_p}
	N^{1-\frac{p}{2}} \hat{\beta}_N(\bm X) \big|\cA_{12} \xrightarrow{D} \frac{1}{2} \left(\frac{1}{p |Z|^{p-2}}\right) +\frac{1}{2}\delta_\infty.
	\end{equation}
	The result in \eqref{eq:threshold_mple_II} follows now from \eqref{eq:threshold_p}, by noting that $\p_\beta(\cs \in A_1\bigcup A_2) \rightarrow \alpha$ and $N^{1-\frac{p}{2}} \hat{\beta}_N(\bm X) \big|\{\cs \in A_3\} \xrightarrow{P} 0$ (by \eqref{odp}). 
\end{itemize} 

%\fi

\chapter{High-Dimensional Logistic Regression with Dependent Observations}\footnotetext{This chapter is a (ongoing) joint work with Sagnik Halder, George Michailidis and Bhaswar B. Bhattacharya}\label{chap:covariateising}

The most popular way of modeling the probability of events with binary outcomes in statistics, is the logistic regression \cite{appliedlogistic}. Given predictors $\bm Z_1,\ldots, \bm Z_N \in \mathbb{R}^d$ and independent response variables $X_1,\ldots,X_N \in \{-1,1\}$, the logistic regression model is given by:
\begin{equation}\label{logisticbasic}
\log \left[\frac{\p(X_i = 1)}{1-\p(X_i=1)}\right] = \bh^\top \bm Z_i\quad(\textrm{for all} ~1\le i\le N)
\end{equation}
where $\bh \in \mathbb{R}^d$ is the vector of regression coefficients. It follows from \eqref{logisticbasic} that the joint distribution of $\bm X := (X_1,\ldots,X_N)$ is given by:
\begin{equation}\label{jointlogistic}
\p(\bm X = \bm x) = \prod_{i=1}^N \frac{1}{1+\exp(-\bh^\top x_i \bm Z_i)} ~\propto~ \exp\left(\frac{1}{2} \bh^\top \sum_{i=1}^N x_i \bm Z_i\right)
\end{equation}
Results on consistency and asymptotic normality for the maximum likelihood estimates of the parameters in logistic regression and generalized linear models were given in \cite{haberman1974,haberman1977}, \cite{monfort} and \cite{kaufmann}. It is well known that in the vanilla logistic regression, the parameter $\bh$ can be estimated within $L^2$ error $O_d(N^{-1/2})$. Classical results in logistic regression mainly deal with the regime when $d$ is fixed. However, in most modern scientific applications, datasets have a large number of features, which necessitates analysis outside the $N \gg d$ regime. Performance of the maximum likelihood estimate for the logistic regression in the regime $N = \Theta(d)$ has been studied recently in \cite{sur1,sur2,sur3}, where it is shown that in this regime, the ML estimate is not even unbiased. Adding a regularizer to the negative log-likelihood makes recovery of the parameter vector possible even when the MLE does not exist due to an inadequate sample size. There has been a significant amount of statistical literature on regularized logistic regression (\cite{florentina}, \cite{kakade}, \cite{vandegeer008}, \cite{salehi}), which often require the parameter to have some structure, for example sparsity of order $o(d)$.  

The independence assumption on the binary response variables is violated in many real-life scenarios. For example, the health status (healthy or ill) of individuals in an epidemic network, which may depend upon a number of personal attributes such as his immunity, age, weight, diet and smoking habit, are highly dependent. The health status of a single individual depends not only on his own health attributes, but also on the health status of other persons in the network he came into contact with. The vanilla logistic regression model \eqref{jointlogistic} can be generalized to capture dependency arising from a network with adjacency matrix $\bm A := ((A_{ij}))_{1\le i,j\le N}$, by introducing a quadratic interaction term in the probability mass function in \eqref{jointlogistic} as follows:
\begin{equation}\label{eq:logisticmodel}
\p_{\beta,\bh} (\bm X = \bm x) ~\propto ~ \exp\left(\bh^\top \sum_{i=1}^N x_i \bm Z_i + \frac{\beta}{2} \bm x^\top \bm A \bm x\right)\quad (\bm x \in \{-1,1\}^N)
\end{equation} 
The model \eqref{eq:logisticmodel} is in fact, an Ising model with varying external magnetic fields (the magnetic field at site $i$ being ${\bh^\top \bm Z_i}$).

 The model \eqref{eq:logisticmodel} has been studied in \cite{cd_ising_I}, where the authors showed that the parameters $\beta, \bh$ can be estimated to within error $O_d(N^{-1/2})$ by the maximum pseudolikelihood approach, under certain assumptions on the underlying network $\bm A$. The norm conditions imposed on $\bm A$ in \cite{cd_ising_I} also appear in \cite{pg_sm}, whose model is a special case of \eqref{eq:logisticmodel} with $d=1$ and $Z_1=\ldots = Z_N$. The dependence of the rate of convergence of the MPLE on $d$ is not made explicit in \cite{cd_ising_I}, which is an issue, if $d$ is allowed to grow with $N$. Keeping this in mind, we use a penalized version of the maximum pseudolikelihood approach described in Chapter \ref{ch:generalmple} to estimate the parameter $(\beta,\bh^\top)$. More precisely, we define
 
 \begin{equation}\label{eq:pseudolikelihood}
 (\hat{\beta},\hat{\bh}^\top) :=  \underset{(\tilde{\beta},\tilde{\bh}^\top) \in \mathbb{R}^{d+1}}{\mathrm{argmin}} L_N(\tilde{\beta},\tilde{\bh}) + \lambda \norm{(\tilde{\beta},\tilde{\bh})}_1 
 \end{equation}
 where $\lambda > 0$ is a tuning parameter to be chosen suitably, and
 \begin{align*}
 L_N(\tilde{\beta},\tilde{\bh}) :&= -\frac{1}{N}\sum_{i=1}^N \log \p_{\tilde{\beta}, \tilde{\bh}} \left(X_i \big| (X_j)_{j\ne i}\right)\\&= \log 2 - \frac{1}{N} \sum_{i=1}^N \left[X_i \left(\tilde{\beta} \sum_{j=1}^N A_{ij} X_j + \tilde{\bh}^\top \bm Z_i\right) - \log \cosh \left(\tilde{\beta} \sum_{j=1}^N A_{ij} X_j + \tilde{\bh}^\top \bm Z_i\right)\right]~.
 \end{align*}
 
 The regularization approach has been used in \cite{ising_nonconcave} and \cite{ravikumarlaff} in the context of structure recovery in Ising models, i.e. learning the matrix $\bm A$, which in their setup, is the unknown parameter.
 
 Consistency of the MLE has been shown to hold in both low dimensions ($d$ is fixed) as well as high dimensions ($N,d\rightarrow\infty$) under the assumption that the true parameter $\bh$ is $s$-sparse. Specifically, in case of high dimensions, sparsity in the ML estimator is induced using an $L^1$-penalized LASSO approach, and the optimal rates of  consistency have been obtained in \cite{vandegeer008}. In this chapter, we show that in the model \eqref{eq:logisticmodel}, if the parameter $\bh$ is sparse, then as long as $d$ grows a little slower than $\sqrt{N}$, the MPLE of $(\beta, \bh^\top)$ converges to the true parameter at rate $O(\sqrt{\log d/N})$.

\section{Main Result}\label{sec:covmain}
In this section, we state the main result of this chapter, which states that as long as the true parameter $\bh$ is sparse, the penalized MPLE \eqref{eq:pseudolikelihood} converges to the true parameter vector at rate $\sqrt{\log d/n}$. We begin with a few notations and assumptions.
\subsection{Notations}
We start by recalling some standard notations for vector and matrix norms from linear algebra.
\begin{enumerate}
	\item For a vector $\bm a := (a_1,\ldots,a_s) \in \mathbb{R}^s$ and $p \in (0,\infty)$,
	\begin{itemize}
		\item $\|\bm a\|_p := \left(\sum_{i=1}^s |a_i|^p\right)^\frac{1}{p}$
		
		\item $\|\bm a\|_\infty := \max_{1\le i\le s} |a_s|$
		
		\item $\|\bm a\|_0 := \sum_{i=1}^s \mathbbm{1}\{a_i \ne 0\}$
	\end{itemize}

Note that $\|\cdot\|_0$ is not a vector norm, since $\|\alpha \bm a\|_0 \ne |\alpha| \|\bm a\|_0$ as long as $\bm a$ is non-zero and $\alpha \ne \pm 1$.

	\item For a matrix $\bm M := ((M_{ij}))_{1\le i \le s,1\le j \le t} \in \mathbb{R}^{s\times t}$, 
	\begin{itemize}
		\item $\|\bm M\|_\infty := \max_{1\le i\le s} \sum_{j=1}^t |M_{ij}|$
		
		\item $\|\bm M\|_1 := \max_{1\le j\le t} \sum_{i=1}^s |M_{ij}|$
		
		\item $\|\bm M\|_2 := \sigma_{\max}(\bm M)$, where $\sigma_{\max}(\bm M)$ denotes the largest singular value of $\bm M$.
		
		\item $\|\bm M\|_F := \sqrt{\sum_{i=1}^s \sum_{j=1}^t M_{ij}^2}$
	\end{itemize}
	For a square matrix $\bm M$, note that $\|\bm M\|_2$ equals the absolute value of the eigenvalue of $\bm M$ with the largest magnitude.
\end{enumerate}

\subsection{Assumptions}
Before stating the main result of this chapter, let us state a few standing assumptions.

\begin{itemize}
	\item \textbf{Assumption 1.}~$\|\bh\|_\infty < \Theta$ and $\|\bm Z_i\|_\infty < M$ for all $1\le i\le N$, for some fixed constants $\Theta,M > 0$.  
	
	\item \textbf{Assumption 2.}~ $\bm A$ is a symmetric matrix with zeros on the diagonal.
	
	\item \textbf{Assumption 3.}~ $\sup_{N \ge 1} \|\bm A\|_\infty \le 1$.
	
	\item \textbf{Assumption 4.}~ $\liminf_{N \rightarrow \infty} \frac{1}{N} \|\bm A\|_F^2 > 0$.
	
	\item \textbf{Assumption 5.}~ $|\beta| < B := 1/4$.
	
		\item \textbf{Assumption 6.} $\liminf_{N\rightarrow\infty} \lambda_{\min}(N^{-1}\mathbf{Z}^\top \mathbf{Z}) > 0$, where $\mathbf{Z} := (\bm Z_1,\ldots, \bm Z_N)^\top$ and $\lambda_{\min}$ denotes the minimum eigenvalue.
		
		\item \textbf{Assumption 7.} $s:= \|(\beta, \bh^\top)\|_0$ is bounded with $N$ and $d$.
	\end{itemize}

It is important to understand that consistency of $(\hat{\beta},\hat{\bh}^\top)$ is not true in the generality. Even when $d=1$ and $\bm Z_1=\ldots =\bm Z_N$, it is shown in \cite{pg_sm} that consistent estimation of both the parameters $\beta$ and $\theta$ is impossible when $\bm A$ is the (scaled) adjacency matrix of a dense Erd\H{o}s-R\'enyi model, which includes the Curie-Weiss model as a special case. Hence, we need more restrictions on $\bm A$ to ensure consistent estimation.  Assumptions 3 and 4 are required in \cite{pg_sm} (see Theorem 1.15) and \cite{cd_ising_I} for $\sqrt{N}$-consistency of the parameters. If we take
$$\bm A := \frac{|V(G)|}{|E(G)|} \cdot \mathcal{A}(G)$$ for a graph $G$ with adjacency matrix $\mathcal{A}(G)$, then Assumption 3 says that the maximum degree of $G$ is of the same order as its average degree, and Assumption 4 says that the average degree of $G$ is bounded. Hence, together they imply that $G$ is a bounded degree graph, which is one of the standing assumptions in \cite{pg_sm} for deriving $\sqrt{N}$-consistency of $(\beta,\theta)$. Moreover, as noted in \cite{cd_ising_I}, the requirement of boundedness of $\|\bm A\|_\infty $ and $\beta$ (Assumptions 3 and 5) is crucial to ensure that the peer effects through $\beta \bm A$ coming from the quadratic dependence term in the probability mass function \eqref{eq:logisticmodel} does not overpower the effect of the signal $\bh$ coming from the linear terms $\bh^\top \bm Z_i$, thereby hindering joint recovery of the correlation term $\beta$ and the signal term $\bh$. A similar logic applies in support of Assumption 1, this time ensuring that the signal term does not dominate the correlation term. Assumption 6 is required crucially in establishing a strong-concavity type condition on the negative log-pseudolikelihood, which is pivotal in ensuring consistency of the estimator. This assumption holds with high probability, if the covariates $\bm Z_1,\ldots,\bm Z_N$ are i.i.d. realizations from a sub-Gaussian distribution on $\mathbb{R}^d$, the minimum eigenvalue of whose covariance matrix is bounded away from $0$ (see Theorem 2.1 in \cite{cd_ising_I}).

Below, we state the main result of this chapter.
\begin{thm}\label{finalmaintheorem}
    Suppose that Assumptions $1-7$ hold. Then there exist constants $\delta > 0$ and $\varepsilon > 0$, such that if if $d^2 \log d \le \varepsilon N$, then by taking $\lambda := \delta \sqrt{\log(d+1)/N}$ in the objective function in \eqref{eq:pseudolikelihood}, we have:
	$$\|(\hat{\beta}-\beta, \hat{\bh}^\top-\bh^\top)\|_2 = O\left(\sqrt{\frac{\log d}{N}}\right)$$
	with probability $1-o(1)$ as $N\rightarrow \infty$ and $d \rightarrow \infty$.
\end{thm}

Below, we provide a sketch of the proof of Theorem \ref{finalmaintheorem}. The actual proof can be found in Section \ref{prfinalmain}.

\subsection{Sketch of the Proof of Theorem \ref{finalmaintheorem}}
The proof of Theorem \ref{finalmaintheorem} proceeds through a number of steps. The first step is to show that if the gradient of $L_N$ at the true parameter $\bg := (\beta,\bh^\top)^\top$ is bounded entrywise by $\lambda$ (which we will refer to as the first order condition), and if $L_N$ satisfies the strong-concavity type condition (which we will refer to as the second order condition):
$$L_N(\hat{\bg}) \ge L_N(\bg) + \nabla L_N(\bg)^\top (\hat{\bg}-\bg) + \alpha \|\hat{\bg}-\bg\|_2^2~,$$  then  $\|\hat{\bg}-\bg\|_2 \lesssim \lambda/\alpha$. This is proved in Lemma \ref{mainlem}. The proof then boils down to deriving $\lambda$ and $\alpha$ which satisfy the first and second order conditions.

The first order condition is proved by applying the method of exchangeable pairs on each element of the gradient of $L_N(\bg)$. More specifically, one starts by showing that all elements of $\nabla L_N(\bg)$ have mean $0$, and then establishes concentration of these elements around their means (which are all $0$), by applying methods from \cite{scthesis}. The high-probability bound $\lambda$ on $\|\nabla L_N(\gamma)\|_\infty$ turns out to be of order $\sqrt{\log (d)/N}$, which gives the rate in Theorem \ref{finalmaintheorem}.

Verifying the second order condition for a constant value of $\alpha$ is more involved. This is essentially equivalent to showing that the lowest eigenvalue of the Hessian $\nabla^2 L_N$ is bounded away from $0$ in a neighborhood of the true parameter $\bg$, stretching from $\bg$ to $\hat{\bg}$. That is why, it is essential to bound the size of this neighborhood, and in Lemma \ref{l1bounded}, we show that the $L^{1}$ radius of this neighborhood, $\|\hat{\bg}-\bg\|_1$, is $O(1)$ as long as $d^2\log d \leq \varepsilon N$ for some sufficiently small $\varepsilon > 0$. 

The remaining part of the proof on bounding the lowest eigenvalue of $\nabla^2 L_N$ away from $0$, involves (in view of a Schur complement argument) showing that the quantity $N^{-1} \|\bm F \bm m\|_2^2$ is bounded away from $0$ with high probability. This is done in two steps. The first step is to show that the expectation of $N^{-1} \|\bm F \bm m\|_2^2$ is bounded away from $0$, and the second step is to prove a concentration of $N^{-1} \|\bm F \bm m\|_2^2$ around its expectation. Interestingly, the concentration step does not follow from the method of exchangeable pairs and the techniques in \cite{scthesis}, but requires a more sophisticated argument on concentration for polynomials in Ising models satisfying the Dobrushin condition (see \cite{radek}). The first step of lower bounding the mean of $N^{-1} \|\bm F \bm m\|_2^2$ is a bit more involved, and requires lower bounding the variance of linear projections of $\bm X$, which involves delicate arguments analogous to those in \cite{cd_ising_estimation}. 

\section{Proof of Theorem \ref{finalmaintheorem}}\label{prfinalmain}
For notational convenience, we will henceforth denote the vector $(\beta, \bh^\top)$ by $\bg^\top$, and $(\hat{\beta}, \hat{\bh}^\top)$ by $\hat{\bg}^\top$. 
Key to the proof of Theorem \ref{finalmaintheorem}, is the following two lemma:
\begin{lem}\label{actuallem}
There exists a constant $C > 0$, such that for every $\delta > 0$ sufficiently large, if we take $\lambda := \delta \sqrt{\log(d+1)/N}$ in the objective function in \eqref{eq:pseudolikelihood}, then we have 
	$$\|\hat{\bg} -\bg\|_2 \leq C\delta \cosh^2\left(B + M(\|\hat{\bg}-\bg\|_1)\right) \sqrt{\frac{\log d}{N}}$$
with probability $1-o(1)$ as $N\rightarrow \infty$ and $d \rightarrow \infty$, as long as $d = o(N)$.
\end{lem}

Theorem \ref{finalmaintheorem} follows immediately from Lemma \ref{actuallem}, if we can ensure that $\|\hat{\bg}-\bg\|_1 = O(1)$ with probability $1-o(1)$.     The following lemma guarantees this, as long as $d^2 \log d \le \varepsilon n$ for some sufficiently small constant $\varepsilon>0$.
\begin{lem}\label{l1bounded}	
For every $\delta > 0$ sufficiently large, if we take $\lambda := \delta \sqrt{\log(d+1)/N}$ in the objective function in \eqref{eq:pseudolikelihood}, then we have 
	$$\min\{\|\hat{\bg}-\bg\|_1~,~1\} = O\left(\delta d \sqrt{\frac{\log d}{N}}\right)~.$$
	with probability $1-o(1)$ as $N\rightarrow \infty$ and $d \rightarrow \infty$, as long as $d = o(N)$.
\end{lem}

It remains to prove Lemmas \ref{actuallem} and \ref{l1bounded}, in order to complete the proof of Theorem \ref{finalmaintheorem}.

\subsection{Proof of Lemma \ref{actuallem}}
Towards proving Lemma \ref{actuallem}, we start with a basic result, which gives the consistency rate of $\bm v := \hat{\bg}-\bg$ under some conditions on the first and the second order derivatives of the function $L_N$. Let us define:
$$S := \{1\le i\le d+1: \gamma_i \neq 0\}~.$$ For any vector $\bm a \in \mathbb{R}^{p}$ and any set $Q \subseteq \{1,\ldots,p\}$, we denote the vector $(a_i)_{i\in Q}$ by $\bm a_Q$.
\begin{lem}\label{mainlem}
	Suppose that $\|\nabla L_N(\bg)\|_{\infty}   \le \lambda/2$ and $L_N(\hat{\bg}) - L_N(\bg) - \nabla L_N(\bg)^\top \bm v \geq \alpha \|\bm v\|_2^2$. Then,
	$$\|\hat{\bg}-\bg\|_2 \leq \frac{3\lambda \sqrt{s}}{2\alpha}~.$$
\end{lem}
\begin{proof}
	By Lemma \ref{lambound}, we have:
	\begin{equation}\label{mainbound}
	\|v_S\|_1 - \|v_{S^c}\|_1 + \frac{\|v\|_1}{2} \geq \frac{\alpha \|v\|_2^2}{\lambda}~.
	\end{equation}
	Using the fact that $\|v\|_1 = \|v_S\|_1 + \|v_{S^c}\|_1$, we have from \ref{mainbound},
	\begin{equation}\label{maineq1}
	\frac{\alpha \|v\|_2^2}{\lambda} \leq \frac{3\|v_S\|_1}{2} - \frac{\|v_{S^c}\|_1}{2} \leq \frac{3\|v_S\|_1}{2} \leq \frac{3\sqrt{s}\|v_S\|_2}{2}\leq \frac{3\sqrt{s}\|v\|_2}{2}~. 
	\end{equation}
	Lemma \ref{mainlem} now follows from \eqref{maineq1}.
\end{proof}

We will refer to the conditions $\|\nabla L_N(\bg)\|_{\infty}   \le \lambda/2$ and $L_N(\hat{\bg}) - L_N(\bg) - \nabla L_N(\bg)^\top \bm v \geq \alpha \|\bm v\|_2^2$ in the statement of Lemma \ref{mainlem} as the first order condition and the second order condition, respectively. Lemma \ref{actuallem} will follow from Lemma \ref{mainlem} if we can verify these conditions for some suitable values of $\lambda$ and $\alpha$. We do this in the next subsection.

\subsubsection{Verifying the First and Second Order Conditions}
In this subsection, we verify the first and second order conditions assumed in the statement of Lemma \ref{mainlem}. The goal is to make suitable choices of $\lambda$ and $\alpha$, such that the hypotheses of Lemma \ref{mainlem} hold with high probability. As usual, we will use the notation
$$m_i(\bs) := \sum_{j=1}^N A_{ij} X_j~.$$  We start by verifying the first order condition.

\begin{lem}\label{devcond}
	\textbf{(First Order Condition)}\hspace{0.2cm}Let $C := (1-4\beta\|\ba\|_2)\big /4 \max \{(\beta+3)^2, M^2(1+\beta)^2\}$ and $\lambda := \delta \sqrt{\log(d+1)/N}$ for some $\delta > 0$. Then, we have:
	$$\bp\left(\norm{\nabla L_N(\bg)}_\infty \leq \frac{\lambda}{2}\right) \geq 1 - 2(d+1)^{1-(C\delta^2/4)}~.$$
\end{lem}
\begin{proof}
	To begin with, note that:
	$$\be\left(X_i - \tanh(\beta m_i(\bs) +\bh^\top \bm Z_i)\Bigg|(X_j)_{j\neq i}, \right) = 0~,$$ and hence,
	\begin{align}\label{centered1}
	\be\left(\frac{\partial L_N}{\partial \beta}\right) &=  -\frac{1}{N} \sum_{i=1}^N \be\left[m_i(\bs)\left(X_i - \tanh(\beta m_i(\bs) +\bh^\top \bm Z_i)\right)\right]\nonumber\\&= -\frac{1}{N} \sum_{i=1}^N \be\left[m_i(\bs)\be\left(X_i - \tanh(\beta m_i(\bs) +\bh^\top \bm Z_i)\Bigg|(X_j)_{j\neq i}\right)\right] = 0~,
	\end{align}
	and
	\begin{align}\label{centered2}
	\be\left(\frac{\partial L_N}{\partial \theta_j}\right) &=  -\frac{1}{N} \sum_{i=1}^N \be\left[Z_{i,j}\left(X_i - \tanh(\beta m_i(\bs) +\bh^\top \bm Z_i)\right)\right]\nonumber\\&= -\frac{1}{N} \sum_{i=1}^N \be\left[Z_{i,j}\be\left(X_i - \tanh(\beta m_i(\bs) +\bh^\top \bm Z_i)\Bigg|(X_\ell)_{\ell\neq i}\right)\right] = 0~,
	\end{align}
	By Lemma 4.4 in \cite{scthesis}, Dobrushin's interdependence matrix for model \eqref{eq:logisticmodel} is $4\beta \ba$. Hence, by \eqref{centered1}, \eqref{centered2}, Lemma \ref{partialbetalips}, Lemma \ref{partialthetalips}, and Theorem 4.3 and Lemma 4.4 in \cite{scthesis}, we have for every $t \geq 0$,
	\begin{equation}\label{conc1}
	\bp\left(\left|\frac{\partial L_N}{\partial \beta}\right| \geq t\right) \leq 2\exp\left(-\frac{N(1-4\beta\|\ba\|_2)t^2}{4(\beta+3)^2}\right)~.
	\end{equation}
	Similarly, for each $j\in [d]$, we have:
	\begin{equation}\label{conc2}
	\bp\left(\left|\frac{\partial L_N}{\partial \theta_j}\right| \geq t\right) \leq 2\exp\left(-\frac{N(1-4\beta\|\ba\|_2)t^2}{4M^2(1+\beta)^2}\right)~.
	\end{equation}
    It thus follows from \eqref{conc1}, \eqref{conc2} and a union bound, that
	\begin{equation}\label{devconc}
	\bp\left(\|\nabla L_N(\bg)\|_\infty \geq t\right) \leq 2(d+1)e^{-CNt^2}~,
	\end{equation}
	where $C$ is defined in the hypothesis of Lemma \ref{devcond}.
	Lemma \ref{devcond} now follows on taking $t = \lambda/2$, with $\lambda$ as in the hypothesis.
\end{proof}

In the following lemma, we verify the second order condition in the hypothesis of Lemma \ref{mainlem}.

\begin{lem}\label{rsccond}
	\textbf{(Second Order Condition)}\hspace{0.2cm}
There exists a constant $C > 0$, such that
	$$L_N(\hat{\bg}) - L_N(\bg) - \nabla L_N(\bg)^\top \bm v \geq \frac{C}{\cosh^2(B+M(\|\bm v\|_1+s\Theta) )}  \norm{\bm v}_2^2$$
	with probability $1-o(1)$.
\end{lem}

\begin{proof}
	By a second order Taylor series expansion, we know that there exists $\underline{\bg} = (\underline{\beta},\underline{\bh}^\top) \in B_{1}(\bg;\|\bm v\|_1)$\footnote{For a vector $\bm a$, positive integer $p$ and positive real $r$, the set $B_p(\bm a;r)$ denotes the open $L^p$ ball of radius $r$ around $u$, i.e. $B_p(\bm a;r) := \{\bm x: \|\bm x-\bm a\|_p < r\}$.}, such that
	\begin{equation}\label{second_taylor}
	L_N(\hat{\bg}) - L_N(\bg) - \nabla L_N(\bg) ^\top \bm v = \frac{1}{2} \bm v^\top \nabla^2 L_N(\underline{\bg}) \bm v = \frac{1}{2N}\sum_{i=1}^N \frac{\bm v^\top \bm U_i \bm U_i^\top \bm v}{\cosh^2(\underline{\beta} m_i(\bs) + \underline{\bh}^\top \bm Z_i)} ~,   
	\end{equation}
	where and $\bm U_i := (m_i(\bs),\bm Z_i^\top)^\top$. Now, note that: $$|\underline{\bh}^\top \bm Z_i| \leq \|\underline{\bh}\|_1 \|\bm Z_i\|_\infty \leq M\left(\|\underline{\bh}-\bh\|_1 + \|\bh\|_1\right) \leq M(\|\bm v\|_1 + s\Theta)$$ Since $\cosh$ is an even function, which is increasing on the positive axis, we have:
	\begin{equation}\label{nnd}
	\frac{1}{2N}\sum_{i=1}^N \frac{\bm v^\top \bm U_i \bm U_i^\top \bm v}{\cosh^2(\underline{\beta} m_i(\bs) + \underline{\bh}^\top \bm Z_i)} \ge  \frac{\bm v^\top \bm G \bm  v}{2\cosh^2(\underline{\beta} m_i(\bs) + \underline{\bh}^\top \bm Z_i)}~,
	\end{equation}
	where 
	
	\begin{equation*}
	\bm G :=  \frac{1}{N}\left(
	\begin{array}{cc}
	\bm m^\top \bm m & \bm m^\top \mathbf{Z}\\
	\mathbf{Z}^{\top} \bm m & \mathbf{Z}^\top \mathbf{Z} 
	\end{array}
	\right)~,\quad\bm m := (m_1(\bm X),\ldots,m_N(\bm X))^\top\quad \textrm{and}\quad \mathbf{Z} = (\bm Z_1,\ldots,\bm Z_N)^\top. 
	\end{equation*}
	 
	Combining \eqref{second_taylor} and \eqref{nnd}, we have:
	\begin{equation}\label{withG}
	L_N(\hat{\bg}) - L_N(\bg) - \nabla L_N(\bg) ^\top \bm v = \frac{1}{2} \bm v^\top \nabla^2 L_N(\underline{\bg}) \bm v \ge \frac{v^\top \bm G v}{2\cosh^2(B + M(\|\bm v\|_1 + s\Theta))} ~. 
	\end{equation}
	
	Lemma \eqref{rsccond} now follows from \eqref{withG} and Lemma \ref{mineiglem}. 
\end{proof}
The proof of Lemma \ref{actuallem} now follows from Lemmas \ref{mainlem}, \ref{devcond} and \ref{rsccond}. \qed

\subsection{Proof of Lemma \ref{l1bounded}}
Define: $\bg_t := t\hat{\bg} + (1-t) \bg$, and $g(t) := (\hat{\bg} - \bg)^\top \nabla L_N(\bg_t)~.$
First, note that:
\begin{equation}\label{step1}
|g(1)-g(0)| = \big|(\hat{\bg} - \bg)^\top (\nabla L_N(\hat{\bg}) - \nabla  L_N(\bg))\big| \le \|\hat{\bg} - \bg\|_2 \cdot \|\nabla L_N(\hat{\bg}) - \nabla  L_N(\bg)\|_2~.
\end{equation}
Next, we have:
\begin{eqnarray*}
	g'(t) &=& (\hat{\bg}-\bg)^\top \nabla^2 L_N(\bg_t) (\hat{\bg}-\bg)\\&=& \frac{1}{N}\sum_{i=1}^N \frac{(\hat{\bg}-\bg)^\top \bm U_i \bm U_i^\top (\hat{\bg}-\bg)}{\cosh^2(\beta_t  m_i(\bs) + \bh_t^\top \bm Z_i)}\quad(\textrm{where}~\bm U_i := (m_i(\bs), \bm Z_i^\top)^\top)\\&\ge& \frac{(\hat{\bg}-\bg)^\top \bm G (\hat{\bg}-\bg)}{\cosh^2(B + M(\|\bg_t - \bg\|_1 + s\Theta))}\\&\ge & \frac{\|\hat{\bg}-\bg\|_2^2 ~\lambda_{\mathrm{min}}(\bm G)}{\cosh^2(B + M(\|\bg_t - \bg\|_1 + s\Theta))}\\&\gtrsim& \frac{\|\hat{\bg}-\bg\|_2^2 }{\cosh^2(B + M(\|t(\hat{\bg} - \bg)\|_1 + s\Theta))}\quad(\textrm{by Lemma} \ref{mineiglem})
\end{eqnarray*}
Hence, we have:
\begin{align}\label{in3}
|g(1)-g(0)|&\ge g(1) - g(0)=\int_0^1 g'(t)\nonumber\\ &\ge \int_0^{\min\{1~,~ \|\hat{\bg}-\bg\|_1^{-1}\}} g'(t)\nonumber\\&\gtrsim \|\hat{\bg} - \bg\|_2^2 ~\min\{1~,~ \|\hat{\bg}-\bg\|_1^{-1}\}
\end{align}
Combining \eqref{step1} and \eqref{in3}, we have:
\begin{equation}\label{in4}
\frac{\|\hat{\bg}-\bg\|_2}{\|\hat{\bg}-\bg\|_1}\cdot  \min\{\|\hat{\bg}-\bg\|_1~,~1\} \le \|\nabla L_n(\hat{\bg}) - \nabla L_n(\bg)\|_2~.
\end{equation}
Using the fact that for every non-zero vector $\bm a \in \mathbb{R}^p$,
$$\frac{\|\bm a\|_2}{\|\bm a\|_1} \ge \frac{1}{\sqrt{p}}$$ we have from \eqref{in4},
\begin{equation}\label{in5}
\min\{\|\hat{\bg}-\bg\|_1~,~1\} \le \sqrt{d+1}\cdot  \|\nabla L_N(\hat{\bg}) - \nabla L_N(\bg)\|_2 \le \sqrt{d+1} \cdot (\|\nabla L_N(\hat{\bg})\|_2 + \|\nabla L_N(\bg)\|_2)
\end{equation}
By Lemma \ref{devcond}, we have with probability at least $1-2(d+1)^{1-(C\delta^2/4)}$ (with $C$ as in the statement of Lemma \ref{devcond}),
$$\|\nabla L_N(\bg)\|_\infty \le \frac{\delta}{2}\sqrt{\log(d+1)/N}\quad \implies \quad \|\nabla L_N(\bg)\|_2 \le \frac{\delta}{2}\sqrt{(d+1)\log(d+1)/N}~.$$
Also, by Lemma \ref{nondiff}, we have:
$$\|\nabla L_n(\hat{\bg})\|_2 \le \lambda \sqrt{d+1}~.$$
The proof of Lemma \ref{l1bounded} is now complete. \qed

This completes the proof of Theorem \ref{finalmaintheorem}.

\section{Positivity of the Hessian $\nabla^2 L_N$}\label{poshes}
In this section, we show that with high probability, the lowest eigenvalue of $\nabla^2 L_N$ is bounded away from $0$. In fact, this is equivalent to showing that the lowest eigenvalue of the following matrix
\begin{equation*}
\bm G :=  \frac{1}{N}\left(
\begin{array}{cc}
\bm m^\top \bm m & \bm m^\top \mathbf{Z}\\
\mathbf{Z}^{\top} \bm m & \mathbf{Z}^\top \mathbf{Z} 
\end{array}
\right) 
\end{equation*}   
 is bounded away from $0$ with high probability. 
 
 \begin{lem}\label{mineiglem}
 	There exists a constant $C>0$ (depending only on $s,\Theta,B$ and $M$), such that $$\p(\lambda_{\min}(\bm G) \ge C) \ge 1-e^{-\Omega(N)}$$ as $N$ and $d \rightarrow \infty$, with $d=o(N)$.
 \end{lem}
%\begin{proof}
	The first step towards proving Lemma \ref{mineiglem} is to observe that:
	$$\det(\bm G-\lambda \bm I) = \left(\frac{1}{N}\|\mf \bm m\|_2^2 - \lambda\right) \cdot \det\left(\frac{1}{n} \bm Z^\top \bm Z - \lambda I\right)~,$$
	where $\bm F := \bm I - \mathbf{Z} (\mathbf{Z}^\top \mathbf{Z})^{-1} \mathbf{Z}^\top$.
	Hence, we have:
	\begin{equation}\label{crucial}
	\lambda_{\min}(\bm G) = \min \left\{\lambda_{\min}\left(\frac{1}{N} \bm Z^\top \bm Z\right)~,~\frac{1}{N}\|\mf \bm m\|_2^2\right\}.
	\end{equation}
	In view of Assumption 6, we only need to show the existence of a constant $C>0$ (depending only on $s,\Theta,B$ and $M$), such that
	\begin{equation}\label{fm}
	\p\left(\frac{1}{N}\|\mf \bm m\|_2^2 \ge C\right) = 1-e^{-\Omega(N)}
	\end{equation}
	in order to complete the proof of Lemma \ref{mineiglem}. We do this in two steps. First, we show that the mean of $N^{-1} \|\bm F \bm m\|_2^2$ is $\Omega(1)$, and then we show that $N^{-1} \|\bm F \bm m\|_2^2$ concentrates around its mean. These steps are implemented in Lemmas \ref{meanlb} and \ref{concmeanlb} respectively.  
	%\end{proof}
\begin{lem}\label{meanlb}
	There exists a universal constant $C > 0$, such that for all $N \ge 1$, we have:
	$$ \mathbb{E}\left(\frac{1}{N}\|\bm F \bm m\|_2^2\right)\ge Ce^{-8 \Theta M s - 4B}~.$$
	\end{lem}

\begin{proof}
For a matrix $\bm M$, we will denote the $i^{\mathrm{th}}$ row of $\bm M$ by $\bm M_i$. Now, we have:
\begin{equation}\label{fs1}
\e \|\bm F \bm m\|_2^2 = \sum_{i=1}^N \e \left([(\bm F \bm A)_i \bm X]^2\right) \ge \sum_{i=1}^N \mathrm{Var}\left((\bm F \bm A)_i \bm X\right) ~.
\end{equation}
The hypothesis of Lemma \ref{lemma10} is satisfied by Assumptions 3 and 5, and hence, by Lemma \ref{lemma10}, we have:
$$\mathrm{Var}\left((\bm F \bm A)_i \bm X\right) \gtrsim \Upsilon^2 \|(\bm F\bm A)_i\|_2^2$$ where $\Upsilon := \min_{1\le i\le N} \mathrm{Var}(X_i|\bm X_{-i})$. It follows from \eqref{fs1} that

\begin{equation}\label{fs2}
\e \|\bm F \bm m\|_2^2  \gtrsim \Upsilon^2 \|\bm F \bm A\|_F^2 \ge  \Upsilon^2 \left(\|\bm A\|_F^2 - d\right)~.
\end{equation}
The last inequality in \eqref{fs2} follows from Lemma \ref{freq}. It follows from \eqref{fs2}, in view of Assumption 4 and the hypothesis $d=o(N)$ of Lemma \ref{mineiglem}, that  
\begin{equation}\label{25conc}
\e \left(\|\bm F \bm m\|_2^2 \right) \gtrsim \Upsilon^2 N~.
\end{equation}
The task now, is to give a lower bound for $\Upsilon^2$. Towards this, note that for every $1\le i\le N$, $\mathrm{Var}(X_i|\bm X_{-i}) = 4p(1-p)$, where $p= \p(X_i =1|\bm X_{-i})$. Now, we have:
$$p = \frac{\exp\left(\bh^\top \bm Z_i + \frac{\beta}{2}\sum_{v \ne j} A_{jv}X_v\right)}{2 \cosh\left(\bh^\top \bm Z_i + \frac{\beta}{2}\sum_{v \ne j} A_{jv}X_v\right)}~.$$
By Lemma \ref{math1} we have:
\begin{equation}\label{minp}
\min\{p, 1-p\} \ge \frac{1}{2} \exp\left(-2\left|\bh^\top \bm Z_i + \frac{\beta}{2}\sum_{v \ne j} A_{jv}X_v\right|\right)
\end{equation}
Now, we have
$$\left|\bh^\top \bm Z_i + \frac{\beta}{2}\sum_{v \ne j} A_{jv}X_v\right| \le |\bh^\top \bm Z_i| + \frac{\beta}{2} \|\bm A\|_\infty \le \Theta M s + \frac{B}{2} $$
and hence, it follows from \eqref{minp} that
$$\min\{p,1-p\} \ge \frac{1}{2} e^{-2\Theta M s -B}~.$$
Hence, for all $1\le i\le N$, we have:
$$\mathrm{Var}(X_i|\bm X_{-i}) \ge e^{-4\Theta M s - 2B}~.$$ Hence, we have:
\begin{equation}\label{25conc2}
\Upsilon^2 \ge e^{-8\Theta M s - 4B}~.
\end{equation}

Lemma \ref{meanlb} now follows from \eqref{25conc} and \eqref{25conc2}.

\end{proof}

\begin{lem}\label{concmeanlb}
	For any $t > 0$, we have:
	$$\bp\left(\|\bm F \bm m\|_2^2 < \mathbb{E}\|\bm F \bm m\|_2^2 - t\right) \le 2 \exp\left(-C\cdot \min\left\{  \frac{t^2}{2N}~,~t\right\}\right)~,$$ where $C$ is a constant, depending only on $\Theta,M,s$ and $B$.
	\end{lem}
\begin{proof}
	Let us denote $\bm F \bm A$ by $\bm W$, and let $\bm H$ be the matrix obtained from $\bm W^\top \bm W$ by zeroing out all its diagonal elements. It follows from Example 2.5 in \cite{radek} that 
	\begin{align}
	\bp\left(\|\bm F \bm m\|_2^2 < \mathbb{E}\left[\|\bm F \bm m\|_2^2\right] - t\right) &= \bp\left(\|\bm X^\top \bm H \bm X\|_2^2 < \mathbb{E}\left[\|\bm X^\top \bm H \bm X\|_2^2\right] - t\right)\nonumber\\&\le 2\exp \left(-c \cdot \min\left\{\frac{t^2}{\|\bm H\|_F^2 + \|\e (\bm H \bm X)\|_2^2}~,~\frac{t}{\|\bm H\|_2} \right\}\right)\\\label{radek1}
	\end{align}
	where $c$ is a constant depending only on $\Theta,M,s$ and $B$ (note that the parameters $\alpha$ and $\rho$ in \cite{radek} can be taken to be $\Theta M s$ and $3/4$ respectively).  Now, note that:
	
	\begin{align}
	\|\bm H\|_F^2 + \|\e (\bm H \bm X)\|_2^2 &\le \|\bm H\|_F^2 + \left(\|\e (\bm W^\top \bm W \bm X)\|_2 + \|\e [(\bm H - \bm W^\top \bm W)\bm X]\|_2\right)^2\nonumber\\ &\le 2\|\bm H\|_F^2 + 2 \e \|[(\bm H - \bm W^\top \bm W)\bm X]\|_2^2 + 2\|\e (\bm W^\top \bm W \bm X)\|_2^2\nonumber\\&= 2\sum_{i\ne j} (\bm W^\top \bm W)_{ij}^2 + 2 \sum_i (\bm W^\top \bm W)_{ii}^2 + 2\|\e (\bm W^\top \bm W \bm X)\|_2^2\nonumber\\\label{eqar1}&= 2\|\bm W^\top \bm W\|_F^2 + 2\|\e (\bm W^\top \bm W \bm X)\|_2^2~.
	\end{align}
	
	We also have $\|\bm H\|_2 \le \|\bm W^\top \bm W\|_2$, since for any vector $\bm a \in \mathbb{R}^N$, we have:
	$$\bm a^\top \bm H \bm a = \bm a^\top \bm W^\top \bm W \bm a - \sum_{i=1}^N (\bm W^\top \bm W)_{ii} a_i^2 \le a^\top \bm W^\top \bm W \bm a~.$$
	
	Hence, it follows from \eqref{radek1} and \eqref{eqar1} that
	
	\begin{align}\label{radek2}
	&\bp\left(\|\bm F \bm m\|_2^2 < \mathbb{E}\left[\|\bm F \bm m\|_2^2\right] - t\right)\nonumber\\ &\le 2\exp \left(-c' \cdot \min\left\{\frac{t^2}{\|\bm W^\top \bm W\|_F^2 + \|\e (\bm W^\top \bm W \bm X)\|_2^2}~,~\frac{t}{\|\bm W^\top \bm W\|_2} \right\}\right)
	\end{align}
	where $c' := c/2$. Next, note that for any two matrices $\bm U$ and $\bm V$ such that $\bm U \bm V$ exists, we have $\|\bm U \bm V\|_F \le \|\bm U\|_2 ~\|\bm V\|_F$. This, together with the fact $\|\bm R\|_2 = \|\bm R^\top\|_2$ for every square matrix $\bm R$, implies that
	\begin{equation}\label{rst1}
		\|\bm W^\top \bm W\|_F^2 \le \|\bm W\|_2^2~\|\bm W\|_F^2~.
		\end{equation}
		Also, by the submultiplicativity of the matrix $L^2$ norm, we have:
		\begin{equation}\label{rst2}
		\|\bm W^\top \bm W\|_2 \le \|\bm W^\top\|_2 \|\bm W\|_2 = \|\bm W\|_2^2
		\end{equation}
		Finally, 
		\begin{equation}\label{rst3}
		\|\e (\bm W^\top \bm W \bm X)\|_2^2 = \|\bm W^\top \e(\bm W\bm X)\|_2^2 \le \|\bm W\|_2^2\cdot \|\e(\bm W\bm X)\|_2^2~.
		\end{equation}
		It follows from \eqref{radek2}, \eqref{rst1}, \eqref{rst2} and \eqref{rst3}, that:
		\begin{equation}\label{prelem}
		\bp\left(\|\bm F \bm m\|_2^2 < \mathbb{E}\left[\|\bm F \bm m\|_2^2\right] - t\right) \le 2\exp\left(-\frac{c'}{\|\bm F \bm A\|_2^2}\cdot \min\left\{  \frac{t^2}{\|\bm F \bm A\|_F^2 + \|\mathbb{E}\left[\bm F \bm m\right] \|_2^2}~,~t\right\}\right)~.
		\end{equation}
		Lemma \ref{concmeanlb} now follows on observing that $\|\bm F \bm A\|_2^2 \le \|\bm F\|_2^2 \|\bm A\|_2^2 \le 1$, $\|\bm F \bm A\|_F^2 \le N\|\bm F \bm A\|_2^2 \le N$, and $\|\mathbb{E}\left[\bm F \bm m\right] \|_2^2 \le \e \|\bm F \bm m\|_2^2 \le \|\bm F \bm A\|_2^2 ~\e\|\bm X\|_2^2 \le N$.
\end{proof}
 Equation \eqref{fm} now follows from Lemmas \ref{meanlb} and \ref{concmeanlb}, on taking $$t= (CN/2) \exp(-8\Theta M s - 4B)$$ in Lemma \ref{concmeanlb}, where $C$ is as in the statement of Lemma \ref{meanlb}. The proof of Lemma \ref{mineiglem} is now complete.

%\appendix

% Changing formatting for appendices
%\newenvironment{appendixf}{}{}
%\titleformat{\chapter}[hang]{\large\center}{APPENDIX}{0 pt}{} % Old {APPENDIX \thechapter}{0 pt}{ : }
%\titlespacing*{\chapter}{0pt}{-33 pt}{6 pt} % The key value here is the -33 pts, I got to it by old fashioned measuring with a ruler....
\begin{appendices}
	
%	\addtocontents{toc}{\protect\setcounter{tocdepth}{-1}} % This is to fix how appendices are shown in ToC
%	\clearpage
%	\chapter{}
%	\addtocontents{toc}{\protect\setcounter{tocdepth}{1}} % This is to bring things back to normal
%	Here, we prove various missing details in the proofs above, and also prove some technical lemmas required in those proofs.

	\chapter{Technical Lemmas from Chapter \ref{curiech}}\label{anatools}

	\section{Special Functions and their Properties}\label{mathfunc} 
	
	In this section, we state a few important properties of some special mathematical functions which arise in our analysis. 
	
	\begin{definition}
		The gamma function $\Gamma: (0,\infty) \mapsto \mathbb{R}$ is defined as:
		$$\Gamma(x) := \int_0^\infty u^{x-1} e^{-u}~du.$$
	\end{definition}
	
	\begin{definition}
		The digamma function $\Gamma: (0,\infty) \mapsto \mathbb{R}$ is defined as:
		$$\psi(x) := \frac{\mathrm d}{\mathrm d x} \log \Gamma(x) = \frac{\Gamma'(x)}{\Gamma(x)}.$$
	\end{definition}
	
	The following standard expansion of the digamma function will be very helpful in our analysis: As  $x \rightarrow \infty$, 
	\begin{equation}\label{digamexp}
	\psi(1+x) = \log x + \frac{1}{2x} + O(x^{-2}). 
	\end{equation}

	\begin{definition}
		For real numbers $x\geq y>0$, the binomial coefficient $x$ choose $y$ is defined as $$\binom{x}{y} := \frac{\Gamma(x+1)}{\Gamma(y+1)\Gamma(x-y+1)}.$$
	\end{definition}
	
	\begin{lem}\label{bindiff}
		Fix $u>0$. Then, for every $x \in (0,u)$, we have
		$$\frac{\mathrm d}{\mathrm d x}\binom{u}{x} = \binom{u}{x}\left[\psi(u-x+1) - \psi(x+1)\right].$$
	\end{lem}
	\begin{proof}
		Let $\iota(x) = \binom{u}{x}$. Then, $\log \iota(x) = \log \Gamma(u+1) - \log \Gamma(x+1) - \log \Gamma(u-x+1)$ and hence,
		\begin{equation}\label{diffelt}
		\frac{\iota'(x)}{\iota(x)} = \frac{\mathrm d}{\mathrm d x}\log \iota(x)= -\psi(x+1) + \psi(u-x+1).
		\end{equation}
		Lemma \ref{bindiff} now follows from \eqref{diffelt}. 
	\end{proof}
	
	\section{Mathematical Approximations}\label{approx} In this section, we give three different types of standard mathematical approximations, which play crucial roles in our analysis.
	
	\begin{lem}[Riemann Approximation]\label{Riemann}
		Let $f:[a,b]\rightarrow \mathbb{R}$ be a differentiable function, and let $a=x_0<x_1<\ldots <x_n=b$. Let $x_s^* \in [x_{s-1},x_s]$ for each $1\leq k\leq n$. Then, we have:
		$$\left|\int_a^b f - \sum_{k=1}^n (x_s - x_{s-1})f(x_s^*)\right| \leq \frac{1}{2}(b-a)\max_{1\leq k\leq n}(x_s - x_{s-1})\sup_{x\in [a,b]}|f'(x)|.$$
	\end{lem}
	\begin{proof}
		Lemma \ref{Riemann} follows from the following string of inequalities:
		\begin{align}
		\left|\int_a^b f - \sum_{s=1}^n (x_s - x_{s-1})f(x_s^*) \right| &= \left|\sum_{s=1}^n \int_{x_{s-1}}^{x_s} (f(x) - f(x_s^*))\mathrm d x \right|\nonumber\\ 
		\label{sq1} &\leq \sum_{s=1}^n \int_{x_{s-1}}^{x_s} \left|f(x) - f(x_s^*)\right|\mathrm d x \\
		\label{sq2} &\leq  \sup_{x\in [a,b]}|f'(x)| \sum_{s=1}^n \int_{x_{s-1}}^{x_s} \left|x-x_s^*\right|\mathrm d x \\
		&= \frac{1}{2}\sup_{x\in [a,b]}|f'(x)| \sum_{s=1}^n \left[(x_s^*-x_{s-1})^2 + (x_s-x_s^*)^2\right]\nonumber\\
		&\leq \frac{1}{2}\sup_{x\in [a,b]}|f'(x)| \sum_{s=1}^n (x_s-x_{s-1})^2 \nonumber\\ 
		&\leq\frac{1}{2}(b-a)\max_{1\leq s \leq n}(x_s - x_{s-1})\sup_{x\in [a,b]}|f'(x)|. \nonumber
		\end{align}
		Note that, in going from \eqref{sq1} to \eqref{sq2}, we used the mean value theorem.
	\end{proof}

	The following lemma gives a Laplace-type approximation of an integral over a shrinking interval. For the classical Laplace approximation, which approximates integrals over fixed intervals, refer to \cite{bruijn,wong}. Even though the proof of Lemma \ref{Laplace} below is exactly similar to that of the classical Laplace approximation, we provide the proof here for the sake of completeness. To this end, for positive sequences $\{a_n\}_{n\geq 1}$ and $\{b_n\}_{n\geq 1}$, $a_n = O_\square(b_n)$ denotes  $a_n \leq C_1(\square) b_n$ and $a_n = \Omega_\square(b_n)$ denotes $a_n \geq C_2(\square) b_n$, for all $n$ large enough and positive constants $C_1(\square), C_2(\square)$, which may depend on the subscripted parameters.  
	
	\begin{lem}[Laplace-Type Approximation-I]\label{Laplace}
		Let $a<b$ be fixed real numbers, $g: [a,b]\mapsto \mathbb{R}$ be a differentiable function on $(a,b)$, and $h_n : [a,b]\mapsto \mathbb{R}$ be a sequence of thrice differentiable functions on $(a,b)$. Suppose that $\{x_n\}$ is a sequence in $(a,b)$ that is bounded away from both $a$ and $b$, satisfying $h_n'(x_n) = 0$ and $h_n''(x_n) < 0$ for all $n$. Suppose further, that for every $a < u<v<b$, $\sup_{x\in [u,v]} |g'(x)| = O_{u,v}(1)$, $\sup_{n\geq 1}\sup_{x\in [u,v]} |h_n^{(3)}(x)|= O_{u,v}(1)$ and $\inf_{x \in [u,v]} |g(x)| = \Omega_{u,v}(1)$. Also, suppose that $\inf_{n\geq 1} |h_n''(x_n)| > 0$. Then, for all $\alpha \in \left(0,\frac{1}{6}\right)$, we have as $n \rightarrow \infty$,
		\begin{equation*}
		\int_{x_n-n^{-\frac{1}{2}+\alpha}}^{x_n+n^{-\frac{1}{2}+\alpha}} g(x) e^{n h_n(x)}\mathrm d x=\sqrt{\dfrac{2\pi}{n\left| h_n''(x_n)\right|}}g(x_n) e^{n h_n(x_n)}\left(1+ O\left(n^{-\frac{1}{2}+3\alpha}\right)\right).
		\end{equation*} 
	\end{lem}
	\begin{proof}
		If we make the change of variables $y = \sqrt{n}(x-x_n)$, we have
		\begin{equation}\label{lapl1}
		\int_{x_n-n^{-\frac{1}{2}+\alpha}}^{x_n+n^{-\frac{1}{2}+\alpha}} g(x) e^{n h_n(x)}\mathrm d x = n^{-\frac{1}{2}} \int_{-n^{\alpha}}^{n^{\alpha}} g(yn^{-\frac{1}{2}}+x_n) e^{n h_n\left(yn^{-\frac{1}{2}}+x_n\right) }~\mathrm d y .
		\end{equation}
		By a Taylor expansion, we have for any sequence $y \in [-n^\alpha, n^\alpha]$,
		\begin{equation}\label{lapl2}
		e^{n h_n\left(yn^{-\frac{1}{2}}+x_n\right) } = \left(1+O\left(n^{3\alpha - \frac{1}{2}}\right)\right)e^{nh_n(x_n) + \frac{y^2}{2}h_n''(x_n)} 
		\end{equation}
		and
		\begin{equation}\label{lapl222}
	g(yn^{-\frac{1}{2}}+x_n) = \left(1+O\left(n^{\alpha - \frac{1}{2}}\right)\right)g(x_n)
		\end{equation}
		Using \eqref{lapl2} and \eqref{lapl222} the right side of \eqref{lapl1} becomes
		\begin{align*}
		&n^{-\frac{1}{2}}\left(1+O\left(n^{3\alpha - \frac{1}{2}}\right)\right) g(x_n)e^{n h_n(x_n)} \int_{-n^\alpha}^{n^\alpha} e^{\frac{y^2}{2}h_n''(x_n)}~\mathrm d y \nonumber\\
		& =  \left(1+O\left(n^{3\alpha - \frac{1}{2}}\right)\right) \sqrt{\dfrac{2\pi}{n\left| h_n''(x_n)\right|}}g(x_n) e^{n h_n(x_n)} \mathbb{P}\left(\left|N\left(0,\frac{1}{|h_n''(x_n)|}\right)\right| \leq n^\alpha\right)\\& =  \left(1+O\left(n^{3\alpha - \frac{1}{2}}\right)\right) \sqrt{\dfrac{2\pi}{n\left| h_n''(x_n)\right|}}g(x_n) e^{n h_n(x_n)}\left(1-O\left(e^{-n^\alpha}\right)\right)\\& =  \left(1+O\left(n^{3\alpha - \frac{1}{2}}\right)\right) \sqrt{\dfrac{2\pi}{n\left| h_n''(x_n)\right|}}g(x_n) e^{n h_n(x_n)}.
		\end{align*}
		The proof of Lemma \ref{Laplace} is now complete.
	\end{proof}
	
	\begin{lem}[Laplace-Type Approximation-II]\label{Laplaceirreg}
		Let $a<b$ be fixed real numbers, $g: [a,b]\mapsto \mathbb{R}$ be a differentiable function on $(a,b)$, and $h_n : [a,b]\mapsto \mathbb{R}$ be a sequence of 5-times differentiable functions on $(a,b)$. Suppose that $\{x_n\}$ is a sequence in $(a,b)$ that is bounded away from both $a$ and $b$, satisfying $h_n'(x_n) = 0$ for all $n\geq 1$. Also, assume that $n^\frac{1}{2}h''_n(x_n)=C_1 + O(n^{-\frac{1}{4}})$, $n^{\frac{1}{4}} h^{(3)}_n(x_n)=C_2 + O(n^{-\frac{1}{4}})$, and $h_n^{(4)}(x_n) = C_3 + O(n^{-\frac{1}{4}})$, where $C_1, C_2$ and $C_3$ are real constants. Suppose further, that for every $a < u<v<b$, $\sup_{x\in [u,v]} |g'(x)| = O_{u,v}(1)$, $\sup_{n\geq 1}\sup_{x\in [u,v]} |h_n^{(5)}(x)|= O_{u,v}(1)$ and $\inf_{x \in [u,v]} |g(x)| = \Omega_{u,v}(1)$. Then, for all $\alpha \in \left(0,\frac{1}{20}\right)$, as $n \rightarrow \infty$,
		\begin{equation*}
		\int_{x_n-n^{-\frac{1}{4}+\alpha}}^{x_n+n^{-\frac{1}{4}+\alpha}} g(x) e^{n h_n(x)}\mathrm d x=n^{-\frac{1}{4}} g(x_n)e^{n h_n(x_n)} \int_{-n^\alpha}^{n^\alpha} e^{\frac{y^2}{2}C_1+\frac{y^3}{6}C_2 + \frac{y^4}{24}C_3}~\mathrm d y \left(1+O\left(n^{5\alpha - \frac{1}{4}}\right)\right).
		\end{equation*}
	\end{lem}
	\begin{proof}
		To begin with, by a change of variables $y = n^{\frac{1}{4}}(x-x_n)$, we have
		\begin{equation}\label{lapl3}
		\int_{x_n-n^{-\frac{1}{4}+\alpha}}^{x_n+n^{-\frac{1}{4}+\alpha}} g(x) e^{n h_n(x)}\mathrm d x = n^{-\frac{1}{4}} \int_{-n^{\alpha}}^{n^{\alpha}} g(yn^{-\frac{1}{4}}+x_n) e^{n h_n\left(yn^{-\frac{1}{4}}+x_n\right) }~\mathrm d y ~
		\end{equation}
		Now, by a Taylor expansion of $n h_n\left(yn^{-\frac{1}{4}}+x_n\right)$ around $x_n$, we have for any sequence $y \in [n^{-\alpha},n^\alpha]$,
		\begin{align}\label{lptaylor}
		\nonumber n h_n\left(yn^{-\frac{1}{4}}+x_n\right) & =  nh_n(x_n)+\frac{n^\frac{1}{2}y^2}{2}h_n''(x_n)+\frac{n^\frac{1}{4}y^3}{6}h_n^{(3)}(x_n)+\frac{y^4}{24}h^{(4)}_n(x_n)+O\left(n^{-\frac{1}{4}}y^5\right)\\
		& = nh_n(x_n)+\frac{y^2}{2}C_1+\frac{y^3}{6}C_2+\frac{y^4}{24}C_3+O\left(n^{5\alpha - \frac{1}{4}}\right).
		\end{align}
		It follows from \eqref{lptaylor}, that
		\begin{align}\label{lapl4}
		e^{n h_n\left(yn^{-\frac{1}{4}}+x_n\right) } & =  \left(1+O\left(n^{5\alpha - \frac{1}{4}}\right)\right)e^{nh_n(x_n) + \frac{y^2}{2}C_1+\frac{y^3}{6}C_2+\frac{y^4}{24}C_3}.
		\end{align}
		Similarly, for any sequence $y \in [-n^\alpha,n^\alpha]$, we have
		\begin{equation}\label{g2}
		g(yn^{-\frac{1}{4}}+x_n) = \left(1+O\left(n^{\alpha - \frac{1}{4}}\right)\right)g(x_n).
		\end{equation}
		Using \eqref{lapl4} and \eqref{g2}, the right side of \eqref{lapl3} becomes
		\begin{align*}
		n^{-\frac{1}{4}} g(x_n)e^{n h_n(x_n)} \int_{-n^\alpha}^{n^\alpha} e^{\frac{y^2}{2}C_1+\frac{y^3}{6}C_2 + \frac{y^4}{24}C_3}~\mathrm d y \left(1+O\left(n^{5\alpha - \frac{1}{4}}\right)\right).
		\end{align*}
		The proof of Lemma \ref{Laplaceirreg} is now complete.
	\end{proof}
	
	\begin{lem}[Stirling's Approximation of the Binomial Coefficient]\label{stir}
		Suppose that $x=x_N$ is a sequence in $(-1,1)$ that is bounded away from both $1$ and $-1$. Then, as $N\rightarrow \infty$, 
		$$\binom{N}{N(1+x)/2} = 2^N\sqrt{\dfrac{2}{\pi N(1-x^2)}} \exp\left(-NI(x)\right)\left(1+O(N^{-1})\right).$$
	\end{lem}
	\begin{proof}
		First, note that by the usual Stirling approximation for the gamma function, we have the following as all of $u$, $v$ and $u-v \rightarrow \infty$,
		\begin{align*}
		\binom{u}{v}&= \frac{\sqrt{2\pi u}\left(\frac{u}{e}\right)^u\left(1+O\left(\frac{1}{u}\right)\right)}{\sqrt{2\pi v}\left(\frac{v}{e}\right)^v\left(1+O\left(\frac{1}{v}\right)\right)\sqrt{2\pi (u-v)}\left(\frac{u-v}{e}\right)^{(u-v)}\left(1+O\left(\frac{1}{ u -v }\right)\right)}\\
		&=	\sqrt{\dfrac{u}{2\pi v (u-v)}} \cdot \dfrac{u^u}{v^v(u-v)^{u-v}}\left(1+O\left(\frac{1}{u}\right)+O\left(\frac{1}{v}\right)+O\left(\frac{1}{u-v}\right)\right).
		\end{align*}
		Substituting $u = N$ and $v = N(1+x)/2$ (the hypothesis of the lemma indeed implies that $u,v$ and $u-v \rightarrow \infty$), we have
		\begin{align*}
		&\binom{N}{N(1+x)/2}\\&= \sqrt{\dfrac{N}{2\pi \frac{N(1+x)}{2}\cdot\frac{N(1-x)}{2}}}\cdot \dfrac{N^N}{\left(\dfrac{N(1+x)}{2}\right)^{N(1+x)/2} \left(\dfrac{N(1-x)}{2}\right)^{N(1-x)/2}}\left(1+O(N^{-1})\right)\\
		&=2^N\sqrt{\dfrac{2}{\pi N(1-x^2)}}\exp\left(-\frac{N(1+x)}{2} \log(1+x)-\frac{N(1-x)}{2}\log(1-x) \right)\left(1+O(N^{-1})\right)\\
		&=2^N\sqrt{\dfrac{2}{\pi N(1-x^2)}} \exp\left(-NI(x)\right)\left(1+O(N^{-1})\right).
		\end{align*}
		This completes the proof of Lemma \ref{stir}.
	\end{proof}

	\section{Properties of the Function $H$ and other Technical Lemmas}\label{sec:technical_lemmas} 
	
	This subsection is devoted to proving several technical lemmas that are used throughout the proofs of our main results. In Section \ref{sec:propH}, we will prove several important properties of the function $H$. Section \ref{sec:techmle} is devoted to proving various technical results related to the ML estimates of $\beta$ and $h$. Finally, we collect the proofs of some other technical lemmas in Section \ref{sec:other_tech}.
	
	\subsection{Properties of the Function $H$}\label{sec:propH} 
	We start by showing that a $p$-strongly critical point arises if and only if $p \geq 4$ is even, and in that case, the only such point is $(\tilde{\beta}_p,0)$ (recall \eqref{eq:betatilde}).
	
	\begin{lem}[Basic properties of the function $H$]\label{derh11}The function $H_{\beta,h,p}$ has the following properties. 
		\begin{enumerate}
			\item[$(1)$] $\sup_{x \in [-1,1]} H_{\beta,h,p}(x) \geq 0$ and equality holds if and only if $(\beta,h) \in  [0,\tilde{\beta}_p]\times \{0\}$.
			\item[$(2)$] Every local maximizer of $H_{\beta,h,p}$ lies in $(-1,1)$. 
			\item[$(3)$] $H_{\beta,h,p}$ can have at most two local maximizers for $p = 3$ and at most three local maximizers for $p \geq 4$. Further, it has three global maximizers if and only if $p\geq 4$ is even, $h=0$ and $\beta = \tilde{\beta}_p$.
		\end{enumerate}
	\end{lem}
	
	\begin{proof}[Proof of $(1)$.]
		First note that $\sup_{x \in [-1,1]} H_{\beta,h,p}(x) \geq H_{\beta,h,p}(0) = 0$. Now, it follows from first principles, that $\lim_{\varepsilon \rightarrow 0} H_{\beta,h,p}(\varepsilon)/\varepsilon = H_{\beta,h,p}'(0) = h$. If $h > 0$, then there exists $0<\varepsilon < 1$ such that $H_{\beta,h,p}(\varepsilon)/\varepsilon > h/2$, and if $h < 0$, then there exists $-1<\varepsilon < 0$ such that $H_{\beta,h,p}(\varepsilon)/\varepsilon < h/2$. In either case, $\sup_{x \in [-1,1]} H_{\beta,h,p}(x) \geq H_{\beta,h,p}(\varepsilon) > \varepsilon h /2 > 0$. Therefore, equality in (1) implies that $h = 0$, and hence, by  the definition in \eqref{eq:betatilde}, we must have $\beta \leq \tilde{\beta}_p$. This proves the ``only if" direction. For the ``if" direction, suppose that $(\beta,h) \in  [0,\tilde{\beta}_p]\times \{0\}$. Consider the case $\beta < \tilde{\beta}_p$ first, so that by  the definition in \eqref{eq:betatilde}, there exists $\beta'  > \beta$ such that $\sup_{x \in [-1,1]} H_{\beta' ,0,p}(x) = 0$. Equality in (1) now follows from: $$0\leq \sup_{x \in [-1,1]} H_{\beta,0,p}(x) = \sup_{x \in [-1,1]} H_{\beta,0,p}(|x|) \leq \sup_{x \in [-1,1]} H_{\beta' ,0,p}(|x|) = \sup_{x \in [-1,1]} H_{\beta' ,0,p}(x) = 0.$$ Finally, let $\beta = \tilde{\beta}_p$, and suppose towards a contradiction, that $H_{\beta,0,p}(x) > 0$ for some $x \in [-1,1]$. Then, $H_{\beta,0,p}(|x|) \geq H_{\beta,0,p}(x) > 0$, and hence, there exists $\beta'  < \beta$ such that $$H_{\beta' ,0,p}(|x|) = H_{\beta,0,p}(|x|) + (\beta' -\beta)|x|^p > 0.$$ This contradicts our previous finding that $\sup_{x \in [-1,1]} H_{\underline{\beta} ,0,p}(x) = 0$ for all $\beta < \tilde{\beta}_p$. The proof of (1) is now complete.
		\medskip
		
		\noindent\emph{Proof of $(2)$.}~Note that $\lim_{x \rightarrow -1^+} H_{\beta,h,p}'(x) = +\infty$ and $\lim_{x \rightarrow 1^-} H_{\beta,h,p}'(x) = -\infty$. Hence, there exists $\varepsilon > 0$, such that $H_{\beta,h,p}$ is strictly increasing on $[-1,-1+\varepsilon]$ and strictly decreasing on $[1-\varepsilon,1]$, showing that none of $-1$ and $1$ can be a local maximizer of $H_{\beta,h,p}$. 
		\medskip
		
		\noindent\emph{Proof of $(3)$.}~Define 
		$$N_{\beta,h,p}(x) := (1-x^2)H_{\beta,h,p}''(x) = \beta p(p-1)x^{p-2}(1-x^2) - 1,$$ 
		for $x \in (-1,1)$. Note that on $(-1,1)$, $N_{\beta,h,p}'(x) = \beta p(p-1)x^{p-3}(p-2-px^2)$ has exactly two roots $\pm \sqrt{1-2/p}$, for $p=3$, and an additional root $0$ for $p\geq 4$. Define: $$K_p := 2\bm{1}\{p =3\} + 3\bm{1}\{p \geq 4\}. $$ Then, by Rolle's theorem, $N_{\beta,h,p}$, and hence, $H_{\beta,h,p}''$ can have at most $K_p+1$ roots on $(-1,1)$. This shows that $H_{\beta,h,p}'$ can have at most $K_p+2$ roots on $(-1,1)$, which by part (2), include all the local maximizers of $H_{\beta,h,p}$. We now claim that for any two local maximizers $a< b$ of $H_{\beta,h,p}$, there exists a root of $H_{\beta,h,p}'$ in $(a,b)$. To see this, note that since $a$ and $b$ are local maximizers of $H_{\beta,h,p}$, by the mean value theorem, there must exist $a_1 < b_1 \in (a,b)$ such that $H_{\beta,h,p}'(a_1) \leq 0$ and $H_{\beta,h,p}'(b_1) \geq 0$. Now, by the intermediate value theorem applied on the continuous function $H_{\beta,h,p}'$, we conclude that there is a $\zeta \in (a_1,b_1)$ such that $H_{\beta,h,p}'(\zeta) = 0$. Hence, if there are $\ell$ local maximizers of $H_{\beta,h,p}$ on $(-1,1)$, then there are at least $2\ell-1$ roots of $H_{\beta,h,p}'$ on $(-1,1)$. Thus, $$2\ell-1 \leq K_p+2,\quad\textrm{i.e.}\quad \ell \leq (K_p+3)/2, $$ which proves the first part of (3). 
		
		To prove the second part of (3), first suppose that $H_{\beta,h,p}$ has three global maximizers. By the first part, $p$ must be at least $4$. We will now show that $p$ is even, by contradiction. If $p$ is odd, then $H_{\beta,h,p}''(x) < 0$ for all $x \leq 0$, and hence, by Rolle's theorem, there can be at most one non-positive root of $H_{\beta,h,p}'$. Now, if $H_{\beta,h,p}'$ has at least four positive roots, then by repeated application of Rolle's theorem, $N_{\beta,h,p}'$ has at least two positive roots. This is a contradiction, since $\sqrt{1-2/p}$ is the only positive root of $N_{\beta,h,p}'$. Hence, $H_{\beta,h,p}'$ can have at most three positive roots. Thus, $H_{\beta,h,p}'$ can have at most four roots, and hence, $H_{\beta,h,p}$ can have at most two local maximizers, a contradiction. Hence, $p$ must be even. 
		
		Next, we show that $h$ must be $0$. If $h > 0$, then $H_{\beta,h,p}(x) < H_{\beta,h,p}(-x)$ for all $x < 0$, and hence, all the three global maximizers of $H_{\beta,h,p}$ must be positive. Thus, $H_{\beta,h,p}'$ has at least $5$ positive roots, which implies that $N_{\beta,h,p}'$ has at least three positive roots, a contradiction. Similarly, if $h<0$, then all the three global maximizers of $H_{\beta,h,p}$ must be negative, and thus, $H_{\beta,h,p}'$ has at least $5$ negative roots, which implies that $N_{\beta,h,p}'$ has at least three negative roots, once again a contradiction. This shows that $h=0$. 
		
		Finally, we show that $\beta = \tilde{\beta}_p$. If $\beta > \tilde{\beta}_p$, then by  the definition in \eqref{eq:betatilde}, $0$ is not a global maximizer of $H_{\beta,h,p}$ and hence, $H_{\beta,h,p}$ being an even function, must have an even number of global maximizers, a contradiction. Therefore, it suffices to assume that $\beta<\tilde{\beta}_p$. We will show that $0$ is the only global maximizer of $H_{\beta,h,p}$, which is enough to complete the proof of the {\it only if} implication. Towards this, suppose that there is a non-zero global maximizer $x^*$ of $H_{\beta,h,p}$. Since $\beta < \tilde{\beta}_p$, we must have $H_{\beta,h,p}(x^*) = 0$, and hence, for every $\beta'  \in (\beta,\tilde{\beta}_p)$, we must have $H_{\beta' ,h,p}(x^*) > 0$, a contradiction to  the definition in \eqref{eq:betatilde}. This completes the proof of the {\it only if} implication. 
		
		For the {\it if} implication, let $\beta_N := \tilde{\beta}_p + \frac{1}{N}$, whence by part (1), $\sup_{x\in [-1,1]} H_{\beta_N,0,p}(x) > 0$ for all $N \geq 1$. Since $H_{\beta_N,0,p}(0) = 0$, for each $N$ there exists $x_N \neq 0$ such that $H_{\beta_N,0,p}(x_N) > 0$. Let $x_{N_k}$ be a convergent subsequence of $x_N$, converging to a point $x^*$. Then, $$\lim_{k \rightarrow \infty} H_{\beta_{N_k},0,p}(x_{N_k}) = H_{\tilde{\beta}_p,0,p}(x^*), $$ and hence, $H_{\tilde{\beta}_p,0,p}(x^*) \geq 0$. However, by part (1), the reverse inequality is true, and hence, $H_{\tilde{\beta}_p,0,p}(x^*) = 0$, and hence, $0, x^*$ and $-x^*$ are all global maximizers of $H_{\tilde{\beta}_p,0,p}$. We will be done, if we can show that $x^* \neq 0$. Towards this, note that since $\lim_{\varepsilon\rightarrow 0} H_{\tilde{\beta}_p,0,p}(\varepsilon)/\varepsilon^2 = -\frac{1}{2}$, there exists $\delta > 0$ such that $H_{\tilde{\beta}_p,0,p}(\varepsilon) < -\varepsilon^2/4$ whenever $|\varepsilon| < \delta$. Suppose that $x^* =0$, i.e. $x_{N_k} \rightarrow 0$ as $k \rightarrow \infty$. Then for all $k$ large enough, we must have
		\begin{equation*}
		H_{\beta_{N_k},0,p}(x_{N_k}) = H_{\tilde{\beta}_p,0,p}(x_{N_k}) + \frac{x_{N_k}^p}{N_k} < -\frac{x_{N_k}^2}{4} + \frac{x_{N_k}^p}{N_k} < 0,
		\end{equation*}
		a contradiction. This shows that $x^* \neq 0$. The proof of $(3)$ and Lemma \ref{derh11} is now complete.
	\end{proof}
	
	\begin{rem}\label{later}
		The argument in the last paragraph of the proof of Lemma \ref{derh11} can be adopted to show that for odd $p$, $H_{\tilde{\beta}_p,0,p}$ has exactly two global maximizers, one at $0$ and the other one positive.
	\end{rem}
	
	We now proceed to describe $p$-special points. To begin with, for convenience in the proof, we introduce the following notation. 
	
	\begin{definition}
		A point $(\beta,h) \in [0,\infty)\times \mathbb{R}$ is said to be $p$-{\it locally special}, if the function $H_{\beta,h,p}$ has a local maximizer $m$ satisfying $H_{\beta,h,p}''(m) = 0$.
	\end{definition}
	
	We will see that every $p$-locally special point is actually $p$-special, and hence, the two notions are identical. In the following lemma, we give exact expressions for $p$-special points.
	\begin{lem}[Description of $p$-special points]\label{derh22}Define $$\check{\beta}_p := \frac{1}{2(p-1)} \left(\frac{p}{p-2}\right)^{\frac{p-2}{2}}\quad\textrm{and}\quad \check{h}_p := \tanh^{-1}\left(\sqrt{\frac{p-2}{p}}\right) - p\check{\beta}_p \left(\frac{p-2}{p}\right)^{\frac{p-1}{2}}. $$ Then, we have the following:
		\begin{enumerate}
			\item[$(1)$] If $p\geq 3$ is odd, then $\left(\check{\beta}_p, \check{h}_p\right)$ is the only $p$-locally special point in $[0,\infty)\times \mathbb{R}$. In this case, $m_*:= \sqrt{1-2/p}$ is the only solution to the equation $H_{\check{\beta}_p, \check{h}_p,p}''(x) =0$. In fact, $m_*$ is a global maximizer of $H_{\check{\beta}_p, \check{h}_p,p}$ satisfying $H_{\check{\beta}_p, \check{h}_p,p}^{(3)}(m_*) =0$ and $H_{\check{\beta}_p, \check{h}_p,p}^{(4)}(m_*) <0$. Further, $m_*$ is the unique stationary point of $H_{\check{\beta}_p, \check{h}_p,p}$.
			\item[$(2)$] If $p\geq 4$ is even, then $\left(\check{\beta}_p, \check{h}_p\right)$ and $\left(\check{\beta}_p, -\check{h}_p\right)$ are the only $p$-locally special points in $[0,\infty)\times \mathbb{R}$. In this case, $m_*(1) := \sqrt{1-2/p}$ and $m_*(-1) := -m_*(1)$ are the only solutions to each of the equations $H_{\check{\beta}_p, i\check{h}_p,p}''(x) =0$ for $i \in \{-1,1\}$. In fact, $m_*(i)$ is a global maximizer of $H_{\check{\beta}_p, i\check{h}_p,p}$ for $i \in \{-1,1\}$ satisfying $$H_{\check{\beta}_p, i\check{h}_p,p}^{(3)}(m_*(i)) =0 \text{ and } H_{\check{\beta}_p, i\check{h}_p,p}^{(4)}(m_*(i)) <0, \quad \text{ for } i \in \{-1,1\}.$$ Further, $m^*(i)$ is the unique global maximizer of $H_{\check{\beta}_p, i\check{h}_p,p}$ for $i \in \{-1,1\}$.
		\end{enumerate}
		Hence, a point $(\beta,h)$ is $p$-locally special if and only if it is $p$-special.
	\end{lem}

	\noindent{\it Proof of Lemma} \ref{derh22}: We start the following proposition: 
	
	\begin{prop}\label{midlem}
		Let $\beta := \check{\beta}_p$, $h \in \mathbb{R}$, and let $y \in (0,1)$ be a local maximum of $H_{\beta,h,p}$, satisfying $H_{\beta,h,p}''(y)=H_{\beta,h,p}^{(3)}(y) = 0$. Then $H_{\beta,h,p}^{(4)}(y) < 0$.
	\end{prop}
	\begin{proof}
		For convenience, we will denote $N_{\beta,h,p} := (1-x^2)H_{\beta,h,p}''(x)$ by $N$ and $H_{\beta,h,p}$ by $H$. Note that $$N''(x) = (1-x^2) H^{(4)}(x) - 4x H^{(3)}(x) - 2 H''(x).$$ 
		By hypothesis, $N''(y) = (1-y^2)H^{(4)}(y)$. Now, $$N''(x) = \beta p(p-1)(p-2)(p-3)x^{p-4} - \beta p^2(p-1)^2 x^{p-2}$$ cannot have any root other than $0$ and $\pm \sqrt{\frac{(p-2)(p-3)}{p(p-1)}}$. But we know from the proof of Lemma \ref{derh22} that $H_{\beta,h,p}''$ cannot have any root other than $\pm \sqrt{1-2/p}$ (note that Proposition \ref{midlem} is not needed to reach this conclusion, and hence, there is no circularity in the argument), and for $p\geq 3$, we have $\frac{(p-2)(p-3)}{p(p-1)} < \frac{p-2}{p}$. Therefore, $y$ is not a root of $N''$, and hence, not a root of $H^{(4)}$. Proposition \ref{midlem} now follows from the standard higher derivative test. 
	\end{proof}
	
	We are now proceed with the proof of Lemma \ref{derh22}. We start by proving that the first coordinate of every $p$-locally special point in $[0,\infty)\times \mathbb{R}$ must be equal to $\check{\beta}_p$. Towards this, we first claim that $H_{\beta,h,p}''(x) < 0$, or equivalently, $N_{\beta,h,p}(x) < 0$ for all $x \in (-1,1)$, if $\beta < \check{\beta}_p$. This will rule out the possibility of $(\beta,h)$ being a candidate for a $p$-locally special point, for $\beta <\check{\beta}_p$. Towards proving this claim, we  can assume that $$\sup_{x \in (-1,1)} N_{\beta,h,p}(x) > -1, $$ since otherwise we would be done. Since $N_{\beta,h,p}(-1) =
	N_{\beta,h,p}(0) = N_{\beta,h,p}(1) = -1$, the function $N_{\beta,h,p}$ attains maximum at some $m \in (-1,1)\setminus\{0\}$, and hence, $m$ is a non-zero solution to the equation $N_{\beta,h,p}'(x) = 0$. Therefore, from the proof of (3) in Lemma \ref{derh11}, that $m \in \{-q,q\}$, where $q := \sqrt{1-2/p}$. Since $N_{\beta,h,p}(q) \geq N_{\beta,h,p}(-q)$, we know for sure that $q$ is a global maximizer of $N_{\beta,h,p}$. Our claim now follows from the observation that $\beta < \check{\beta}_p \implies N_{\beta,h,p}(q) < 0$. 
	
	Now, we are going to rule out the possibility $\beta > \check{\beta}_p$, as well. Suppose that $\beta > \check{\beta}_p$, and let $m_*$ be a local maximizer of $H_{\beta,h,p}$ satisfying $H_{\beta,h,p}''(m_*) = 0$,
	i.e. $N_{\beta,h,p}(m_*) = 0$. Now, $N_{\beta,h,p}(0) = -1 \implies m_* \neq 0$. Next, since $\beta > \check{\beta}_p$, it follows that $N_{\beta,h,p}(q) > 0$, and hence, $m_* \neq q$. If $p$ is even, then $N_{\beta,h,p}(-q) = N_{\beta,h,p}(q) > 0$, and if $p$ is odd, then $N_{\beta,h,p}(x) < -1$ for all $x < 0$. Thus, in either case, $m_* \neq -q$. All these show that $N_{\beta,h,p}'(m_*) \neq 0$. Suppose that $N_{\beta,h,p}'(m_*) > 0$. Since $N_{\beta,h,p}(m_*) = 0$, there exists $\varepsilon > 0$ such that $N_{\beta,h,p}(x) > 0$ for all $x \in (m_*,m_*+\varepsilon)$ and $N_{\beta,h,p}(x) < 0$ for all $x \in (m_*,m_*-\varepsilon)$. Thus, $H_{\beta,h,p}''(x) > 0$ for all $x \in (m_*,m_*+\varepsilon)$ and $H_{\beta,h,p}''(x) < 0$ for all $x \in (m_*-\varepsilon,m_*)$. Since $H_{\beta,h,p}'(m_*) = 0$, we must have $$H_{\beta,h,p}'(x) > 0\quad\textrm{for all}~x \in (m_*-\varepsilon,m_*+\varepsilon)\setminus \{m_*\}. $$ This implies that $H_{\beta,h,p}$ is strictly increasing on $[m_*,m_*+\varepsilon)$, contradicting that $m_*$ is a local maximizer of $H_{\beta,h,p}$. Similarly, if $N_{\beta,h,p}'(m_*) < 0$, then there exists $\varepsilon > 0$ such that $H_{\beta,h,p}'(x) < 0$ for all $x \in (m_*-\varepsilon,m_*+\varepsilon)\setminus \{m_*\}$, and so, $H_{\beta,h,p}(x)$ is strictly decreasing on $(m_*-\varepsilon,m_*]$, contradicting once again, that $m_*$ is a local maximizer of $H_{\beta,h,p}$. We have thus proved our claim, that the first coordinate of every $p$-special point in $[0,\infty)\times \mathbb{R}$ must be equal to $\check{\beta}_p$. 
	In what follows, let $\beta := \check{\beta}_p$.
	\medskip
	
	\noindent\emph{Proof of $(1)$.} Let $p \geq 3$ be odd and let $m_*$ be any solution to the equation $H_{\beta,h,p}''(x) = 0$, or equivalently, to the equation $N_{\beta,h,p}(x) = 0$. Since $N_{\beta,h,p}(x) \leq -1$ for all $x \leq 0$, it follows that $m_* \in (0,1)$. Now, we already know that the only positive root of $N_{\beta,h,p}'$ is $q:=\sqrt{1-2/p}$, and since $N_{\beta,h,p}(q) = 0$, by Rolle's theorem, $N_{\beta,h,p}$ cannot have any positive root other than $q$. Thus, $m_* = q$ is the only root of $H_{\beta,h,p}''$.  Since $N_{\beta,h,p}(m_*) = N_{\beta,h,p}'(m_*) = 0$, we have
	$$H_{\beta,h,p}^{(3)}(m_*) = \frac{N_{\beta,h,p}'(m_*)(1-m_*^2) + 2m_*N_{\beta,h,p}(m_*)}{(1-m_*^2)^2} = 0.$$ Now, $m_*$ is a stationary point of $H_{\beta,h,p}$, i.e. $H_{\beta,h,p}'(m_*) = 0$ if and only if $h = \check{h}_p$. Hence, $(\check{\beta}_p,\check{h}_p)$ is the only candidate for being a $p$-locally special point in $[0,\infty)\times \mathbb{R}$. Let $h := \check{h}_p$ throughout the rest of the proof of (a). Since $H_{\beta,h,p}'(m_*) = 0$ and $m_*$ is the only root of $H_{\beta,h,p}''$, by Rolle's theorem, $H_{\beta,h,p}'$ cannot have any root other than $m_*$. This implies that the sign of $H_{\beta,h,p}'$ remains constant on each of the intervals $(-1,m^*)$ and $(m^*,1)$. Since $$\lim_{x \rightarrow -1^+} H_{\beta,h,p}' (x) = +\infty\quad\textrm{and}\quad\lim_{x \rightarrow 1^-} H_{\beta,h,p}' (x) = -\infty, $$ we conclude that $H_{\beta,h,p}' > 0$ on $(-1,m^*)$ and $H_{\beta,h,p}' < 0$ on $(m^*,1)$, thereby showing that $m^*$ is a global maximizer, and also the unique stationary point of $H_{\beta,h,p}$, and verifying that $(\check{\beta}_p,\check{h}_p)$ is actually a $p$-special point. The result in part (1) now follows from Proposition \ref{midlem}.
	\medskip
	
	\noindent\textit{Proof of} (2): Let $p \geq 4$ be even. Since $m_*(1)$ and $m_*(-1)$ are the only non-zero roots of $N_{\beta,h,p}'$, and they are also roots of $N_{\beta,h,p}$, by Rolle's theorem, they are the only roots of $N_{\beta,h,p}$, as well. Hence, the only roots of $H_{\beta,h,p}''$ are $m_*(1)$ and $m_*(-1)$, and so, $H_{\beta,h,p}^{(3)}(m_*(1)) = H_{\beta,h,p}^{(3)}(m_*(-1)) = 0$. 
	
	For $i \in \{-1,1\}$, note that $m_*(i)$ is a stationary point of $H_{\beta,h,p}$, i.e. $H_{\beta,h,p}'(m_*(i)) = 0$, if and only if $h = i\check{h}_p$. Hence, $(\check{\beta}_p,\check{h}_p)$ and $(\check{\beta}_p,-\check{h}_p)$ are the only candidates for being $p$-locally special points in $[0,\infty)\times \mathbb{R}$. Let $h := \check{h}_p$ throughout the rest of the proof of (2). Since $H_{\beta,ih,p}'(m_*(i)) = 0$ and $m_*(i)$ is the only root of $H_{\beta,ih,p}''$ with sign $i$, by Rolle's theorem, $H_{\beta,ih,p}'$ cannot have $0$ or any point with sign $i$ as a root, other than $m_*(i)$. This implies that the sign of $H_{\beta,h,p}'$ remains constant on each of the intervals $[0,m_*(1))$ and $(m_*(1),1)$, and the sign of $H_{\beta,-h,p}'$ remains constant on each of the intervals $(-1,m_*(-1))$ and $(m_*(-1),0]$. Since $$\lim_{x \rightarrow -1^+} H_{\beta,\pm h,p}' (x) = +\infty\quad \textrm{and}\quad\lim_{x \rightarrow 1^-} H_{\beta,\pm h,p}' (x) = -\infty, $$ we conclude that $H_{\beta,h,p}' < 0$ on $(m_*(1),1)$ and $H_{\beta,-h,p}' > 0$ on $(-1,m_*(-1))$. Now, note that 
	\begin{align*}
	h = \tanh^{-1}\left(\sqrt{\frac{p-2}{p}}\right) - \check{\beta}_p p \left(\frac{p-2}{p}\right)^{\frac{(p-1)}{2}} & =\left[\sum_{k=0}^\infty \frac{\left(\sqrt{\frac{p-2}{p}}\right)^{2k+1}}{2k+1}\right] - \frac{p}{2(p-1)} \sqrt{\frac{p-2}{p}}\\
	&\geq  \sqrt{\frac{p-2}{p}} - \frac{p}{2(p-1)} \sqrt{\frac{p-2}{p}} \\ 
	& = \frac{p-2}{2(p-1)}\sqrt{\frac{p-2}{p}} > 0. 
	\end{align*}     
	Hence, $H_{\beta,h,p}'(0) = h > 0$ and $H_{\beta,-h,p}'(0) = -h < 0$. Consequently, $H_{\beta,h,p}' > 0$ on $[0,m_*(1))$ and $H_{\beta,-h,p}' < 0$ on $(m_*(-1),0]$. Thus, $m_*(i)$ is the unique global maximizer of $H_{\beta,ih,p}$ over the interval $i[0,1] := \{ix: x \in [0,1]\}$. Now, it is easy to see that $H_{\beta,ih,p}(x) < H_{\beta,ih,p}(-x)$ for all $x \in [-1,1]\setminus i[0,1]$. This shows that $m_*(i)$ is the unique global maximizer of $H_{\beta,ih,p}$ over $[-1,1]$. Part (2) now follows from Proposition \ref{midlem}, and the proof of Lemma \ref{derh22} is now complete. \qed \\

	Next, we give a description of $p$-weakly critical points that is, points $(\beta,h)$ for which the function $H_{\beta,h,p}$ has exactly two global maximizers). Note that we already have a full characterization of $p$-strongly critical points (that is, points $(\beta,h)$ for which the function $H_{\beta,h,p}$ has exactly three global maximizers) by part (3) of Lemma \ref{derh11}. To elaborate, we know that there cannot be any $p$-strongly critical point if $p$ is odd, and if $p\geq 4$ is even, then $(\tilde{\beta}_p,0)$ is the only $p$-strongly critical point. In the following lemma, we show that the set of all $p$-critical points is a one-dimensional continuous curve in the plane $[0,\infty) \times \mathbb{R}$. We also prove some other interesting properties of this curve, for instance, the only limit point(s) of the curve which is (are) outside it, is (are) the $p$-special point(s).  
	
	\begin{lem}[Description of $p$-weakly critical points]\label{derh33}
		For every $p \geq 3$, $\check{\beta}_p < \tilde{\beta}_p$, and the set $\cp^+$ can be characterized as follows.
		\begin{enumerate}
			\item[$(1)$] For every even $p \geq 4$, there exists a continuous function $\varphi_p: (\check{\beta}_p,\infty)\mapsto [0,\infty)$ which is strictly decreasing on $(\check{\beta}_p,\tilde{\beta}_p)$ and vanishing on $[\tilde{\beta}_p,\infty)$, such that 
			$$\mathscr{C}_p^+ = \left\{(\beta, \pm \varphi_p(\beta)): \beta \in (\check{\beta}_p,\infty)\setminus \{\tilde{\beta}_p\}\right\}.$$ 
			\item[$(2)$] For every odd $p \geq 3$, there exists a strictly decreasing, continuous function $\varphi_p: (\check{\beta}_p,\infty)\mapsto \mathbb{R}$ satisfying $\varphi_p(\tilde{\beta}_p) = 0$ and $\lim_{\beta \rightarrow \infty}\varphi_p(\beta) = -\infty$, such that
			$$\mathscr{C}_p^+ = \left\{(\beta, \varphi_p(\beta)): \beta \in (\check{\beta}_p,\infty)\right\}.$$ 
		\end{enumerate}
		In both cases, $\lim_{\beta \rightarrow \check{\beta}_p^+} \varphi_p(\beta) = \tanh^{-1}(m_*) - p\check{\beta}_p m_*^{p-1}$, where $m_* := \sqrt{\frac{p-2}{p}}$.
	\end{lem}
	\begin{proof}
		First, we prove that $\check{\beta}_p < \tilde{\beta}_p$ for all $p \geq 3$. Since $$\sup_{x\in [-1,1]} H_{\beta,0,p+1}(x) = \sup_{x\in [0,1]} H_{\beta,0,p+1}(x) \leq \sup_{x\in [0,1]} H_{\beta,0,p}(x) = \sup_{x\in [-1,1]} H_{\beta,0,p}(x),$$ it follows that $\tilde{\beta}_{p+1} \geq \tilde{\beta}_p$, i.e. $\tilde{\beta}_p$ is increasing in $p$. Therefore, $\tilde{\beta}_p \geq \tilde{\beta}_2  = \frac{1}{2}$ for all $p \geq 3$. First note that $\check{\beta}_3 = \frac{\sqrt{3}}{4} < \frac{1}{2}$. Next, note that for $p \geq 4$,
		$$\check{\beta}_p = \frac{1}{2(p-1)} \left(1+\frac{2}{p-2}\right)^{\frac{p-2}{2} } \leq \frac{e}{2(p-1)} \leq  \frac{e}{6} < \frac{1}{2}.$$ Hence, $\check{\beta}_p < \frac{1}{2} \leq \tilde{\beta}_p$ for all $p \geq 3$. 
		
		Next, we show that $\mathscr{C}_p^+ \subseteq (\check{\beta}_p,\infty) \times \mathbb{R}$. Towards this, first let $\beta < \check{\beta}_p$ and $h\in \mathbb{R}$. It follows from the proof of Lemma \ref{derh22}, that $H_{\beta,h,p}'' < 0$ on $[-1,1]$, so $H_{\beta,h,p}$ is strictly concave on $[-1,1]$, and hence, can have at most one global maximum. Therefore, $(\beta, h)\notin \mathscr{C}_p^+$. Now, let $\beta = \check{\beta}_p$ and $h \in \mathbb{R}$. From the proof of Lemma \ref{derh22}, we know that $H_{\beta,h,p}''$ cannot have any root on $[-1,1]$ other than possibly $\pm \sqrt{1-2/p}$. Since $H_{\beta,h,p}''(-1) = H_{\beta,h,p}''(1)=-\infty$, $H_{\beta,h,p}''(0) = -1$ and $H_{\beta,h,p}''$ is continuous, $H_{\beta,h,p}''(x) < 0$ for all $x \in [-1,1]\setminus \{\pm\sqrt{1-2/p} \}$. This shows that $H_{\beta,h,p}'$ is strictly decreasing on $[-1,1]$, and hence, $H_{\beta,h,p}$ can have at most one stationary point. Consequently, $(\beta, h)\notin \mathscr{C}_p^+$, proving our claim that $\mathscr{C}_p^+ \subseteq (\check{\beta}_p,\infty) \times \mathbb{R}$. We now consider the cases of even and odd $p$ separately.
		
		\medskip
		\noindent\textit{Proof of} (1): Let $p \geq 4$ be even. Since $x \mapsto \beta x^p - I(x)$ is an even function, the set $\mathscr{C}_p^+$ is symmetric about the line $h=0$, i.e. $(\beta,h) \in \cp^+ \implies (\beta,-h) \in \cp^+$. Next, we show that for every $\beta > \check{\beta}_p$, there exists at most one $h \geq 0$ such that $(\beta,h) \in \cp^+$. Suppose towards a contradiction, that there exists $\beta > \check{\beta}_p$ and $h_2>h_1\geq 0$, such that both $(\beta,h_1)$ and $(\beta,h_2) \in \cp^+$. Letting $m_* := \sqrt{1-2/p}$, it follows that $H_{\beta,h,p}''(m_*) > 0$ for all $h \in \mathbb{R}$. Recalling that $H_{\beta,h,p}''$ can have at most two roots in $[0,1]$, and using the facts $$H_{\beta,h,p}''(0)=-1, H_{\beta,h,p}''(1) = -\infty, $$ it follows that there exist $0<a_1<m_*<a_2<1$, such that $H_{\beta,h,p}''<0$ on $[0,a_1)$, $H_{\beta,h,p}''(a_1)=0$, $H_{\beta,h,p}''>0$ on $(a_1,a_2)$, $H_{\beta,h,p}''(a_2) = 0$ and $H_{\beta,h,p}''< 0$ on $(a_2,1]$. This shows that $H_{\beta,h,p}'$ is strictly decreasing on $[0,a_1]$, strictly increasing on $[a_1,a_2]$ and strictly decreasing on $[a_2,1]$. 
		
		First assume that $h_1 > 0$, whence the two global maximizers $m_1(h_i) < m_2(h_i)$ of $H_{\beta,h_i,p}$ must be positive roots of $H_{\beta,h_i,p}'$ for $i \in \{1,2\}$. Note that the monotonicity pattern of the function $H_{\beta,h_i,p}'$ implies that $m_1(h_i) \in (0,a_1)$ and $m_2(h_i) \in (a_2,1)$. Hence, $H_{\beta,h_i,p}'(a_1) < 0$ and $H_{\beta,h_i,p}'(a_2) > 0$, and by the intermediate value theorem, there exists $m(h_i) \in (a_1,a_2)$ such that $$H_{\beta,h_i,p}'(m(h_i)) = 0. $$ Observe that $H_{\beta,h_i,p}'$ is positive on $[0,m_1(h_i))$, negative on $(m_1(h_i),m(h_i))$, positive on $(m(h_i),m_2(h_i))$ and negative on $(m_2(h_i),1]$. Since $h_2 > h_1$, it follows that $H_{\beta,h_2,p}'>0$ on $[0,m_1(h_1)]$ and on $[m(h_1),m_2(h_1)]$. However, since $m_1(h_2), m(h_2)$ and $m_2(h_2)$ are roots of $H_{\beta,h_2,p}'$ on $(0,a_1), (a_1,a_2)$ and $(a_2,1)$ respectively, it follows that $m_1(h_1)<m_1(h_2)$, $m(h_2)<m(h_1)$ and $m_2(h_1)<m_2(h_2)$. 
		Combining all these, gives 
		\begin{equation}\label{intim1}
		\int_{m_1(h_1)}^{m(h_1)} H_{\beta,h_1,p}'(t) \mathrm d t < \int_{m_1(h_2)}^{m(h_2)} H_{\beta,h_1,p}' (t) \mathrm d t < \int_{m_1(h_2)}^{m(h_2)} H_{\beta,h_2,p}' (t) \mathrm d t
		\end{equation}
		and 
		\begin{equation}\label{intim2}
		\int_{m(h_1)}^{m_2(h_1)} H_{\beta,h_1,p}' (t) \mathrm d t < \int_{m(h_1)}^{m_2(h_1)} H_{\beta,h_2,p}' (t) \mathrm d t < \int_{m(h_2)}^{m_2(h_2)} H_{\beta,h_2,p}' (t) \mathrm d t
		\end{equation}
		Adding \eqref{intim1} and \eqref{intim2}, we have
		\begin{equation}\label{intim3}
		\int_{m_1(h_1)}^{m_2(h_1)} H_{\beta,h_1,p}' (t) \mathrm d t < \int_{m_1(h_2)}^{m_2(h_2)} H_{\beta,h_2,p}' (t) \mathrm d t .
		\end{equation}
		This is a contradiction, since both sides of \eqref{intim3} are $0$.
		
		Therefore, it must be that $h_1 = 0$. In this case, the global maximizers $m_1(h_1) < m_2(h_1)$ of $H_{\beta,h_1,p}$ satisfy $m_1(h_1) = - m_2(h_1)$. Since $H_{\beta,h_1,p}'$ vanishes at $0$, it must be negative on $(0,a_1]$. Hence, $m_2(h_1) \in (a_2,1)$. This shows that $H_{\beta,h_1,p}'(a_2) > 0$, and hence, there exists $m(h_1) \in (a_1,a_2)$ such that $H_{\beta,h_1,p}'(m(h_1)) = 0$. Observe that $H_{\beta,h_1,p}'$ is negative on $(0,m(h_1))$, positive on $(m(h_1),m_2(h_1))$ and negative on $(m_2(h_1),1)$. Therefore, since $h_2>h_1$, $H_{\beta,h_2,p}'>0$ on $[m(h_1),m_2(h_1)]$. Since $m(h_2)$ and $m_2(h_2)$ are roots of $H_{\beta,h_2,p}'$ on $(a_1,a_2)$ and $(a_2,1)$ respectively, we must have $m(h_2)<m(h_1)$ and $m_2(h_1) < m_2(h_2)$. Hence, we have
		\begin{equation}\label{intim11}
		\int_{0}^{m(h_1)} H_{\beta,h_1,p}' (t) \mathrm d t < \int_{m_1(h_2)}^{m(h_2)} H_{\beta,h_1,p}' (t) \mathrm d t < \int_{m_1(h_2)}^{m(h_2)} H_{\beta,h_2,p}' (t) \mathrm d t 
		\end{equation}
		and 	
		\begin{equation}\label{intim21}
		\int_{m(h_1)}^{m_2(h_1)} H_{\beta,h_1,p}' (t) \mathrm d t < \int_{m(h_1)}^{m_2(h_1)} H_{\beta,h_2,p}' (t) \mathrm d t  < \int_{m(h_2)}^{m_2(h_2)} H_{\beta,h_2,p}' (t) \mathrm d t 
		\end{equation}
		Adding \eqref{intim11} and \eqref{intim21}, gives 
		\begin{equation}\label{intim31}
		\int_{0}^{m_2(h_1)} H_{\beta,h_1,p}' (t) \mathrm d t < \int_{m_1(h_2)}^{m_2(h_2)} H_{\beta,h_2,p}'(t) \mathrm d t .
		\end{equation}
		Once again, this is a contradiction, since the right side of \eqref{intim31} is $0$, whereas the left side of \eqref{intim31} is non-negative. This completes the proof of our claim that for every $\beta > \check{\beta}_p$, there exists at most one $h \geq 0$ such that $(\beta,h) \in \cp^+$.

		We now show that for all $\beta \in (\check{\beta}_p,\infty)\setminus \{\tilde{\beta}_p\}$, there exists at least one $h \geq 0$ such that $(\beta,h) \in \cp^+$. First, suppose that $\beta > \tilde{\beta}_p$. In this case, $\sup_{x\in [-1,1]} H_{\beta,0,p}(x) > 0$ by  the definition in \eqref{eq:betatilde}, and hence, $H_{\beta,0,p}$ has a non-zero global maximizer $m_*$. Since $H_{\beta,0,p}$ is an even function, $-m_*$ is also a global maximizer. It now follows from part (3) of Lemma \ref{derh11}, that $H_{\beta,0,p}$ has exactly two global maximizers, and hence, $(\beta,0) \in \cp^+$.
		
		Next, let $\beta \in (\check{\beta}_p,\tilde{\beta}_p)$. Recall that the function $H_{\beta,0,p}'$ is continuous and strictly decreasing on each of the intervals $[0,a_1]$ and $[a_2,1)$. Hence, the functions $$\psi_1 := H_{\beta,0,p}'\Big|_{[0,a_1]} \quad  \text{and} \quad \psi_2 := H_{\beta,0,p}'\Big|_{[a_2,1)}$$ are invertible, and by Proposition 2.1 in \cite{irish}, the functions $\psi_1^{-1}$ and $\psi_2^{-1}$ are continuous. Hence, the function $\Lambda: [H_{\beta,0,p}'(a_1), \min\{0,H_{\beta,0,p}'(a_2) \} ] \rightarrow \mathbb{R}$ defined as:
		$$\Lambda(h) := \int_{\psi_1^{-1}(h)}^{\psi_2^{-1}(h)} H_{\beta,-h,p}' (t) \mathrm d t = \int_{\psi_1^{-1}(h)}^{\psi_2^{-1}(h)} H_{\beta,0,p}' (t) \mathrm d t + h\left(\psi_1^{-1}(h)-\psi_2^{-1}(h)\right)$$ is continuous. Since the function $t \mapsto H_{\beta,0,p}'(t) - H_{\beta,0,p}'(a_1)$ is strictly positive on the interval $( a_1, \psi_2^{-1}(H_{\beta,0,p}'(a_1)) )$ (because it is strictly increasing on $[a_1,a_2]$, strictly decreasing on $[a_2,1)$, and vanishes at the endpoints $a_1$ and $\psi_2^{-1}(H_{\beta,0,p}'(a_1))$ of the interval), 
		\begin{equation}\label{hshiftintm1}
		\Lambda(H_{\beta,0,p}'(a_1)) = \int_{a_1}^{\psi_2^{-1}(H_{\beta,0,p}'(a_1))} \left(H_{\beta,0,p}'(t) - H_{\beta,0,p}'(a_1)\right) \mathrm dt > 0.
		\end{equation}
		Next, suppose that $H_{\beta,0,p}'(a_2) \leq 0$. Since the function $t \mapsto H_{\beta,0,p}'(t) - H_{\beta,0,p}'(a_2)$ is strictly negative on the interval $(\psi_1^{-1}(H_{\beta,0,p}'(a_2)), a_2 )$ (because it is strictly decreasing on $[0,a_1]$, strictly increasing on $[a_1,a_2]$, and vanishes at the endpoints $\psi_1^{-1}(H_{\beta,0,p}'(a_2))$ and $a_2$ of the interval), 
		\begin{equation}\label{hshiftintm2}
		\Lambda(H_{\beta,0,p}'(a_2)) = \int_{\psi_1^{-1}(H_{\beta,0,p}'(a_2))}^{a_2} \left(H_{\beta,0,p}'(t) - H_{\beta,0,p}'(a_2)\right) \mathrm dt < 0.
		\end{equation}
		Finally, suppose that $H_{\beta,0,p}'(a_2) > 0$. Then we have
		\begin{equation}\label{hshiftintm3}
		\Lambda(0) = \int_{0}^{\psi_2^{-1}(0)} H_{\beta,0,p}' (t) \mathrm d t = H_{\beta,0,p}(\psi_2^{-1}(0)) < 0.
		\end{equation}
		The last inequality in \eqref{hshiftintm3} follows from the facts that $\psi_2^{-1}(0) > 0$ and $\beta < \tilde{\beta}_p$. 
		
		Using \eqref{hshiftintm1}, \eqref{hshiftintm2}, \eqref{hshiftintm3} and the intermediate value theorem, we conclude that there exists $h(\beta) \in (H_{\beta,0,p}'(a_1), \min\{0,H_{\beta,0,p}'(a_2) \})$ such that $\Lambda(h(\beta)) = 0$, i.e. 
		\begin{equation}\label{eqgl}
		H_{\beta,-h(\beta),p}(\psi_1^{-1}(h(\beta))) = H_{\beta,-h(\beta),p}(\psi_2^{-1}(h(\beta))).
		\end{equation}
		Now, $\psi_1^{-1}(h(\beta)) \in (0,a_1)$ and $\psi_2^{-1}(h(\beta)) \in (a_2,1)$, and hence, $H_{\beta,-h(\beta),p}'$ is strictly decreasing on some open neighborhoods of $\psi_1^{-1}(h(\beta))$ and $\psi_2^{-1}(h(\beta))$. Since $H_{\beta,-h(\beta),p}'(\psi_1^{-1}(h(\beta))) = H_{\beta,-h(\beta),p}'(\psi_2^{-1}(h(\beta))) = 0$, the points $\psi_1^{-1}(h(\beta))$ and $\psi_2^{-1}(h(\beta))$ are local maximizers of $H_{\beta,-h(\beta),p}$. Since $-h(\beta)>0$, any global maximizer of $H_{\beta,-h(\beta),p}$ must be a positive root of $H_{\beta,-h(\beta),p}'$, and further, it cannot lie on the interval $[a_1,a_2]$, since $H_{\beta,-h(\beta),p}'$ is strictly increasing on this interval. Hence, one of $\psi_1^{-1}(h(\beta))$ and $\psi_2^{-1}(h(\beta))$ must be a global maximizer of $H_{\beta,-h(\beta),p}$, and by \eqref{eqgl}, both must be global maximizers of $H_{\beta,-h(\beta),p}$. By part (3) of Lemma \ref{derh11}, these are the only global maximizers of $H_{\beta,-h(\beta),p}$, and hence, $(\beta,-h(\beta)) \in \cp^+$.

		Next, if $\beta = \tilde{\beta}_p$, then $H_{\beta,0,p}$ has three global maximizers, so $(\beta,0) \notin \cp^+$. One of these global maximizers is $0$ and the other two are negative of one another. It follows from the argument used in proving the uniqueness of $h$ under the case $h_1 = 0$, that $$\int_{m_1(h)}^{m_2(h)} H_{\beta,h,p}'(t) \mathrm d t > 0,$$ for every $h > 0$, where $m_2(h) > m_1(h) > 0$ are possible global maximizers of $H_{\beta,h,p}$ (see inequality \eqref{intim31}), which is a contradiction. Hence, $$\cp^+ \subseteq \left(\{\tilde{\beta}_p\} \times \mathbb{R}\right)^c. $$ At this point, we completed proving that for every $\beta \in (\check{\beta}_p,\infty)\setminus \{\tilde{\beta}_p\}$, there exists unique $h \geq 0$ such that $(\beta,h) \in \cp^+$, and further, there exists no such $h$ for $\beta = \tilde{\beta}_p$. Denote by $\varphi_p(\beta)$, this unique $h$ corresponding to $\beta \in (\check{\beta}_p,\infty)\setminus \{\tilde{\beta}_p\}$. Our proof so far, also reveals that $\varphi_p(\beta) = 0$ for $\beta > \tilde{\beta}_p$ and $\varphi_p(\beta) > 0$ for $\beta \in (\check{\beta}_p,\tilde{\beta}_p)$. Define $\varphi_p(\tilde{\beta}_p) = 0$ for the sake of completing its definition on the whole of $(\check{\beta}_p,\infty)$. 
		
		We now show that $\varphi_p$ is strictly decreasing on $(\check{\beta}_p,\tilde{\beta}_p)$. Towards this, take $\check{\beta}_p < \beta_1<\beta_2<\tilde{\beta}_p$. Let $h_1 := \varphi_p(\beta_1)$ and $h_2 := \varphi_p(\beta_2)$ (we already know from the proof of the existence part, that $h_1$ and $h_2$ are positive), and suppose towards a contradiction, that $h_1 \leq h_2$. Then, $H_{\beta_1,h_1,p}' < H_{\beta_2,h_2,p}'$ on $(0,1]$. Let $m_{11}<m_{13}$ be the global maximizers of $H_{\beta_1,h_1,p}$ and $m_{21}<m_{23}$ be the global maximizers of $H_{\beta_2,h_2,p}$. Also, let $m_{12} \in (m_{11},m_{13})$ and $m_{22}\in (m_{21},m_{23})$ be local minimizers of $H_{\beta_1,h_1,p}$ and $H_{\beta_2,h_2,p}$, respectively. We have already shown that for $i \in \{1,2\}$, the function $H_{\beta_i,h_i,p}'$ is positive on $[0,m_{i1})$, negative on $(m_{i1},m_{i2})$, positive on $(m_{i2},m_{i3})$ and negative on $(m_{i3},1)$. Since $H_{\beta_2,h_2,p}' > 0$ on $[0,m_{11}]$, we must have $m_{21} > m_{11}$. On the other hand, we have $m_{21} < m_* := \sqrt{1-2/p} < m_{13}$. This, combined with the fact that $H_{\beta_2,h_2,p}' > 0$ on $[m_{12},m_{13}]$, implies that $m_{21}<m_{12}$. Next, since $H_{\beta_1,h_1,p}' < 0$ on $[m_{21},m_{22}]$ and $H_{\beta_1,h_1,p}'(m_{12}) = 0$, it follows that $m_{22} < m_{12}$. Finally, since $H_{\beta_1,h_1,p}' <0$ on $[m_{23},1)$, we must have $m_{13}<m_{23}$. Hence, we have
		$$m_{11} < m_{21} <m_{22}<m_{12}<m_{13}<m_{23}.$$
		Using this and proceeding exactly as in the proof of the uniqueness of $h$, we have
		$$\int_{m_{11}}^{m_{12}} H_{\beta_1,h_1,p}' (t) \mathrm d t < \int_{m_{21}}^{m_{22}} H_{\beta_2,h_2,p}' (t) \mathrm d t \quad\textrm{and}\quad\int_{m_{12}}^{m_{13}} H_{\beta_1,h_1,p}' (t) \mathrm d t < \int_{m_{22}}^{m_{23}} H_{\beta_2,h_2,p}' (t) \mathrm d t .$$
		Adding the above two inequalities, we have
		$$\int_{m_{11}}^{m_{13}} H_{\beta_1,h_1,p}' (t) \mathrm d t < \int_{m_{21}}^{m_{23}} H_{\beta_2,h_2,p}' (t) \mathrm d t ,$$ which is a contradiction once again, since both sides of the above inequality are $0$. Hence, we must have $h_1 > h_2$, showing that $\varphi_p$ is strictly decreasing on $(\check{\beta}_p,\tilde{\beta}_p)$.

		Next, we show that $\varphi_p$ is continuous on $(\check{\beta}_p,\tilde{\beta}_p]$. Towards this, first take $\beta \in (\check{\beta}_p,\tilde{\beta}_p)$, and let $\{\beta_n\}_{n\geq 1}$ be a monotonic sequence in $(\check{\beta}_p,\tilde{\beta}_p)$ converging to $\beta$. Since $\varphi_p$ is decreasing on $(\check{\beta}_p,\tilde{\beta}_p)$, it follows that $\varphi_p(\beta_n)$ is monotonic as well (the direction of monotonicity being opposite to that of $\beta_n$). Moreover, $\varphi_p(\beta_n)$ is bounded between $\varphi_p(\beta_1)$ and $\varphi_p(\beta)$. Hence, $\lim_{n \rightarrow \infty} \varphi_p(\beta_n)$ exists, which we call $h$. Let $m_1(n) <m_2(n)$ denote the global maximizers of $H_{\beta_n,\varphi_p(\beta_n),p}$. Choose a subsequence $n_k$ such that $m_1(n_k) \rightarrow m_1$ and $m_2(n_k)\rightarrow m_2$ for some $m_1,m_2 \in [-1,1]$. Since $$H_{\beta_{n_k},\varphi_p(\beta_{n_k}),p}(m_i(n_k)) \geq H_{\beta_{n_k},\varphi_p(\beta_{n_k}),p}(x)\quad \textrm{for all}~ x \in [-1,1]~\textrm{and}~i \in \{1,2\}, $$ taking limit as $k \rightarrow \infty$ on both sides, we have $H_{\beta,h,p}(m_i) \geq H_{\beta,h,p}(x)$ for all $x \in [-1,1]$ and $i \in \{1,2\}$, showing that $m_1$ and $m_2$ are global maximizers of $H_{\beta,h,p}$. We now show that $m_1 < m_2$. Since $\beta_n \rightarrow \beta > \check{\beta}_p$, there exists $\underline{\beta} > \check{\beta}_p$ such that $\beta_n >  \underline{\beta}$ for all large $n$. If $a_1(\underline{\beta}) < a_2(\underline{\beta})$ are the positive roots of $H_{\underline{\beta},0,p}''$, then $H_{\beta_n,0,p}'' > 0$ on $[a_1(\underline{\beta}), a_2(\underline{\beta})]$ for all large $n$, and hence, $m_1(n) < a_1(\underline{\beta})$ and $m_2(n) > a_2(\underline{\beta})$ for all large $n$. This shows that $$m_1\leq a_1(\underline{\beta}) < a_2(\underline{\beta})\leq m_2$$ and hence, $m_1 < m_2$. Thus $H_{\beta,h,p}$ has at least two global maximizers. But $\beta \neq \tilde{\beta}_p$, and $H_{\beta,h,p}$ must therefore have exactly two global maximizers, showing that $(\beta,h) \in \cp^+$. Since $h \geq 0$, by the uniqueness property, we must have $h = \varphi_p(\beta)$. Hence, $\lim_{n\rightarrow \infty} \varphi_p(\beta_n) = \varphi_p(\beta)$, showing that $\varphi_p$ is continuous on $(\check{\beta}_p,\tilde{\beta}_p)$. 
		
		To show that $\lim_{\beta \rightarrow (\tilde{\beta}_p)^{-}} \varphi_p(\beta) = 0$, take a sequence $\beta_n \in (\check{\beta}_p,\tilde{\beta}_p)$ increasing to $\tilde{\beta}_p$, whence $\varphi_p(\beta_n)$ decreases to some $h \geq 0$. By the same arguments as before, it follows that $H_{\tilde{\beta}_p,h,p}$ has at least two global maximizers. If $h > 0$, then $H_{\tilde{\beta}_p,h,p}$ will have exactly two global maximizers. Therefore $(\tilde{\beta}_p,h) \in \cp^+$, contradicting our finding that $\cp^+ \subseteq (\{\tilde{\beta}_p\} \times \mathbb{R})^c$. This shows that $h=0$, completing the proof of (1). 
		
		\vskip 10pt
		\noindent\textit{Proof of} (2): Let $p \geq 3$ be odd. In this case, $H_{\beta,0,p}''<0$ on $[-1,0]$ for all $\beta \geq 0$. Let $\beta > \check{\beta}_p.$ Once again, $H_{\beta,0,p}''$ can have at most two positive roots, which, together with the facts $H_{\beta,0,p}''(m_*) > 0$ and $H_{\beta,0,p}''(1) = -\infty$, imply the existence of $0<a_1<m_*<a_2<1$, such that $H_{\beta,0,p}'' < 0$ on $[-1,a_1)\bigcup (a_2,1]$ and $H_{\beta,0,p}''> 0$ on $(a_1,a_2)$. One can now follow the proof of (a) modulo obvious modifications, to show that there exists at most one $h \in \mathbb{R}$ such that $(\beta,h) \in \cp^+$.
		
		To show the existence of at least one such $h \in \mathbb{R}$, one can once again essentially follow the proof of (a) modulo a couple of minor modifications. To be specific, if we modify the definition of $\psi_1$ to $H_{\beta,0,p}'\big|_{(-1,a_1]}$, and change the domain of $\Lambda$ to $[H_{\beta,0,p}'(a_1), H_{\beta,0,p}'(a_2)]$, then by following the proof of (a), we can show the existence of $h(\beta) \in (H_{\beta,0,p}'(a_1), H_{\beta,0,p}'(a_2))$ such that $(\beta,-h(\beta)) \in \cp^+$. If we denote the unique $h$ corresponding to each $\beta > \check{\beta}_p$ such that $(\beta,h) \in \cp^+$ by $\varphi_p(\beta)$, then continuity and the strict decreasing nature of $\varphi_p$ once again follow from the proof of (a).
		
		Next, it follows from Remark \ref{later}, that $\varphi_p(\tilde{\beta}_p) = 0$. We now show that $\lim_{\beta \rightarrow \infty} \varphi_p(\beta) = -\infty$. Towards this, note that the monotonicity pattern of $H_{\beta,\varphi_p(\beta),p}'$ for $\beta > \check{\beta}_p$ implies that $H_{\beta,\varphi_p(\beta),p}$ has exactly two local maximizers $m_1(\beta) \in (-1,a_1(\beta))$ and $m_2(\beta) \in (a_2(\beta),1)$, where $a_1(\beta)$ and $a_2(\beta)$ are the inflection points of $H_{\beta,\varphi_p(\beta),p}$, satisfying $0<a_1(\beta) < m_* < a_2(\beta)<1$ for all $\beta > \check{\beta}_p$. Hence, $m_1(\beta)$ and $m_2(\beta)$ are global maximizers of $H_{\beta,\varphi_p(\beta),p}$. Let $\beta > \tilde{\beta}_p$, whence the strictly decreasing nature of $\varphi_p$ implies that $\varphi_p(\beta) < 0$. Since $H_{\beta,\varphi_p(\beta),p}'(-1) = \infty$ and $H_{\beta,\varphi_p(\beta),p}'(0) = \varphi_p(\beta) < 0$, the intermediate value theorem implies that $m_1(\beta) < 0$. Hence, $$\beta (m_1(\beta))^p - I(m_1(\beta)) < 0,\quad\textrm{that is,}\quad H_{\beta,\varphi_p(\beta),p}(m_1(\beta)) < \varphi_p(\beta) m_1(\beta). $$
		Now, since $$H_{\beta,\varphi_p(\beta),p}(m_1(\beta)) = H_{\beta,\varphi_p(\beta),p}(m_2(\beta)) = \beta (m_2(\beta))^p + \varphi_p(\beta)m_2(\beta)- I(m_2(\beta)), $$ 
		we have $\beta (m_2(\beta))^p + \varphi_p(\beta)m_2(\beta)- I(m_2(\beta))  < \varphi_p(\beta)m_1(\beta)$. This implies, 
		\begin{equation}\label{divide}
		-2\varphi_p(\beta) > \varphi_p(\beta)(m_1(\beta) - m_2(\beta)) > \beta (m_2(\beta))^p - I(m_2(\beta)) \geq \beta m_*^p - I(m_2(\beta)).
		\end{equation}
		The proof of our claim now follows from \eqref{divide} since $\lim_{\beta\rightarrow \infty} \beta m_*^p - I(m_2(\beta)) = \infty$. This completes the proof of part (2).
		
		Finally, we prove that $\lim_{\beta \rightarrow \check{\beta}_p^+} \varphi_p(\beta) = \tanh^{-1}(m_*) - p\check{\beta}_p m_*^{p-1}$, where $m_* := \sqrt{1-2/p}$. Towards this, let $0< \varepsilon < \tilde{\beta}_p-\check{\beta}_p$ be given, and take any $$\beta \in \left(\check{\beta}_p,\check{\beta}_p + \frac{\varepsilon}{2p(p-1)}\right).$$ As before, let $0<a_1<a_2<1$ be the points such that $H_{\beta,0,p}''<0$ on $[0,a_1)\bigcup(a_2,1]$ and $H_{\beta,0,p}''>0$ on $(a_1,a_2)$. Since $H_{\check{\beta}_p,0,p}'' \leq 0$ on $[0,1]$, it follows that $H_{\beta,0,p}'' \leq (\beta-\check{\beta}_p)p(p-1) < \varepsilon/2$ on $[0,1]$. Hence, for every $h \in \mathbb{R}$, we have
		\begin{equation}\label{lengthinc}
		H_{\beta,h,p}'(a_2) - H_{\beta,h,p}'(a_1) = \int_{a_1}^{a_2} H_{\beta,0,p}'' (t) \mathrm d t \leq \varepsilon (a_2-a_1)/2 <\varepsilon/2. 
		\end{equation}
		Since $H_{\beta,0,p}''(m_*) > 0$, we must have $m_* \in (a_1,a_2)$. If $m_1 < m_2$ are the two global maximizers of $H_{\beta,\varphi_p(\beta),p}$, then $m_1 \in (0,a_1)$ and $m_2 \in (a_2,1)$. Since $H_{\beta,\varphi_p(\beta),p}'$ is strictly decreasing on each of the intervals $[0,a_1]$ and $[a_2,1)$, we must have $H_{\beta,\varphi_p(\beta),p}'(a_1) < 0$ and $H_{\beta,\varphi_p(\beta),p}'(a_2) > 0$. Hence, there exists $a_3 \in (a_1,a_2)$ such that $H_{\beta,\varphi_p(\beta),p}'(a_3) = 0$. Now, since $H_{\beta,\varphi_p(\beta),p}'$ is increasing on $[a_1,a_2]$, we have from \eqref{lengthinc}, 
		$$\big|H_{\beta,\varphi_p(\beta),p}'(a_3) - H_{\beta,\varphi_p(\beta),p}'(m_*)\big| \leq H_{\beta,\varphi_p(\beta),p}'(a_2) - H_{\beta,\varphi_p(\beta),p}'(a_1) < \varepsilon/2,$$
		and hence, $\big|H_{\beta,\varphi_p(\beta),p}'(m_*)\big| = \big|\tanh^{-1}(m_*) - p\beta m_*^{p-1} - \varphi_p(\beta)\big| < \varepsilon/2$. Now, $\big|p\beta m_*^{p-1} - p\check{\beta}_p m_*^{p-1}\big| \leq p(\beta - \check{\beta}_p) < \varepsilon/2$. By triangle inequality, we thus have
		\begin{align}\label{triang}
		&\big|\tanh^{-1}(m_*) - p\check{\beta}_p m_*^{p-1} - \varphi_p(\beta)\big|\nonumber\\ & \leq \big|\tanh^{-1}(m_*) - p\beta m_*^{p-1} - \varphi_p(\beta)\big| + \big|p\beta m_*^{p-1} - p\check{\beta}_p m_*^{p-1}\big| \nonumber \\ 
		& < \varepsilon.
		\end{align}
		Our claim now follows from \eqref{triang}. The proof of (2) and Lemma \ref{derh33} is now complete.
	\end{proof}
	
	Now, we will prove some properties of the function $H$, when the underlying parameter $(\beta,h)$ is perturbed to $(\beta_N,h_N)$, where $(\beta_N,h_N) \rightarrow (\beta,h)$, as $N \rightarrow \infty$. Investigating the properties of the function $H_{\beta_N,h_N,p}$ is especially important, since our analysis hinges more upon these perturbed functions, rather than the original function $H_{\beta,h,p}$.
	
	\begin{lem}\label{derh44}
		Suppose that $(\beta_N,h_N) \in [0,\infty)\times \mathbb{R}$ is a sequence converging to a point $(\beta,h) \in [0,\infty)\times \mathbb{R}$. Then, we have the following:
		\begin{enumerate}
			\item[$(1)$] Suppose that $(\beta,h)$ is a $p$-regular point, and let $m_*$ be the global maximizer of $H_{\beta,h,p}$. Then, for any sequence $(\beta_N,h_N) \in [0,\infty)\times \mathbb{R}$ converging to $(\beta,h)$, the function $H_{\beta_N,h_N,p}$ will have unique global maximizer $m_*(N)$ for all large $N$, and $m_*(N) \rightarrow m_*$ as $N \rightarrow \infty$. 
			\item[$(2)$]Let $m$ be a local maximizer of the function $H_{\beta,h,p}$, where the point $(\beta,h)$ is not $p$-special. Suppose that $(\beta_N,h_N) \in [0,\infty)\times \mathbb{R}$ is a sequence converging to $(\beta,h)$. Then for all large $N$, the function $H_{\beta_N,h_N,p}$ will have a local maximizer $m(N)$, such that $m(N) \rightarrow m$ as $N \rightarrow \infty$. Further, if $A\subseteq [-1,1]$ is a closed interval such that $m \in \textrm{int}(A)$ and $H_{\beta,h,p}(m) > H_{\beta,h,p}(x)$ for all $x \in A\setminus \{m\}$, then there exists $N_0\geq 1$, such that for all $N \geq N_0$, we have $H_N(m(N)) > H_N(x)$ for all $x \in A \setminus \{m(N)\}$.
		\end{enumerate}
	\end{lem}
	\begin{proof}[Proof of $(1)$] The set $\cR_p$ of all $p$-regular points is an open subset of $[0,\infty)\times \mathbb{R}$. To see this, note that $\cR_p^c$ is given by $\cp \bigcup \{(\check{\beta}_p,\check{h}_p)\}$ if $p$ is odd, and by $\cp \bigcup \{(\check{\beta}_p,\check{h}_p), (\check{\beta}_p,-\check{h}_p)\}$ if $p$ is even. By Lemma \ref{derh33}, $\cR_p^c$ is a closed set in either case. Hence, the function $H_{\beta_N,h_N,p}$ will have unique global maximizer $m_*(N)$ for all large $N$. 
		
		To show that $m_*(N) \rightarrow m_*$, let $\{N_k\}_{k \geq 1}$ be a subsequence of the natural numbers. Then, $\{N_k\}_{k \geq 1}$ will have a further subsequence $\{N_{k_\ell}\}_{\ell \geq 1}$, such that $m_*(N_{k_\ell})$ converges to some $m' \in [-1,1]$. Since $H_{\beta_{N_{k_\ell}},h_{N_{k_\ell}},p}\left(m_*(N_{k_\ell})\right) \geq H_{\beta_{N_{k_\ell}},h_{N_{k_\ell}},p}(x)$ for all $x \in [-1,1]$, by taking limit as $\ell \rightarrow \infty$ on both sides, we have $H_{\beta,h,p}(m') \geq H_{\beta,h,p}(x)$ for all $x \in [-1,1]$, showing that $m'$ is a global maximizer of $H_{\beta,h,p}$. Since $m_*$ is the unique global maximizer of $H_{\beta,h,p}$, it follows that $m' = m_*$, completing the proof of (1).  
		
		\medskip
		
		\noindent\noindent\textit{Proof of} (2):~Let us denote $H_{\beta,h,p}$ by $H$ and $H_{\beta_N,h_N,p}$ by $H_N$. It is easy to show that there exists $M \geq 1$ odd, and points $-1=a_0<a_1<\ldots<a_M=1$, such that $H'$ is strictly decreasing on $[a_{2i},a_{2i+1}]$ and strictly increasing on $[a_{2i+1},a_{2i+2}]$ for all $0\leq i \leq \frac{M-1}{2}$. Hence, the local maximizer $m$ of $H$ lies in $(a_{2i},a_{2i+1})$ for some $0\leq i \leq \frac{M-1}{2}$. Since $H'(a_{2i}) > 0$ and $H'(a_{2i+1})<0$, we also have $H_N'(a_{2i}) > 0$ and $H_N'(a_{2i+1})<0$ for all large $N$, and hence $H_N'$ has a root $m(N) \in (a_{2i},a_{2i+1})$ for all large $N$. 
		
		Let us now show that $m(N) \rightarrow m$. Towards this, let $\{N_k\}_{k\geq 1}$ be a subsequence of the natural numbers, whence there is a further subsequence $\{N_{k_\ell}\}_{\ell \geq 1}$ of $\{N_k\}_{k\geq 1}$, such that $m(N_{k_\ell}) \rightarrow m'$ for some $m' \in [a_{2i},a_{2i+1}]$. Since $H_{N_{k_\ell}}'(m(N_{k_\ell})) = 0$ for all $\ell \geq 1$, we have $H'(m') = 0$. But the strict decreasing nature of $H'$ on $[a_{2i},a_{2i+1}]$ implies that $m$ is the only root of $H'$ on this interval, and hence, $m'=m$. This shows that $m(N) \rightarrow m$. 
		
		Next, we show that $m(N)$ is a local maximizer of $H_N$ for all $N$ sufficiently large. For this, we prove something stronger than needed, because this will be useful in proving the last statement of (2). Since $H''(m) < 0$, there exists $\varepsilon > 0$ such that $[m-\varepsilon,m+\varepsilon] \subset (a_{2i},a_{2i+1})$ and $H''<0$ on $[m-\varepsilon,m+\varepsilon]$. If $m_0 \in [m-\varepsilon,m+\varepsilon]$ is such that $H''(m_0) = \sup_{x \in [m-\varepsilon,m+\varepsilon]} H''(x) < 0$, then since $H_N''$ converges to $H''$ uniformly on $(-1,1)$, $$\sup_{x\in [m-\varepsilon,m+\varepsilon]}H_N''(x) < H''(m_0)/2\quad\textrm{for all large}~N.$$ In particular, since $m(N) \in [m-\varepsilon,m+\varepsilon]$ for all large $N$, we have $H_N''(m(N)) < 0$ for all large $N$, showing that $m(N)$ is a local maximizer of $H_N$ for all large $N$. Also, since $H_N'(m(N)) = 0$ and $\sup_{x\in [m-\varepsilon,m+\varepsilon]}H_N''(x) < 0$ for all large $N$, we must have $$H_N(m(N)) > H_N(x)\quad\textrm{for all}~ x \in [m-\varepsilon,m+\varepsilon]\setminus\{m(N)\},\quad \textrm{for all large} N.$$   Finally, suppose that $A \subseteq [-1,1]$ is a closed interval such that $m \in \textrm{int}(A)$ and $H(m) > H(x)$ for all $x \in A\setminus \{m\}$. By Lemma \ref{mest}, there exists $\varepsilon' >0 $ such that for all $0 < \delta \leq \varepsilon'$, $\sup_{x \in A\setminus (m-\delta,m+\delta)} H(x) = H(m\pm \delta)$. Let $\alpha = \min\{\varepsilon,\varepsilon'\}$. Then, 
		\begin{equation}\label{fs}
		H_N(m(N)) > H_N(x) \quad\textrm{for all} ~x \in [m-\alpha,m+\alpha]\setminus \{m(N)\},\quad\textrm{for all large}~ N,
		\end{equation}
		and $\sup_{x \in A\setminus (m-\alpha,m+\alpha)} H(x) = H(m\pm \alpha) < H(m)$ (since $H'(m) = 0$ and $H''< 0$ on $[m-\alpha,m+\alpha]$). Hence, 
		\begin{equation}\label{ss}
		\sup_{x \in A\setminus (m-\alpha,m+\alpha)} H_N(x) < H_N(m(N))\quad\textrm{for all large}~ N. 
		\end{equation}
		The proof of (2) now follows from \eqref{fs} and \eqref{ss}, and the proof of Lemma \ref{derh44} is now complete. 
	\end{proof}
	
	\subsection{Technical Properties of the ML Estimates }\label{sec:techmle}
	In this subsection, we prove some technical properties related to $\hat{\beta}_{N}$ and $\hat{h}_{N}$. We begin with a result which says that the functions $u_{N,p}$ and $u_{N,1}$ appearing in the left-hand sides of equations \eqref{eqmle} and \eqref{eqmleh} are strictly increasing in both $\beta$ and $h$. This result is particularly important in the proofs of the results in Section \ref{sec:mle_beta_h_I}.
	
	\begin{lem}\label{increasing}
		For every fixed $h$, the function $\beta \mapsto F_N(\beta,h,p)$ is strictly convex, and for every fixed $\beta$, the function $h \mapsto F_N(\beta,h,p)$ is strictly convex. Consequently, the maps $u_{N,1}$ and $u_{N,p}$ defined in \eqref{eqmle} and \eqref{eqmleh} respectively, are strictly increasing in both $\beta$ and $h$.
	\end{lem}
	
	\begin{proof}
		Let $\psi_N(\beta,h) := F_N(\beta,h,p) + N\log 2 = \log \sum_{\bs \in \sa_N} e^{N\beta \overline{X}^p_N + Nh\os}$. Then for every $\beta_1, \beta_2, h$ and $\lambda \in (0,1)$, we have by H\"older's inequality,
		\begin{align*}
		\psi_N(\lambda \beta_1 + (1-\lambda)\beta_2,h)& =  \log \sum_{\bs \in \sa_N} e^{N\lambda (\beta_1 \overline{X}^p_N + h\os)} e^{N(1-\lambda)(\beta_2 \overline{X}^p_N + h \os)}\\& <  \log \left[\left(\sum_{\bs \in \sa_N} e^{N\beta_1 \overline{X}^p_N + Nh\os} \right)^\lambda \left(\sum_{\bs \in \sa_N} e^{N\beta_2 \overline{X}^p_N + Nh\os} \right)^{1-\lambda}\right]\\& =  \lambda \psi_N(\beta_1,h) + (1-\lambda)\psi_N(\beta_2,h).
		\end{align*}
		Similarly, for every $h_1,h_2,\beta$ and $\lambda \in (0,1)$, we have by H\"older's inequality,
		\begin{align*}
		\psi_N(\beta,\lambda h_1 + (1-\lambda)h_2)& =  \log \sum_{\bs \in \sa_N} e^{N\lambda (\beta \overline{X}^p_N + h_1\os)} e^{N(1-\lambda)(\beta \overline{X}^p_N + h_2 \os)}\\& <  \log \left[\left(\sum_{\bs \in \sa_N} e^{N\beta \overline{X}^p_N + Nh_1\os} \right)^\lambda \left(\sum_{\bs \in \sa_N} e^{N\beta \overline{X}^p_N + Nh_2\os} \right)^{1-\lambda}\right]\\& =  \lambda \psi_N(\beta,h_1) + (1-\lambda)\psi_N(\beta,h_2).
		\end{align*}
		This shows strict convexity of the functions $\beta \mapsto F_N(\beta,h,p)$ and $h \mapsto F_N(\beta,h,p)$. Now, note that $$\frac{\partial}{\partial \beta} F_N(\beta,h,p) = Nu_{N,p}(\beta,h,p)\quad\textrm{and}\quad \frac{\partial}{\partial h} F_N(\beta,h,p) = Nu_{N,1}(\beta,h,p).$$ Lemma \ref{increasing} now follows from the fact that the first derivative of a differentiable, strictly convex function is strictly increasing.
	\end{proof}
	
	In the following Lemma, we show that for fixed $\beta$, the ML Estimate  of $h$ exists, and for fixed $h$, the ML Estimate  of $\beta$ exists, asymptotically almost surely. However, if $p$ is even, then the joint ML Estimate  of $(\beta,h)$ does not exist.
	
	\begin{lem}\label{mleexist}
		Fix $N \geq 1$. Then $\hat{\beta}_N$ and $\hat{h}_N$ exist in $[-\infty,\infty]$ and are unique.  Further, $\hat{h}_N $ exists in $\mathbb{R}$ if and only if $|\os| < 1$. For odd $p$, $\hat{\beta}_N$ exists in $\mathbb{R}$ if and only if $|\os| < 1$, and for even $p$, $\hat{\beta}_N$ exists in $\mathbb{R}$ if and only if $\os \notin \{-1,0,1\}$. Hence, $$\lim_{N \rightarrow \infty} \p_{\beta,h,p} \left(\hat{\beta}_N~\textrm{exists in}~\mathbb{R}\right) = \lim_{N \rightarrow \infty} \p_{\beta,h,p} \left(\hat{h}_N~\textrm{exists in}~\mathbb{R}\right) =1.$$ However, if $p$ is even, then for all  $N \geq 1$ and all $\bs \in \sa_N$, the joint ML Estimate  of $(\beta,h)$ does not exist.
	\end{lem}
	\begin{proof}
		The log-likelihood function is given by
		$$\ell_p(\beta,h|\bs) = - N \log 2 + \beta N\overline{X}^p_N + hN\os - F_N(\beta,h,p).$$ By Lemma \ref{increasing}, the functions $\beta \mapsto F_N(\beta,h,p)$ and $h \mapsto F_N(\beta,h,p)$ are strictly convex, and hence, the functions $\beta \mapsto \ell_p(\beta,h|\bs)$ and $h \mapsto \ell_p(\beta,h|\bs)$ are strictly concave. Consequently, $\beta \mapsto \ell_p(\beta,h|\bs)$ attains maximum at $\hat{\beta} \in \mathbb{R}$ if and only if $\frac{\partial}{\partial \beta} \ell_p(\beta,h|\bs)\Big|_{\beta = \hat{\beta}} = 0$, and $h \mapsto \ell_p(\beta,h|\bs)$ attains maximum at $\hat{h} \in \mathbb{R}$ if and only if $\frac{\partial}{\partial h} \ell_p(\beta,h|\bs)\Big|_{h = \hat{h}} = 0$. In those cases, $\hat{\beta}$ and $\hat{h}$ are the unique maximizers of $\ell_p(\beta,h|\bs)$ over $\beta \in \mathbb{R}$ and $h \in \mathbb{R}$, respectively. Now, the equations $\frac{\partial}{\partial \beta} \ell_p(\beta,h|\bs) = 0$ and $\frac{\partial}{\partial h} \ell_p(\beta,h|\bs) = 0$ are (respectively) equivalent to the equations
		\begin{equation}\label{mleeq1}
		\frac{\partial}{\partial \beta} F_N(\beta,h,p) = N\overline{X}^p_N\quad\textrm{and}\quad \frac{\partial}{\partial h} F_N(\beta,h,p) = N\os.
		\end{equation}
		One can easily show that
		\begin{equation}\label{hinm}
		\lim_{h \rightarrow \infty} \frac{\partial}{\partial h} F_N(\beta,h,p) = N\quad\textrm{and}\quad \lim_{h \rightarrow -\infty} \frac{\partial}{\partial h} F_N(\beta,h,p) = -N.
		\end{equation}
		Similarly, if $p$ is odd, we have
		\begin{equation}\label{binm1}
		\lim_{\beta \rightarrow \infty} \frac{\partial}{\partial \beta} F_N(\beta,h,p) = N\quad\textrm{and}\quad \lim_{\beta \rightarrow -\infty} \frac{\partial}{\partial \beta} F_N(\beta,h,p) = -N.
		\end{equation}
		Finally, for even $p$, we have
		\begin{equation}\label{binm2}
		\lim_{\beta \rightarrow \infty} \frac{\partial}{\partial \beta} F_N(\beta,h,p) = N\quad\textrm{and}\quad \lim_{\beta \rightarrow -\infty} \frac{\partial}{\partial \beta} F_N(\beta,h,p) = 0.
		\end{equation}
		The existence and uniqueness of $\hat{h}_N$ and $\hat{\beta}_N$ in $[-\infty,\infty]$, and the necessary and sufficient conditions about the existence of $\hat{h}_N$ and $\hat{\beta}_N$ in $\mathbb{R}$ now follow from \eqref{mleeq1}, \eqref{hinm}, \eqref{binm1} and \eqref{binm2}, since the functions $h \mapsto \frac{\partial}{\partial h} F_N(\beta,h,p)$ and $\beta \mapsto \frac{\partial}{\partial \beta} F_N(\beta,h,p)$ are strictly increasing and continuous. 
		
		Next, we show that the ML estimates are real valued with probability (under $\beta,h$) going to $1$. Towards this, first note that under $\p_{\beta,h,p} $, $\os$ converges weakly to a discrete measure supported on the set of all global maximizers of $H_{\beta,h,p}$ (see Theorem \ref{cltintr1}). Since $-1$ and $1$ are not global maximizers of $H_{\beta,h,p}$, it follows that $\p_{\beta,h,p} (|\os| = 1) \rightarrow 0$ as $N \rightarrow \infty$. If $h\neq 0$, then $0$ is not a global maximizer of $H_{\beta,h,p}$, so $\p_{\beta,h,p} (\os \in \{-1,0,1\}) \rightarrow 0$ as $N \rightarrow \infty$. Therefore, assume that $h=0$. By Stirling-type bounds, $$\p_{\beta,0,p}(\os = 0) = \frac{1}{2^N}\binom{N}{\frac{N}{2}}Z_{N}(\beta,0,p)^{-1}\bm{1}\{N~\textrm{is even}\} \leq \frac{e}{\pi \sqrt{N}},$$ 
		where the last inequality uses the fact that $F_N(\beta,0,p) \geq 0$ for all $\beta \geq 0$. Hence, $\p_{\beta,0,p}(\os = 0) \rightarrow 0$ as $N \rightarrow \infty$, completing the proof of the finiteness of $\hat{h}_N$ and $\hat{\beta}_N$ for all $\beta,h,p$.

		Finally, let $p$ be even and $N \geq 1$. If the joint ML Estimate  of $(\beta,h)$ exists, then by \eqref{eqmle} and \eqref{eqmleh}, we must have
		$\e_{\beta,h,p}  \overline{X}^p_N = \left(\e_{\beta,h,p}  \os\right)^p$.
		Since each of the measures $\p_{\beta,h,p} $ has support $\mathcal{M}_N$ and the function $x \mapsto x^p$ is non-affine, convex on $\mathcal{M}_N$, we arrive at a contradiction to Jensen's inequality.
	\end{proof}

	\subsection{Other Technical Lemmas}\label{sec:other_tech}
	In this subsection, we collect the proofs of the remaining technical lemmas, which are used in the proofs of the main results in various places. We start with a result that gives implicit expressions for the partial derivatives of any stationary point of $H_{\beta,h,p}$ with respect to $\beta$ and $h$. This result is required in the proof of Theorem \ref{cltun}.

	\begin{lem}\label{derh2}
		Let $m = m(\beta,h,p)$ satisfy the implicit relation $H_{\beta,h,p}'(m) = 0$, and suppose that $H_{\beta,h,p}''(m) \neq 0$. Then, the partial derivatives of $m$ with respect to $\beta$ and $h$ are given by:
		\begin{equation}\label{stpar}
		\frac{\partial m}{\partial \beta} = -\frac{pm^{p-1}}{H_{\beta,h,p}''(m)}\quad\quad\textrm{and}\quad\quad \frac{\partial m}{\partial h} = -\frac{1}{H_{\beta,h,p}''(m)}.
		\end{equation}
		Moreover, $\big|\frac{\partial^2 m}{\partial \beta^2}\big| < \infty$ and $\big|\frac{\partial^2 m}{\partial h^2}\big| < \infty$, if $H_{\beta,h,p}''(m) \neq 0$. 
	\end{lem}
	
	\begin{proof} 
		Differentiating both sides of the identity $\beta p m^{p-1} +h - \tanh^{-1}(m) = 0$ with respect to $\beta$ and  $h$ separately, we get the following two first order partial differential equations, respectively:
		\begin{equation}\label{diffbeta}
		pm^{p-1} + \beta p(p-1)m^{p-2}\frac{\partial m}{\partial \beta} - \frac{1}{1-m^2} \frac{\partial m}{\partial \beta} = 0,\quad\textrm{that is,}\quad pm^{p-1} + H_{\beta,h,p}''(m)\frac{\partial m}{\partial \beta} = 0~;
		\end{equation}
		\begin{equation}\label{diffh}
		\beta p(p-1)m^{p-2}\frac{\partial m}{\partial h} + 1 - \frac{1}{1-m^2} \frac{\partial m}{\partial h} = 0,\quad\textrm{that is,}\quad 1+H_{\beta,h,p}''(m)\frac{\partial m}{\partial h} = 0~; 
		\end{equation}
		The expressions in \eqref{stpar} follow from \eqref{diffbeta} and \eqref{diffh}. Another implicit differentiation of \eqref{diffbeta} with respect to $\beta$ and \eqref{diffh} with respect to $h$ yields the following two second order partial differential equations, respectively:
		\begin{equation}\label{diffbeta22}
		2p(p-1) m^{p-2} \frac{\partial m}{\partial \beta} + H_{\beta,h,p}^{(3)}(m)\left(\frac{\partial m}{\partial \beta}\right)^2 + H_{\beta,h,p}''(m)\frac{\partial^2 m}{\partial \beta^2} = 0;
		\end{equation}
		\begin{equation}\label{diffh33}
		H_{\beta,h,p}^{(3)}(m)\left(\frac{\partial m}{\partial h}\right)^2 + H_{\beta,h,p}''(m)\frac{\partial^2 m}{\partial h^2} = 0;
		\end{equation}
		The finiteness of the second order partial derivatives of $m$ as long as $H_{\beta,h,p}''(m) \neq 0$, now follow from the fact that $H_{\beta,h,p}''(m)$ is the coefficient of $\frac{\partial^2 m}{\partial \beta^2}$ and $\frac{\partial^2 m}{\partial h^2}$ in the differential equations \eqref{diffbeta22} and \eqref{diffh33}.
	\end{proof}
	
	We now derive some important properties of the function $\zeta$ defined in \eqref{xidef}. The following lemma is used in the proof of Lemma \ref{ex}.
	\begin{lem}\label{xiact}
		For any sequence $x \in (-1,1)$ that is bounded away from both $1$ and $-1$, we have
		$$\zeta(x) = \sqrt{\frac{2}{\pi N (1-x^2)}} e^{NH_N(x)}\left(1+O(N^{-1})\right).$$
	\end{lem}
	\begin{proof}
		The proof of Lemma \ref{xiact} follows immediately from Lemma \ref{stir}.
	\end{proof}

	Now, we bound the derivative of the function $\zeta$ in a neighborhood of the point $m_*(N)$. This result appears in the proof of Lemma \ref{ex}.

	\begin{lem}\label{imp1}
		For every $\alpha \geq 0$ and $p$-regular point $(\beta,h)$, we have the following bound:
		$$\sup_{x\in A_{N,\alpha}} |\zeta'(x)| = \zeta(m_*(N))O\left(N^{\frac{1}{2}+\alpha}\right),$$
		where $A_{N,\alpha} := \left(m_*(N)-N^{-\frac{1}{2} + \alpha},m_*(N) +N^{-\frac{1}{2} + \alpha}\right)$ and $m_*(N)$ is the global maximizer of $H_N$. 
	\end{lem}

	\begin{proof}[Proof of Lemma \ref{imp1}] We begin with the following lemma: 
		\begin{lem}\label{xiprime}
			For any sequence $x \in (-1,1)$ that is bounded away from both $1$ and $-1$, we have
			$$\zeta'(x) =\zeta(x)\left(NH_N'(x) + \frac{x}{1-x^2} + O(N^{-1})\right).$$
		\end{lem}
		
		\begin{proof}
			By Lemma \ref{bindiff} and \eqref{digamexp}, we have	
			\begin{align}\label{ff1}
			& \dfrac{\mathrm d}{\mathrm d x} \binom{N}{N(1+x)/2}\nonumber\\ & = \frac{N}{2}\binom{N}{N(1+x)/2}\left[\psi\left(1-\frac{Nx}{2}+\frac{N}{2}\right)-\psi\left(1+\frac{Nx}{2}+\frac{N}{2}\right)\right]\nonumber\\
			& = \frac{N}{2}\binom{N}{N(1+x)/2} \left(\log\left(\frac{ N}{2}(1-x)\right)-\log\left(\frac{N}{2}(1+x)\right)+\frac{1}{N(1-x)}-\frac{1}{N(1+x)}\right)\nonumber\\ &+ \frac{N}{2}\binom{N}{N(1+x)/2} O(N^{-2})\nonumber\\
			& = \binom{N}{N(1+x)/2}\left[-N \tanh^{-1}(x)+\frac{x}{1-x^2}+O(N^{-1})\right].
			\end{align} 
			We thus have by the product rule of differential calculus and \eqref{ff1},
			\begin{align*}
			\zeta'(x) & =  \zeta(x)(N\beta_N px^{p-1} + Nh_N) + \exp\left\{N(\beta_N x^p + h_Nx - \log 2)\right\}\frac{\mathrm d}{\mathrm d x} \binom{N}{N(1+x)/2}\\& =  \zeta(x)(N\beta_N px^{p-1} + Nh_N) + \zeta(x)\left[-N \tanh^{-1}(x)+\frac{x}{1-x^2}+O(N^{-1})\right]\\& =  \zeta(x)\left(NH_N'(x) + \frac{x}{1-x^2} + O(N^{-1})\right),
			\end{align*}
			completing the proof of Lemma \ref{xiprime}.
		\end{proof}
		
		Now, we proceed with the proof of Lemma \ref{imp1}. First note that, since $H_N'(m_*(N)) = 0$, we have by the mean value theorem,
		\begin{equation}\label{imp1eq1}
		\sup_{x \in A_{N,\alpha}}\big|H_N'(x)\big| \leq \sup_{x \in A_{N,\alpha}}\big|x-m_*(N)\big|\sup_{x\in A_{N,\alpha}}|H_N''(x)|=O\left(N^{-\frac{1}{2}+\alpha}\right).
		\end{equation}
		It follows from \eqref{imp1eq1} and Lemma \ref{xiprime} that
		\begin{equation}\label{imp1eq2}
		\sup_{x \in A_{N,\alpha}}|\zeta'(x)| \leq O\left(N^{\frac{1}{2}+\alpha}\right)\sup_{x \in A_{N,\alpha}}\zeta(x).
		\end{equation}
		Now, Lemma \ref{xiact} implies that
		\begin{equation}\label{imp1eq3}
		\sup_{x \in A_{N,\alpha}} \zeta(x) \leq \left(1+O(N^{-1})\right) \zeta(m_*(N))\sup_{x\in A_{N,\alpha}} \sqrt{\frac{1-m_*(N)^2}{1-x^2}} = \zeta(m_*(N))O(1).
		\end{equation}
		Lemma \ref{imp1} now follows from \eqref{imp1eq2} and \eqref{imp1eq3}.
	\end{proof}
	
	Lemma \ref{imp1} has an analogous version for $p$-special points $(\beta,h)$, which is stated below. In this case, the bound on $\zeta'$ is better, and holds on a slightly larger region, too. 
	
	\begin{lem}\label{xiprimeirreg}
		Let $m_*(N)$ be the unique global maximizer of the function $H_N := H_{\beta_N,h_N,p}$, where $(\beta_N,h_N) := \left(\beta + \bar{\beta}N^{-3/4}, h +\bar{h}N^{-3/4}\right)$ for some $\bar{\beta},\bar{h}\in \mathbb{R}$, and $(\beta,h)$ is a $p$-special point. Then, for all $\alpha \geq 0$,
		$$\sup_{x\in \mathcal{A}_{N,\alpha}} |\zeta'(x)| = \zeta(m_*(N))O\left(N^{\frac{1}{4}+3\alpha}\right)$$
		where $\mathcal{A}_{N,\alpha} := \left(m_*(N) - N^{-\frac{1}{4}+\alpha}, m_*(N) + N^{-\frac{1}{4} + \alpha}\right)$.
	\end{lem}
	
	\begin{proof}
		The proof of Lemma \ref{xiprimeirreg} is similar to that of Lemma \ref{imp1}, the only difference being a change in the estimate of $\sup_{x \in \mathcal{A}_{N,\alpha}} |H_N'(x)|$ from the estimate in \eqref{imp1eq1}. Note that
		\begin{align*}
		\sup_{x\in \mathcal{A}_{N,\alpha}} |H_N''(x)| &= \sup_{x\in \mathcal{A}_{N,\alpha}} |H''(x)| + O\left(N^{-\frac{3}{4}}\right)\\  
		&\leq \sup_{x \in \mathcal{A}_{N,\alpha}} \tfrac{1}{2}(x-m_*)^2\sup_{x\in \mathcal{I}(\mathcal{A}_{N,\alpha}\cup \{m_*\})}H^{(4)} (x) + O\left(N^{-\frac{3}{4}}\right) = O\left(N^{-\frac{1}{2}+ 2\alpha}\right),
		\end{align*}
		where $m_*$ denotes the global maximizer of $H_{\beta,h,p}$ and for a set $A \subseteq \mathbb{R}$, $\mathcal{I}(A)$ denotes the smallest interval containing $A$. The last equality follows from the observation $$\sup_{x\in \mathcal{A}_{N,\alpha}} |x-m_*| \leq \sup_{x\in \mathcal{A}_{N,\alpha}} |x-m_*(N)| + |m_*(N) - m_*| \leq N^{-\frac{1}{4} + \alpha} + O\left(N^{-\frac{1}{4}}\right) = O\left(N^{-\frac{1}{4} + \alpha}\right),$$ by Lemma \ref{mnminusm}. Following \eqref{imp1eq1}, we have
		$$\sup_{x\in \mathcal{A}_{N,\alpha}} |H_N'(x)| = O\left(N^{-\frac{3}{4} + 3\alpha}\right).$$ The rest of the proof is exactly same as that of Lemma \ref{imp1}.
	\end{proof}
	
	In the next lemma, we prove an asymptotic expansion of a local maximum value of the perturbed function $H_N$, around the corresponding local maximum value of the original function $H$. This is required in the proof of Lemma \ref{redunsup}.  
	
	\begin{lem}\label{bdhbest}
		Let $m$ be a local maximizer of $H := H_{\beta,p,h}$. Let $\beta_N := \beta + \bar{\beta}x_N$ and $h_N := h + \bar{h} y_N$ for some fixed constants $\beta,\bar{\beta},h,\bar{h}$ and sequences $x_N, y_N \rightarrow 0$. Suppose that the point $(\beta,h)$ is not $p$-special. Let $m(N)$ denote the local maximizer of $H_{\beta_N,h_N,p}$ converging to $m$. Then we have as $N \rightarrow \infty$,
		$$H_{\beta_N,h_N,p}(m(N)) = H(m) + \bar{\beta} x_N m^p + \bar{h} y_N m + O\left((x_N + y_N)^2\right).$$
	\end{lem}
	\begin{proof}
		For any sequence $({\beta}_N', {h}_N')  \rightarrow (\beta,h)$, let us denote by $m({\beta}_N', {h}_N' ,p)$ the local maximum of $H_{{\beta}_N',{h}_N',p}$ converging to $m$. In particular, $m(\beta_N,h_N,p) = m(N)$ and $m(\beta,h,p) = m$. By a simple application of Taylor's theorem and Lemma \ref{derh2}, we have
		\begin{align}\label{ty1}
		m(N) - m & = m(\beta_N,h_N,p) - m(\beta_N,h,p) + m(\beta_N,h,p) - m(\beta,h,p)\nonumber\\
		&= -\frac{\bar{h} y_N}{H_{\beta_N,h,p}''(m(\beta_N,h,p))} -  \frac{\bar{\beta}pm^{p-1} x_N}{H''(m)} + O\left(x_N^2 + y_N^2\right)\nonumber\\
		&= O(x_N + y_N).
		\end{align}
		By another application of Taylor's theorem, we have
		\begin{align}\label{ty2}
		H_{\beta_N,h_N,p}(m(N)) - H(m) &  = H_{\beta_N,h_N,p}(m(N)) - H_{\beta_N,h_N,p}(m) + H_{\beta_N,h_N,p}(m) - H(m)\nonumber\\
		&= O\left((m(N)-m)^2\right) + \bar{\beta} x_N m^p + \bar{h}y_N m.
		\end{align}
		Lemma \ref{bdhbest} now follows from \eqref{ty1} and \eqref{ty2}.
	\end{proof}
	
	The following lemma provides estimates of the first four derivatives of the function $H$ at the maximizer $m_*(N)$ for a perturbation of a $p$-special point. This key result is used in the proof of Lemma \ref{ex2}. 
	
	\begin{lem}\label{mnminusm}
		Let $(\beta,h)$ be a $p$-special point and $(\beta_N,h_N) := (\beta+\bar{\beta}N^{-\frac{3}{4}},h+\bar{h}N^{-\frac{3}{4}} )$ for some $\bar{\beta},\bar{h} \in \mathbb{R}$. If $m_*$ and $m_*(N)$ denote the unique global maximizers of $H := H_{\beta,h,p}$ and $H_N := H_{\beta_N,h_N,p}$ respectively, then we have the following:
		\begin{align} 
		\label{id1} N^\frac{1}{4}(m_*(N) - m_*) & = -\left(\frac{6(\bar{\beta}pm_*^{p-1} + \bar{h})}{H^{(4)}(m_*)}\right)^\frac{1}{3} + O\left(N^{-\frac{1}{4}}\right), \\
		\label{id2}	N^\frac{1}{2}H''(m_*(N)) & = \frac{1}{2}\left(6(\bar{\beta}pm_*^{p-1} + \bar{h})\right)^\frac{2}{3}\left(H^{(4)}(m_*)\right)^{\frac{1}{3}} + O\left(N^{-\frac{1}{4}}\right), \\
		\label{id3} N^\frac{1}{4}H^{(3)}(m_*(N)) & = -\left(6(\bar{\beta}pm_*^{p-1} + \bar{h})\right)^\frac{1}{3}\left(H^{(4)}(m_*)\right)^\frac{2}{3} + O\left(N^{-\frac{1}{4}}\right), \\
		\label{id4}
		H^{(4)}(m_*(N)) & = H^{(4)}(m_*) + O\left(N^{-\frac{1}{4}}\right).
		\end{align} 
	\end{lem} 
	
	\begin{proof}
		Let us start by noting that $$H'(m_*(N)) = H_N'(m_*(N)) - (\bar{\beta}pm_*(N)^{p-1} + \bar{h})N^{-\frac{3}{4}} = -(\bar{\beta} p m_*(N)^{p-1}+ \bar{h})N^{-\frac{3}{4}}.$$ On the other hand, by a Taylor expansion of $H'$ around $m_*$ and using the fact $H'(m_*) = H''(m_*) = H^{(3)}(m_*) = 0$ (see Lemma \ref{derh22}), we have
		$$H'(m_*(N)) = \tfrac{1}{6}(m_*(N)-m_*)^3 H^{(4)}(\zeta_N),$$ where $\zeta_N$ lies between $m_*(N)$ and $m_*$. Hence, $$N^\frac{3}{4} (m_*(N)-m_*)^3 = -\frac{6(\bar{\beta}pm_*(N)^{p-1} + \bar{h})}{H^{(4)}(\zeta_N)}.$$	
		Now, it follows from the proof of Lemma \ref{derh44}, part (1), that $m_*(N) \rightarrow m_*$, and hence, $\zeta_N \rightarrow m_*$. This implies that 
		\begin{equation}\label{help1}
		\lim_{N \rightarrow \infty} N^{\tfrac{1}{4}} (m_*(N) - m_*) = -\left(\frac{6(\bar{\beta}pm_*^{p-1} + \bar{h})}{H^{(4)}(m_*)}\right)^\frac{1}{3} .
		\end{equation}
		By a $5$-term Taylor expansion of $H'(m_*(N))$ around $m_*$, one obtains
		\begin{equation}\label{addtayl}
		\tfrac{1}{6}(m_*(N)-m_*)^3 H^{(4)}(m_*) + \tfrac{1}{24}(m_*(N)-m_*)^4 H^{(5)}(\zeta_N') = -(\bar{\beta} p m_*(N)^{p-1}+ \bar{h})N^{-\frac{3}{4}}.
		\end{equation}
		for some sequence $\zeta_N'$ lying between $m_*(N)$ and $m_*$. From \eqref{addtayl} and \eqref{help1}, we have
		\begin{align}\label{addtaylor}
		N^\frac{3}{4} (m_*(N) - m_*)^3 & =  -\frac{6(\bar{\beta} p m_*(N)^{p-1}+ \bar{h})}{H^{(4)}(m_*)} - \frac{N^\frac{3}{4}(m_*(N) - m_*)^4 H^{(5)}(\zeta_N')}{4H^{(4)}(m_*)}\nonumber\\& =  -\frac{6(\bar{\beta} p m_*^{p-1}+ \bar{h})}{H^{(4)}(m_*)} + O\left(N^{-\frac{1}{4}}\right).
		\end{align}
		\eqref{id1} now follows from \eqref{addtaylor}, and \eqref{id2}, \eqref{id3}, \eqref{id4} follow by substituting \eqref{id1} into the following expansions
		\begin{equation*}
		H''(m_*(N)) = \tfrac{1}{2}\left(m_*(N) - m_*\right)^2 H^{(4)}(m_*) + O\left((m_*(N)-m_*)^3\right),
		\end{equation*}
		\begin{equation*}
		H^{(3)}(m_*(N)) = \left(m_*(N) - m_*\right)H^{(4)}(m_*) + O\left((m_*(N)-m_*)^2\right),
		\end{equation*}
		and 	$H^{(4)}(m_*(N))  = H^{(4)}(m_*) + O(m_*(N) - m_*)$. 
	\end{proof}
	
	The final lemma shows that if a function has non-vanishing curvature at a unique point of maxima, then for every sufficiently small open interval $I$ around that point of maxima, it attains its maximum on $I^c$ at either of the endpoints of $I$. This fact is used in the proofs of Lemmas \ref{conc} and \ref{lem:multiple_max}. 
	
	\begin{lem}\label{mest} Let $A \subseteq [-1,1]$ be a closed interval. Suppose that $f: A \mapsto \mathbb{R}$ is continuous on $A$ and twice continuously differentiable on $\mathrm{int}(A)$. Suppose that there exists $x_*\in \mathrm{int}(A)$ such that $f(x_*) > f(x)$ for all $x \in A\setminus \{x_*\}$, and $f''(x_*) < 0$. Then, there exists $\eta > 0$ such that for all $0<\varepsilon \leq \eta$, $f$ attains maximum on the set $A\setminus (x_*-\varepsilon,x_*+\varepsilon)$ at either $x_*-\varepsilon$ or $x_*+\varepsilon$. 
	\end{lem}
	\begin{proof}
		Since $f''$ is continuous on $\mathrm{int}(A)$ and negative at $x_*$, there exists $\delta > 0$ such that $f''(x)<0$ for all $x \in (x_*-\delta,x_*+\delta)$. Hence, $f'$ is strictly decreasing on $(x_*-\delta,x_*+\delta)$. Since $f'(x_*) = 0$, we have $f'(x) > 0$ for all $x \in (x_*-\delta,x_*)$ and $f'(x) < 0$ for all $x \in (x_*,x_*+\delta)$. Hence, $f$ is strictly increasing on $(x_*-\delta,x_*]$ and strictly decreasing on $[x_*,x_*+\delta)$. 
		
		Suppose now, towards a contradiction, that the lemma is not true. Then, there is a sequence $\varepsilon_n \rightarrow 0$ such that neither $x_*-\eps_n$ nor $x_*+\eps_n$ is a point of maximum of $f$ on $A\setminus (x_*-\eps_n,x_*+\eps_n)$. Let $x_n \in A\setminus [x_*-\eps_n,x_*+\eps_n]$ be such that $f(x_n) = \sup_{x \in A\setminus (x_*-\eps_n,x_*+\eps_n)} f(x)$, which exists by the continuity of $f$ and compactness of the set $A\setminus (x_*-\eps_n,x_*+\eps_n)$. Since $f(x_*-\eps_n) \leq f(x_n) \leq f(x_*)$ for all $n$, and $f$ is continuous, it follows that $f(x_n) \rightarrow f(x_*)$. If $x_{n_k}$ is a convergent subsequence of $x_n$ converging to some $y \in A$, then by continuity of $f$, we have $f(y) = f(x_*)$. This implies that $y = x_*$. Therefore, there exists $k$ such that $x_{n_k} \in (x_*-\delta,x_*+\delta)\setminus \{x_*\}$ and $\eps_{n_k} < \delta$. For this $k$, we have $f(x_{n_k})< \max\{f(x_*-\eps_{n_k}), f(x_*+\eps_{n_k})\}$. This contradicts the fact that $x_{n_k}$ maximizes $f$ on the set $A\setminus (x_*-\eps_{n_k},x_*+\eps_{n_k})$, completing the proof of Lemma \ref{mest}.
	\end{proof}

\chapter{Properties of the Curie-Weiss Threshold} 
\label{cwtp}

Here, we will prove various properties of the Curie-Weiss threshold $\bw=\beta_{\mathrm{ER}}^*(p, 1)$ (recall \eqref{hypergraph_random}). 

\begin{lem} The Curie-Weiss threshold $\bw$ has the following properties: 
	\begin{itemize}
		\item[(1)] $\lim_{p \rightarrow \infty} \bw = \log 2$. 
		
		\item[(2)] The sequence $\{\bw\}_{p\geq 2}$ is strictly increasing. 
		
		\item[(3)]  $\beta_{\mathrm{CW}}^*(2) = 0.5$. 
	\end{itemize} 
\end{lem}

\begin{proof} Define the function $g_{\beta,p}(t) := \beta t^p - I(t).$ Since $g_{\beta,p}$(1)$ = \beta - \log 2$, recalling \eqref{hypergraph_random},  it immediately follows that $\bw \leq \log 2$. Now, take any $\beta < \log 2$. Note that $g_{\beta,2}$(1)$ < 0$ and the function $t \mapsto g_{\beta,2}(t)$ is continuous at $1$. Therefore, there exists $r \in (0,1)$, such that $g_{\beta,2}(t) < 0$ for all $t \in [r,1]$. Clearly, $g_{\beta,p}(t) \leq g_{\beta,2}(t) <0$ for all $p\geq 2$ and $t \in [r,1]$. Now, note that for all $t \in [0,1)$,
	\begin{equation}\label{conx}
	g_{\beta,p}''(t) = \beta p (p-1) t^{p-2} - \frac{1}{1-t^2} \leq  \beta p (p-1) t^{p-2} - 1.
	\end{equation}
	Since $\lim_{p\rightarrow \infty} p(p-1) r^{p-2} = 0$, there exists $p(\beta) \geq 2$, such that $g_{\beta,p}''(r) < 0$ for all $p \geq p(\beta)$. Hence, $g_{\beta,p}''(t) < 0$ for all $t \in [0,r]$ and $p \geq p(\beta)$. This, together with the fact that $g_{\beta,p}'(0) = 0$, implies that $g_{\beta,p}$ is strictly decreasing on $[0,r]$ for all $p \geq p(\beta)$. Moreover, because $g_{\beta,p}(0) = 0$, it follows that $g_{\beta,p}(t) \leq 0$ for all $t \in [0,r]$ and $p \geq p(\beta)$. Hence, $g_{\beta,p}(t) \leq 0$ for all $t \in [0,1]$ and $p \geq p(\beta)$, i.e. $\bw \geq \beta$ for all $p \geq p(\beta)$. This shows that $\bw \rightarrow \log 2$, as $p \rightarrow \infty$.
	
	Next, we show that the sequence $\{\bw\}_{p\geq 2}$ is strictly increasing. Towards this, take any $p\geq 3$. It follows from \cite[Lemma F.1]{mlepaper}, that there exists $a \in (0,1)$, such that $0$ and $a$ are both global maximizers of $g_{\bw,p}$. In particular, $g_{\bw,p}(a) = 0$, and hence, $g_{\bw,p-1}(a) > 0$. The function $\beta \mapsto g_{\beta,p-1}(a)$ being continuous, there exists $\beta < \bw$, such that $g_{\beta,p-1}(a) > 0$. Hence, $\beta_{\mathrm{CW}}^*(p-1) \leq \beta$, establishing that $\beta_{\mathrm{CW}}^*(p-1) < \bw$.
	
	Finally, we show that $\beta_{\mathrm{CW}}^*(2) = 0.5$. By \eqref{conx}, $\sup_{t\in [0,1]} g_{\beta,p}''(t) <0$ for all $\beta < 0.5$. This, coupled with the facts that $g_{\beta,p}'$ and $g_{\beta,p}$ vanish at $0$, implies that $\sup_{t \in [0,1]} g_{\beta,p}(t) = 0$ for all $\beta < 0.5$. Hence, $\beta_{\mathrm{CW}}^*(2) \geq 0.5$. On the other hand, for any $\beta > 0.5$, $g_{\beta,p}''(0) = 2\beta - 1 > 0$ and hence, by continuity of the function $g_{\beta,p}''$ at $0$, there exists $\varepsilon > 0$, such that $\inf_{t\in [0,\varepsilon]} g_{\beta,p}''(t) > 0$. Once again, since $g_{\beta,p}'$ and $g_{\beta,p}$ vanish at $0$,	we have $g_{\beta,p}(t) > 0$ for all $t \in (0,\varepsilon]$. This shows that $\beta_{\mathrm{CW}}^*(2) \leq 0.5$. 
\end{proof}

\chapter{Technical Lemmas from Chapter \ref{chap:covariateising}}
In this section, we prove some technical results, which are the key ingredients behind the Proof of Theorem \ref{finalmaintheorem}.
\begin{lem}\label{lambound}
	Suppose that $\|\nabla L_n(\bg)\|_{\infty}   \le \lambda/2$. Then, we have:
	$$\lambda\left(\|\bm v_S\|_1 - \|\bm v_{S^c}\|_1 + \frac{\|\bm v\|_1}{2}\right) \geq L_n(\hat{\bg}) - L_N(\bg) - \nabla L_N(\bg)^\top \bm v~.$$	
\end{lem}
\begin{proof}
	Denote the quantity $L_N(\hat{\bg}) - L_N(\bg) - \nabla L_N(\bg)^\top \bm v$ by $\Delta$. Now, it follows from the definition of $\hat{\bg}$, that 
	$$L_N(\hat{\bg}) + \lambda \|\hat{\bg}\|_1 \leq L_N(\bg) + \lambda \|\bg\|_1~.$$
	Hence, we have:
	\begin{equation}\label{firsteq}
	\lambda(\|\bg\|_1 - \|\hat{\bg}\|_1) \geq L_N(\hat{\bg}) - L_N(\bg) = \nabla L_N(\bg)^\top \bm v + \Delta \geq -\|\nabla L_N(\bg)\|_\infty \|\bm v\|_1 + \Delta \geq -\frac{\lambda \|\bm v\|_1}{2} + \Delta~. 
	\end{equation}
	\begin{gather*}
	\|\hat{\bg}\|_1 = \|\bg + \bm v\|_1 = \|\bg_{S} + \bm v_S + \bm v_{S^c}\|_1 = \|\bg_{S} + \bm v_S\|_1 + \|\bm v_{S^c}\|_1 \nonumber\\
	\ge \|\bg_{S}\|_1 - \|\bm v_S\|_1 + \|\bm v_{S^c}\|_1 = \|\bg\|_1 - \|\bm v_S\|_1 + \|\bm v_{S^c}\|_1
	\end{gather*}
	which implies that 
	\begin{equation}\label{secstep2}
	\|\bg\|_1 - \|\hat{\bg}\|_1 \le \|\bm v_S\|_1 - \|\bm v_{S^c}\|_1   .  
	\end{equation}
	Lemma \ref{lambound} now follows from \eqref{firsteq} and \eqref{secstep2}.	
\end{proof}

\begin{lem}\label{partialbetalips}
	Let $f(\bs) := \partial L_N(\bg) /\partial \beta$, and suppose that $\bs,\bs' \in \{-1,1\}^n$ are two vectors differing in at most one coordinate. Then $|f(\bs) - f(\bs')| \leq (2\beta + 6)/N$.
\end{lem}
\begin{proof}
	To begin with, note that:
	$$f(\bs) = -\frac{1}{N} \sum_{i=1}^N m_i(\bs)\left[X_i - \tanh(\beta m_i(\bs) +\bh^\top \bm Z_i)\right]~,$$ and hence, for any two $\bs,\bs' \in \{-1,1\}^N$, we have:
	\begin{align}\label{lips1}
	&|f(\bs) - f(\bs')|\nonumber\\ &\leq \frac{1}{N} \left|\sum_{i=1}^N \left[m_i(\bs) X_i - m_i(\bs') X_i'\right]\right|\nonumber\\ &+ \frac{1}{N} \left|\sum_{i=1}^N  \left[m_i(\bs) \tanh(\beta m_i(\bs) +\bh^\top \bm Z_i) - m_i(\bs') \tanh(\beta m_i(\bs') +\bh^\top \bm Z_i)\right]\right|~.
	\end{align}
	Now, assume that $\bs$ and $\bs'$ differ only in the $k^{\textrm{th}}$ coordinate, for some $k \in [N]$. Towards bounding the first term in the right side of \eqref{lips1}, note that
	\begin{align}\label{lipsfirst}
	\left|\sum_{i=1}^n \left[m_i(\bs) X_i - m_i(\bs') X_i'\right]\right| &= \left|\sum_{(i,j)\in [n]^2} A_{ij} (X_i X_j - X_i' X_j') \right|\nonumber\\&= \left|\sum_{i=1}^n A_{ik} (X_i X_k - X_i'X_k') + \sum_{j=1}^n A_{kj} (X_k X_j - X_k'X_j')\right|\nonumber\\&= 2\left|\sum_{i=1}^n A_{ik}X_i(X_k-X_k') \right|\nonumber\\&\leq 4\left|\sum_{i=1}^n A_{ik}\right| \leq 4\|\bm A\|_1 \leq 4~.
	\end{align}
	Next, we proceed to bound the second term in the right side of \eqref{lips1}. Towards this, first note that $m_i(\bs) - m_i(\bs') = A_{ik}(X_k - X_k')$, and hence, we have:
	\begin{align}\label{lipssecond}
	&\left|\sum_{i=1}^N  \left[m_i(\bs) \tanh(\beta m_i(\bs) +\bh^\top \bm Z_i) - m_i(\bs') \tanh(\beta m_i(\bs') +\bh^\top \bm Z_i)\right]\right|\nonumber\\&\leq \left|\sum_{i=1}^N m_i(\bs) \left[ \tanh(\beta m_i(\bs) +\bh^\top \bm Z_i) -  \tanh(\beta m_i(\bs') +\bh^\top \bm Z_i)\right]\right|\nonumber\\ &+ \left|\sum_{i=1}^N A_{ik}(X_k-X_k')  \tanh(\beta m_i(\bs') +\bh^\top \bm Z_i) \right|\nonumber\\&\leq \beta \left|\sum_{i=1}^N m_i(\bs) \left(m_i(\bs) - m_i(\bs')\right)\right| + 2\sum_{i=1}^N |A_{ik}|\nonumber\\&\leq 2\beta \sum_{i=1}^n |A_{ik}||m_i(\bs)| + 2\|\bm A\|_1\nonumber\\&\leq 2\beta \|\bm A\|_1\|\bm A\|_\infty + 2\|\bm A\|_1 = 2\beta + 2~.
	\end{align}
	Lemma \ref{partialbetalips} now follows from \eqref{lipsfirst} and \eqref{lipssecond}.
\end{proof}

\begin{lem}\label{partialthetalips}
	Let $f_j(\bs) := \partial L_N(\bg) /\partial \theta_j$ for all $j \in [d]$, and suppose that $\bs,\bs' \in \{-1,1\}^N$ are two vectors differing in at most one coordinate. Then $|f_j(\bs) - f_j(\bs')| \leq 2M(1+\beta)/n$.
\end{lem}
\begin{proof}
	To begin with, note that:
	$$f_j(\bs) = -\frac{1}{N} \sum_{i=1}^N Z_{i,j}\left[X_i - \tanh(\beta m_i(\bs) +\bh^\top \bm Z_i)\right]~,$$ and hence, for any two $\bs,\bs' \in \{-1,1\}^N$ differing in the $k^{\textrm{th}}$ coordinate only, we have:
	\begin{align*}
	|f_j(\bs) - f_j(\bs')| &\leq \frac{1}{N} \left|\sum_{i=1}^N Z_{i,j}(X_i-X_i')\right|\\ &+ \frac{1}{N} \left|\sum_{i=1}^N  Z_{i,j} \left[\tanh(\beta m_i(\bs) +\bh^\top \bm Z_i) - \tanh(\beta m_i(\bs') +\bh^\top \bm Z_i)\right]\right|\nonumber\\&= \frac{2|Z_{k,j}|}{N} + \frac{\beta}{N}\left|\sum_{i=1}^n Z_{i,j} (m_i(\bs) - m_i(\bs'))\right|\nonumber\\&= \frac{2|Z_{k,j}|}{N} + \frac{\beta}{N}\left|\sum_{i=1}^n Z_{i,j} A_{ik}(X_k-X_k')\right|\nonumber\\&\leq \frac{2M}{N} + \frac{2M\beta\|\bm A\|_1}{N}  \leq \frac{2M(1+\beta)}{N}~,\nonumber
	\end{align*}
	as desired.
\end{proof}

	\begin{lem}\label{nondiff}
		We have:
		$$\|\nabla L_n(\hat{\bg})\|_\infty \le \lambda~.$$
	\end{lem}
	\begin{proof}
		Fix $1\le i\le d$, and define the univariate function:
		$$f(x) := L_n(\hat{\gamma}_1,\ldots,\hat{\gamma}_{i-1},x,\hat{\gamma}_{i+1},\ldots,\hat{\gamma}_{d+1})~.$$
		Note that $f'(\hat{\gamma}_i) = \nabla_i L_n(\hat{\bg})$. By definition of $\hat{\bg}$, we have:
		$$f(\hat{\gamma}_i) + \lambda |\hat{\gamma}_i| \le f(x) + \lambda |x|\quad\textrm{i.e.}\quad f(x)-f(\hat{\gamma}_i) \ge \lambda(|\hat{\gamma}_i| - |x|)$$ for all $x$. 
		
		First, suppose that $\hat{\gamma}_i = 0$. Therefore, for all $x > 0$, we have:
		$$\frac{f(x) - f(0)}{x} \ge -\lambda \frac{|x|}{x} = -\lambda~,$$
		and for all $x<0$, we have:
		$$\frac{f(x) - f(0)}{x} \le -\lambda \frac{|x|}{x} = \lambda~.$$
		This shows that the left derivative of $f$ at $0$ is bounded above by $\lambda$, whereas the right derivative of $f$ at $0$ is bounded below by $-\lambda$. Since $f'$ exists, this implies that $|f'(0)| \le \lambda$. 
		
		Next, suppose that $\hat{\gamma}_i > 0$. Then, for all $x > \hat{\gamma}_i$, we have:
		$$\frac{f(x)-f(\hat{\gamma}_i)}{x-\hat{\gamma}_i} \ge \lambda \frac{|\hat{\gamma}_i|-|x|}{x-\hat{\gamma}_i} = -\lambda~,$$
		and for all $0<x<\hat{\gamma}_i$, we have:
		$$\frac{f(x)-f(\hat{\gamma}_i)}{x-\hat{\gamma}_i} \le \lambda \frac{|\hat{\gamma}_i|-|x|}{x-\hat{\gamma}_i} = -\lambda~,$$ showing that $f'(\hat{\gamma}_i) = -\lambda$.
		
		Finally, suppose that $\hat{\gamma}_i < 0$. Then, for all $0>x > \hat{\gamma}_i$, we have:
		$$\frac{f(x)-f(\hat{\gamma}_i)}{x-\hat{\gamma}_i} \ge \lambda \frac{|\hat{\gamma}_i|-|x|}{x-\hat{\gamma}_i} = \lambda~,$$
		and for all $x<\hat{\gamma}_i$, we have:
		$$\frac{f(x)-f(\hat{\gamma}_i)}{x-\hat{\gamma}_i} \le \lambda \frac{|\hat{\gamma}_i|-|x|}{x-\hat{\gamma}_i} = \lambda~,$$ showing that $f'(\hat{\gamma}_i) = \lambda$. This completes the proof of Lemma \ref{nondiff}.
	\end{proof}
	
	Below, we prove Lemmas 14 and 4 in \cite{cd_ising_estimation} for our model \eqref{eq:logisticmodel}. The following definition will be heavily used in the proofs.
	
	\begin{definition}
		Suppose that $\boldsymbol{\sigma} \in \{-1,1\}^N$ is a sample from the Ising model:
		\begin{equation}\label{mgamma}
		\p_{\beta,\bm h}(\boldsymbol{\sigma}) ~\propto~ \exp\left(\sum_{i=1}^N h_i\sigma_i + \frac{\beta}{2} \boldsymbol{\sigma}^\top \bm D \boldsymbol \sigma\right)
		\end{equation}
		where $\bm D$ is a symmetric matrix with zeros on the diagonal, and $\sup_N \|\bm D\|_\infty \le R$. Also, suppose that with probability $1$, 
		$$\min_{1\le i \le N} \mathrm{Var}(\sigma_i|\boldsymbol{\sigma}_{-i}) \ge \Upsilon~.$$ Then, we refer to the model \eqref{mgamma} as an $(R,\Upsilon)$- Ising model.
	\end{definition}

	\begin{lem}\label{lemma14}
		Let $\bs \in \{-1,1\}^N$ be a sample from an $(R,\Upsilon)$- Ising model for some $R < 1/(4\beta)$. Then, for each $i \in [N]$, we have:
		$$W_1\left(\bp_{\bm X_{-i}|X_i=1}~,~\bp_{\bm X_{-i}|X_i=-1} \right) \le \frac{8\beta R}{1-4\beta R}$$ where $W_1$ is the $L^1$-Wasserstein distance, namely
		$$W_1(\mu,\nu) := \min_{\pi \in \mathscr{C}_{\mu,\nu}} \e_{(\bm U,\bm V)\sim \pi} \|\bm U-\bm V\|_1~,$$ $\mathscr{C}_{\mu,\nu}$ denoting the set of all couplings of the probability measures $\mu$ and $\nu$.
	\end{lem} 
	\begin{proof}
		It follows from Lemma 4.9 in \cite{kostis3}, that there exists a coupling $\pi$ of the conditional measures $\bp_{\bm X_{-i}|X_i=1}$ and $\bp_{\bm X_{-i}|X_i=-1}$, such that:
		\begin{equation}\label{is}
		\mathbb{E}_{(\bm U,\bm V)\sim \pi} \left[d_H(\bm U,\bm V)\right] \le \frac{\alpha}{1-\alpha}~,
		\end{equation}
		where $d_H$ denotes the Hamming distance and $\alpha$ denotes the Dobrushin coefficient. By Lemma 4.4 in \cite{scthesis}, we know that Dobrushin's interdependence matrix is given by $4\beta\ba$. Hence, it follows from \eqref{is} that (Dobrushin's coefficient in Theorem 2.3 in \cite{kostis3} is given by $\alpha = 4\beta\|\ba\|_2$; see Theorem 4.3 in \cite{scthesis}) 
		\begin{equation}\label{is1}
		W_1\left(\bp_{\bm X_{-i}|X_i=1}~,~\bp_{\bm X_{-i}|X_i=-1} \right) \le 2~\mathbb{E}_{(\bm U,\bm V)\sim \pi} \left[d_H(\bm U,\bm V)\right] \le \frac{8\beta R }{1-4\beta R}~.
		\end{equation}
		Lemma \ref{lemma14} follows from \eqref{is1}.
	\end{proof}
	
	\begin{lem}\label{combinatorial1}
		Let $\bs\in \{-1,1\}^N$ be a sample from an $(R,\Upsilon)$- Ising model, and fix $\eta \in (0,R)$. Then there exist subsets $I_1,\ldots,I_\ell\subseteq [N]$ with $\ell \lesssim R^2 \log N/\eta^2$ such that:
		\begin{enumerate}
			\item For all $1\le i \le N$,
			$$|\{j \in \ell: i \in I_j \}| = \lceil \eta \ell / 8R\rceil$$
			\item For all $1 \le j \le \ell$, the conditional distribution of $\bm X_{I_j}$ given $\bm X_{-I_j} := (X_u)_{u\in [N]\setminus I_j}$ is an $(\eta,\Upsilon)$- Ising model.
		\end{enumerate}
	Furthermore, for any non-negative vector $\bm a \in \mathbb{R}^N$, there exists $j \in \ell$ such that
	$$\sum_{i \in I_j} a_i \ge \frac{\eta}{8R}\sum_{i=1}^N a_i~.$$
	\end{lem}
	
	\begin{proof}
		Suppose that $\bm X$ comes from the model:
		$$\p_{\beta,\bm h}(\bm X) ~\propto~\exp\left(\sum_{i=1}^N h_i X_i + \frac{\beta}{2} \bm X^\top \bm A \bm X\right)~.$$
		We apply Lemma 17 in \cite{cd_ising_estimation} on the matrix $\bm D := \bm A/R$, where $\bm A$ is the interaction matrix corresponding to the distribution of $\bm X$. Note that $\bm D$ satisfies the hypotheses of Lemma 17 in \cite{cd_ising_estimation}. Define $\eta' := \eta/R$. This ensures that $\eta' \in (0,1)$. By Lemma 17 in \cite{cd_ising_estimation}, there exist subsets $I_1,\ldots, I_\ell \subseteq [N]$ with $\ell \lesssim  R^2\log N/\eta^2$, such that for all $1\le i \le N$,  $|\{j \in \ell: i\in I_j\}| = \lceil \eta \ell/8R\rceil$, and for all $j \in \ell$, $\|\bm D_{I_j \times I_j}\|_\infty \le \eta' \implies \|\bm A_{I_j \times I_j}\|_\infty \le \eta$~\footnote{For a matrix $\bm M \in \mathbb{R}^{s\times t}$ and for sets $S\subseteq \{1,\ldots,s\}, T\subseteq \{1,\ldots,t\}$, we define $\bm M_{S\times T} := ((M_{ij}))_{i\in S, j \in T} \in \mathbb{R}^{|S|\times |T|}$}. Now, we have:
		\begin{eqnarray}\label{lastpiece}
		\frac{\mathbb{P}(\bm X_{I_j} = \bm y|\bm X_{-I_j} = \bm x_{-I_j})}{\mathbb{P}(\bm X_{I_j} = \bm y'|\bm X_{-I_j} = \bm x_{-I_j})} &=& \frac{  \exp\left(\frac{\beta}{2} \bm y_{I_j}^\top \bm A|_{I_j\times I_j} \bm y_{I_j} + \frac{\beta}{2}\sum_{u\in I_j, v\notin I_j} \bm A_{uv} y_u x_v  + \sum_{i\in I_j} h_i y_i \right)}{\exp\left(\frac{\beta}{2} \bm y_{I_j}'^\top \bm A|_{I_j\times I_j} \bm y_{I_j}' + \frac{\beta}{2}\sum_{u\in I_j, v\notin I_j} \bm A_{uv} y_u' x_v  +  \sum_{i\in I_j} h_i y_i' \right)}~,
		\end{eqnarray}
		which represents an Ising model $\mu$ on $\{-1,1\}^{|I_j|}$ with infinity norm of the interaction matrix being bounded above by $\eta$, and the external magnetic field term at site $i$ being $$h_i' = \frac{\beta}{2} \sum_{v \notin I_j} A_{iv} x_v + h_i~.$$
		The next step is to show that if $\bm Y \sim \mu$, then $\min_{1\le u\le |I_j|} \mathrm{Var}(Y_u|\bm Y_{-u}) \ge \Upsilon$. Towards this, note that for any $1\le u \le |I_j|$, we have:
		$$\p(Y_u = 1| \bm Y_{-u} = \bm y_{-u}) = \p(X_{I_j^u} = 1| \bm X_{I_j\setminus \{I_j^u\}} = \bm y_{-u}, \bm X_{-I_j} = \bm x_{-I_j})$$ where $I_j^u$ is the $u^{\mathrm{th}}$ smallest element of $I_j$. Hence,
		$$\mathrm{Var}(Y_u|\bm Y_{-u} = \bm y_{-u}) = \mathrm{Var}(X_{I_j^u} | \bm X_{I_j\setminus \{I_j^u\}} = y_{-u}, \bm X_{-I_j} = \bm x_{-I_j}) \ge  \Upsilon~,$$ thereby proving our claim.
		
		To conclude the last part of Lemma \ref{combinatorial1}, let us pick $j \in [\ell]$ uniformly at random, and note that for any vector $\bm a$, we have:
		$$\mathbb{E} \left(\sum_{i \in I_j} a_i \right) = \sum_{i=1}^N a_i \mathbb{P}(I_j \ni i) = \sum_{i=1}^N a_i \frac{\lceil \eta \ell/8R\rceil}{\ell}  \ge \frac{\eta}{8R} \sum_{i=1}^N a_i~.$$
		This means that there is a fixed sample point $j$, such that $$\sum_{i \in I_j} a_i \ge \frac{\eta}{8R} \sum_{i=1}^N a_i~.$$ This completes the proof of Lemma \ref{combinatorial1}.
	\end{proof}

	We now show that Lemma 10 in \cite{cd_ising_estimation} holds for our model \eqref{eq:logisticmodel} too.

	\begin{lem}\label{lemma10}
		Let $\bs$ be a sample from an $(R,\Upsilon)$- Ising model for some $R < 1/(4\beta)$. Then for any vector $\bm a \in \mathbb{R}^N$,
		$$\mathrm{Var}(\bm a^\top \bs) \gtrsim \frac{\|\bm a\|_2^2\Upsilon^2}{R}~.$$  
	\end{lem}  
	
	\begin{proof}
		 First, consider the case $R \le \Upsilon/(16\beta) \le 1/(16\beta)$. In this case, $1-4R\beta \ge 3/4$, so 
		\begin{equation}\label{med1}
		\frac{4R\beta}{1-4R\beta} \le \frac{\Upsilon}{3}
		\end{equation}
	Now, it follows from Lemma \ref{lemma14} that for every $i$,
	\begin{align}
	\sum_{j\in [N]\setminus \{i\}} \left|\e \left(X_j|X_i=1\right) - \e\left(X_j|X_i=-1\right)\right| &\le \sum_{j\in [N]\setminus \{i\}} W_1\left(\p_{X_j|X_i=1}, \p_{X_j|X_i=-1})\right)\nonumber\\ &\le W_1\left(\bp_{\bm X_{-i}|X_i=1}~,~\bp_{\bm X_{-i}|X_i=-1} \right)\le \frac{8\beta R}{1-4\beta R}\nonumber\\\label{med2}
	\end{align}
	
	Next, note that:
	\begin{align}
	\mathrm{Cov}(X_i,X_j) &= \e[(X_i-\e X_i) X_j] = \e\left[(X_i-\e X_i)\e(X_j|X_i)\right]\nonumber\\&= \p(X_i=1) (1-\e X_i) \e(X_j|X_i=1) - \p(X_i=-1) (1+\e X_i) \e(X_j|X_i=-1)\\&= \frac{1+\e X_i}{2} (1-\e X_i) \e(X_j|X_i=1) - \frac{1-\e X_i}{2}(1+\e X_i) \e(X_j|X_i=-1)\nonumber\\&= \frac{1}{2} [1-(\e X_i)^2]\left[\e(X_j|X_i=1) - \e(X_j|X_i=-1)\right]\nonumber\\\label{med3}&\le \frac{1}{2} \left[\e(X_j|X_i=1) - \e(X_j|X_i=-1)\right]
	\end{align}
	
	Combining \eqref{med1}, \eqref{med2} and \eqref{med3}, we have for every $i$,
	\begin{equation}\label{med4}
	\sum_{j\in [N]\setminus \{i\}} \left|\mathrm{Cov}(X_i,X_j)\right| \le \frac{4\beta R}{1-4\beta R} \le \frac{\Upsilon}{3}~.
	\end{equation}
	
	Hence, we have by \eqref{med4},
	\begin{align}
	\mathrm{Var}(\bm a^\top \bs) &\ge \sum_{i=1}^N a_i^2 \mathrm{Var}(X_i) - \sum_{i\ne j} \left|a_i a_j \mathrm{Cov}(X_i,X_j)\right|\nonumber\\&\ge \sum_{i=1}^N a_i^2 \mathrm{Var}(X_i) - \sum_{i\ne j} \frac{(a_i^2+a_j^2)\left|\mathrm{Cov}(X_i,X_j)\right|}{2}\nonumber\\&= \sum_{i=1}^N a_i^2\left[\mathrm{Var}(X_i) - \sum_{j\in [N]\setminus \{i\}} \left|\mathrm{Cov}(X_i,X_j)\right| \right]\nonumber\\&\ge\sum_{i=1}^N a_i^2 \left(\Upsilon - \frac{\Upsilon}{3}\right) = \label{c111}\frac{2\Upsilon}{3} \|\bm a\|_2^2\\ &\ge \label{c222}\frac{\|\bm a\|_2^2 \Upsilon^2}{R}\cdot \frac{2R}{3} \gtrsim \frac{\|\bm a\|_2^2 \Upsilon^2}{R}~.
	\end{align}
	
		Now  consider the case $R > \Upsilon/(16\beta)$. By Lemma \ref{combinatorial1}, we choose a subset $I$ of $[N]$ such that conditioned on $\bs_{-I}$, $\bs_I$ is a $(\Upsilon/(16\beta),\Upsilon)$- Ising model, and
		\begin{equation}\label{rtour}
		\|\bm a_I\|_2^2 \ge \frac{\Upsilon}{128R\beta}\|\bm a\|_2^2~.
		\end{equation}
		Hence, we have from \eqref{c111} and \eqref{rtour},
		\begin{equation}\label{pl}
		\mathrm{Var}\left(\bm a^\top \bs | \bs_{-I}\right) = \mathrm{Var}\left(\bm a_I^\top \bs_I | \bs_{-I}\right) \ge \frac{2\Upsilon \|\bm a_I\|_2^2}{3} \ge \frac{\Upsilon^2 \|\bm a\|_2^2}{192 R\beta}~.
		\end{equation}
		Lemma \ref{lemma10} now follows from \eqref{pl} on observing that $\mathrm{Var}(\bm a^\top \bs) \ge \mathbb{E} \left[\mathrm{Var}\left(\bm a^\top \bs | \bs_{-I}\right)\right]$.
	\end{proof}
	
	\begin{lem}\label{freq}
		We have:
		$$\|\bm F \bm A\|_F^2 \ge \|\bm A\|_F^2 - d~.$$
		\end{lem}
	\begin{proof}
		For a matrix $\bm M$, let $\sigma_i(\bm M)$ denote the $i^{\mathrm{th}}$ largest singular value of $\bm M$. Then, by Theorem 2 in \cite{wangfr}, we have:
		\begin{equation}\label{wangeq}
		\sum_{t=1}^N \sigma_t^2(\bm F\bm A) \ge \sum_{t=1}^N \sigma_t^2(\bm F)~ \sigma_{N-t+1}^2(\bm A)
		\end{equation}
		Since $\bm F$ is idempotent with trace $N-d$, it follows that $\sigma_1(\bm F)=\ldots =\sigma_{N-d}(\bm F) = 1$ and $\sigma_{N-d+1}(\bm F)=\ldots =\sigma_{N}(\bm F) = 0$. Hence, we have from \eqref{wangeq},
		\begin{equation}\label{wangeq2}
		\|\bm F \bm A\|_F^2 = \sum_{t=1}^N \sigma_t^2(\bm F\bm A) \ge \sum_{t=1}^{N-d} \sigma_{N-t+1}^2(\bm A) = \|\bm A\|_F^2 - \sum_{i=1}^d \sigma_i^2(\bm A)~.
		\end{equation}
		Lemma \ref{freq} follows from \eqref{wangeq2} on observing that $\|\bm A\|_2 \le 1$, and hence, $\sigma_i^2(\bm A) \le 1$ for all $1\le i\le N$.
	\end{proof}
	
	\begin{lem}\label{math1}
		The following inequality holds for all $x \in \mathbb{R}$:
		$$\frac{e^x}{2\cosh(x)} \ge \frac{1}{2} e^{-2|x|}~.$$
	\end{lem}
\begin{proof}
	First, suppose that $x<0$. Then, $e^{2x}+1 < 2$, and so, we have:
	$$\frac{1}{e^x+e^{-x}} > \frac{1}{2} e^x \implies \frac{e^x}{2\cosh(x)} = \frac{e^x}{e^x+e^{-x}} > \frac{1}{2} e^{2x} = \frac{1}{2} e^{-2|x|}~.$$
	Next, let $x\ge 0$. Then, $e^x \ge e^{-x}$ and hence, 
	$$\frac{e^x}{2\cosh(x)} \ge \frac{1}{2} \ge \frac{1}{2} e^{-2|x|}~.$$ Lemma \eqref{math1} follows. 
\end{proof}

%	\addcontentsline{toc}{chapter}{APPENDIX} % You have to write here everything you want the ToC to show including \thechapter if you want numbering
%	\addtocontents{toc}{\protect\setcounter{tocdepth}{-1}} % This is so sections in the appendix are not shown in the ToC
%	
%	%\subfile{Appendix}
%	
%	\addtocontents{toc}{\protect\setcounter{tocdepth}{1}} % This is to bring things back to normal
\end{appendices}

\iffalse
\begin{appendices}
	\addtocontents{toc}{\protect\setcounter{tocdepth}{0}}
	\chapter{First appendix}
	\section{First section}
	\section{Second section}
	
	\chapter{Second appendix}
	\section{First section}
	\section{Second section}
\end{appendices}
\fi

\end{document}